\documentclass[12pt,oneside,reqno,english]{amsart}

\usepackage[T1]{fontenc}
\usepackage[latin9]{inputenc}
\usepackage{geometry}
\usepackage{color}
\usepackage{mathrsfs}
\usepackage{amsbsy}
\usepackage{amstext}
\usepackage{amsthm}
\usepackage{amssymb}
\usepackage{stmaryrd}
\usepackage{graphicx}
\usepackage{setspace}
\makeatletter
\numberwithin{equation}{section}
\numberwithin{figure}{section}
\theoremstyle{remark}
\newtheorem*{acknowledgement*}{\protect\acknowledgementname}
\theoremstyle{plain}
\newtheorem{thm}{\protect\theoremname}[section]
\theoremstyle{definition}
\newtheorem{xca}[thm]{\protect\exercisename}
\theoremstyle{remark}
\newtheorem{rem}[thm]{\protect\remarkname}
\theoremstyle{plain}
\newtheorem{lem}[thm]{\protect\lemmaname}
\theoremstyle{remark}
\newtheorem{notation}[thm]{\protect\notationname}
\theoremstyle{plain}
\newtheorem{prop}[thm]{\protect\propositionname}
\theoremstyle{definition}
\newtheorem*{example*}{\protect\examplename}
\theoremstyle{definition}
\newtheorem{example}[thm]{\protect\examplename}
\theoremstyle{definition}
\newtheorem{defn}[thm]{\protect\definitionname}
\theoremstyle{plain}
\newtheorem{cor}[thm]{\protect\corollaryname}

\usepackage{lmodern}
\DeclareFontFamily{OMX}{MnSymbolE}{}
\DeclareFontShape{OMX}{MnSymbolE}{m}{n}{
    <-6>  MnSymbolE5
   <6-7>  MnSymbolE6
   <7-8>  MnSymbolE7
   <8-9>  MnSymbolE8
   <9-10> MnSymbolE9
  <10-12> MnSymbolE10
  <12->   MnSymbolE12}{}
\DeclareSymbolFont{mnlargesymbols}{OMX}{MnSymbolE}{m}{n}
\SetSymbolFont{mnlargesymbols}{bold}{OMX}{MnSymbolE}{b}{n}
\DeclareMathDelimiter{\llangle}{\mathopen}{mnlargesymbols}{'164}{mnlargesymbols}{'164}
\DeclareMathDelimiter{\rrangle}{\mathclose}{mnlargesymbols}{'171}{mnlargesymbols}{'171}

\makeatother

\usepackage{babel}
\providecommand{\acknowledgementname}{Acknowledgement}
\providecommand{\corollaryname}{Corollary}
\providecommand{\definitionname}{Definition}
\providecommand{\examplename}{Example}
\providecommand{\exercisename}{Exercise}
\providecommand{\lemmaname}{Lemma}
\providecommand{\notationname}{Notation}
\providecommand{\propositionname}{Proposition}
\providecommand{\remarkname}{Remark}
\providecommand{\theoremname}{Theorem}

\begin{document}
\title[Generalized Dirichlet and Thomson Principles]{Generalized Dirichlet and Thomson Principles and Their Applications}
\author{Insuk Seo }
\address{Department of Mathematical Sciences and R.I.M, Seoul National University,
Republic of Korea. }
\email{insuk.seo@snu.ac.kr}

\maketitle
\tableofcontents

\section*{Introduction}

This lecture note is intended to introduce the recently-developed
potential theory for the non-reversible Markov processes and to explain
applications of this new theory to the study of metastability of huge
stochastic interacting systems.

Regarding irreducible Markov processes, it is well-known that the
distribution of the process at time $t>0$ converges to its unique
invariant measure as $t\rightarrow\infty$, regardless of its starting
distribution, and this asymptotic behavior is called the mixing property
of Markov processes. The speed of this convergence is one of the main
concerns in the study of Markov processes, as it is related to a multitude
of important problems such as the performance of Markov chain Monte
Carlo algorithm, equilibration of non-equilibrium physical systems,
and metastability of random dynamics.

In the study of the mixing property of Markov processes, one of the
most useful tools is potential theory, especially the quantity called
\textit{capacity} with respect to the Markov process under consideration.
Capacity is measured for two disjoint subsets of the state space of
the Markov process, and it is inversely related to how well the corresponding
Markov process commutes between these two disjoint sets. Since the
convergence explained above will take a long time if the Markov process
cannot quickly commute between two large (with respect to the invariant
measure) se\textcolor{black}{ts, ca}pacity is a useful notion in the
analysis of mixing properties.

Classic potential theory is developed only when the underlying Markov
process is reversible with respect to its invariant measure, and has
been widely used in the study of the mixing property of Markov processes
(e.g., \cite{B-denH} or \cite[Chapters 9, 10]{L-P-W}). In the potential
theory of reversible Markov processes, the so-called Dirichlet and
Thomson principles provide a robust way of estimating the capacity
via construction of a test function or a test flow.

Potential theory for non-reversible processes has been developed very
recently. In particular, \cite{G-L} and \cite{S} established the
Dirichlet and Thomson principles for non-reversible Markov processes,
respectively. These formulae are far more involved than the corresponding
principles for the reversible processes, and technical difficulties
arise in the application of these principles. To minimize these technical
issues, a more generalized version of the Dirichlet and Thomson principles
were developed in \cite{L-Mariani-Seo,Seo NRZRP}. In the first part
of the current note, we give a comprehensive review on these recent
developments in the potential theory of non-reversible Markov processes
based on \cite{G-L,Seo NRZRP,S}.

In the second and third parts of this note, we explain two applications
of the recently-developed potential theory to the study of metastability.
The metastability is a ubiquitous phenomenon appearing when a Markov
process possesses a poor mixing property because of the existence
of multiple locally stable sets, or metastable sets. For example,
metastability occurs for the models such as
\begin{itemize}
\item small random perturbations of dynamical systems (e.g., \cite{B-E-G-K,B-G-K,F-W,L-Mariani-Seo,L-Seo NR2,Lee-Seo1,Lee-Seo2,L-M,R-S}),
\item interacting particle systems with condensing phenomena (e.g., \cite{B-L ZRP,B-D-G,Kim,Kim-Seo1,L-M-Seo1,L-M-Seo2,O-R,Seo NRZRP}),
and
\item stochastic spin systems in the low-temperature regime (e.g., \cite{A-C,BA-C,B-B-I,B-denH,B-denH-N,B-denH-S,B-M,C-G-O-V,C-O,L-L,L-L 19,L-Seo Potts,Lee,N-S Ising1,N-S Ising2,N-Z}).
\end{itemize}
Readers are referred to monographs \cite{B-denH,O-V} for more comprehensive
discussions regarding the mathematical study of metastability.

The potential theory plays a crucial role in the rigorous analysis
of metastability. In particular, two representative ways of quantitatively
analyzing the metastable behavior are the Eyring--Kramers law \cite{Eyr,Kra}
and Markov chain model reduction \cite{B-L TM,B-L TM2,B-L MG,L-M-Seo2}.

The Eyring--Kramers law describes the precise asymptotics of the
mean transition time from a metastable set to other metastable sets.
Since such a transition between metastable sets is the signature behavior
of metastability, the Eyring--Kramers law is clearly a crucial problem.
A robust methodology to prove the Eyring--Kramers law based on the
potential theory (known as the potential-theoretic approach) is developed
in \cite{B-E-G-K}. We refer to the monograph \cite{B-denH} for a
comprehensive review on this approach. \textit{In Part 2, we derive
the Eyring--Kramers law for a stochastic spin system known as the
Ising model on a large, finite two-dimensional lattice without external
field as an application of the potential theory explained in Part
1. }This part is largely based on the recent article \cite{Kim-Seo Potts}.
We remark that the article \cite{Kim-Seo Potts} addresses more general
situations. This article not only considers the Ising model on a two-dimensional
lattice but also the Potts model (which is a generalization of the
Ising model) on two- and three-dimensional lattices. In particular,
the three-dimensional model is more cumbersome for carrying out rigorous
analyses. Moreover, this article not only concerns the Eyring--Kramers
law but also the precise analyses of the energy landscape and the
typical path of transitions. In this note, we only focus on the Eyring--Kramers
law for the two-dimensional model to convey the overall idea. For
interested readers, we refer to the article \cite{Kim-Seo Potts}
for more comprehensive results.

If there are several metastable sets and the transitions between them
take place successively, it is tempting to analyze these successive
transitions all at once. A natural way of carrying this out is to
approximately describe, after a suitable time-rescaling, the successive
transitions between metastable sets as a Markov chain whose state
space consists of metastable sets of the original Markov process.
This methodology for describing the metastable behavior is a special
case of the Markov chain model reduction. A robust methodology for
the verification of this Markov chain model reduction based on potential
theory has been developed in \cite{B-L TM,B-L TM2,B-L MG}, and this
method is called the martingale approach. \textit{In Part 3, we combine
this approach and the potential theory for non-reversible processes
to analyze the metastable behavior of non-reversible zero-range processes.}
This part is largely based on the recent article \cite{Seo NRZRP}.
For conciseness of the discussion, we only consider the asymmetric
nearest neighbor random walk on a cycle, but the discussion given
here can be applied to the general model; we refer to \cite{Seo NRZRP}
for the interested readers.
\begin{acknowledgement*}
This lecture note is written with the support of the Sangsan Lecture
Note fund of the Research Institute of Mathematics of the Seoul National
University. The contents of the lecture note have been developed with
the support of the National Research Foundation of Korea (NRF) grant
funded by the Korean government (MSIT) (No. 2017R1A5A1015626 and No.
2018R1C1B6006896). The author thanks Seonwoo Kim and Jungkyoung Lee
for careful reading of the early version of the note and for helping
to clarify the presentation.
\end{acknowledgement*}
\newpage

\part{Potential Theory}

In the first part, we review the potential theory of continuous-time
Markov processes, introduce the Dirichlet and Thomson principles,
and then finally explain the generalized Dirichlet and Thomson principles
developed in \cite{Seo NRZRP}. Although we explain the whole theory
in the context of continuous-time Markov processes for the convenience
of the discussion, the corresponding results are also valid for discrete-time
Markov chains or diffusion processes. For the discussion of diffusion
processes, we refer to \cite{L-Mariani-Seo}.

\section{Potential Theory of Markov Processes}

\subsection{Markov processes}

We start by introducing several relevant notions regarding a continuous-time
Markov process $(X(t))_{t\ge0}$ on a finite set $\mathcal{H}$.

\subsubsection*{Continuous-time Markov processes}

For $x\in\mathcal{H}$, we denote by $\mathbb{P}_{x}$ the law of
the process $X(\cdot)$ starting from $x$, and by $\mathbb{E}_{x}$
the expectation with respect to $\mathbb{P}_{x}$. We assume that
the process $X(\cdot)$ is irreducible, in the sense that for all
$x,\,y\in\mathcal{H}$\footnote{In this lecture note, writing $a,\,b\in A$ always implies that $a$
and $b$ are \textit{different} elements of a set $A$.},
\[
\mathbb{P}_{x}[X(t)=y\text{ for some }t>0]=1\;.
\]

We denote by $r:\mathcal{H}\times\mathcal{H}\rightarrow[0,\,\infty)$
the jump rate of the Markov process $X(\cdot)$. Namely, for $x,\,y\in S$,
the quantity $r(x,\,y)\ge0$ represents the rate of the jump from
$x$ to $y$ for the Markov process $X(\cdot)$. For convenience,
we set $r(x,\,x)=0$ for all $x\in\mathcal{H}$. Denote by
\begin{equation}
\lambda(x)=\sum_{y\in\mathcal{H}}r(x,\,y)\;\;\;\;;\;x\in\mathcal{H}\label{holding}
\end{equation}
the holding rate of the process $X(\cdot)$ at $x$. Then, the dynamics
$X(\cdot)$ can be described as follows: if $X(t)=x$, then the process
waits for an exponential time of mean $\lambda(x)^{-1}$. Then, it
jumps to $y\in\mathcal{H}$ with probability $r(x,\,y)/\lambda(x)$.

\subsubsection*{Embedded chain}

We denote by $(\widehat{X}(n))_{n\in\mathbb{Z}^{+}}$ where $\mathbb{Z}^{+}=\mathbb{Z}\cap[0,\,\infty)$
the discrete-time Markov chain with jump probability $p(x,\,y)=r(x,\,y)/\lambda(x)$.
This chain is referred to as the embedded chain of $X(\cdot)$, and
represents the jumping dynamics (irrespective of the exponential waiting
time between successive jumps) of $X(\cdot)$. For $x\in\mathcal{H}$,
denote by $\widehat{\mathbb{P}}_{x}$ the law of the embedded chain
$\widehat{X}(\cdot)$ starting from $x$, and by $\widehat{\mathbb{E}}_{x}$
the expectation with respect to $\widehat{\mathbb{P}}_{x}$.

\subsubsection*{Invariant measure and reversibility}

By irreducibility of the process $X(\cdot)$, there exists a unique
probability distribution $\mu(\cdot)$ on $\mathcal{H}$ that satisfies
\begin{equation}
\sum_{x\in\mathcal{H}}\mu(x)r(x,\,y)=\sum_{x\in\mathcal{H}}\mu(y)r(y,\,x)\;.\label{1inv1}
\end{equation}
One can readily infer from the irreducibility that
\begin{equation}
\mu(x)>0\text{ for all }x\in\mathcal{H}\;.\label{inv2}
\end{equation}

\begin{xca}
Suppose that the Markov process $X(\cdot)$ is irreducible. Prove
that there exists a unique probability distribution $\mu(\cdot)$
on $\mathcal{H}$ satisfying (\ref{1inv1}). Then, prove that this
unique $\mu(\cdot)$ satisfies (\ref{inv2}).
\end{xca}

The distribution $\mu(\cdot)$ is called the\textit{ invariant (or
stationary) distribution} since the marginal distribution of the process
$X(\cdot)$ at any later time $t>0$ is $\mu$, provided $X(0)$ is
distributed according to $\mu$. We say that the process $X(\cdot)$
is \textit{reversible} if the following detailed balance condition
holds:
\begin{equation}
\mu(x)r(x,\,y)=\mu(y)r(y,\,x)\;\;\;\text{for all }x,\,y\in\mathcal{H}\;.\label{1detbal}
\end{equation}
Note that (\ref{1detbal}) immediately implies (\ref{1inv1}). Such
a process is called reversible since the time-reversed process has
the same law with the original process. If the process $X(\cdot)$
is not reversible, it is called a non-reversible or irreversible process.

In addition, we can readily check that a measure $M(\cdot)$ on $\mathcal{H}$
given by
\begin{equation}
M(x)=\lambda(x)\mu(x)\;\;\;\;;\;x\in\mathcal{H}\label{inv_emb}
\end{equation}
is an invariant measure (not necessarily a probability measure) for
the embedded chain $\widehat{X}(\cdot)$. Moreover, the chain $\widehat{X}(\cdot)$
is reversible, i.e., $M(x)p(x,\,y)=M(y)p(y,\,x)$ for all $x,\,y\in\mathcal{H}$,
if and only if the original process $X(\cdot)$ is reversible.

\subsubsection*{Generator and Dirichlet form}

The generator $\mathscr{L}$ associated with the process $X(\cdot)$
is an operator acting on each function $f:\mathcal{H}\rightarrow\mathbb{R}$
in a way that
\[
(\mathscr{L}f)(x)=\sum_{y\in\mathcal{H}}r(x,\,y)(f(y)-f(x))\;\;\;;\;x\in\mathcal{H}\;.
\]
Namely, $\mathscr{L}f$ is another real function on $\mathcal{H}$.
We denote by $L^{2}(\mu)$ the $L^{2}$ space of real functions on
$\mathcal{H}$ with respect to the measure $\mu$. Since $\mathcal{H}$
is a finite set, the space $L^{2}(\mu)$ is merely a collection of
all real functions on $\mathcal{H}$. \footnote{Of course, this is no longer true if we consider the diffusion case.}
Denote by $\left\langle \cdot,\,\cdot\right\rangle _{\mu}$ the inner
product on $L^{2}(\mu)$, i.e., for $f,\,g:\mathcal{H}\rightarrow\mathbb{R}$,
\[
\left\langle f,\,g\right\rangle _{\mu}=\sum_{x\in\mathcal{H}}f(x)g(x)\mu(x)\;.
\]
The Dirichlet form associated to the process $X(\cdot)$ is defined
by, for $f:\mathcal{H}\rightarrow\mathbb{R}$,
\begin{equation}
\mathscr{D}(f)=\left\langle f,\,-\mathscr{L}f\right\rangle _{\mu}\;.\label{e_Df}
\end{equation}
This plays an important role in the potential theory. By the summation
of parts and (\ref{1inv1}), we can write
\begin{equation}
\mathscr{D}(f)=\frac{1}{2}\sum_{x\in\mathcal{H}}\sum_{y\in\mathcal{H}}\mu(x)r(x,\,y)[f(y)-f(x)]^{2}\;.\label{DF}
\end{equation}

We note that the analyses of the reversible process are far more convenient
than those of the non-reversible one, mainly because the operator
$\mathscr{L}$ is self-adjoint in the space $L^{2}(\mu)$ in the sense
that, for all $f,\,g:\mathcal{H}\rightarrow\mathbb{R}$,
\[
\left\langle f,\,\mathscr{L}g\right\rangle _{\mu}=\left\langle \mathscr{L}f,\,g\right\rangle _{\mu}\;.
\]
By the summation by parts and (\ref{1detbal}), we can check that
both sides of the previous identity equal
\[
-\frac{1}{2}\sum_{x\in\mathcal{H}}\sum_{y\in\mathcal{H}}\mu(x)r(x,\,y)[f(y)-f(x)][g(y)-g(x)]\;.
\]

\subsubsection*{Adjoint process}

For the non-reversible case, we define the \textit{adjoint process}
$(X^{\dagger}(t))_{t\ge0}$, which is another continuous-time Markov
process on $\mathcal{H}$ with rate
\[
r^{\dagger}(x,\,y)=\frac{\mu(y)r(y,\,x)}{\mu(x)}\;\;\;;\;x,\,y\in\mathcal{H}\;.
\]
We shall denote by $\mathbb{P}_{x}^{\dagger}$ the law of the adjoint
process $X^{\dagger}(\cdot)$ starting from $x$, and by $\mathbb{E}_{x}^{\dagger}$
the expectation with respect to $\mathbb{P}_{x}^{\dagger}$.

The process $X^{\dagger}(\cdot)$ is a time-reversed process of $X(\cdot)$,
and we can notice from (\ref{1detbal}) that $X^{\dagger}(\cdot)$
is defined by the same law with $X(\cdot)$ in the reversible case;
hence the time-reversing does not change the law. We define the generator
for the adjoint process $X^{\dagger}(\cdot)$ as, for $f:\mathcal{H}\rightarrow\mathbb{R}$,
\[
(\mathscr{L}^{\dagger}f)(x)=\sum_{y\in\mathcal{H}}r^{\dagger}(x,\,y)(f(y)-f(x))\;\;\;;\;x\in\mathcal{H}\;.
\]
The importance of the adjoint process in the context of the potential
theory follows from the fact that $\mathscr{L}^{\dagger}$ is indeed
the adjoint operator of $\mathscr{L}$ in the sense that, for all
$f,\,g:\mathcal{H}\rightarrow\mathbb{R}$,
\begin{equation}
\left\langle f,\,\mathscr{L}g\right\rangle _{\mu}=\left\langle \mathscr{L}^{\dagger}f,\,g\right\rangle _{\mu}\;.\label{adjrel}
\end{equation}

\begin{xca}
\label{ex21}
\begin{enumerate}
\item Verify (\ref{adjrel}).
\item Prove that $\left\langle f,\,\mathscr{L}g\right\rangle _{\mu}=0$
if $f$ is a constant function. In particular, for any $g:\mathcal{H}\rightarrow\mathbb{R}$,
we have
\[
\sum_{x\in\mathcal{H}}\mu(x)(\mathscr{L}g)(x)=0\;.
\]
\end{enumerate}
\end{xca}

\begin{rem}
\label{remgf}Inserting $g=-f$ at (\ref{adjrel}), we can observe
that the Dirichlet form for the adjoint process is also given as $\mathscr{D}(\cdot)$.
\end{rem}

We can also consider the embedded chain of the adjoint process. Write
$\widehat{X}^{\dagger}(\cdot)$ the embedded chain with respect to
the process $X^{\dagger}(\cdot)$. One can readily verify that the
jump rate $p^{\dagger}(\cdot,\,\cdot)$ of the chain $\widehat{X}^{\dagger}(\cdot)$
is given by
\begin{equation}
p^{\dagger}(x,\,y)=\frac{M(y)p(y,\,x)}{M(x)}\;\;\;\;;\;x,\,y\in\mathcal{H}\;,\label{e_225}
\end{equation}
and furthermore $M(\cdot)$ is again the invariant measure for the
process $\widehat{X}^{\dagger}(\cdot)$. Similarly, we denote by $\widehat{\mathbb{P}}_{x}^{\dagger}$
the law of the process $\widehat{X}^{\dagger}(\cdot)$ starting at
$x\in\mathcal{H}$, and by $\mathbb{\widehat{E}}_{x}^{\dagger}$ the
expectation with respect to $\widehat{\mathbb{P}}_{x}^{\dagger}$.

\subsection{\label{sec_pot}Equilibrium potential and capacity}

Two crucial notions in the potential theory of Markov processes are
the equilibrium potential and the capacity. In this section, we define
these objects and review their elementary properties.

\subsubsection*{Equilibrium potential}

For $\mathcal{A}\subset\mathcal{H}$, we denote by $\tau_{\mathcal{A}}$
the hitting time of the set $\mathcal{A}$:
\[
\tau_{\mathcal{A}}=\inf\{t\ge0:X(t)\in\mathcal{A}\}\;.
\]
For two non-empty and disjoint subsets $\mathcal{A}$ and $\mathcal{B}$
of $\mathcal{H}$, we define the equilibrium potential between $\mathcal{A}$
and $\mathcal{B}$ with respect to the process $X(\cdot)$ as a function
$h_{\mathcal{A},\,\mathcal{B}}:\mathcal{H}\rightarrow[0,\,1]$ defined
by
\[
h_{\mathcal{A},\,\mathcal{B}}(x)=\mathbb{P}_{x}[\tau_{\mathcal{A}}<\tau_{\mathcal{B}}]\;\;\;\;;\;x\in\mathcal{H}\;.
\]
By definition, it is clear that
\begin{equation}
h_{\mathcal{B},\,\mathcal{A}}=1-h_{\mathcal{A},\,\mathcal{B}}\;.\label{hba}
\end{equation}
The following lemma gives the basic properties of the equilibrium
potential $h_{\mathcal{A},\,\mathcal{B}}$.
\begin{lem}
\label{lem21}For two non-empty and disjoint subsets $\mathcal{A}$
and $\mathcal{B}$ of $\mathcal{H}$, the equilibrium potential $h_{\mathcal{A},\,\mathcal{B}}$
satisfies
\begin{equation}
\begin{cases}
h_{\mathcal{A},\,\mathcal{B}}\equiv1 & \text{on }\mathcal{A}\;,\\
h_{\mathcal{A},\,\mathcal{B}}\equiv0 & \text{on }\mathcal{B}\;,\;\text{and}\\
\mathscr{L}h_{\mathcal{A},\,\mathcal{B}}\equiv0 & \text{on }(\mathcal{A}\cup\mathcal{B})^{c}=\mathcal{H}\setminus(\mathcal{A}\cup\mathcal{B})\;.
\end{cases}\label{e_eq_prop}
\end{equation}
\end{lem}

\begin{proof}
The first two properties are evident from the definition of $h_{\mathcal{A},\,\mathcal{B}}$.
Let us focus on the last one. Fix $x\in(\mathcal{A}\cup\mathcal{B})^{c}$.
Then, since the process $X(\cdot)$ starting at $x$ jumps to $y$
with probability $r(x,\,y)/\lambda(x)$, by the Markov property we
can write
\[
h_{\mathcal{A},\,\mathcal{B}}(x)=\mathbb{P}_{x}[\tau_{\mathcal{A}}<\tau_{\mathcal{B}}]=\sum_{y\in\mathcal{H}}\frac{r(x,\,y)}{\lambda(x)}\mathbb{P}_{y}[\tau_{\mathcal{A}}<\tau_{\mathcal{B}}]=\sum_{y\in\mathcal{H}}\frac{r(x,\,y)}{\lambda(x)}h_{\mathcal{A},\,\mathcal{B}}(y)\;.
\]
Multiplying both sides by $\lambda(x)$ and reorganizing give us $\mathscr{L}h_{\mathcal{A},\,\mathcal{B}}(x)=0$.
\end{proof}
\begin{rem}
\label{rem22}Of course, we can define the equilibrium potential $h_{\mathcal{A},\,\mathcal{B}}^{\dagger}:\mathcal{H}\rightarrow[0,\,1]$
with respect to the adjoint process $X^{\dagger}(\cdot)$. Then, an
analogue of Lemma \ref{lem21} holds for $h_{\mathcal{A},\,\mathcal{B}}^{\dagger}$.
It suffices to replace the last property of (\ref{e_eq_prop}) with
$\mathscr{L}^{\dagger}h_{\mathcal{A},\,\mathcal{B}}^{\dagger}\equiv0$
on $(\mathcal{A}\cup\mathcal{B})^{c}$.
\end{rem}

\subsubsection*{Capacity}

For two non-empty and disjoint subsets $\mathcal{A}$ and $\mathcal{B}$
of $\mathcal{H}$, we define the capacity between $\mathcal{A}$ and
$\mathcal{B}$ with respect to the process $X(\cdot)$ as
\begin{equation}
\textup{cap}(\mathcal{A},\,\mathcal{B})=\mathscr{D}(h_{\mathcal{A},\,\mathcal{B}})\;.\label{e_capAB}
\end{equation}
By the expression (\ref{DF}) of the Dirichlet form and (\ref{hba}),
it holds that
\begin{equation}
\textup{cap}(\mathcal{A},\,\mathcal{B})=\mathscr{D}(h_{\mathcal{A},\,\mathcal{B}})=\mathscr{D}(h_{\mathcal{B},\,\mathcal{A}})=\textup{cap}(\mathcal{B},\,\mathcal{A})\;.\label{e_cap_synm}
\end{equation}

\begin{notation}
\label{not_single}If $\mathcal{A}=\{a\}$ or $\mathcal{B}=\{b\}$
(or both), we simply write $a$ or $b$ instead of $\{a\}$ or $\{b\}$,
respectively, in the subscript of $h_{\mathcal{A},\,\mathcal{B}}$
and $\textup{cap}(\mathcal{A},\,\mathcal{B})$. For instance, if $\mathcal{A}=\{a\}$
and $\mathcal{B}=\{b\}$, we write $h_{a,\,b}$ and $\textup{cap}(a,\,b)$,
instead of $h_{\{a\},\,\{b\}}$ and $\textup{cap}(\{a\},\,\{b\})$,
respectively.
\end{notation}

\begin{xca}
\label{ex26}Let $\mathcal{H}=\mathbb{T}_{N}=\mathbb{Z}/(N\mathbb{Z})$($=\mathbb{Z}_{N}$)
be a discrete torus of length $N$ (i.e., a cycle of length $N$).
Define a rate as
\[
r(x,\,y)=\begin{cases}
p & \text{if }x-y\equiv1\;(\text{mod }N)\;,\\
1-p & \text{if }x-y\equiv-1\;(\text{mod }N)\;,\\
0 & \text{otherwise\;,}
\end{cases}
\]
for some $p\in[0,\,1]$. For the Markov process $X(\cdot)$ on $\mathbb{T}_{N}$
with rate $r(\cdot,\,\cdot)$, answer the following questions.
\begin{enumerate}
\item Prove that the uniform measure $\mu(\cdot)$ on $\mathbb{T}_{N}$,
namely,
\[
\mu(x)=\frac{1}{N}\;\;\;\;\text{for all }x\in\mathbb{T}_{N}\;,
\]
is the unique invariant measure for the process $X(\cdot)$, and moreover
that the process $X(\cdot)$ is reversible if and only if $p=1/2$.
\item For $x,\,y\in\mathbb{T}_{N}$, compute $\textup{cap}(x,\,y)$. (cf.
Notation \ref{not_single})
\item For any non-empty and disjoint subsets $\mathcal{A}$ and $\mathcal{B}$
of $\mathbb{T}_{N}$, compute $\textup{cap}(\mathcal{A},\,\mathcal{B})$.
\end{enumerate}
\end{xca}

Next, we introduce an alternative expression for the capacity that
turns out to play an important role in using the capacity in various
instances. We write $\tau_{\mathcal{A}}^{+}$ for the return time
to the set $\mathcal{A}$:
\[
\tau_{\mathcal{A}}^{+}=\inf\{t>0:X(t)\in\mathcal{A}\text{ and }X(s)\neq X(0)\text{ for some }s\in[0,\,t]\}\;.
\]
Namely, this time expresses the first time at which $X(t)$ arrives
at $\mathcal{A}$ after leaving its initial location. In particular,
if the process starts from $x\notin\mathcal{A}$, we have $\tau_{\mathcal{A}}^{+}=\tau_{\mathcal{A}}$.
Recall the measure $M(\cdot)$ from (\ref{inv_emb}).
\begin{lem}
\label{lem_cap_alt}For two non-empty and disjoint subsets $\mathcal{A}$
and $\mathcal{B}$ of $\mathcal{H}$, it holds that
\[
\textup{cap}(\mathcal{A},\,\mathcal{B})=\sum_{x\in\mathcal{A}}M(x)\mathbb{P}_{x}[\tau_{\mathcal{B}}<\tau_{\mathcal{A}}^{+}]\;.
\]
\end{lem}

\begin{proof}
By (\ref{e_Df}) and (\ref{e_capAB}), we can write
\[
\textup{cap}(\mathcal{A},\,\mathcal{B})=\left\langle h_{\mathcal{A},\,\mathcal{B}},\,-\mathscr{L}h_{\mathcal{A},\,\mathcal{B}}\right\rangle _{\mu}=\sum_{x\in\mathcal{H}}h_{\mathcal{A},\,\mathcal{B}}(x)\,(-\mathscr{L}h_{\mathcal{A},\,\mathcal{B}})(x)\,\mu(x)\;.
\]
By (\ref{e_eq_prop}), we have $h_{\mathcal{A},\,\mathcal{B}}(x)\,(\mathscr{L}h_{\mathcal{A},\,\mathcal{B}})(x)=0$
for all $x\notin\mathcal{A}$, and thus we can write
\[
\textup{cap}(\mathcal{A},\,\mathcal{B})=\sum_{x\in\mathcal{A}}(-\mathscr{L}h_{\mathcal{A},\,\mathcal{B}})(x)\,\mu(x)\;.
\]
Note that we used the fact that $h_{\mathcal{A},\,\mathcal{B}}\equiv1$
on $\mathcal{A}$. By the definition of the generator and (\ref{hba}),
we can further write
\begin{align}
\textup{cap}(\mathcal{A},\,\mathcal{B}) & =\sum_{x\in\mathcal{A}}\sum_{y\in\mathcal{H}}\mu(x)r(x,\,y)[h_{\mathcal{A},\,\mathcal{B}}(x)-h_{\mathcal{A},\,\mathcal{B}}(y)]\nonumber \\
 & =\sum_{x\in\mathcal{A}}\sum_{y\in\mathcal{H}}\mu(x)r(x,\,y)[1-h_{\mathcal{A},\,\mathcal{B}}(y)]=\sum_{x\in\mathcal{A}}\sum_{y\in\mathcal{H}}\mu(x)r(x,\,y)h_{\mathcal{B},\,\mathcal{A}}(y)\;.\label{e_l241}
\end{align}
On the other hand, by the Markov property, for $x\in\mathcal{A}$
we have
\begin{equation}
\mathbb{P}_{x}[\tau_{\mathcal{B}}<\tau_{\mathcal{A}}^{+}]=\sum_{y\in\mathcal{H}}p(x,\,y)\mathbb{P}_{y}[\tau_{\mathcal{B}}<\tau_{\mathcal{A}}]=\sum_{y\in\mathcal{H}}p(x,\,y)h_{\mathcal{B},\,\mathcal{A}}(y)\;.\label{e_l242}
\end{equation}
Since $\mu(x)r(x,\,y)=M(x)p(x,\,y)$, we can complete the proof from
(\ref{e_l241}) and (\ref{e_l242}).
\end{proof}
The capacity with respect to the adjoint process $X^{\dagger}(\cdot)$
is given by (cf. Remark \ref{remgf})
\begin{equation}
\textup{cap}^{\dagger}(\mathcal{A},\,\mathcal{B})=\mathscr{D}(h_{\mathcal{A},\,\mathcal{B}}^{\dagger})\;.\label{e_CapABda}
\end{equation}
Then, by the same reasoning as above, it holds that $\textup{cap}^{\dagger}(\mathcal{A},\,\mathcal{B})=\textup{cap}^{\dagger}(\mathcal{B},\,\mathcal{A})$.

Now, we give two important properties of the capacity based on Lemma
\ref{lem_cap_alt}. The first is a somewhat unexpected property in
view of the definitions (\ref{e_capAB}) and (\ref{e_CapABda}) of
capacities.
\begin{prop}
\label{prop_cap_sym}For two non-empty and disjoint subsets $\mathcal{A}$
and $\mathcal{B}$ of $\mathcal{H}$, it holds that
\[
\textup{cap}(\mathcal{A},\,\mathcal{B})=\textup{cap}^{\dagger}(\mathcal{A},\,\mathcal{B})\;.
\]
\end{prop}

\begin{proof}
We first claim that, for all $x\in\mathcal{A}$ and $y\in\mathcal{B}$,
\begin{equation}
M(x)\mathbb{P}_{x}\left[\tau_{\mathcal{B}}<\tau_{\mathcal{A}}^{+},\,\tau_{\mathcal{B}}=\tau_{y}\right]=M(y)\mathbb{P}_{y}^{\dagger}\left[\tau_{\mathcal{A}}<\tau_{\mathcal{B}}^{+},\,\tau_{\mathcal{A}}=\tau_{x}\right]\;.\label{e2261}
\end{equation}
To prove this, we write the left-hand side as
\begin{equation}
\sum_{T=1}^{\infty}\sum_{(\omega_{t})_{t=0}^{T}:\omega_{0}=x,\,\omega_{T}=y}M(x)\prod_{t=0}^{T-1}p(\omega_{t},\,\omega_{t+1})\;,\label{e226}
\end{equation}
where the summation is carried out for the paths $(\omega_{t})_{t=0}^{T}$
such that $p(\omega_{t},\,\omega_{t+1})>0$ for all $t\in\llbracket0,\,T-1\rrbracket$\footnote{Here, for integers $a$ and $b$, $\llbracket a,\,b\rrbracket$ denotes
$[a,\,b]\cap\mathbb{Z}$.} and $\omega_{t}\notin\mathcal{A}\cup\mathcal{B}$ for all $t\in\llbracket1,\,T-1\rrbracket$.
By (\ref{e_225}), we have
\[
M(x)\prod_{t=0}^{T-1}p(\omega_{t},\,\omega_{t+1})=M(y)\prod_{t=0}^{T-1}p^{\dagger}(\omega_{t+1},\,\omega_{t})\;.
\]
Therefore, we can rewrite (by reversing the path) (\ref{e226}) as
\[
\sum_{T=1}^{\infty}\sum_{(\omega_{t})_{t=0}^{T}:\omega_{0}=y,\,\omega_{T}=x}M(y)\prod_{t=0}^{T-1}p^{\dagger}(\omega_{t},\,\omega_{t+1})\;,
\]
where the summation is carried out for the paths $(\omega_{t})_{t=0}^{T}$
such that $p^{\dagger}(\omega_{t},\,\omega_{t+1})>0$ for all $t\in\llbracket0,\,T-1\rrbracket$
and $\omega_{t}\notin\mathcal{A}\cup\mathcal{B}$ for all $t\in\llbracket1,\,T-1\rrbracket$.
By the same reasoning as above, this corresponds to the right-hand
side of (\ref{e2261}). Hence, we have proved (\ref{e2261}).

Therefore, by Lemma \ref{lem_cap_alt},
\begin{align*}
\textup{cap}(\mathcal{A},\,\mathcal{B}) & =\sum_{x\in\mathcal{A}}\sum_{y\in\mathcal{B}}M(x)\mathbb{P}_{x}\left[\tau_{\mathcal{B}}<\tau_{\mathcal{A}}^{+},\,\tau_{B}=\tau_{y}\right]\\
 & =\sum_{x\in\mathcal{A}}\sum_{y\in\mathcal{B}}M(y)\mathbb{P}_{y}^{\dagger}\left[\tau_{\mathcal{A}}<\tau_{\mathcal{B}}^{+},\,\tau_{\mathcal{A}}=\tau_{x}\right]\\
 & =\sum_{y\in\mathcal{B}}M(y)\mathbb{P}_{y}^{\dagger}\left[\tau_{\mathcal{A}}<\tau_{\mathcal{B}}^{+}\right]=\textup{cap}^{\dagger}(\mathcal{B},\,\mathcal{A})\;.
\end{align*}
Now, it suffices to recall (\ref{e_cap_synm}).
\end{proof}
\begin{prop}
\label{prop:inc}Suppose that $\mathcal{A}'$ and $\mathcal{B}'$
are non-empty disjoint subsets of $\mathcal{H}$. Let $\mathcal{A}$
and $\mathcal{B}$ be non-empty subsets of $\mathcal{A}'$ and $\mathcal{B}'$,
respectively. Then, it holds that
\begin{equation}
\textup{cap}(\mathcal{A},\,\mathcal{B})\le\textup{cap}(\mathcal{A}',\,\mathcal{B}')\;.\label{eq:inc}
\end{equation}
\end{prop}

\begin{proof}
It suffices to prove that, the capacity is monotone in the second
argument, i.e.,
\begin{equation}
\textup{cap}(\mathcal{A},\,\mathcal{B})\le\textup{cap}(\mathcal{A},\,\mathcal{B}')\;,\label{inc1}
\end{equation}
since by the symmetry (\ref{e_cap_synm}), we can proceed as
\[
\textup{cap}(\mathcal{A},\,\mathcal{B})\le\textup{cap}(\mathcal{A},\,\mathcal{B}')=\textup{cap}(\mathcal{B}',\,\mathcal{A})\le\textup{cap}(\mathcal{B}',\,\mathcal{A}')=\textup{cap}(\mathcal{A}',\,\mathcal{B}')\;,
\]
provided that we have (\ref{inc1}).

Now, let us prove (\ref{inc1}). By Lemma \ref{lem_cap_alt}, it suffices
to prove
\[
\sum_{x\in\mathcal{A}}M(x)\mathbb{P}_{x}[\tau_{\mathcal{B}}<\tau_{\mathcal{A}}^{+}]\le\sum_{x\in\mathcal{A}}M(x)\mathbb{P}_{x}[\tau_{\mathcal{B}'}<\tau_{\mathcal{A}}^{+}]\;.
\]
Since $\mathcal{B}\subset\mathcal{B}'$, we trivially have $\mathbb{P}_{x}[\tau_{\mathcal{B}}<\tau_{\mathcal{A}}^{+}]\le\mathbb{P}_{x}[\tau_{\mathcal{B}'}<\tau_{\mathcal{A}}^{+}]$.
\end{proof}
In the investigation of the mixing property of Markov processes, use
of the capacity defined above is crucial, and its (more of less accurate)
estimation is required. The definition of the capacity given above
is easy to understand, but it is not suitable for the estimation.
Instead, the variational expression known as the Dirichlet and Thomson
principles are typically used in the estimation of the capacity. The
remainder of Part 1 is devoted to explain this strategy.

To explore this advanced strategy to estimate the capacity, we need
to reinterpret the capacity in the context of flow structure explained
below. We refer to \cite{G-L,S,L-Seo NR1} for more comprehensive
discussions on the flow structure of Markov processes, and to \cite{L-Mariani-Seo}
for the flow structure of diffusion processes.

\subsection{\label{sec13}Flow structure for reversible case}

Since the flow structure is clearer when the Markov process $X(\cdot)$
is reversible, we start with this case. The general case will be treated
in the next subsection.

Let us assume throughout this subsection that $X(\cdot)$ is reversible,
i.e., (\ref{1detbal}) holds.

For $x,\,y\in\mathcal{H}$, we write $x\sim y$ if $r(x,\,y)>0$.
Since $r(x,\,y)>0$ if and only if $r(y,\,x)>0$, we observe that
$x\sim y$ if and only if $y\sim x$. Then, we define the set of directed
edges by
\begin{equation}
\mathfrak{E}=\{(x,\,y)\in\mathcal{H}\times\mathcal{H}:x\sim y\}\;.\label{edge}
\end{equation}
Note that $(x,\,y)\in\mathfrak{E}$ if and only if $(y,\,x)\in\mathfrak{E}$
by the previous remark.

A function $\phi:\mathfrak{E}\rightarrow\mathbb{R}$ is called a flow
if it is anti-symmetric, in the sense that
\[
\phi(x,\,y)=-\phi(y,\,x)\text{ for all }x,\,y\in\mathcal{H}\;.
\]
Here, $\phi(x,\,y)$ is indeed a shorthand of $\phi((x,\,y))$. This
is called flow, since the quantity $\phi(x,\,y)$ represents the flux
of the flow from site $x$ to $y$ (and hence should be $-\phi(y,\,x)$).

The \textit{divergence} of the flow $\phi$ at site $x$ is defined
by
\[
(\textup{div}\,\phi)(x)=\sum_{y:x\sim y}\phi(x,\,y)\;,
\]
and represents the amount of the net flow coming from $x$. For $\mathcal{A}\subset\mathcal{H}$,
define
\[
(\textup{div}\,\phi)(\mathcal{A})=\sum_{x\in\mathcal{A}}(\textup{div}\,\phi)(x)\;.
\]
A flow $\phi$ is called \textit{divergence-free }at $x\in\mathcal{H}$
if $(\textup{div\,}\phi)(x)=0$, and is called divergence-free on
$\mathcal{A}\subset\mathcal{H}$ if $(\textup{div\,}\phi)(x)=0$ for
all $x\in\mathcal{A}$.

Now, we define an $L^{2}$-structure on the space of flows. Define
the \textit{conductance} between the sites as
\begin{equation}
c(x,\,y)=\mu(x)r(x,\,y)\;\;\;\;;\;x,\,y\in\mathcal{H}\;,\label{cxy}
\end{equation}
so that $c(x,\,y)=c(y,\,x)$ by (\ref{1detbal}). Denote by $\mathfrak{F}$
the space of flows. For $\phi\in\mathfrak{F}$ and $\psi\in\mathfrak{F}$,
define an inner product
\begin{equation}
\left\langle \phi,\,\psi\right\rangle _{\mathfrak{F}}=\frac{1}{2}\sum_{(x,\,y)\in\mathfrak{E}}\frac{\phi(x,\,y)\psi(x,\,y)}{c(x,\,y)}\;.\label{FN}
\end{equation}
The \textit{flow norm} of a flow $\phi$ is naturally defined by $\Vert\phi\Vert_{\mathfrak{F}}=\left\langle \phi,\,\phi\right\rangle _{\mathfrak{F}}^{1/2}$.
\begin{example*}
For $f:\mathcal{H}\rightarrow\mathbb{R}$, we define a flow $\Psi_{f}$
as
\begin{equation}
\Psi_{f}(x,\,y)=c(x,\,y)[f(y)-f(x)]\;\;\;\;;\;(x,\,y)\in\mathfrak{E}\;.\label{psif}
\end{equation}
The anti-symmetry, i.e., $\Psi_{f}(x,\,y)=-\Psi_{f}(y,\,x)$, is a
consequence of (\ref{1detbal}). A crucial feature of this flow is
the fact that
\begin{equation}
\Vert\Psi_{f}\Vert_{\mathfrak{F}}^{2}=\mathscr{D}(f)\;,\label{psif1}
\end{equation}
which follows from (\ref{DF}), (\ref{cxy}), and (\ref{FN}). Thus,
for any two disjoint and non-empty subsets $\mathcal{A}$ and $\mathcal{B}$
of $\mathcal{H}$, we have
\begin{equation}
\Vert\Psi_{h_{\mathcal{A},\,\mathcal{B}}}\Vert_{\mathfrak{F}}^{2}=\textup{cap}(\mathcal{A},\,\mathcal{B})\;.\label{psif2}
\end{equation}
This fact will be critically used later to derive the Thomson principle.
\end{example*}
Now, we can observe the following elementary properties.
\begin{prop}
\label{pflow-rev}With the notations as above, the followings hold.
\begin{enumerate}
\item For all $f:\mathcal{H}\rightarrow\mathbb{R}$ and $x\in\mathcal{H}$,
\[
(\textup{div}\,\Psi_{f})(x)=\mu(x)\,(\mathscr{L}f)(x)\;.
\]
In particular, for two disjoint non-empty subsets $\mathcal{A},\,\mathcal{B}$
of $\mathcal{H}$, the flow $\Psi_{h_{\mathcal{A},\mathcal{\,B}}}$
is divergence-free on $(\mathcal{A}\cup\mathcal{B})^{c}$.
\item For all $f:\mathcal{H}\rightarrow\mathbb{R}$ and $\phi\in\mathfrak{F}$,
\[
\left\langle \Psi_{f},\,\phi\right\rangle _{\mathfrak{F}}=-\sum_{x\in\mathcal{H}}f(x)\,(\mbox{\textup{div}}\,\phi)(x)\;.
\]
\end{enumerate}
\end{prop}

\begin{proof}
The proof follows from elementary computations. For the first assertion
of (1),
\begin{align*}
(\textup{div}\,\Psi_{f})(x) & =\sum_{y:x\sim y}\Psi_{f}(x,\,y)=\sum_{y\in\mathcal{H}}c(x,\,y)[f(y)-f(x)]\\
 & =\mu(x)\sum_{y\in\mathcal{H}}r(x,\,y)[f(y)-f(x)]=\mu(x)\,(\mathscr{L}f)(x)\;,
\end{align*}
where the second equality holds since for $y$ such that $x\not\sim y$,
we have $c(x,\,y)=0$. The second assertion of (1) follows directly
from (\ref{e_eq_prop}).

For (2), by the definition of $\Psi_{f}$,
\begin{align*}
\left\langle \Psi_{f},\,\phi\right\rangle _{\mathfrak{F}} & =\frac{1}{2}\sum_{(x,\,y)\in\mathfrak{E}}\phi(x,\,y)[f(y)-f(x)]=-\sum_{x\in\mathcal{H}}\sum_{y:y\sim x}f(x)\phi(x,\,y)\\
 & =-\sum_{x\in\mathcal{H}}f(x)\,(\mbox{\textup{div}}\,\phi)(x)\;.
\end{align*}
\end{proof}

\subsection{\label{sec14}Flow structure for non-reversible case}

Now, we turn to the general case that is developed in \cite{G-L}.
We say that $x\sim y$ if $r(x,\,y)+r(y,\,x)>0$. Similarly as before,
$x\sim y$ if and only if $y\sim x$. With this modified equivalence
relationship, we define $\mathfrak{E}$ as in (\ref{edge}), and then
the flow is defined as anti-symmetric functions on $\mathfrak{E}$.
The divergence is also defined in an identical manner.

The difference now appears at the inner product structure. Recall
(\ref{cxy}) and define
\[
c^{s}(x,\,y)=\frac{1}{2}[c(x,\,y)+c(y,\,x)]=\frac{1}{2}[\mu(x)r(x,\,y)+\mu(y)r(y,\,x)]\;,
\]
so that $c^{s}(x,\,y)=c^{s}(y,\,x)$. Then, the inner product is defined
by
\begin{equation}
\left\langle \phi,\,\psi\right\rangle _{\mathfrak{F}}=\frac{1}{2}\sum_{(x,\,y)\in\mathfrak{E}}\frac{\phi(x,\,y)\psi(x,\,y)}{c^{s}(x,\,y)}\;.\label{FN-1}
\end{equation}
Note that this definition is in accordance with (\ref{FN}) in the
reversible case. Then, the flow norm is again defined as $\Vert\phi\Vert_{\mathfrak{F}}=\left\langle \phi,\,\phi\right\rangle _{\mathfrak{F}}^{1/2}$.
\begin{example*}
For $f:\mathcal{H}\rightarrow\mathbb{R}$, define three flows as
\begin{align}
\Phi_{f}(x,\,y) & =f(y)c(y,\,x)-f(x)c(x,\,y)\;,\nonumber \\
\Phi_{f}^{*}(x,\,y) & =f(y)c(x,\,y)-f(x)c(y,\,x)\;,\label{flow}\\
\Psi_{f}(x,\,y) & =c^{s}(x,\,y)\left[f(y)-f(x)\right]=(1/2)(\Phi_{f}+\Phi_{f}^{*})(x,\,y)\;.\nonumber
\end{align}
Note that the definition of $\Psi_{f}$ is in accordance with (\ref{psif}),
and moreover we have $\Phi_{f}=\Phi_{f}^{*}=\Psi_{f}$ in the reversible
case. We remark that the relations (\ref{psif1}) and (\ref{psif2})
are still in force in this case. However, unlike the reversible case,
the expression (\ref{psif2}) for the capacity is not sufficient to
derive the Dirichlet and Thomson principles, and hence the flows $\Phi_{f}$
and $\Phi_{f}^{*}$ have to be crucially used.
\end{example*}
We conclude this subsection with the following proposition, which
summarizes several elementary properties that will be useful later.
\begin{prop}
\label{pflow_nonrev}With the notations as above, the followings hold.
\begin{enumerate}
\item For all $f:\mathcal{H}\rightarrow\mathbb{R}$ and $x\in\mathcal{H}$,
\[
(\textup{div}\,\Phi_{f})(x)=\mu(x)\,(\mathscr{L}^{\dagger}f)(x)\;\;\mbox{and\;\;}(\textup{div}\,\Phi_{f}^{*})(x)=\mu(x)\,(\mathscr{L}f)(x)\;.
\]
In particular, for two disjoint non-empty subsets $\mathcal{A},\,\mathcal{B}$
of $\mathcal{H}$, the flows $\Phi_{h_{\mathcal{A},\mathcal{\,B}}^{\dagger}}$
and $\Phi_{h_{\mathcal{A},\,\mathcal{B}}}^{*}$ are divergence-free
on $(\mathcal{A}\cup\mathcal{B})^{c}$.
\item For all $f:\mathcal{H}\rightarrow\mathbb{R}$ and $\phi\in\mathfrak{F}$,
\[
\left\langle \Psi_{f},\,\phi\right\rangle _{\mathfrak{F}}=-\sum_{x\in\mathcal{H}}f(x)\,(\mbox{\textup{div}\,}\phi)(x)\;.
\]
\item For all $f,\,g:\mathcal{H}\rightarrow\mathbb{R}$,
\[
\left\langle \Psi_{f},\,\Phi_{g}\right\rangle _{\mathfrak{F}}=\left\langle -\mathscr{L}f,\,g\right\rangle _{\mu}\;\;\text{and}\;\;\left\langle \Psi_{f},\,\Phi_{g}^{*}\right\rangle _{\mathfrak{F}}=\left\langle -\mathscr{L}^{\dagger}f,\,g\right\rangle _{\mu}\;.
\]
\end{enumerate}
\end{prop}

\begin{proof}
Proofs of (1) and (2) are similar to those of Proposition \ref{pflow-rev}
and hence are left to the readers. For (3), we first consider $\left\langle \Psi_{f},\,\Phi_{g}\right\rangle _{\mathfrak{F}}$.
By part (2), we can write
\[
\left\langle \Psi_{f},\,\Phi_{g}\right\rangle _{\mathfrak{F}}=-\sum_{x\in\mathcal{H}}f(x)\,(\text{div}\,\Phi_{g})(x)\;.
\]
Applying part (1), we get
\[
\left\langle \Psi_{f},\,\Phi_{g}\right\rangle _{\mathfrak{F}}=-\sum_{x\in\mathcal{H}}f(x)\,[\mu(x)\,(\mathscr{L}^{\dagger}g)(x)]=\left\langle f,\,-\mathscr{L}^{\dagger}g\right\rangle _{\mu}\;.
\]
Now the proof is completed by recalling (\ref{adjrel}). The proof
for $\left\langle \Psi_{f},\,\Phi_{g}^{*}\right\rangle _{\mathfrak{F}}$
is identical and will be omitted.
\end{proof}

\subsection{Application of potential theory: an example}

Before proceeding further regarding variational expression of the
capacity, we explain an application of the potential theory in the
estimate of expected hitting time or related quantities (see discussions
after Proposition \ref{p_mht}).

We fix two non-empty and disjoint subsets $\mathcal{A}$ and $\mathcal{B}$
of $\mathcal{H}$ throughout this subsection.\textbf{ }We define the
so-called equilibrium measure between $\mathcal{A}$ and $\mathcal{B}$
on $\mathcal{A}$ with respect to the process $X(\cdot)$ as
\[
\nu_{\mathcal{A},\,\mathcal{B}}(x)=\frac{M(x)\,\mathbb{P}_{x}[\tau_{\mathcal{B}}<\tau_{\mathcal{A}}^{+}]}{\textup{cap}(\mathcal{A},\,\mathcal{B})}\;\;\;\;;\;x\in\mathcal{A}\;.
\]
By Lemma \ref{lem_cap_alt}, $\nu_{\mathcal{A},\,\mathcal{B}}(\cdot)$
is a probability measure on $\mathcal{A}$. Similarly, we can define
the equilibrium measure $\nu_{\mathcal{A},\,\mathcal{B}}^{\dagger}(\cdot)$
with respect to the adjoint process $X^{\dagger}(\cdot)$:
\begin{equation}
\nu_{\mathcal{A},\,\mathcal{B}}^{\dagger}(x)=\frac{M(x)\,\mathbb{P}_{x}^{\dagger}[\tau_{\mathcal{B}}<\tau_{\mathcal{A}}^{+}]}{\textup{cap}(\mathcal{A},\,\mathcal{B})}\;\;\;\;;\;x\in\mathcal{A}\;,\label{e_emdag}
\end{equation}
where $M(x)$ and $\textup{cap}(\mathcal{A},\,\mathcal{B})$ are not
changed since $M(\cdot)$ is still the invariant measure for the embedded
chain of the adjoint process and since Proposition \ref{prop_cap_sym},
respectively.
\begin{rem}
\label{rem28}
\begin{enumerate}
\item Define the boundary $\partial\mathcal{A}$ of $\mathcal{A}$ as
\[
\partial\mathcal{A}=\{x\in\mathcal{A}:r(x,\,y)>0\text{ for some }y\notin\mathcal{A}\}\;.
\]
Note that we have $\mathbb{P}_{x}[\tau_{\mathcal{B}}<\tau_{\mathcal{A}}^{+}]=0$
for $x\in\mathcal{A}\setminus\partial\mathcal{A}$. Hence, the measure
$\nu_{\mathcal{A},\,\mathcal{B}}$ (as well as $\nu_{\mathcal{A},\,\mathcal{B}}^{\dagger}$)
is concentrated on the boundary $\partial\mathcal{A}$.
\item If $\mathcal{A}=\{a\}$ is a singleton, the measure $\nu_{\mathcal{A},\,\mathcal{B}}$
is merely the Dirac measure on $\{a\}$.
\end{enumerate}
\end{rem}

For a probability measure $\pi$ on $\mathcal{H}$, denote by $\mathbb{P}_{\pi}$
the law of the process $X(\cdot)$ when $X(0)$ is distributed according
to $\pi$, and by $\mathbb{E}_{\pi}$ the associated expectation.
The following proposition is the main result of the current subsection.
\begin{prop}
\label{p_mht}For any $f:\mathcal{H}\rightarrow\mathbb{R}$, we have
that
\begin{equation}
\mathbb{E}_{\nu_{\mathcal{A},\,\mathcal{B}}^{\dagger}}\left[\int_{0}^{\tau_{\mathcal{B}}}f(X(t))dt\right]=\frac{\left\langle f,\,h_{\mathcal{A},\,\mathcal{B}}^{\dagger}\right\rangle _{\mu}}{\textup{cap}(\mathcal{A},\,\mathcal{B})}\;.\label{e291}
\end{equation}
\end{prop}

Before proving this proposition, we explain several direct applications
of this proposition. First, we take $f\equiv1$ to deduce
\begin{equation}
\mathbb{E}_{\nu_{\mathcal{A},\,\mathcal{B}}^{\dagger}}\left[\tau_{\mathcal{B}}\right]=\frac{\sum_{x\in\mathcal{H}}h_{\mathcal{A},\,\mathcal{B}}^{\dagger}(x)\mu(x)}{\textup{cap}(\mathcal{A},\,\mathcal{B})}\;.\label{exp hit}
\end{equation}
Moreover, by taking $\mathcal{A}=\{z\}$, the left-hand side becomes
the mean hitting time $\mathbb{E}_{z}\left[\tau_{\mathcal{B}}\right]$
(cf. Remark \ref{rem28}-(2)), and thus we obtain
\begin{equation}
\mathbb{E}_{z}\left[\tau_{\mathcal{B}}\right]=\frac{\sum_{x\in\mathcal{H}}h_{z,\,\mathcal{B}}^{\dagger}(x)\mu(x)}{\textup{cap}(z,\,\mathcal{B})}\;.\label{exphit2}
\end{equation}
 Note from (\ref{e_eq_prop}) that
\begin{equation}
\sum_{x\in\mathcal{H}}h_{\mathcal{A},\,\mathcal{B}}^{\dagger}(x)\mu(x)=\mu(\mathcal{A})+\sum_{x\in(\mathcal{A}\cup\mathcal{B})^{c}}h_{\mathcal{A},\,\mathcal{B}}^{\dagger}(x)\mu(x)\le\mu(\mathcal{A})+\mu((\mathcal{A}\cup\mathcal{B})^{c})=1-\mu(\mathcal{B})\;.\label{ebn}
\end{equation}
Hence, by deriving a lower bound on $\textup{cap}(\mathcal{A},\,\mathcal{B})$,
we can obtain an upper bound on the expectation $\mathbb{E}_{\nu_{\mathcal{A},\,\mathcal{B}}^{\dagger}}\left[\tau_{\mathcal{B}}\right]$
of the hitting time from (\ref{exp hit}). In the next two sections,
we will discuss how to get a lower and an upper bound on $\textup{cap}(\mathcal{A},\,\mathcal{B})$.
Of course, in the real application, we may need more refined estimates
than (\ref{ebn}) by studying the equilibrium potential.

Next, by taking $f=\mathbf{1}_{\mathcal{C}}$ for some $\mathcal{C}\subset\mathcal{H}\setminus\mathcal{B}$,
the previous proposition becomes
\[
\mathbb{E}_{\nu_{\mathcal{A},\,\mathcal{B}}^{\dagger}}\left[\int_{0}^{\tau_{\mathcal{B}}}\mathbf{1}_{\mathcal{C}}(X(t))dt\right]=\frac{\sum_{x\in\mathcal{\mathcal{C}}}h_{\mathcal{A},\,\mathcal{B}}^{\dagger}(x)\mu(x)}{\textup{cap}(\mathcal{A},\,\mathcal{B})}\;.
\]
The left-hand side now measures the amount of time the process spends
on $\mathcal{C}$ before arriving at $\mathcal{B}$. For this setting,
the numerator of the right-hand side can be trivially bounded from
above by $\mu(\mathcal{C})$ since $h_{\mathcal{A},\,\mathcal{B}}^{\dagger}\le1$.
Now, let us return to Proposition \ref{p_mht}. The following is from
the arguments given in \cite[Proof of Proposition 6.10]{B-L TM} and
\cite[Proof of Proposition A.2]{B-L TM2}.
\begin{proof}[Proof of Proposition \ref{p_mht}]
It suffices to prove the proposition when $f=\mathbf{1}_{\{z\}}$
for all $z\in\mathcal{H}$. Let us fix $z\in\mathcal{H}$. If $z\in\mathcal{B}$,
both sides of (\ref{e291}) are trivially $0$, and hence we can assume
$z\in\mathcal{H}\setminus\mathcal{B}$.

Since the embedded chain is obtained from the original Markov process
via the time changing, we have
\[
h_{\mathcal{A},\,\mathcal{B}}^{\dagger}(z)=\mathbb{P}_{z}^{\dagger}[\tau_{\mathcal{A}}<\tau_{\mathcal{B}}]=\widehat{\mathbb{P}}_{z}^{\dagger}[\widehat{\tau}_{\mathcal{A}}<\widehat{\tau}_{\mathcal{B}}]\;,
\]
where the hitting time and the return time appearing on the right-hand
side are computed with respect to the process $\widehat{X}^{\dagger}(\cdot)$.
Write
\[
L_{\mathcal{A},\,\mathcal{B}}=\sup\{n\ge0:\widehat{X}^{\dagger}(n)\in\mathcal{A}\text{ and }n<\widehat{\tau}_{\mathcal{B}}\}\;,
\]
where we use the convention that $\sup\emptyset=-\infty$. With these
notations, we can rewrite $h_{\mathcal{A},\,\mathcal{B}}^{\dagger}(z)$
as
\begin{align}
h_{\mathcal{A},\,\mathcal{B}}^{\dagger}(z) & =\widehat{\mathbb{P}}_{z}^{\dagger}[\widehat{\tau}_{\mathcal{A}}<\widehat{\tau}_{\mathcal{B}}]=\sum_{n=0}^{\infty}\widehat{\mathbb{P}}_{z}^{\dagger}[L_{\mathcal{A},\,\mathcal{B}}=n]\nonumber \\
 & =\sum_{n=0}^{\infty}\sum_{y\in\mathcal{A}}\widehat{\mathbb{P}}_{z}^{\dagger}[\widehat{X}^{\dagger}(n)=y,\,n<\widehat{\tau}_{\mathcal{B}}]\,\widehat{\mathbb{P}}_{y}^{\dagger}[\widehat{\tau}_{\mathcal{B}}<\widehat{\tau}_{\mathcal{A}}^{+}]\nonumber \\
 & =\sum_{y\in\mathcal{A}}\left[\widehat{\mathbb{P}}_{y}^{\dagger}[\widehat{\tau}_{\mathcal{B}}<\widehat{\tau}_{\mathcal{A}}^{+}]\sum_{n=0}^{\infty}\widehat{\mathbb{P}}_{z}^{\dagger}[\widehat{X}^{\dagger}(n)=y,\,n<\widehat{\tau}_{\mathcal{B}}]\right]\;,\label{e_222}
\end{align}
where the second equality follows from the Markov property.

For $u,\,v\in\mathcal{H}\setminus\mathcal{B}$ and $n\ge0$, denote
by $P(u,\,v;n,\,\mathcal{B})$ the collection of paths $(w_{0},\,w_{1},\,\cdots,\,w_{n})$
such that $w_{0}=u$, $w_{n}=v$, and $w_{i}\notin\mathcal{B}$ for
all $0\le i\le n$. Note that
\begin{equation}
(w_{0},\,w_{1},\,\cdots,\,w_{n})\in P(u,\,v;n,\,\mathcal{B})\;\;\;\text{if and only if \;\;}(w_{n},\,w{}_{n-1},\,\cdots,\,w_{0})\in P(v,\,u;n,\,\mathcal{B})\;.\label{e_rev}
\end{equation}
With this notation (noting that we assumed $z\notin\mathcal{B}$),
we can write
\begin{equation}
M(z)\widehat{\mathbb{P}}_{z}^{\dagger}[\widehat{X}^{\dagger}(n)=y,\,n<\widehat{\tau}_{\mathcal{B}}]=\sum_{(w_{0},\,w_{1},\,\cdots,\,w_{n})\in P(z,\,y;n,\,\mathcal{B})}\sum_{i=0}^{n-1}M(w_{i})p^{\dagger}(w_{i},\,w_{i+1})\;.\label{e224}
\end{equation}
Hence, we can deduce from (\ref{e_225}), (\ref{e_rev}), and (\ref{e224})
that if $y\in A$, then
\begin{align*}
M(z)\widehat{\mathbb{P}}_{z}^{\dagger}[\widehat{X}^{\dagger}(n)=y,\,n<\widehat{\tau}_{\mathcal{B}}] & =\sum_{(w_{0},\,w_{1},\,\cdots,\,w_{n})\in P(z,\,y;n,\,\mathcal{B})}\sum_{i=0}^{n-1}M(w_{i+1})p(w_{i+1},\,w_{i})\\
 & =\sum_{(w_{0}',\,w_{1}',\,\cdots,\,w_{n}')\in P(y,\,z;n,\,\mathcal{B})}\sum_{i=0}^{n-1}M(w_{i}')p(w_{i}',\,w_{i+1}')\\
 & =M(y)\widehat{\mathbb{P}}_{y}[\widehat{X}(n)=z,\,n<\widehat{\tau}_{\mathcal{B}}]\;,
\end{align*}
Inserting this into (\ref{e_222}), we get
\begin{align}
h_{\mathcal{A},\,\mathcal{B}}^{\dagger}(z) & =\sum_{y\in\mathcal{A}}\left[\frac{M(y)}{M(z)}\widehat{\mathbb{P}}_{y}^{\dagger}[\widehat{\tau}_{\mathcal{B}}<\widehat{\tau}_{\mathcal{A}}^{+}]\sum_{n=0}^{\infty}\widehat{\mathbb{P}}_{y}[\widehat{X}(n)=z,\,n<\widehat{\tau}_{\mathcal{B}}]\right]\nonumber \\
 & =\frac{\textup{cap}(\mathcal{A},\,\mathcal{B})}{M(z)}\sum_{y\in\mathcal{A}}\left[\nu_{\mathcal{A},\,\mathcal{B}}^{\dagger}(y)\sum_{n=0}^{\widehat{\tau}_{\mathcal{B}}-1}\widehat{\mathbb{P}}_{y}[\widehat{X}(n)=z]\right]\nonumber \\
 & =\frac{\textup{cap}(\mathcal{A},\,\mathcal{B})}{M(z)}\widehat{\mathbb{E}}_{\nu_{\mathcal{A},\,\mathcal{B}}^{\dagger}}\left[\sum_{n=0}^{\widehat{\tau}_{\mathcal{B}}-1}\mathbf{1}\{\widehat{X}(n)=z\}\right]\;,\label{e225}
\end{align}
where the second equality follows from the explicit formula (\ref{e_emdag}),
while the last equality follows from the Fubini theorem. Since if
the original chain $X(\cdot)$ arrives at $z$, then it spends mean
$\lambda(z)^{-1}$ exponential random time there, and hence we can
conclude that
\[
\mathbb{E}_{\nu_{\mathcal{A},\,\mathcal{B}}^{\dagger}}\left[\int_{0}^{\tau_{B}}\mathbf{1}_{\{z\}}(X(t))dt\right]=\frac{1}{\lambda(z)}\widehat{\mathbb{E}}_{\nu_{\mathcal{A},\,\mathcal{B}}^{\dagger}}\left[\sum_{n=0}^{\widehat{\tau}_{\mathcal{B}}-1}\mathbf{1}\{\widehat{X}(n)=z\}\right]\;.
\]
Inserting this to (\ref{e_225}), we get
\[
\mathbb{E}_{\nu_{\mathcal{A},\,\mathcal{B}}^{\dagger}}\left[\int_{0}^{\tau_{B}}\mathbf{1}_{\{z\}}(X(t))dt\right]=\frac{\mu(z)h_{\mathcal{A},\,\mathcal{B}}^{\dagger}(z)}{\textup{cap}(\mathcal{A},\,\mathcal{B})}\;.
\]
This completes the proof of proposition when $f=\mathbf{1}_{\{z\}}$,
and we are done.
\end{proof}
\begin{rem}
One can expect that a quantity such as $\mathbb{E}_{a}\left[\tau_{\mathcal{B}}\right]$
is closely related with the mixing of the Markov process $X(\cdot)$.
This relation has been explained in \cite[Chapters 9 and 10]{L-P-W}.
The potential-theoretic notions are closely connected with the mixing
of Markov chains.
\end{rem}

\subsection{Bound on equilibrium potential via capacities}

Let us again fix two non-empty and disjoint subsets $\mathcal{A}$
and $\mathcal{B}$ of $\mathcal{H}$. We know that $h_{\mathcal{A},\,\mathcal{B}}\equiv1$
on $\mathcal{A}$ and $h_{\mathcal{A},\,\mathcal{B}}\equiv0$ on $\mathcal{B}$,
but the value of $h_{\mathcal{A},\,\mathcal{B}}$ on $(\mathcal{A}\cup\mathcal{B})^{c}$
is described only in terms of the Laplace equation (cf. (\ref{e_eq_prop})),
and hence the exact value is almost impossible to compute in most
applications. However, in many instances, we need to bound the value
of $h_{\mathcal{A},\,\mathcal{B}}$ on $(\mathcal{A}\cup\mathcal{B})^{c}$
to carry out an estimation. For example, with such a bound, we can
carry out a much better estimate in (\ref{ebn}).

In this subsection, we present the following useful upper bound on
the value of $h_{\mathcal{A},\,\mathcal{B}}$ on $(\mathcal{A}\cup\mathcal{B})^{c}$
in terms of capacities. This bound will be frequently used in various
situations. The following proof is an excerpt from \cite[Section 3]{L-L}.
\begin{prop}
\label{p_eqp_bd}We have that
\[
h_{\mathcal{A},\,\mathcal{B}}(x)\le\frac{\textup{cap}(x,\,\mathcal{A})}{\textup{cap}(x,\,\mathcal{A\cup\mathcal{B}})}\;\;\;\text{for all }x\in(\mathcal{A}\cup\mathcal{B})^{c}.
\]
\end{prop}

\begin{proof}
Fix $x\in(\mathcal{A}\cup\mathcal{B})^{c}$. By the strong Markov
property, we can write
\begin{align*}
\mathbb{P}_{x}[\tau_{\mathcal{A}}<\tau_{\mathcal{B}}] & =\mathbb{P}_{x}[\tau_{x}^{+}<\tau_{\mathcal{\mathcal{A}\cup B}},\,\tau_{\mathcal{A}}<\tau_{\mathcal{B}}]+\mathbb{P}_{x}[\tau_{x}^{+}>\tau_{\mathcal{\mathcal{A}\cup B}},\,\tau_{\mathcal{A}}<\tau_{\mathcal{B}}]\\
 & =\mathbb{P}_{x}[\tau_{x}^{+}<\tau_{\mathcal{\mathcal{A}\cup B}}]\mathbb{P}_{x}[\tau_{\mathcal{A}}<\tau_{\mathcal{B}}]+\mathbb{P}_{x}[\tau_{\mathcal{A}}<\tau_{\mathcal{B}}<\tau_{x}^{+}]\;.
\end{align*}
Therefore, we have that
\[
\mathbb{P}_{x}[\tau_{\mathcal{A}}<\tau_{\mathcal{B}}]=\frac{\mathbb{P}_{x}[\tau_{\mathcal{A}}<\tau_{\mathcal{B}}<\tau_{x}^{+}]}{\mathbb{P}_{x}[\tau_{x}^{+}>\tau_{\mathcal{\mathcal{A}\cup B}}]}\le\frac{\mathbb{P}_{x}[\tau_{x}^{+}>\tau_{\mathcal{\mathcal{A}}}]}{\mathbb{P}_{x}[\tau_{x}^{+}>\tau_{\mathcal{\mathcal{A}\cup B}}]}\;.
\]
The proof is completed since by Lemma \ref{lem_cap_alt},
\begin{align*}
\textup{cap}(x,\,\mathcal{A}) & =M(x)\mathbb{P}_{x}[\tau_{x}^{+}>\tau_{\mathcal{\mathcal{A}}}]\;\;\;\text{and}\\
\textup{cap}(x,\,\mathcal{A}\cup\mathcal{B}) & =M(x)\mathbb{P}_{x}[\tau_{x}^{+}>\tau_{\mathcal{\mathcal{A}\cup B}}]\;.
\end{align*}
\end{proof}
We note that, in view of Proposition \ref{prop_cap_sym}, the same
result holds for $h_{\mathcal{A},\,\mathcal{B}}^{\dagger}(x)$ in
place of $h_{\mathcal{A},\,\mathcal{B}}(x)$. In addition, the bound
obtained in the previous proposition is particularly useful since
there are numerous robust tools to estimate capacities. We discuss
such robust tools in the following sections.

\section{Dirichlet and Thomson Principles}

In the application of the potential theory, it is important to (more
or less precisely) estimate the capacity. Classic tools for this purpose
are the Dirichlet and Thomson principles that we introduce in this
section.

Let us fix two disjoint and non-empty subsets $\mathcal{A}$ and $\mathcal{B}$
throughout the section. Then, we explain strategies to estimate the
capacity $\textup{cap}(\mathcal{A},\,\mathcal{B})$.

\subsection{Spaces of functions and flows}

To explain the variational principles for capacities, we need to define
classes of functions and flows as follows:
\begin{itemize}
\item For real numbers $a$ and $b$, denote by $\mathfrak{C}_{a,\,b}(\mathcal{A},\,\mathcal{B})$
the set of all real-valued functions $f$ on $\mathcal{H}$ satisfying
$f|_{\mathcal{A}}\equiv a$ and $f|_{\mathcal{B}}\equiv b$, i.e.,
\[
\mathfrak{C}_{a,\,b}(\mathcal{A},\,\mathcal{B})=\left\{ f:\mathcal{H}\rightarrow\mathbb{R}:f(x)=a,\,\forall x\in\mathcal{A}\;\mbox{and}\;f(x)=b,\,\forall x\in\mathcal{B}\right\} \;.
\]
\item For $a\in\mathbb{R}$, let $\mathfrak{U}_{a}(\mathcal{A},\,\mathcal{B})$
be the set of all flows $\phi\in\mathfrak{F}$ which are divergence
free on $(\mathcal{A}\cup\mathcal{B})^{c}$, i.e.,
\[
(\textup{div}\,\phi)(x)=0\;\mbox{for all }x\in(\mathcal{A}\cup\mathcal{B})^{c}\;,
\]
and satisfy
\[
(\textup{div}\,\phi)(\mathcal{A})=-(\textup{div}\,\phi)(\mathcal{B})=a\;.
\]
In particular, a flow belonging to $\mathfrak{U}_{1}$ is called a
\textit{unit flow}.
\end{itemize}
\begin{example}
\label{ex_31}The equilibrium potential $h_{\mathcal{A},\,\mathcal{B}}$
belongs to the class $\mathfrak{C}_{1,\,0}(\mathcal{A},\,\mathcal{B})$.
\end{example}

\begin{xca}
\label{ex_32}
\begin{enumerate}
\item Suppose that the process $X(\cdot)$ is reversible. Prove that  the
flow
\begin{equation}
\psi_{\mathcal{A},\,\mathcal{B}}:=-\frac{1}{\textup{cap}(\mathcal{A},\,\mathcal{B})}\Psi_{h_{\mathcal{A},\,\mathcal{B}}}\label{psiAB}
\end{equation}
is a unit flow between $\mathcal{A}$ and $\mathcal{B}$. (Hint: Proposition
\ref{pflow-rev}-(1))
\item Suppose that the process $X(\cdot)$ is non-reversible. Prove that
the flows
\[
\phi_{\mathcal{A},\,\mathcal{B}}:=-\frac{1}{\textup{cap}(\mathcal{A},\,\mathcal{B})}\Phi_{h_{\mathcal{A},\,\mathcal{B}}^{\dagger}}\;\;\;\;\text{and\;\;\;\;}\phi_{\mathcal{A},\,\mathcal{B}}^{*}:=-\frac{1}{\textup{cap}(\mathcal{A},\,\mathcal{B})}\Phi_{h_{\mathcal{A},\,\mathcal{B}}}^{*}
\]
are unit flows between $\mathcal{A}$ and $\mathcal{B}$. (Hint: Proposition
\ref{pflow_nonrev}-(1))
\end{enumerate}
\end{xca}

\subsection{Dirichlet and Thomson principles: reversible case}

We begin with the Dirichlet and Thomson principles for reversible
Markov processes. Hence, in this subsection, we temporarily assume
that the process $X(\cdot)$ is reversible.

The Dirichlet principle provides a minimization problem for the capacity.
\begin{thm}[Dirichlet principle for reversible Markov processes]
\label{t_DP_rev}We have that
\[
\textup{cap}(\mathcal{A},\,\mathcal{B})=\inf_{f\in\mathfrak{C}_{1,\,0}(\mathcal{A},\,\mathcal{B})}\mathscr{D}(f)\;,
\]
and the unique minimizer is given by $f=h_{\mathcal{A},\,\mathcal{B}}$.\footnote{Note that $h_{\mathcal{A},\,\mathcal{B}}\in\mathfrak{C}_{1,\,0}(\mathcal{A},\,\mathcal{B})$
as we observed in Example \ref{ex_31}.}
\end{thm}

\begin{proof}
Let $f\in\mathfrak{C}_{1,\,0}(\mathcal{A},\,\mathcal{B})$. Then,
write $g=f-h_{\mathcal{A},\,\mathcal{B}}$ so that $g\in\mathfrak{C}_{0,\,0}(\mathcal{A},\,\mathcal{B})$.
Then,
\begin{align*}
\mathscr{D}(f) & =\left\langle h_{\mathcal{A},\,\mathcal{B}}+g,\,-\mathscr{L}(h_{\mathcal{A},\,\mathcal{B}}+g)\right\rangle _{\mu}\\
 & =\mathscr{D}(h_{\mathcal{A},\,\mathcal{B}})+\mathscr{D}(g)-2\left\langle g,\,\mathscr{L}h_{\mathcal{A},\,\mathcal{B}}\right\rangle _{\mu}\;,
\end{align*}
where at the second equality we used the reversibility which implies
the self-adjointness of $\mathscr{L}$. Since $g\equiv0$ on $\mathcal{A}\cup\mathcal{B}$,
and since $\mathscr{L}h_{\mathcal{A},\,\mathcal{B}}\equiv0$ on $(\mathcal{A}\cup\mathcal{B})^{c}$,
we get $\left\langle g,\,\mathscr{L}h_{\mathcal{A},\,\mathcal{B}}\right\rangle _{\mu}=0$.
Therefore,
\begin{equation}
\mathscr{D}(f)=\mathscr{D}(h_{\mathcal{A},\,\mathcal{B}})+\mathscr{D}(g)\ge\mathscr{D}(h_{\mathcal{A},\,\mathcal{B}})\;,\label{e_dr}
\end{equation}
and the equality holds only when $\mathscr{D}(g)=0$, i.e., when $g$
is a constant function. Since $g\in\mathfrak{C}_{0,\,0}(\mathcal{A},\,\mathcal{B})$,
$g$ must be the zero function to obtain the equality in (\ref{e_dr}).
This completes the proof.
\end{proof}
On the other hand, the Thomson principle provides a maximization problem
for the capacity.
\begin{thm}[Thomson principle for reversible Markov processes]
\label{t_TP_rev}We have that
\[
\textup{cap}(\mathcal{A},\,\mathcal{B})=\sup_{\phi\in\mathfrak{U}_{1}(\mathcal{A},\,\mathcal{B})}\frac{1}{\Vert\phi\Vert_{\mathfrak{F}}^{2}}\;,
\]
and the unique maximizer is given by $\phi=\psi_{\mathcal{A},\,\mathcal{B}}$
(cf. (\ref{psiAB})).
\end{thm}

\begin{proof}
Let $\phi\in\mathfrak{U}_{1}(\mathcal{A},\,\mathcal{B})$. By Proposition
\ref{pflow-rev}, we have
\[
\left\langle \Psi_{h_{\mathcal{A},\,\mathcal{B}}},\,\phi\right\rangle _{\mathfrak{F}}=-\sum_{x\in\mathcal{H}}h_{\mathcal{A},\,\mathcal{B}}(x)\,(\mbox{\textup{div}}\,\phi)(x)=-\sum_{x\in\mathcal{A}}h_{\mathcal{A},\,\mathcal{B}}(x)\,(\mbox{\textup{div}}\,\phi)(x)\;,
\]
where the second equality holds since $h_{\mathcal{A},\,\mathcal{B}}\equiv0$
on $\mathcal{B}$ and $\mbox{\textup{div}}\,\phi\equiv0$ on $(\mathcal{A}\cup\mathcal{B})^{c}$.
Since $h_{\mathcal{A},\,\mathcal{B}}\equiv1$ on $\mathcal{A}$ and
since $(\mbox{\textup{div}}\,\phi)(\mathcal{A})=1$, we can conclude
that
\[
\left\langle \Psi_{h_{\mathcal{A},\,\mathcal{B}}},\,\phi\right\rangle _{\mathfrak{F}}=-1\;.
\]
By the Cauchy--Schwarz inequality and (\ref{psif2}),
\[
1=\left\langle \Psi_{h_{\mathcal{A},\,\mathcal{B}}},\,\phi\right\rangle _{\mathfrak{F}}^{2}\le\Vert\Psi_{h_{\mathcal{A},\,\mathcal{B}}}\Vert_{\mathfrak{F}}^{2}\,\Vert\phi\Vert_{\mathfrak{F}}^{2}=\textup{cap}(\mathcal{A},\,\mathcal{B})\,\Vert\phi\Vert_{\mathfrak{F}}^{2}\;.
\]
This proves $\textup{cap}(\mathcal{A},\,\mathcal{B})\ge\frac{1}{\Vert\phi\Vert_{\mathfrak{F}}^{2}}$.
Since the equality of the previous Cauchy--Schwarz inequality holds
only when $\phi=c\Psi_{h_{\mathcal{A},\,\mathcal{B}}}$ for some $c\in\mathbb{R}$,
we must have $\phi=\psi_{\mathcal{A},\,\mathcal{B}}$ since $\phi$
is a unit flow.
\end{proof}
\begin{rem}
\label{rem:DTP}At this point, it is now clear how to use the Dirichlet
and Thomson principles to estimate the capacity. If we take any test
function $f\in\mathfrak{C}_{1,\,0}(\mathcal{A},\,\mathcal{B})$ and
any test flow $\phi\in\mathfrak{U}_{1}(\mathcal{A},\,\mathcal{B})$,
we can deduce from Theorems \ref{t_DP_rev} and \ref{t_TP_rev} that
\[
\frac{1}{\Vert\phi\Vert_{\mathfrak{F}}^{2}}\le\textup{cap}(\mathcal{A},\,\mathcal{B})\le\mathscr{D}(f)\;.
\]
If one wants these lower and upper bounds to be sharp, it is necessary
to take $f$ and $\phi$ as objects close to the genuine optimizers,
namely, as $f\approx h_{\mathcal{A},\,\mathcal{B}}$ and $\phi\approx\psi_{\mathcal{A},\,\mathcal{B}}$.
For a concrete example of such a construction, we refer to \cite{L-M-T}.

We note that there is no special technical difficulty in finding such
a test function. On the other hand, constructing an appropriate test
flow is fundamentally more difficult, since the object that we constructed
as a test flow must satisfy the divergence-free condition on $(\mathcal{A}\cup\mathcal{B})^{c}$,
and there is no trivial way of defining such an object. This issue
will be discussed in more detail in the next section.
\end{rem}

\begin{rem}
\label{rem_alt}In the reversible case, there is an alternative way,
based on a Cauchy--Schwarz-type argument, of obtaining a lower bound
for the capacity without relying on the Thomson principle. More precisely,
if we are able to prove that $\mathscr{D}(f)$ is bounded below by
a constant $c$ for all $f\in\mathfrak{C}_{1,\,0}(\mathcal{A},\,\mathcal{B})$
via the Cauchy--Schwarz inequalities, then by the Dirichlet principle
we have the lower bound $\mathscr{D}(f)\ge c$. This bound can be
sharp if we apply the inequalities in a careful manner. We refer to
\cite{B-L ZRP,B-D-G,B-E-G-K} for examples of this method. This method
is difficult to use when the underlying energy landscape is complicated.
\end{rem}

\subsection{Dirichlet and Thomson principles: non-reversible case}

The Dirichlet and Thomson principles were known only for the reversible
case, but recently the corresponding principles for the non-reversible
case have been revealed. The following theorem is a summary of these
results. We no longer assume that the process $X(\cdot)$ is reversible.
\begin{thm}
\label{t_dt}The following variational expressions for the capacity
hold.
\begin{enumerate}
\item It holds that
\begin{align}
\textup{cap}(\mathcal{A},\,\mathcal{B}) & =\inf_{f\in\mathfrak{C}_{1,\,0}(\mathcal{A},\,\mathcal{B}),\;\phi\in\mathfrak{U}_{0}(\mathcal{A},\,\mathcal{B})}\,||\Phi_{f}-\phi||^{2}\;,\label{DP}
\end{align}
and the unique minimizer is given by
\begin{equation}
(f,\,\phi)=\left(\frac{h_{\mathcal{A},\,\mathcal{B}}+h_{\mathcal{A},\,\mathcal{B}}^{\dagger}}{2},\,-\frac{\Phi_{h_{\mathcal{A},\mathcal{\,B}}^{\dagger}}-\Phi_{h_{\mathcal{A},\mathcal{\,B}}}^{*}}{2}\right)\;,\label{DPo}
\end{equation}
\item It holds that
\begin{equation}
\textup{cap}(\mathcal{A},\,\mathcal{B})=\sup_{g\in\mathfrak{C}_{0,\,0}(\mathcal{A},\,\mathcal{B}),\;\psi\in\mathfrak{U}_{1}(\mathcal{A},\,\mathcal{B})}\,\frac{1}{||\Phi_{g}-\psi||^{2}}\;.\label{TP}
\end{equation}
and the unique maximizer is given by
\begin{equation}
(g,\,\psi)=\left(\frac{h_{\mathcal{A},\mathcal{\,B}}^{\dagger}-h_{\mathcal{A},\mathcal{\,B}}}{2\textup{cap}(\mathcal{A},\,\mathcal{B})},\,-\frac{\Phi_{h_{\mathcal{A},\mathcal{\,B}}^{\dagger}}+\Phi_{h_{\mathcal{A},\,\mathcal{B}}}^{*}}{2\textup{cap}(\mathcal{A},\,\mathcal{B})}\right)\;.\label{TPo}
\end{equation}
\end{enumerate}
\end{thm}

In the previous theorem, the Dirichlet principle (\ref{DP}) and the
Thomson principle (\ref{TP}) were established in \cite{G-L} and
\cite{S}, respectively. Note also that
\[
-\frac{\Phi_{h_{\mathcal{A},\mathcal{\,B}}^{\dagger}}-\Phi_{h_{\mathcal{A},\mathcal{\,B}}}^{*}}{2}\in\mathfrak{U}_{0}(\mathcal{A},\,\mathcal{B})\;\;\;\;\text{and\;\;\;\;}-\frac{\Phi_{h_{\mathcal{A},\mathcal{\,B}}^{\dagger}}+\Phi_{h_{\mathcal{A},\,\mathcal{B}}}^{*}}{2\textup{cap}(\mathcal{A},\,\mathcal{B})}\in\mathfrak{U}_{1}(\mathcal{A},\,\mathcal{B})
\]
follows from Example \ref{ex_32}. We now turn to the proof.
\begin{proof}
Let $f\in\mathfrak{C}_{a,\,0}(\mathcal{A},\,\mathcal{B})$. By Proposition
\ref{pflow_nonrev}-(3), we have that
\[
\left\langle \Psi_{h_{\mathcal{A},\,\mathcal{B}}},\,\Phi_{f}\right\rangle _{\mathfrak{F}}=\left\langle -\mathscr{L}h_{\mathcal{A},\,\mathcal{B}},\,f\right\rangle _{\mu}=\sum_{x\in\mathcal{H}}f(x)(-\mathscr{L}h_{\mathcal{A},\,\mathcal{B}})(x)\mu(x)\;.
\]
Since $-\mathscr{L}h_{\mathcal{A},\,\mathcal{B}}\equiv0$ on $(\mathcal{A}\cup\mathcal{B})^{c}$
and $f=ah_{\mathcal{A},\,\mathcal{B}}$ on $\mathcal{A}\cup\mathcal{B}$,
we can conclude from the previous identity that
\begin{align}
\left\langle \Psi_{h_{\mathcal{A},\,\mathcal{B}}},\,\Phi_{f}\right\rangle _{\mathfrak{F}} & =a\sum_{x\in\mathcal{\mathcal{A}}}h_{\mathcal{A},\,\mathcal{B}}(x)(-\mathscr{L}h_{\mathcal{A},\,\mathcal{B}})(x)\mu(x)\nonumber \\
 & =a\mathscr{D}(h_{\mathcal{A},\,\mathcal{B}})=a\,\textup{cap}(\mathcal{A},\,\mathcal{B})\;.\label{edt1}
\end{align}
Let $\phi\in\mathfrak{U}_{a}(\mathcal{A},\,\mathcal{B})$. Then, by
Proposition \ref{pflow_nonrev}-(2),
\[
\left\langle \Psi_{h_{\mathcal{A},\,\mathcal{B}}},\,\phi\right\rangle _{\mathfrak{F}}=-\sum_{x\in\mathcal{H}}h_{\mathcal{A},\,\mathcal{B}}(x)\,(\mbox{\textup{div}\,}\phi)(x)\;.
\]
Since $\textup{div}\,\phi\equiv0$ on $(\mathcal{A}\cup\mathcal{B})^{c}$
and $h_{\mathcal{A},\,\mathcal{B}}=\mathbf{1}_{\mathcal{A}}$ on $\mathcal{A}\cup\mathcal{B}$,
\begin{equation}
\left\langle \Psi_{h_{\mathcal{A},\,\mathcal{B}}},\,\phi\right\rangle _{\mathfrak{F}}=-\sum_{x\in\mathcal{A}}(\mbox{\textup{div}\,}\phi)(x)=-(\mbox{\textup{div}\,}\phi)(\mathcal{A})=-a\;,\label{edt2}
\end{equation}
where the last equality follows from $\phi\in\mathfrak{U}_{a}(\mathcal{A},\,\mathcal{B})$.

Now, we first look at (1). If $f\in\mathfrak{C}_{1,\,0}(\mathcal{A},\,\mathcal{B})$
and $\phi\in\mathfrak{U}_{0}(\mathcal{A},\,\mathcal{B})$, then by
(\ref{edt1}) and (\ref{edt2}),
\[
\left\langle \Psi_{h_{\mathcal{A},\,\mathcal{B}}},\,\Phi_{f}-\phi\right\rangle _{\mathfrak{F}}=\textup{cap}(\mathcal{A},\,\mathcal{B})\;.
\]
By the Cauchy--Schwarz inequality,
\[
\textup{cap}(\mathcal{A},\,\mathcal{B})^{2}=\left\langle \Psi_{h_{\mathcal{A},\,\mathcal{B}}},\,\Phi_{f}-\phi\right\rangle _{\mathfrak{F}}^{2}\le\Vert\Psi_{h_{\mathcal{A},\,\mathcal{B}}}\Vert_{\mathfrak{F}}^{2}\,\Vert\Phi_{f}-\phi\Vert_{\mathfrak{F}}^{2}=\textup{cap}(\mathcal{A},\,\mathcal{B})\,\Vert\Phi_{f}-\phi\Vert_{\mathfrak{F}}^{2}\;.
\]
Therefore, we get $\Vert\Phi_{f}-\phi\Vert_{\mathfrak{F}}^{2}\ge\textup{cap}(\mathcal{A},\,\mathcal{B})$.
The equality holds only when $\Phi_{f}-\phi=c\Psi_{h_{\mathcal{A},\,\mathcal{B}}}$
for some $c\in\mathbb{R}$. By carefully analyzing this restriction,
we can conclude that equality holds only for (\ref{DPo}).

Next, we consider (2). If $g\in\mathfrak{C}_{0,\,0}(\mathcal{A},\,\mathcal{B})$
and $\psi\in\mathfrak{U}_{1}(\mathcal{A},\,\mathcal{B})$, then again
by (\ref{edt1}) and (\ref{edt2}), we get
\[
\left\langle \Psi_{h_{\mathcal{A},\,\mathcal{B}}},\,\Phi_{g}-\psi\right\rangle _{\mathfrak{F}}=-1\;.
\]
Thus, by the Cauchy--Schwarz inequality,
\[
1=\left\langle \Psi_{h_{\mathcal{A},\,\mathcal{B}}},\,\Phi_{g}-\psi\right\rangle _{\mathfrak{F}}^{2}\le\Vert\Psi_{h_{\mathcal{A},\,\mathcal{B}}}\Vert_{\mathfrak{F}}^{2}\,\Vert\Phi_{g}-\psi\Vert_{\mathfrak{F}}^{2}=\textup{cap}(\mathcal{A},\,\mathcal{B})\,\Vert\Phi_{g}-\psi\Vert_{\mathfrak{F}}^{2}\;.
\]
Hence, we get $\textup{cap}(\mathcal{A},\,\mathcal{B})\ge\Vert\Phi_{g}-\psi\Vert_{\mathfrak{F}}^{-2}$.
One can also readily check that the equality holds only for the selection
(\ref{TPo}).
\end{proof}
Now, Theorem \ref{t_dt} can be used to estimate the capacity in the
non-reversible case in the same manner as Remark \ref{rem:DTP}. We
note that now divergence-free test flows are needed for both upper
and lower bounds, and thus we must address this technical issue directly
to use these principles. Note that, in the non-reversible case, an
argument such as Remark \ref{rem_alt} does not exist.
\begin{rem}
In fact, Proposition \ref{prop:inc} for the reversible case is a
consequence of the Dirichlet principle (Theorem \ref{t_DP_rev}),
since we have $\mathfrak{C}_{1,\,0}(\mathcal{A}',\,\mathcal{B}')\subset\mathfrak{C}_{1,\,0}(\mathcal{A},\,\mathcal{B})$
if $\mathcal{A}\subset\mathcal{A}'$ and $\mathcal{B}\subset\mathcal{B}'$.
On the other hand, for the non-reversible case, we do not have such
a simple argument since it holds that $\mathfrak{U}_{0}(\mathcal{A},\,\mathcal{B})\subset\mathfrak{U}_{0}(\mathcal{A}',\,\mathcal{B}')$
instead of $\mathfrak{U}_{0}(\mathcal{A}',\,\mathcal{B}')\subset\mathfrak{U}_{0}(\mathcal{A},\,\mathcal{B})$.
\end{rem}

\subsection{Comparison result for capacity}

One can observe from the Dirichlet and Thomson principles that the
capacity estimates of non-reversible processes are far more complicated
than those of reversible processes. Hence, if one only needs a rough
capacity estimate of a non-reversible process, it would be very handy
if a comparison result between the capacity of a reversible process
and that of a non-reversible one exists. In this section, we provide
such a result based on the Dirichlet principle. This comparison result
will be used in Part 3.

Define a symmetrized rate as
\[
r^{s}(x,\,y)=\frac{1}{2\mu(x)}[\mu(x)r(x,\,y)+\mu(y)r(y,\,x)]\;\;\;\;;\;x,\,y\in\mathcal{H}\;,
\]
and let $(X^{s}(t))_{t\ge0}$ be a continuous-time Markov process
on $\mathcal{H}$ with rate $r^{s}(\cdot,\,\cdot)$. One can observe
now that the following detailed balance condition holds:
\[
\mu(x)r^{s}(x,\,y)=\mu(y)r^{s}(y,\,x)\;.
\]
Hence, $\mu(\cdot)$ is the invariant measure for the process $X^{s}(\cdot)$,
and furthermore $X^{s}(\cdot)$ is a reversible process.

We write $h_{\mathcal{A},\,\mathcal{B}}^{s}$ and $\textup{cap}^{s}(\mathcal{A},\,\mathcal{B})$
the equilibrium potential and the capacity, respectively, with respect
to the process $X^{s}(\cdot)$, for two disjoint and non-empty subsets
$\mathcal{A}$ and $\mathcal{B}$ of $\mathcal{H}$. One can easily
check that the Dirichlet form of this symmetrized process is still
$\mathscr{D}(\cdot)$ (cf. Remark \ref{remgf}).

Since $X^{s}(\cdot)$ is reversible, it could be much simpler to estimate
$\textup{cap}^{s}(\mathcal{A},\,\mathcal{B})$ than to estimate $\textup{cap}(\mathcal{A},\,\mathcal{B})$.
The purpose of this subsection is to compare these two capacities.

Firstly, we can show that the symmetrized capacity is always smaller.
\begin{prop}
\label{p_capcomp1}For any two disjoint and non-empty subsets $\mathcal{A}$
and $\mathcal{B}$ of $\mathcal{H}$, it holds that
\[
\textup{cap}^{s}(\mathcal{A},\,\mathcal{B})\le\textup{cap}(\mathcal{A},\,\mathcal{B})\;.
\]
\end{prop}

\begin{proof}
Since $h_{\mathcal{A},\,\mathcal{B}}\in\mathfrak{C}_{1,\,0}(\mathcal{A},\,\mathcal{B})$,
by the Dirichlet principle for reversible processes (Theorem \ref{t_DP_rev}),
\[
\textup{cap}^{s}(\mathcal{A},\,\mathcal{B})=\inf_{f\in\mathfrak{C}_{1,\,0}(\mathcal{A},\,\mathcal{B})}\mathscr{D}(f)\le\mathscr{D}(h_{\mathcal{A},\,\mathcal{B}})=\textup{cap}(\mathcal{A},\,\mathcal{B})\;.
\]
\end{proof}
We next investigate the opposite bound. To this end, we have to introduce
the sector condition.
\begin{defn}
\label{sector}A Markov process $X(\cdot)$ is said to satisfy the
sector condition with constant $C_{0}>0$ if
\begin{equation}
\left\langle f,\,-\mathscr{L}g\right\rangle _{\mu}^{2}\le C_{0}\mathscr{D}(f)\mathscr{D}(g)\label{e_sect}
\end{equation}
for all $f,\,g:\mathcal{H}\rightarrow\mathbb{R}$.
\end{defn}

Heuristically, this is called the sector condition since the eigenvalues
of $\mathscr{L}$ satisfying (\ref{e_sect2}) are located on a certain
sector at the complex plane originating from $0$. In this sense,
one regards a Markov process with the sector condition as a process
which is not far from reversibility. A huge class of Markov processes
under consideration satisfies the sector condition. We shall check,
for instance, whether the non-reversible zero-range process considered
in Part \ref{pt4} satisfies the sector condition (cf. Proposition
\ref{psector}).
\begin{xca}
If $X(\cdot)$ is reversible, prove that one can write
\begin{equation}
\left\langle f,\,-\mathscr{L}g\right\rangle _{\mu}=\frac{1}{2}\sum_{x\in\mathcal{H}}\sum_{y\in\mathcal{H}}\mu(x)r(x,\,y)(g(y)-g(x))(f(y)-f(x))\;,\label{eq:fgmu}
\end{equation}
and therefore $X(\cdot)$ satisfies the sector condition with constant
$1$.
\end{xca}

\begin{rem}
\label{rmk_sector}Of course, if $X(\cdot)$ is non-reversible, the
expression (\ref{eq:fgmu}) does not hold, and therefore checking
the inequality (\ref{e_sect}) is not trivial at all. To check (\ref{e_sect}),
one usually proves inequality of the form
\begin{equation}
\left\langle f,\,-\mathscr{L}g\right\rangle _{\mu}\le C_{1}\mathscr{D}(f)+C_{2}\mathscr{D}(g)\label{e_sect2}
\end{equation}
for some constant $C_{1},\,C_{2}>0$ for all $f,\,g:\mathcal{H}\rightarrow\mathbb{R}$.
We first note that the inequality (\ref{e_sect}) is trivial if $f$
or $g$ is a constant function (cf. Exercise \ref{ex21}). Otherwise,
inserting $f:=\sqrt{C_{2}\mathscr{D}(g)}f$ and $g:=\sqrt{C_{1}\mathscr{D}(f)}g$
to (\ref{e_sect2}), we get
\[
\sqrt{C_{1}C_{2}\mathscr{D}(f)\mathscr{D}(g)}\left\langle f,\,-\mathscr{L}g\right\rangle _{\mu}\le2C_{1}C_{2}\mathscr{D}(f)\mathscr{D}(g)\;.
\]
Therefore, we can conclude that $X(\cdot)$ satisfies the sector condition
with constant $4C_{1}C_{2}$.
\end{rem}

Now, we are ready to establish the opposite bound of the one established
in Proposition \ref{p_capcomp1}.
\begin{prop}
\label{p_capcomp2}Suppose that a Markov process $X(\cdot)$ satisfies
the sector condition with constant $C_{0}>0$. Then, we have that
\[
\textup{cap}(\mathcal{A},\,\mathcal{B})\le C_{0}\,\textup{cap}^{s}(\mathcal{A},\,\mathcal{B})\;.
\]
\end{prop}

\begin{proof}
We may assume that $\textup{cap}(\mathcal{A},\,\mathcal{B})>0$, as
otherwise the inequality is trivial. We first note that
\[
\textup{cap}(\mathcal{A},\,\mathcal{B})=\left\langle h_{\mathcal{A},\,\mathcal{B}},\,-\mathscr{L}h_{\mathcal{A},\,\mathcal{B}}\right\rangle _{\mu}=\left\langle h_{\mathcal{A},\,\mathcal{B}}^{s},\,-\mathscr{L}h_{\mathcal{A},\,\mathcal{B}}\right\rangle _{\mu}\;,
\]
where the second equality holds since $\mathscr{L}h_{\mathcal{A},\,\mathcal{B}}=0$
on $(\mathcal{A}\cup\mathcal{B})^{c}$ and $h_{\mathcal{A},\,\mathcal{B}}=h_{\mathcal{A},\,\mathcal{B}}^{s}$
on $\mathcal{A}\cup\mathcal{B}$. Therefore, by the sector condition,
\[
\textup{cap}(\mathcal{A},\,\mathcal{B})^{2}\le C_{0}\mathscr{D}(h_{\mathcal{A},\,\mathcal{B}}^{s})\mathscr{D}(h_{\mathcal{A},\,\mathcal{B}})=C_{0}\textup{cap}^{s}(\mathcal{A},\,\mathcal{B})\textup{cap}(\mathcal{A},\,\mathcal{B})\;.
\]
Dividing both sides by $\textup{cap}(\mathcal{A},\,\mathcal{B})>0$
completes the proof.
\end{proof}

\section{Generalized Dirichlet--Thomson Principles}

Let us fix two disjoint and non-empty subsets $\mathcal{A}$ and $\mathcal{B}$
of $\mathcal{H}$. In the previous section, we explain a general strategy
to estimate or bound the capacity $\textup{cap}(\mathcal{A},\,\mathcal{B})$
based on the Dirichlet and Thomson principles. To apply this strategy,
one has to construct suitable test functions or test flows. As we
have mentioned earlier, the Thomson principle for the reversible case
and the Dirichlet and Thomson principles for the non-reversible case
require us to construct a test flow which must be divergence-free
on $(\mathcal{A}\cup\mathcal{B})^{c}$ which is a major technical
problem in applications of these method. In this section, we introduce
alternative variational principles that do not require us to construct
a divergence-free flow, and hence are suitable for many applications.

\subsection{Reversible case}

Let us start by considering the reversible case. Hence, we assume
in this subsection that the process $X(\cdot)$ is reversible. We
also emphasize that we do not need to develop a generalized Dirichlet
principle, since the Dirichlet principle for reversible Markov processes
is not involved with the flow structure.

The generalized Thomson principle is given as follows. We write $\mathfrak{F}_{0}$
the collection of non-zero flows, i.e.,
\[
\mathfrak{F}_{0}=\{\phi\in\mathfrak{F}:\Vert\phi\Vert_{\mathfrak{F}}>0\}\;.
\]

\begin{thm}[Generalized Thomson principle: reversible case]
\label{t_GTP}It holds that
\begin{equation}
\mathrm{cap}(\mathcal{A},\,\mathcal{B})=\sup_{\phi\in\mathfrak{F}_{0}}\frac{1}{\|\phi\|_{\mathfrak{F}}^{2}}\bigg[\,\sum_{\sigma\in\mathcal{H}}h_{\mathcal{A},\,\mathcal{B}}(x)\,(\mathrm{div}\,\phi)(x)\,\bigg]^{2}\;.\label{genThomeq}
\end{equation}
Moreover, the optimizers are given by $\phi=c\Psi_{h_{\mathcal{A},\,\mathcal{B}}}$,
$c\neq0$.
\end{thm}

\begin{proof}
By Proposition \ref{pflow-rev}-(2), we have that
\[
\left\langle \Psi_{h_{\mathcal{A},\,\mathcal{B}}},\,\phi\right\rangle _{\mathfrak{F}}=-\sum_{x\in\mathcal{H}}h_{\mathcal{A},\,\mathcal{B}}(x)\,(\mbox{\textup{div}}\,\phi)(x)\;.
\]
Thus, by the Cauchy--Schwarz inequality, it holds that
\begin{align*}
\bigg[\sum_{x\in\mathcal{H}}h_{\mathcal{A},\,\mathcal{B}}(x)\,(\mbox{\textup{div}}\,\phi)(x)\bigg]^{2} & =\left\langle \Psi_{h_{\mathcal{A},\,\mathcal{B}}},\,\phi\right\rangle _{\mathfrak{F}}^{2}\\
 & \le\|\Psi_{h_{\mathcal{A},\,\mathcal{B}}}\|_{\mathfrak{F}}^{2}\,\Vert\phi\Vert_{\mathfrak{F}}^{2}=\mathrm{cap}(\mathcal{A},\,\mathcal{B})\,\Vert\phi\Vert_{\mathfrak{F}}^{2}\;.
\end{align*}
Hence, it holds that
\[
\mathrm{cap}(\mathcal{A},\,\mathcal{B})\ge\frac{1}{\|\phi\|_{\mathfrak{F}}^{2}}\bigg[\sum_{x\in\mathcal{H}}h_{\mathcal{A},\,\mathcal{B}}(x)\,(\mbox{\textup{div}}\,\phi)(x)\bigg]^{2}\;.
\]
From the Cauchy--Schwarz inequality, it is clear that the equality
holds only for $\phi=c\Psi_{h_{\mathcal{A},\,\mathcal{B}}}$, $c\neq0$.
\end{proof}
The advantage of this generalized Thomson principle is very clear.
We no longer impose the divergence-free condition on test flows, and
hence any flow approximating $h_{\mathcal{A},\,\mathcal{B}}$ (e.g.,
$\Psi_{h_{\mathcal{A},\,\mathcal{B}}}$) can be used as a test flow.
For instance, if one constructed a test function $f$ approximating
$h_{\mathcal{A},\,\mathcal{B}}$ and obtained an upper bound on the
capacity by injecting this test function $f$ to the Dirichlet principle,
then one can also use $\Psi_{f}$ as the test flow in this generalized
Thomson principle. If we encounter a technical issue in a certain
region, we can modify the flow accordingly in this region to obtain
a test flow. This idea was used in \cite{Kim-Seo Potts} to analyze
the metastability of Ising and Potts models on large, fixed lattices
without external fields. For this model, the energy landscape is extremely
complex, and it is very difficult to construct a divergence-free flow.
We explain a special case of this result in Part \ref{pt2}.

Clearly, the crucial disadvantage of the generalized Thomson principle
is the appearance of the equilibrium potential in the variational
principle. Hence, this generalized version turns the difficulty stemming
from the divergence-free restriction to the difficulty of handling
the equilibrium potential. Of course, Proposition \ref{p_eqp_bd}
plays an important role in controlling the equilibrium potential.

\subsection{Non-reversible case }

Now, we no longer assume that the process $X(\cdot)$ is reversible.
Then, the variational problem becomes more complicated.
\begin{thm}
\label{t_genDTP_nonrev}The followings hold.
\begin{enumerate}
\item (Generalized Dirichlet principle) We have that \textbf{
\begin{equation}
\textup{cap}(\mathcal{A},\,\mathcal{B})=\inf_{f\in\mathfrak{C}_{1,\,0}(\mathcal{A},\,\mathcal{B}),\,\phi\in\mathfrak{F}}\left\{ ||\Phi_{f}-\phi||^{2}-2\sum_{x\in\mathcal{H}}h_{\mathcal{A},\,\mathcal{B}}(x)\,(\mbox{\textup{div}\,}\phi)(x)\right\} \;,\label{DPG}
\end{equation}
}and (\ref{DPo}) is a minimizer.
\item (Generalized Thomson principle) We have that
\begin{equation}
\textup{cap}(\mathcal{A},\,\mathcal{B})=\sup_{g\in\mathfrak{C}_{0,\,0}(\mathcal{A},\,\mathcal{B}),\,\psi\in\mathfrak{F}_{0}}\frac{1}{||\Phi_{g}-\psi||^{2}}\left[\sum_{x\in\mathcal{H}}h_{\mathcal{A},\,\mathcal{B}}(x)\,(\mbox{\textup{div}\,}\psi)(x)\right]^{2}\;,\label{TPG}
\end{equation}
and the constant multiples of (\ref{TPo}), i.e.,
\begin{equation}
(g,\,\psi)=\left(c\frac{h_{\mathcal{A},\mathcal{\,B}}^{\dagger}-h_{\mathcal{A},\mathcal{\,B}}}{2\textup{cap}(\mathcal{A},\,\mathcal{B})},\,-c\frac{\Phi_{h_{\mathcal{A},\mathcal{\,B}}^{\dagger}}+\Phi_{h_{\mathcal{A},\,\mathcal{B}}}^{*}}{2\textup{cap}(\mathcal{A},\,\mathcal{B})}\right)\;,\;\;c\neq0\label{TPo-gen}
\end{equation}
are maximizers.
\end{enumerate}
\end{thm}

\begin{proof}
In the proof of Theorem \ref{t_dt}, we showed that for $f\in\mathfrak{C}_{a,\,0}(\mathcal{A},\,\mathcal{B})$,
\begin{equation}
\left\langle \Psi_{h_{\mathcal{A},\,\mathcal{B}}},\,\Phi_{f}\right\rangle _{\mathfrak{F}}=a\,\textup{cap}(\mathcal{A},\,\mathcal{B})\;.\label{ob1}
\end{equation}
On the other hand, by Proposition \ref{pflow_nonrev}-(2), we have
\begin{equation}
\left\langle \Psi_{h_{\mathcal{A},\,\mathcal{B}}},\,\phi\right\rangle _{\mathfrak{F}}=-\sum_{x\in\mathcal{H}}h_{\mathcal{A},\,\mathcal{B}}(x)\,(\mbox{\textup{div}\,}\phi)(x)\;.\label{ob2}
\end{equation}

For part (1), let $f\in\mathfrak{C}_{a,\,0}(\mathcal{A},\,\mathcal{B})$.
Then, by (\ref{ob1}) and (\ref{ob2}),
\begin{equation}
\left\langle \Psi_{h_{\mathcal{A},\,\mathcal{B}}},\,\Phi_{f}-\phi\right\rangle _{\mathfrak{F}}=\textup{cap}(\mathcal{\mathcal{A}},\,\mathcal{B})+\sum_{x\in\mathcal{H}}h_{\mathcal{A},\,\mathcal{B}}(x)\,(\mbox{\textup{div}\,}\phi)(x)\;.\label{e451}
\end{equation}
Furthermore, by the Cauchy--Schwarz inequality and (\ref{psif2})
(which still holds for non-reversible processes)
\begin{equation}
\left\langle \Psi_{h_{\mathcal{A},\,\mathcal{B}}},\,\Phi_{f}-\phi\right\rangle _{\mathfrak{F}}^{2}\le\Vert\Psi_{h_{\mathcal{A},\,\mathcal{B}}}\Vert_{\mathfrak{F}}^{2}\,\Vert\Phi_{f}-\phi\Vert_{\mathfrak{F}}^{2}=\textup{cap}(\mathcal{\mathcal{A}},\,\mathcal{B})\,\Vert\Phi_{f}-\phi\Vert_{\mathfrak{F}}^{2}\;.\label{e452}
\end{equation}
By (\ref{e451}) and (\ref{e452}),
\[
\textup{cap}(\mathcal{\mathcal{A}},\,\mathcal{B})\,\Vert\Phi_{f}-\phi\Vert_{\mathfrak{F}}^{2}\ge\textup{cap}(\mathcal{\mathcal{A}},\,\mathcal{B})^{2}+2\textup{cap}(\mathcal{\mathcal{A}},\,\mathcal{B})\sum_{x\in\mathcal{H}}h_{\mathcal{A},\,\mathcal{B}}(x)\,(\mbox{\textup{div}\,}\phi)(x)\;.
\]
Thus, part (1) is proved if we check that the equality holds for (\ref{DPo}).

The proof of part (2) is similar. For $g\in\mathfrak{C}_{0,\,0}(\mathcal{A},\,\mathcal{B})$,
again by (\ref{ob1}) and (\ref{ob2}), we have
\begin{align*}
\left\langle \Psi_{h_{\mathcal{A},\,\mathcal{B}}},\,\Phi_{g}-\psi\right\rangle _{\mathfrak{F}} & =\sum_{x\in\mathcal{H}}h_{\mathcal{A},\,\mathcal{B}}(x)\,(\mbox{\textup{div}\,}\psi)(x)\;.
\end{align*}
Hence, by the computations as before, the proof of part (2) is completed.
\end{proof}
\begin{rem}
We did not attempt to characterize all the optimizers in the previous
principles.
\end{rem}

When we use these principles, it is important to control terms of
the form
\[
\sum_{x\in\mathcal{H}}h_{\mathcal{A},\,\mathcal{B}}(x)\,(\mbox{\textup{div}\,}\phi)(x)\;.
\]
For the Thomson principle, we used $\psi$, instead of $\phi$, to
denote the test flow, but in what follows we denote by $\phi$ the
flow for the Thomson principle as well for convenience.

Indeed, this is trade-off in order to avoid the construction of a
divergence-free flow. By the property of the equilibrium potential
(cf. (\ref{e_eq_prop})), this summation can be decomposed into
\[
(\mbox{\textup{div}\,}\phi)(\mathcal{A})+\sum_{x\in(\mathcal{A}\cup\mathcal{B})^{c}}h_{\mathcal{A},\,\mathcal{B}}(x)\,(\mbox{\textup{div}\,}\phi)(x)\;.
\]
If we take the test function and flow as a good approximation of the
optimizers (\ref{DPo}) and (\ref{TPo}), we have $(\mbox{\textup{div}\,}\phi)(\mathcal{A})\simeq0$
for the Dirichlet principle and $(\mbox{\textup{div}\,}\phi)(\mathcal{A})\simeq1$
for the Thomson principle. Since $\phi$ can be approximately divergence-free
on $(\mathcal{A}\cup\mathcal{B})^{c}$, we also have
\[
\sum_{x\in(\mathcal{A}\cup\mathcal{B})^{c}}h_{\mathcal{A},\,\mathcal{B}}(x)\,(\mbox{\textup{div}\,}\phi)(x)\simeq0\;.
\]
Since the equilibrium potential is trivially bounded by $1$, we may
hope
\[
\sum_{x\in(\mathcal{A}\cup\mathcal{B})^{c}}|(\mbox{\textup{div}\,}\phi)(x)|\simeq0\;,
\]
but in general it may not be true (since there are too many elements
in $(\mathcal{A}\cup\mathcal{B})^{c}$). Instead, we need to decompose
$(\mathcal{A}\cup\mathcal{B})^{c}$ into two regions $\mathcal{C}_{1}$
and $\mathcal{C}_{2}$ so that
\[
\sum_{x\in\mathcal{C}_{1}}|(\mbox{\textup{div}\,}\phi)(x)|\simeq0\;,
\]
but on $\mathcal{C}_{2}$ the summation is small because $h_{\mathcal{A},\,\mathcal{B}}$
is small. To prove that $h_{\mathcal{A},\,\mathcal{B}}$ is sufficiently
small, Proposition \ref{p_eqp_bd} can be useful.

\section{\label{sec_col}Collapsed Processes}

In this section, we introduce the notion known as the collapsed process,
which is essentially obtained by contracting a subset $\mathcal{E}\subset\mathcal{H}$
to a single point $\mathfrak{e}$. This process was introduced in
\cite{G-L} to study the Dirichlet principle for non-reversible processes.
Moreover, in \cite{L-Seo NR1}, it is observed that the collapsed
process is a crucial notion (along with the capacity) in the precise
estimate of the so-called mean jump rate, which is key to the martingale
approach of metastability (cf. \cite{B-L TM,B-L TM2,B-L MG}).

In this section, we fix a set $\mathcal{E}\subset\mathcal{H}$. We
note that the contents of the current subsection are from \cite[Section 8]{L-Seo NR1}.

\subsection{Definition of collapsed process}

As mentioned earlier, our aim is collapsing a set $\mathcal{E}$ into
a single point $\mathfrak{e}$. To this end, let us first define the
state space $\overline{\mathcal{H}}=(\mathcal{H}\setminus\mathcal{E})\cup\{\mathfrak{e}\}$.
Then, (recalling that $\mu(\cdot)$ is the invariant measure for the
process $X(\cdot)$) define a rate $\overline{r}:\overline{\mathcal{H}}\times\overline{\mathcal{H}}\rightarrow[0,\,\infty)$
as
\begin{equation}
\begin{cases}
\overline{r}(x,\,y)=r(x,\,y) & \text{for }x,\,y\in\mathcal{H}\setminus\mathcal{E}\;,\\
\overline{r}(x,\,\mathfrak{e})=\sum_{\boldsymbol{z}\in\mathcal{E}}r(x,\,z) & \text{for }x\in\mathcal{H}\setminus\mathcal{E}\;,\\
\overline{r}(\mathfrak{e},\,y)=\frac{1}{\mu(\mathcal{E})}\sum_{\boldsymbol{z}\in\mathcal{E}}\mu(z)r(z,\,y) & \text{for }y\in\mathcal{H}\setminus\mathcal{E}\;.
\end{cases}\label{eq:coll rate}
\end{equation}
The collapsed process is defined as a continuous-time Markov process
$(\overline{X}(t))_{t\ge0}$ on $\overline{\mathcal{H}}$ with rate
$\overline{r}(\cdot,\,\cdot)$.

Denote by $\overline{\mathbb{P}}_{x}$ the law of $\overline{X}(\cdot)$
starting from $x$, and by $\overline{\mathscr{L}}$ and $\overline{\mathscr{D}}(\cdot)$
the generator and the Dirichlet form corresponding to the collapsed
process $\overline{X}(\cdot)$, respectively. Define a probability
measure $\overline{\mu}(\cdot)$ on $\overline{\mathcal{H}}$ as
\begin{equation}
\begin{cases}
\overline{\mu}(x)\;=\;\mu(x) & \text{if }x\in\mathcal{H}\setminus\mathcal{E}\;,\\
\overline{\mu}(\mathfrak{e})\;=\;\mu(\mathcal{E})\;.
\end{cases}\label{eq:coll_inv}
\end{equation}

\begin{xca}
\label{ex_collinv}Answer the following questions.
\begin{enumerate}
\item Prove that the measure $\overline{\mu}(\cdot)$ is the invariant measure
for the process $\overline{X}(\cdot)$.
\item Prove that the process $\overline{X}(\cdot)$ is reversible if the
process $X(\cdot)$ is reversible. Is the converse true?
\end{enumerate}
\end{xca}

\subsection{Flow space of collapsed process}

Next, we investigate the flow structure with respect to the collapsed
process $\overline{X}(\cdot)$. For $x,\,y\in\mathcal{H}$, we defined
the conductance between $x$ and $y$ with respect to the original
process $X(\cdot)$ as (cf. (\ref{cxy}))
\[
c(x,\,y)=\mu(x)r(x,\,y)\;.
\]
Similarly, for $x,\,y\in\overline{\mathcal{H}}$, we define the conductance
with respect to the collapsed process $\overline{X}(\cdot)$ as
\[
\overline{c}(x,y)=\overline{\mu}(x)\overline{r}(x,\,y)\;.
\]
Then, by (\ref{eq:coll rate}) and (\ref{eq:coll_inv}), this conductance
$\overline{c}(\cdot,\,\cdot)$ can be rewritten as
\begin{equation}
\begin{cases}
\overline{c}(x,\,y)=c(x,\,y) & \text{for }x,\,y\in\mathcal{H}\setminus\mathcal{E}\;,\\
\overline{c}(x,\,\mathfrak{e})=\sum_{z\in\mathcal{E}}c(x,\,z) & \text{for }x\in\mathcal{H}\setminus\mathcal{E}\;,\\
\overline{c}(\mathfrak{e},\,y)=\sum_{z\in\mathcal{E}}c(z,\,y) & \text{for }y\in\mathcal{H}\setminus\mathcal{E}\;.
\end{cases}\label{eq:coll conduct}
\end{equation}
Define the symmetrized conductance as
\[
\overline{c}^{s}(x,\,y)=\frac{1}{2}[\overline{c}(x,\,y)+\overline{c}(y,\,x)]\;\;\;\;;\;x,\,y\in\overline{\mathcal{H}}\;.
\]

For $x,\,y\in\overline{\mathcal{H}}$, we write $x\sim y$ if $\overline{c}^{s}(x,\,y)>0$.
Since $\overline{c}^{s}(x,\,y)=\overline{c}^{s}(y,\,x)$, we observe
that $x\sim y$ if and only if $y\sim x$. Then, the set of directed
edges are defined by
\begin{equation}
\overline{\mathfrak{E}}=\{(x,\,y)\in\mathcal{H}\times\mathcal{H}:x\sim y\}\;.\label{edge-1}
\end{equation}
As before, we can define a flow structure on the set $\overline{\mathfrak{F}}$
of flows on $\overline{\mathfrak{E}}$ which are anti-symmetric functions
on $\overline{\mathfrak{E}}$. Then, we can induce the Hilbert space
structure on $\overline{\mathfrak{F}}$, as we did in Sections \ref{sec13}
and \ref{sec14}. Denote the corresponding inner product and the flow
norm by $\left\langle \cdot,\,\cdot\right\rangle _{\mathcal{\overline{\mathfrak{F}}}}$
and $\Vert\cdot\Vert_{\overline{\mathfrak{F}}}$, respectively. In
particular, we can write
\begin{align*}
\left\langle \phi,\,\psi\right\rangle _{\mathcal{\overline{\mathfrak{F}}}} & =\frac{1}{2}\sum_{(x,\,y)\in\mathfrak{\overline{E}}}\frac{\phi(x,\,y)\psi(x,\,y)}{\overline{c}^{s}(x,\,y)}\;,\;\text{and}\\
\Vert\phi\Vert_{\overline{\mathfrak{F}}}^{2} & =\frac{1}{2}\sum_{(x,\,y)\in\overline{\mathfrak{E}}}\frac{\phi(x,\,y)^{2}}{\overline{c}^{s}(x,\,y)}\;.
\end{align*}

For each flow $\phi\in\mathfrak{F}$, define the collapsed flow $\overline{\phi}\in\overline{\mathfrak{F}}$
by
\begin{equation}
\begin{cases}
\overline{\phi}(x,\,y)=\phi(x,\,y) & \text{for }x,\,y\in\mathcal{H}\setminus\mathcal{E}\;,\\
\overline{\phi}(x,\,\mathfrak{e})=\sum_{z\in\mathcal{E}}\phi(x,\,z) & \text{for }x\in\mathcal{H}\setminus\mathcal{E}\;,\\
\overline{\phi}(\mathfrak{e},\,y)=\sum_{z\in\mathcal{E}}\phi(z,\,y) & \text{for }y\in\mathcal{H}\setminus\mathcal{E}\;.
\end{cases}\label{eq:coll_flow}
\end{equation}

\begin{xca}
Prove that
\begin{equation}
\begin{cases}
(\mbox{div}\,\overline{\phi})(x)\;=\;(\mbox{div }\phi)(x) & \text{for }x\in\mathcal{H}\setminus\mathcal{E}\;,\\
(\mbox{div}\,\overline{\phi})(\mathfrak{e})\;=\;(\mbox{div }\phi)(\mathcal{E})\;.
\end{cases}\label{eq:coll_div}
\end{equation}
\end{xca}

The following contraction property of the flow norm is useful later.
\begin{lem}
\label{lem:coll_flow}For all $\phi\in\mathfrak{F}$ and its collapsed
flow $\overline{\phi}\in\overline{\mathfrak{F}}$, it holds that
\[
\Vert\overline{\phi}\Vert_{\overline{\mathfrak{F}}}\le\Vert\phi\Vert_{\mathfrak{F}}\;.
\]
Moreover, the equality holds if and only if
\begin{equation}
\begin{cases}
\phi(x,\,y)=0 & \text{if }x,\,y\in\mathcal{E}\;,\text{ and}\\
\frac{\phi(x,\,y)}{c^{s}(x,\,y)}=\frac{\phi(x',\,y)}{c^{s}(x',\,y)} & \text{if }y\in\mathcal{E}\text{ and }x,\,x'\in\mathcal{H}\setminus\mathcal{E\text{ satisfies }}x\sim y\text{ and }x'\sim y\;.
\end{cases}\label{eq:eqcon}
\end{equation}
\end{lem}

\begin{proof}
Decompose the flow norm of the flow $\phi$ as
\[
\Vert\phi\Vert_{\mathfrak{F}}^{2}=\frac{A_{1}}{2}+A_{2}+\frac{A_{3}}{2}\;,
\]
where
\begin{align*}
A_{1}= & \sum_{(x,\,y)\in\mathfrak{E}:x,\,y\in\mathcal{H}\setminus\mathcal{E}}\frac{\phi(x,\,y)^{2}}{c^{s}(x,\,y)}\;,\\
A_{2}= & \sum_{(x,\,y)\in\mathfrak{E}:x\in\mathcal{H}\setminus\mathcal{E},\,y\in\mathcal{E}}\frac{\phi(x,\,y)^{2}}{c^{s}(x,\,y)}\;,\;\text{and}\\
A_{3}= & \sum_{(x,\,y)\in\mathfrak{E}:x,\,y\in\mathcal{E}}\frac{\phi(x,\,y)^{2}}{c^{s}(x,\,y)}\;.
\end{align*}
Then, decompose the flow norm of the collapsed flow $\overline{\phi}$
as
\[
\Vert\overline{\phi}\Vert_{\mathfrak{\overline{F}}}^{2}=\frac{\overline{A}_{1}}{2}+\overline{A}_{2}\;,
\]
 where
\begin{eqnarray*}
\overline{A}_{1} & = & \sum_{(x,\,y)\in\mathfrak{\overline{E}}:x,\,y\in\mathcal{H}\setminus\mathcal{E}}\frac{\overline{\phi}(x,\,y)^{2}}{\overline{c}^{s}(x,\,y)}\;\;\;\;\text{and}\\
\overline{A}_{2} & = & \sum_{x\in\overline{\mathcal{H}}:(x,\,\mathfrak{e})\in\mathfrak{\overline{E}}}\frac{\overline{\phi}(x,\,\mathfrak{e})^{2}}{\overline{c}^{s}(x,\,\mathfrak{e})}\;.
\end{eqnarray*}
By (\ref{eq:coll conduct}) and (\ref{eq:coll_flow}), we immediately
have that $A_{1}=\overline{A}_{1}$.

Therefore, it suffices to prove $A_{2}\ge\overline{A}_{2}$. For each
$x\in\mathcal{H}\setminus\mathcal{E}$ adjacent to at least one point
of $\mathcal{E}$, by (\ref{eq:coll conduct}), (\ref{eq:coll_flow}),
and the Cauchy--Schwarz inequality, we obtain
\begin{eqnarray*}
\sum_{y\in\mathcal{E}:(x,\,y)\in\mathfrak{E}}\frac{\phi(x,\,y)^{2}}{c^{s}(x,\,y)} & \ge & \frac{\left[\sum_{y\in\mathcal{E}:(x,\,y)\in\mathfrak{E}}\phi(x,\,y)\right]^{2}}{\sum_{y\in\mathcal{E}:(x,\,y)\in\mathfrak{E}}c^{s}(x,\,y)}=\frac{\overline{\phi}(x,\,\mathfrak{e})^{2}}{\overline{c}^{s}(x,\,\mathfrak{e})}\;.
\end{eqnarray*}
By adding this inequality over $x\in\mathcal{H}\setminus\mathcal{E}$,
we obtain $A_{2}\ge\overline{A}_{2}$, and the proof is completed.
\end{proof}
\begin{xca}
\label{ex:eqcon}For $f:\mathcal{H}\rightarrow\mathbb{R}$ which is
constant over $\mathcal{E}$, prove that the flow $\Psi_{f}$ satisfies
the equality condition (\ref{eq:eqcon}).
\end{xca}

If a function $f:\mathcal{H}\rightarrow\mathbb{R}$ is constant over
$\mathcal{E}$, we define a collapsed function $\overline{f}:\overline{\mathcal{H}}\rightarrow\mathbb{R}$
as
\begin{equation}
\begin{cases}
\overline{f}(x)=f(x)\;\;\text{if }x\in\mathcal{H}\setminus\mathcal{E}\;,\\
\overline{f}(\mathfrak{e})=\mbox{the constant value of \ensuremath{f} on }\mathcal{E}\;.
\end{cases}\label{col_f}
\end{equation}

\begin{lem}
\label{lem_coll_dr}Suppose that the functions $f,\,g:\mathcal{H}\rightarrow\mathbb{R}$
are constant over $\mathcal{E}$, and let $\overline{f},\,\overline{g}:\overline{\mathcal{H}}\rightarrow\mathbb{R}$
be the collapsed function of $f,\,g$ (cf. (\ref{col_f})), respectively.
Then, we have
\begin{equation}
\left\langle g,\,-\mathscr{L}f\right\rangle _{\mu}=\left\langle \overline{g},\,-\overline{\mathscr{L}}\overline{f}\right\rangle _{\overline{\mu}}\;.\label{e_GLF}
\end{equation}
In particular, we have
\begin{equation}
\overline{\mathscr{D}}(\overline{f})=\mathscr{D}(f)\;.\label{e_DF}
\end{equation}
\end{lem}

\begin{proof}
Since $f$ is constant over $\mathcal{E}$, we can write
\begin{align}
 & \left\langle g,\,-\mathscr{L}f\right\rangle _{\mu}\nonumber \\
= & \frac{1}{2}\left[\sum_{x\in\mathcal{H}\setminus\mathcal{E}}\sum_{y\in\mathcal{H}\setminus\mathcal{E}}+\sum_{x\in\mathcal{H}\setminus\mathcal{E}}\sum_{y\in\mathcal{E}}+\sum_{x\in\mathcal{E}}\sum_{y\in\mathcal{H}\setminus\mathcal{E}}\right]\mu(x)r(x,\,y)[f(y)-f(x)]g(x)\;.\label{0sm}
\end{align}
Note that the first summation is equal to
\begin{equation}
\sum_{x\in\mathcal{H}\setminus\mathcal{E}}\sum_{y\in\mathcal{H}\setminus\mathcal{E}}\overline{\mu}(x)\overline{r}(x,\,y)[\overline{f}(y)-\overline{f}(x)]\overline{g}(x)\;,\label{eq:1sm}
\end{equation}
since $\mu=\overline{\mu}$, $r=\overline{r}$, $f=\overline{f}$,
and $g=\overline{g}$ on $\mathcal{H}\setminus\mathcal{E}$. On the
other hand, we have $f(y)=\overline{f}(\mathfrak{e})$ for all $y\in\mathcal{E}$,
and thus the second summation is equal to
\begin{align}
 & \sum_{x\in\mathcal{H}\setminus\mathcal{E}}\sum_{y\in\mathcal{E}}\overline{\mu}(x)r(x,\,y)[\overline{f}(\mathfrak{e})-\overline{f}(x)]\overline{g}(x)\nonumber \\
= & \sum_{x\in\mathcal{H}\setminus\mathcal{E}}\overline{\mu}(x)\overline{r}(x,\,\mathfrak{e})[\overline{f}(\mathfrak{e})-\overline{f}(x)]\overline{g}(x)\;,\label{2sm}
\end{align}
where the equality follows from the second line of (\ref{eq:coll rate}).
Finally, a similar computation yields that the third summation is
equal to
\begin{align}
 & \sum_{x\in\mathcal{E}}\sum_{y\in\mathcal{H}\setminus\mathcal{E}}\mu(x)r(x,\,y)[\overline{f}(y)-\overline{f}(\mathfrak{e})]\overline{g}(\mathfrak{e})\nonumber \\
= & \sum_{y\in\mathcal{H}\setminus\mathcal{E}}\overline{\mu}(\mathfrak{e})\overline{r}(\mathfrak{e},\,y)[\overline{f}(y)-\overline{f}(\mathfrak{e})]\overline{g}(\mathfrak{e})\;,\label{3sm}
\end{align}
where the equality follows from the third line of (\ref{eq:coll rate})
and (\ref{eq:coll_inv}). By inserting (\ref{eq:1sm}), (\ref{2sm})
and (\ref{3sm}) into (\ref{0sm}), we can conclude that
\begin{align*}
 & \left\langle g,\,-\mathscr{L}f\right\rangle _{\mu}\\
= & \frac{1}{2}\left[\sum_{x\in\mathcal{H}\setminus\mathcal{E}}\sum_{y\in\mathcal{H}\setminus\mathcal{E}}+\sum_{x\in\mathcal{H}\setminus\mathcal{E}}\sum_{y\in\{\mathfrak{e}\}}+\sum_{x\in\{\mathfrak{e}\}}\sum_{y\in\mathcal{H}\setminus\mathcal{E}}\right]\overline{\mu}(x)\overline{r}(x,\,y)[\overline{f}(y)-\overline{f}(x)]\overline{g}(x)\\
= & \left\langle \overline{g},\,-\overline{\mathscr{L}}\overline{f}\right\rangle _{\overline{\mu}}\;,
\end{align*}
and the proof of (\ref{e_GLF}) is completed. Now, (\ref{e_Df}) follows
from (\ref{e_GLF}) by inserting $g=f$.
\end{proof}
For a function $g:\overline{\mathcal{H}}\rightarrow\mathbb{R}$, define
$\overline{\Phi}_{g},\,\overline{\Phi}_{g}^{*}$ and $\overline{\Psi}_{g}$
as, for $x,\,y\in\overline{\mathcal{H}}$,
\begin{eqnarray}
\overline{\Phi}_{g}(x,\,y) & = & g(y)\overline{c}(y,\,x)-g(x)\overline{c}(x,\,y)\;,\label{eq: col  phi_f}\\
\overline{\Phi}_{g}^{*}(x,\,y) & = & g(y)\overline{c}(x,\,y)-g(x)\overline{c}(y,\,x)\;,\label{eq: col phi*_f}\\
\overline{\Psi}_{g}(x,\,y) & = & \overline{c}^{s}(x,\,y)(g(y)-g(x))\;.\label{eq: col psi_f}
\end{eqnarray}

\begin{lem}
\label{lem:coll_Phi_f}Suppose that the function $f:\mathcal{H}\rightarrow\mathbb{R}$
is constant over $\mathcal{E}$, and let $\overline{f}:\overline{\mathcal{H}}\rightarrow\mathbb{R}$
be the collapsed function of $f$ (cf. (\ref{col_f})). Then, the
flow $\overline{\Phi_{f}}$, which is the collapsed flow of $\Phi_{f}$
defined in (\ref{flow}), coincides with the flow $\overline{\Phi}_{\overline{f}}$.
Similarly, we have that
\begin{equation}
\overline{\Phi_{f}^{*}}=\overline{\Phi}_{\overline{f}}^{*}\;\;\;\text{and\;\;\;}\overline{\Psi_{f}}=\overline{\Psi}_{\overline{f}}\;.\label{ex44}
\end{equation}
\end{lem}

\begin{proof}
We only prove that two flows $\overline{\Phi_{f}}$ and $\overline{\Phi}_{\overline{f}}$
coincide, and leave the proof for the other two as exercise, since
the proofs are quite similar.

Since $\overline{\Phi_{f}}(x,\,y)=\overline{\Phi}_{\overline{f}}(x,\,y)$
for $x,\,y\in\mathcal{H}\setminus\mathcal{E}$ holds trivially from
the definitions, it suffices to prove that $\overline{\Phi_{f}}(x,\,\mathfrak{e})=\overline{\Phi}_{\overline{f}}(x,\,\mathfrak{e})$
for $x\in\mathcal{H}\setminus\mathcal{E}$. This can be verified by
\begin{eqnarray*}
\overline{\Phi_{f}}(x,\,\mathfrak{e})=\sum_{\boldsymbol{z}\in\mathcal{E}}\Phi_{f}(x,\,z) & = & \sum_{z\in\mathcal{E}}\left[f(z)c(z,\,x)-f(x)c(x,\,z)\right]\\
 & = & \overline{f}(\mathfrak{e})\overline{c}(\mathfrak{e},\,x)-\overline{f}(x)\overline{c}(x,\,\mathfrak{e})\\
 & = & \overline{\Phi}_{\overline{f}}(x,\,\mathfrak{e})\;.
\end{eqnarray*}
\end{proof}
\begin{xca}
\label{ex47}Prove (\ref{ex44}).
\end{xca}

\subsection{Capacity and sector condition of collapsed process}

For two non-empty and disjoint subsets $\mathcal{A}$ and $\mathcal{B}$
of $\overline{\mathcal{H}},$we denote by $\overline{h}_{\mathcal{A},\,\mathcal{B}}:\mathcal{H}\rightarrow\mathbb{R}$
the equilibrium potential between $\mathcal{A}$ and $\mathcal{B}$,
and we denote by $\overline{\textup{cap}}(\mathcal{A},\,\mathcal{B})$
and $\overline{\textup{cap}}^{s}(\mathcal{A},\,\mathcal{B})$ the
capacity between $\mathcal{A}$ and $\mathcal{B}$ with respect to
the collapsed process $\overline{X}(\cdot)$ and the symmetrized process
$\overline{X}^{s}(\cdot)$ of $\overline{X}(\cdot)$ (which is a Markov
process on $\overline{\mathcal{H}}$ associated with the generator
$\frac{1}{2}(\overline{\mathscr{L}}+\overline{\mathscr{L}}^{\dagger})$,
where $\overline{\mathscr{L}}^{\dagger}$ is the adjoint generator
of $\overline{\mathscr{L}}$), respectively. In general, for $\mathcal{A},\,\mathcal{B}\subset\mathcal{H}\setminus\mathcal{E}$,
it is difficult to compare $\overline{\textup{cap}}(\mathcal{A},\,\mathcal{B})$
and $\textup{cap}(\mathcal{A},\,\mathcal{B})$.
\begin{xca}
Suppose that $\mathcal{A}$ and $\mathcal{B}$ are two non-empty and
disjoint subsets of $\mathcal{H}\setminus\mathcal{E}$. Then, can
you prove either $\overline{\textup{cap}}(\mathcal{A},\,\mathcal{B})\le\textup{cap}(\mathcal{A},\,\mathcal{B})$
or $\textup{cap}(\mathcal{A},\,\mathcal{B})\le\overline{\textup{cap}}(\mathcal{A},\,\mathcal{B})$?
\end{xca}

However, we have the following identity, which is useful in later
discussions.
\begin{lem}
\label{lem:coll_CAP}For any non-empty $\mathcal{A}\subset\mathcal{H}\setminus\mathcal{E}$,
we have
\[
\overline{\mbox{\textup{cap}}}(\mathfrak{e},\,\mathcal{A})=\textup{cap}(\mathcal{E},\,\mathcal{A})\;.
\]
\end{lem}

\begin{proof}
Recall that $h_{\mathcal{E},\,\mathcal{A}}(\cdot)$ denotes the equilibrium
potential between $\mathcal{E}$ and $\mathcal{A}$. Since the behaviors
of the processes $X(\cdot)$ and $\overline{X}(\cdot)$ are identical
on $\mathcal{H}\setminus\mathcal{E}$, we immediately have that
\[
h_{\mathcal{E},\,\mathcal{A}}(x)=\overline{h}_{\mathfrak{e},\,\mathcal{A}}(x)\text{ for all }x\in\mathcal{H}\setminus\mathcal{E}\;.
\]
Since $h_{\mathcal{E},\,\mathcal{A}}\equiv1$ on $\mathcal{E}$ and
$\overline{h}_{\mathfrak{e},\,\mathcal{A}}(\mathfrak{e})=1$, we can
conclude that $\overline{h}_{\mathfrak{e},\,\mathcal{A}}(\cdot)$
is the collapsed function of $h_{\mathcal{E},\,\mathcal{A}}(\cdot)$,
i.e,.,
\[
\overline{h}_{\mathfrak{e},\,\mathcal{A}}=\overline{h_{\mathcal{E},\,\mathcal{A}}}\;.
\]
Therefore, by Lemma \ref{lem_coll_dr}, we can conclude that
\[
\overline{\mbox{\textup{cap}}}(\mathfrak{e},\,\mathcal{A})=\overline{\mathscr{D}}(\overline{h}_{\mathfrak{e},\,\mathcal{A}})=\overline{\mathscr{D}}(\overline{h_{\mathcal{E},\,\mathcal{A}}})=\mathscr{D}(h_{\mathcal{E},\,\mathcal{A}})=\textup{cap}(\mathcal{E},\,\mathcal{A})\;.
\]
\end{proof}
Next, we assert that the sector condition of the original process
is inherited by the collapsed process.
\begin{lem}
\label{lem:col_sector}Suppose that the process $X(\cdot)$ satisfies
the sector condition with a constant $C>0$ (cf. Definition \ref{sector}).
Then, the process $\overline{X}(\cdot)$ also satisfies the sector
condition with the same constant $C$. In particular, it holds for
any two non-empty and disjoint subsets $\mathcal{A}$ and $\mathcal{B}$
of $\overline{\mathcal{H}}$ that
\begin{equation}
\overline{\textup{cap}}^{s}(\mathcal{A},\,\mathcal{B})\le\overline{\textup{cap}}(\mathcal{A},\,\mathcal{B})\le C\,\overline{\textup{cap}}^{s}(\mathcal{A},\,\mathcal{B})\;.\label{compcap}
\end{equation}
\end{lem}

\begin{proof}
For two functions $f,\,g:\overline{\mathcal{H}}\rightarrow\mathbb{R}$,
define their extended functions $F,\,G:\mathcal{H}\rightarrow\mathbb{R}$
as
\begin{eqnarray*}
F(x) & = & \begin{cases}
f(x) & \mbox{if }x\in\mathcal{H\setminus\mathcal{E}}\;,\\
f(\mathfrak{e}) & \mbox{if }x\in\mathcal{E}\;,
\end{cases}\;\;\;\mbox{and\;\;\;}G(x)=\begin{cases}
g(x) & \mbox{if }x\in H\setminus\mathcal{E}\;,\\
g(\mathfrak{e}) & \mbox{if }x\in\mathcal{E}\;.
\end{cases}
\end{eqnarray*}
so that
\begin{equation}
\overline{F}=f\;\;\;\;\text{and\;\;\;\;}\overline{G}=g\;.\label{e_FFGG}
\end{equation}
Hence, by (\ref{e_FFGG}), Lemma \ref{lem_coll_dr}, and the sector
condition of $X(\cdot)$,
\[
\left\langle f,-\overline{\mathscr{L}}g\right\rangle _{\overline{\mu}}=\left\langle F,-\mathscr{L}G\right\rangle _{\mu}\le C\mathscr{D}(F)\mathscr{D}(G)=C\overline{\mathscr{D}}(f)\overline{\mathscr{D}}(g)\;,
\]
and hence the process $\overline{X}(\cdot)$ also satisfies the sector
condition with a constant $C>0$. Now, (\ref{compcap}) is clear from
Propositions \ref{p_capcomp1} and \ref{p_capcomp2}.
\end{proof}
\newpage

\part{\label{pt2}Two-dimensional Ising Model without External Field}

In this second part of the lecture note, as an application of the
general theory developed so far, we thoroughly analyze the metastable
behavior of the Ising model on large but fixed lattice boxes. In particular,
we focus on the model without an external field, which posed a longstanding
mathematical challenge because of the complexity of the energy landscape.
The dynamics is reversible, and the analysis is based on the Dirichlet
principle (Theorem \ref{t_DP_rev}) and the generalized Thomson principle
(Theorem \ref{t_GTP}).

The contents of the current part is based on \cite{Kim-Seo Potts}
which considered more complex models, namely the Potts model and the
model in three-dimensional boxes. We did not investigate these models
in this note, since the two-dimensional Ising model is enough to deliver
the core of our idea.

\section{Ising Model on Two-dimensional Lattice}

\subsection{\label{sec21}Model }

In this subsection, we introduce the model and review its basic features.

\subsubsection*{Ising model }

For two positive integers $K,\,L$, we write
\begin{equation}
\Lambda=\mathbb{T}_{K}\times\mathbb{T}_{L}\;,\label{box_torus}
\end{equation}
where $\mathbb{T}_{k}=\mathbb{Z}/(k\mathbb{Z})$ is the discrete one-dimensional
torus. For the convenience of the discussion, we assume that $K\le L$
and moreover $K\ge5$.

We will consider the spin system on $\Lambda$; hence, we consider
a spin system on the box with periodic boundary conditions. The model
that we consider in this second part is defined now.
\begin{defn}[Ising model on $\Lambda$ without external field]
\begin{itemize}
\item Denote by $\Omega=\{+,\,-\}$ the set of spins and by $\mathcal{X}=\Omega^{\Lambda}$
the space of spin configurations on the box $\Lambda$. A configuration
$\sigma\in\mathcal{X}$ is written as $\sigma=(\sigma(x))_{x\in\Lambda}$
where $\sigma(x)\in\Omega$ denotes the spin of $\sigma$ at site
$x\in\Lambda$.
\item For $x,\,y\in\Lambda$, let us write $x\sim y$ if they are neighboring
sites in $\Lambda$, that is, $\Vert x-y\Vert=1$, where $\Vert\cdot\Vert$
denotes the Euclidean distance in $\Lambda$ where the periodic boundary
condition has to be taken into account.
\item Define the Hamiltonian $H:\mathcal{X}\rightarrow\mathbb{R}$ as
\begin{equation}
H(\sigma)=\sum_{x\sim y}\mathbf{1}\{\sigma(x)\ne\sigma(y)\}\;\;\;\;;\;\sigma\in\mathcal{X}\;.\label{ham}
\end{equation}
Note that there is no external field in this Hamiltonian; only the
spin--spin interaction is considered.
\item Denote by $\mu_{\beta}(\cdot)$ the Gibbs measure on $\mathcal{X}$
associated to the Hamiltonian $H$ at inverse temperature $\beta>0$,
i.e.,
\begin{equation}
\mu_{\beta}(\sigma)=\frac{1}{Z_{\beta}}e^{-\beta H(\sigma)}\;\;\;\;;\;\sigma\in\mathcal{X}\;,\label{mudef}
\end{equation}
where $Z_{\beta}$ is the partition function defined by
\begin{equation}
Z_{\beta}=\sum_{\sigma\in\mathcal{X}}e^{-\beta H(\sigma)}\;.\label{partition}
\end{equation}
\end{itemize}
The spin system on $\Lambda$ corresponding to the probability measure
$\mu_{\beta}(\cdot)$ on $\mathcal{X}_{d}$ is called the Ising model.
\end{defn}

\subsubsection*{Ground states}

We denote by $\boxplus\in\mathcal{X}$ (resp. $\boxminus\in\mathcal{X}$)
the configuration such that all spins are $+$ (resp. $-$), i.e.,
$\boxplus(x)=+$ (resp. $\boxminus(x)=-$) for all $x\in\Lambda$.
We write
\begin{equation}
\mathcal{S}=\{\boxplus,\,\boxminus\}\subset\mathcal{X}\;.\label{e_Sd}
\end{equation}
Note that the Hamiltonian $H(\cdot)$ attains its minimum value $0$
(only) at $\mathcal{S}$. Hence, $\boxplus$ and $\boxminus$ are
the ground states of the model. Based on this observation, we obtain
the following characterization of the partition function $Z_{\beta}$
defined in (\ref{partition}), as well as the Gibbs measure $\mu_{\beta}$
as $\beta\rightarrow\infty$.
\begin{prop}
\label{p_inv}The following hold:
\begin{enumerate}
\item The partition function satisfies the asymptotics
\begin{equation}
Z_{\beta}=2+O(e^{-2\beta})\;.\label{Zbest}
\end{equation}
\item We have
\[
\lim_{\beta\rightarrow\infty}\mu_{\beta}(\boxplus)=\lim_{\beta\rightarrow\infty}\mu_{\beta}(\boxminus)=\frac{1}{2}\;,\;\;\;\text{and thus\;}\lim_{\beta\rightarrow\infty}\mu_{\beta}(\mathcal{S})=1\;.
\]
\end{enumerate}
\end{prop}

\begin{proof}
We can readily observe that $H(\sigma)\ge2$ for $\sigma\notin\mathcal{S}$.
The estimate (\ref{Zbest}) comes directly from this observation along
with the expression (\ref{partition}). Part (2) of the theorem follows
directly from part (1) and the expression (\ref{mudef}) of $\mu_{\beta}$.
\end{proof}

\subsubsection*{Continuous-time Metropolis dynamics }

We now define a continuous-time Metropolis-type Glauber dynamics which
is a standard heat-bath dynamics in the study of the Ising model (cf.
\cite{N-S Ising1}).\textbf{ }For $x\in\Lambda$, we denote by $\sigma^{x}\in\mathcal{X}$
the configuration obtained from $\sigma$ by flipping the spin at
site $x$.
\begin{defn}
The continuous-time Metropolis dynamics is defined as a continuous
time Markov process $\{\sigma_{\beta}(t)\}_{t\ge0}$ on $\mathcal{X}$
with transition rates
\begin{equation}
c_{\beta}(\sigma,\,\zeta)=\begin{cases}
e^{-\beta[H(\zeta)-H(\sigma)]_{+}} & \text{if }\zeta=\sigma^{x}\ne\sigma\text{ for some }x\in\Lambda\;,\\
0 & \text{otherwise}\;,
\end{cases}\label{e_rate}
\end{equation}
where $[a]_{+}=\max\{a,\,0\}$.
\end{defn}

For $\sigma,\,\zeta\in\mathcal{X}$, we write $\sigma\sim\zeta$ if
$c_{\beta}(\sigma,\,\zeta)>0$, i.e., if $\zeta$ is obtained from
$\sigma$ by flipping the spin at a site (or vice versa). Note that
the relationship $\sigma\sim\zeta$ does not depend on $\beta$. Moreover,
the following detailed balance condition holds:
\begin{equation}
\mu_{\beta}(\sigma)\,c_{\beta}(\sigma,\,\zeta)=\mu_{\beta}(\zeta)\,c_{\beta}(\zeta,\,\sigma)=\begin{cases}
\min\{\mu_{\beta}(\sigma),\,\mu_{\beta}(\zeta)\} & \text{if }\sigma\sim\zeta\;,\\
0 & \text{otherwise\;.}
\end{cases}\label{muprop}
\end{equation}
Consequently, $\mu_{\beta}(\cdot)$ is the unique\footnote{It is clear that the Markov process $\sigma_{\beta}(\cdot)$ is irreducible.}
invariant measure for the Markov process $\sigma_{\beta}(\cdot)$,
and furthermore $\sigma_{\beta}(\cdot)$ is reversible with respect
to $\mu_{\beta}(\cdot)$. We denote by $\mathbb{P}_{\sigma}^{\beta}$
the law of the process $\sigma_{\beta}(\cdot)$ starting from $\sigma$,
and by $\mathbb{E}_{\sigma}^{\beta}$ the associated expectation.

\subsubsection*{Metastability of the model}

The primary concern in this second part is the metastable behavior
of the process $\sigma_{\beta}(\cdot)$ defined above when $\beta$
is large. More precisely, by the expression (\ref{e_rate}) of the
jump rate, we can see that the dynamics $\sigma_{\beta}(\cdot)$ tends
to lower the energy (for large $\beta$) since it jumps to a configuration
with higher energy with exponentially small rate. Hence, in view of
Proposition \ref{p_inv}, the process $\sigma_{\beta}(\cdot)$ starting
from a configuration $\boxplus$ may tend to stay in some neighborhood
of $\boxplus$ for a long time. However, by the irreducibility of
the process $\sigma_{\beta}(\cdot)$, it will eventually make a transition
to $\boxminus$. Similar behavior is expected to occur when the process
starts from $\boxminus$. Hence, such rare transitions between $\boxplus$
and $\boxminus$ will take place successively. This type of behavior
is the metastable behavior of the process $\sigma_{\beta}(\cdot)$.
In this part, we wish to quantitatively analyze this behavior to a
precise level. For instance, we will give precise asymptotic of the
mean transition time from $\boxplus$ to $\boxminus$ in the very
low temperature regime, i.e., when $\beta\rightarrow\infty$.

\subsection{\label{sec22}Main results }

We now explain the main results regarding the metastability of the
stochastic Ising model.

\subsubsection*{Energy barrier between ground states}

We first explain the energy barrier between $\boxplus$ and $\boxminus$.
\begin{itemize}
\item A sequence of configurations $(\omega_{t})_{t=0}^{T}=(\omega_{0},\,\omega_{1},\,\dots,\,\omega_{T})\subseteq\mathcal{X}$
for some $T\ge0$ is called a \textit{path} if $\omega_{t}\sim\omega_{t+1}$
for all $t\in\llbracket0,\,T-1\rrbracket$. A path $(\omega_{t})_{t=0}^{T}$
is a path connecting two configurations $\sigma$ and $\zeta$ in
$\mathcal{X}$ if $\omega_{0}=\sigma$ and $\omega_{T}=\zeta$ or
vice versa.
\item The \textsl{communication height} between two configurations $\sigma,\,\zeta\in\mathcal{X}$
is defined by
\[
\Phi(\sigma,\,\zeta)=\min_{(\omega_{t})_{t=0}^{T}}\,\max_{t\in\llbracket0,\,T\rrbracket}H(\omega_{t})\;,
\]
where the minimum is taken over all paths connecting $\sigma$ and
$\zeta$.
\item The \textit{energy barrier} between ground states is defined by
\[
\Gamma=\Gamma(K,\,L):=\Phi(\text{\ensuremath{\boxplus,\,\boxminus}})=\Phi(\boxminus,\,\text{\ensuremath{\boxplus}})\;,
\]
where the last equality holds from the symmetry of the model.
\end{itemize}
The following result has been verified in \cite{N-Z}. We note that
we have assumed $L\ge K\ge5$.
\begin{thm}
\label{t_energy barrier}The energy barrier is given by $\Gamma=2K+2$.
\end{thm}

The proof of this theorem is given in \cite{N-Z} based on combinatorial
arguments. We do not give the proof of this in the current note in
order to focus more on the role of potential theory in the analysis
of the current model.

\subsubsection*{Eyring--Kramers law}
\begin{notation}
In the current part, a collection $(a_{\beta}=a_{\beta}(K,\,L))_{\beta>0}$
of real numbers is written as $a_{\beta}=o_{\beta}(1)$ if $\lim_{\beta\rightarrow\infty}a_{\beta}=0$
for all $K$ and $L$.
\end{notation}

By Theorem \ref{t_energy barrier} and the large deviation principle,
one can deduce (cf. \cite{N-Z}) the following estimate of the mean
transition time $\mathbb{E}_{\boxplus}^{\beta}[\tau_{\boxminus}]$
and $\mathbb{E}_{\boxminus}^{\beta}[\tau_{\boxplus}]$:
\begin{equation}
\lim_{\beta\rightarrow\infty}\frac{1}{\beta}\log\mathbb{E}_{\boxplus}^{\beta}[\tau_{\boxminus}]=\lim_{\beta\rightarrow\infty}\frac{1}{\beta}\log\mathbb{E}_{\boxminus}^{\beta}[\tau_{\boxplus}]=\Gamma\;.\label{e_ldp}
\end{equation}
Note that $\tau_{\boxplus}$ and $\tau_{\boxminus}$ represent hitting
time of the set $\{\boxplus\}$ and $\{\boxminus\}$, respectively.

Along with the potential theory explained in the first part, we can
derive the precise sub-exponential prefactor of the previous large-deviation
estimate to get sharp asymptotics of the mean transition time.
\begin{thm}
\label{t_EK}There exists a constant $\kappa=\kappa(K,\,L)>0$ such
that
\begin{equation}
\mathbb{E}_{\boxplus}^{\beta}[\tau_{\boxminus}]=\mathbb{E}_{\boxminus}^{\beta}[\tau_{\boxplus}]=(1+o_{\beta}(1))\,\kappa e^{\Gamma\beta}\;.\label{eq:EK}
\end{equation}
Moreover, the constant $\kappa$ satisfies
\begin{align}
\lim_{K\rightarrow\infty}\kappa(K,\,L) & =\begin{cases}
1/4 & \text{if }K<L\;,\\
1/8 & \text{if }K=L\;.
\end{cases}\label{gam2}
\end{align}
\end{thm}

Precise asymptotics such as (\ref{eq:EK}) are called the \textit{Eyring--Kramers
law} (cf. \cite{Berg1} for more detail) for the Metropolis dynamics
$\sigma_{\beta}(\cdot)$. The constant $\kappa$ is explained more
precisely later. Although we have not provided the formula for the
constant $\kappa$ at this point, there exists a complicated but explicit
expression for this constant (cf. (\ref{e_def_kappa}), Proposition
\ref{p_e2est} and Remark \ref{r_e2est}).

This theorem is the main result for the current part. The proof is
divided into several stages. Firstly, in Section \ref{sec3}, we use
the potential theory to reduce the proof of Theorem \ref{t_EK} to
a capacity estimate. To estimate the capacity to a precise level,
we need a much more accurate understanding of the energy landscape
than that needed to derive (\ref{e_ldp}). This analysis of the energy
landscape is carried out in Sections \ref{sec5}-\ref{sec7}. Then,
the capacity estimate will be carried out in Sections \ref{sec9}
and \ref{sec10} based on the Dirichlet principle and the generalized
Thomson principle, respectively.
\begin{rem}
The followings are some comments on Theorem \ref{t_EK}.
\begin{enumerate}
\item If $K<L$, there is only one direction for the transition between
ground states, whereas if $K=L$, there are two possible directions.
This is the reason for the dependency in the asymptotics of $\kappa$
on the relation between $K$ and $L$.
\item The constant $\Gamma$ is model-independent, in the sense that it
will be the same for other Glauber dynamics. However, the constant
$\kappa$ is model-dependent. For other Glauber dynamics, this constant
may be different.
\end{enumerate}
\end{rem}

\section{\label{sec3}Application of Potential-Theoretic Approach}

The proof of Theorem \ref{t_EK} is based on the potential-theoretic
arguments developed in \cite{B-E-G-K} and accurate analyses of the
energy landscape. In this section, based on the argument developed
in \cite{B-E-G-K} along with the Dirichlet and the generalized Thomson
principle (cf. Theorem \ref{t_GTP}) for reversible Markov processes,
we reduce the proof of Theorem \ref{t_EK} to constructions of a test
function and a test flow in Propositions \ref{p_Capup} and \ref{p_Caplow},
respectively.

\subsection{\label{sec3.1}Main capacity estimate}

We first introduce the potential-theoretic notions. These notions
are introduced in Section \ref{sec_pot}, but we rename these objects
in the context of the Ising model.
\begin{itemize}
\item The \textit{Dirichlet form} $\mathscr{D}_{\beta}(\cdot)$ associated
with the reversible process $\sigma_{\beta}(\cdot)$ is given by,
for $f:\mathcal{X}\rightarrow\mathbb{R}$,
\begin{equation}
\mathscr{D}_{\beta}(f)=\frac{1}{2}\sum_{\sigma,\,\zeta\in\mathcal{X}}\mu_{\beta}(\sigma)\,c_{\beta}(\sigma,\,\zeta)\,\{f(\zeta)-f(\sigma)\}^{2}\;.\label{e_Diri}
\end{equation}
\item Let $\mathcal{P}$ and $\mathcal{Q}$ be disjoint and non-empty subsets
of $\mathcal{X}$. The \textit{equilibrium potential} between $\mathcal{P}$
and $\mathcal{Q}$ is the function $h_{\mathcal{P},\,\mathcal{Q}}^{\beta}:\mathcal{X}\rightarrow\mathbb{R}$
defined by
\begin{equation}
h_{\mathcal{P},\,\mathcal{Q}}^{\beta}(\sigma)=\mathbb{P}_{\sigma}^{\beta}\,[\tau_{\mathcal{P}}<\tau_{\mathcal{Q}}]\;,\label{e_eqp}
\end{equation}
and the \textit{capacity} between $\mathcal{P}$ and $\mathcal{Q}$
is defined by
\begin{equation}
\mathrm{cap}_{\beta}(\mathcal{P},\,\mathcal{Q})=\mathscr{D}_{\beta}(h_{\mathcal{P},\,\mathcal{Q}}^{\beta})\;.\label{e_Capdef}
\end{equation}
\end{itemize}
The following theorem is the main capacity estimate.
\begin{thm}
\label{t_Cap}We have that
\begin{equation}
\mathrm{cap}_{\beta}(\boxplus,\,\boxminus)=\frac{1+o_{\beta}(1)}{2\kappa}\,e^{-\Gamma\beta}\;,\label{e_Cap}
\end{equation}
where $\kappa$ is the constant appearing in Theorem \ref{t_EK}.
\end{thm}

Before proceeding to the proof of Theorem \ref{t_Cap}, we first explain
the proof of Theorem \ref{t_EK} by assuming Theorem \ref{t_Cap}.
\begin{proof}[Proof of Theorem \ref{t_EK}]
Since $\mathbb{E}_{\boxplus}^{\beta}\,[\tau_{\boxminus}]=\mathbb{E}_{\boxminus}^{\beta}\,[\tau_{\boxplus}]$
by symmetry, we only focus on the estimate of $\mathbb{E}_{\boxplus}^{\beta}\,[\tau_{\boxminus}]$.
By Proposition \ref{p_mht}, (or more precisely, by (\ref{exphit2})),
we have
\begin{equation}
\mathbb{E}_{\boxplus}^{\beta}\,[\tau_{\boxminus}]=\frac{1}{\mathrm{cap}_{\beta}(\boxplus,\,\boxminus)}\sum_{\sigma\in\mathcal{X}}\mu_{\beta}(\sigma)\,h_{\boxplus,\,\boxminus}(\sigma)\;.\label{e_mjr}
\end{equation}
By Proposition \ref{p_inv} and the fact that $h_{\boxplus,\,\boxminus}(\boxplus)=1$
and $h_{\boxplus,\,\boxminus}(\boxminus)\equiv0$, we rewrite the
last summation as
\[
\frac{1}{2}+o_{\beta}(1)+\sum_{\sigma\in\mathcal{X}\setminus\mathcal{S}}\mu_{\beta}(\sigma)\,h_{\boxplus,\,\boxminus}(\sigma)\;.
\]
Since $|h_{\boxplus,\,\boxminus}|\le1$, again by Proposition \ref{p_inv},
we have
\[
\Big|\,\sum_{\sigma\in\mathcal{X}\setminus\mathcal{S}}\mu_{\beta}(\sigma)\,h_{\boxplus,\,\boxminus}(\sigma)\,\Big|\le\mu_{\beta}(\mathcal{X}\setminus\mathcal{S})=o_{\beta}(1)\;.
\]
In summary, we obtain
\[
\sum_{\sigma\in\mathcal{X}}\mu_{\beta}(\sigma)\,h_{\boxplus,\,\boxminus}(\sigma)=\frac{1}{2}+o_{\beta}(1)\;.
\]
Now, inserting this and Theorem \ref{t_Cap} to (\ref{e_mjr}), we
can complete the proof.
\end{proof}

\subsection{\label{sec4.1}The constant $\kappa$}

To explain the main result for the capacity estimate, we first have
to introduce the bulk constant $\mathfrak{b}$ and the edge constant
$\mathfrak{e}$. The reason for the choice of the words ``bulk''
and ``edge'' will become clear as we analyze the energy landscape
more deeply (cf. Remark \ref{r_testf2}).

Firstly, the bulk constant $\mathfrak{b}$ is defined explicitly as
\begin{equation}
\text{\ensuremath{\mathfrak{b}}}=\begin{cases}
\frac{(K+2)(L-4)}{4KL} & \text{if }K<L\;,\\
\frac{(K+2)(L-4)}{8KL} & \text{if }K=L\;.
\end{cases}\label{e_b2}
\end{equation}
On the other hand, we do not provide a precise definition of the edge
constant $\mathfrak{e}$ at this point. This is a complicated constant
defined in (\ref{e_e2def}) which satisfies (cf. Proposition \ref{p_e2est})
\begin{equation}
0<\mathfrak{e}\le\frac{1}{L}\;.\label{e_e2}
\end{equation}
We stress that these constants depend on $K$ and $L$ even though
the dependency is not highlighted in the notation.

Now, we define the constant $\kappa$ as
\begin{equation}
\kappa=\mathfrak{b}+2\mathfrak{e}\;.\label{e_def_kappa}
\end{equation}
We note that the bulk constant $\mathfrak{b}$ is the constant associated
to the bulk part of the transition between $\boxplus$ and $\boxminus$,
while the edge constant $\mathfrak{e}$ is related to the edge behavior
of the transition. Since there are two edge parts (around $\boxplus$
and around $\boxminus$), the constant $2$ has been multiplied in
front of $\mathfrak{e}$ in (\ref{e_def_kappa}). Moreover, one can
readily observe that, when $K$ (and hence $L$) is large, the edge
constant $\mathfrak{e}$ is much smaller than $\mathfrak{b}$. Hence,
the bulk effect dominates the edge effect. We also note that (\ref{gam2})
follows directly from (\ref{e_b2}) and (\ref{e_e2}).

\subsection{\label{sec4.2}Capacity estimate}

The upper bound estimate is based on the Dirichlet principle for reversible
Markov processes (Theorem \ref{t_DP_rev}). To use this principle,
we will prove the following proposition.
\begin{prop}
\label{p_Capup}There exists a function $f_{0}:\mathcal{X}\rightarrow\mathbb{R}$
such that $f_{0}\in\mathfrak{C}_{1,\,0}(\{\boxplus\},\,\{\boxminus\})$
and that
\begin{equation}
\mathscr{D}_{\beta}(f_{0})=\frac{1+o_{\beta}(1)}{2\kappa}\,e^{-\Gamma\beta}\;.\label{e_Capup}
\end{equation}
\end{prop}

Finding the test function $f_{0}$ requires a deep insight into the
energy landscape, as well as the typical patterns of the Metropolis
dynamics in a suitable neighborhood of saddle configurations. We construct
this test function and prove Proposition \ref{p_Capup} in Section
\ref{sec9}.

To explain the lower bound of the capacity, we use the generalized
Thomson principle (Theorem \ref{t_GTP}). For convenience, we write
the flow norm associated with the process $\sigma_{\beta}(\cdot)$
as $\|\cdot\|_{\beta}$. We shall prove the following proposition
later to establish the lower bound of the capacity.
\begin{prop}
\label{p_Caplow}There exists a flow $\psi_{0}$ such that
\begin{equation}
\|\psi_{0}\|_{\beta}^{2}=(2+o_{\beta}(1))\,\kappa\,e^{\Gamma\beta}\text{\;\;\;\; and\;\;\;\;}\sum_{\sigma\in\mathcal{X}}h_{\boxplus,\,\boxminus}^{\beta}(\sigma)\,(\mathrm{div}\,\psi_{0})(\sigma)=1+o_{\beta}(1)\;.\label{e_Caplow1}
\end{equation}
\end{prop}

We construct the test flow $\psi_{0}$ in Section \ref{sec10} (cf.
Definition \ref{d_testfl2}), and then verify in the same section
that our test flow $\psi_{0}$ indeed satisfies (\ref{e_Caplow1}).

We now prove Theorem \ref{t_Cap} by assuming Propositions \ref{p_Capup}
and \ref{p_Caplow}.
\begin{proof}[Proof of Theorem \ref{t_Cap}]
By Theorem \ref{t_DP_rev} and Proposition \ref{p_Capup}, we get
\begin{equation}
\mathrm{cap}_{\beta}(\boxplus,\,\boxminus)\le\mathscr{D}_{\beta}(f_{0})=\frac{1+o_{\beta}(1)}{2\kappa}\,e^{-\Gamma\beta}\;.\label{e_up}
\end{equation}
On the other hand, by Theorem \ref{t_GTP} and Proposition \ref{p_Caplow},
we obtain
\begin{equation}
\mathrm{cap}_{\beta}(\boxplus,\,\boxminus)\ge\frac{1}{\|\psi_{0}\|_{\beta}^{2}}\,\Big[\,\sum_{\sigma\in\mathcal{X}}h_{\boxplus,\,\boxminus}^{\beta}(\sigma)\,(\mathrm{div}\,\psi_{0})(\sigma)\,\Big]^{2}=\frac{1+o_{\beta}(1)}{2\kappa}\,e^{-\Gamma\beta}\;.\label{e_low}
\end{equation}
The proof is completed by (\ref{e_up}) and (\ref{e_low}).
\end{proof}
Hence, to prove Theorem \ref{t_EK}, it only remains to prove Propositions
\ref{p_Capup} and \ref{p_Caplow}. The proof is given in the remainder
of Part \ref{pt2}.

\section{\label{sec5}Neighborhood of Configurations}

For $c\in\mathbb{R}$, a path $(\omega_{t})_{t=0}^{T}$ in $\mathcal{X}$
is called a \textit{$c$-path} if we have $H(\omega_{t})\le c$ for
all $t\in\llbracket0,\,T\rrbracket$. Heuristically, if two configurations
are connected by a $(\Gamma-1)$-path, in a suitable sense, these
two configurations are indistinguishable in the transition scale $e^{\beta\Gamma}$,
since $\sigma_{\beta}(\cdot)$ commutes them in a shorter scale. Moreover,
if two configurations are not connected by a $\Gamma$-path, the process
$\sigma_{\beta}(\cdot)$ cannot commute these two configurations in
the transition scale $e^{\beta\Gamma}$. The following definition
of neighborhoods is inspired from these observations.
\begin{defn}[Neighborhood of configurations]
\label{d_nbd}
\begin{enumerate}
\item For $\sigma\in\mathcal{X}$, the \textit{neighborhood} $\mathcal{N}(\sigma)$
and the \textit{extended neighborhood} $\widehat{\mathcal{N}}(\sigma)$
are defined as
\begin{align*}
\mathcal{N}(\sigma) & =\{\zeta\in\mathcal{X}:\exists\,\text{a }(\Gamma-1)\text{-path }(\omega_{t})_{t=0}^{T}\text{ connecting }\sigma\text{ and }\zeta\}\;\text{and}\\
\widehat{\mathcal{N}}(\sigma) & =\{\zeta\in\mathcal{X}:\exists\,\text{a }\Gamma\text{-path }(\omega_{t})_{t=0}^{T}\text{ connecting }\sigma\text{ and }\zeta\}\;.
\end{align*}
If $H(\sigma)>\Gamma-1$ (resp. $H(\sigma)>\Gamma$), we set $\mathcal{N}(\sigma)=\emptyset$
(resp. $\widehat{\mathcal{N}}(\sigma)=\emptyset$).
\item For $\mathcal{P}\subseteq\mathcal{X}$, we define
\[
\mathcal{N}(\mathcal{P})=\bigcup_{\sigma\in\mathcal{P}}\mathcal{N}(\sigma)\;\;\;\;\text{and\;\;\;\;}\widehat{\mathcal{N}}(\mathcal{P})=\bigcup_{\sigma\in\mathcal{P}}\mathcal{\widehat{\mathcal{N}}}(\sigma)\;.
\]
\item A path $(\omega_{t})_{t=0}^{T}$ is said to be a path in $\mathcal{A}\subset\mathcal{X}$
if $\omega_{t}\in\mathcal{A}$ for all $t\in\llbracket0,\,T\rrbracket.$
For $\mathcal{Q}\subset\mathcal{X}$ and $\sigma\in\mathcal{X\setminus\mathcal{Q}}$,
we define
\[
\widehat{\mathcal{N}}(\sigma\,;\,\mathcal{Q})=\{\zeta\in\mathcal{X}:\exists\,\text{a }\Gamma\text{-path in }\mathcal{X}\setminus\mathcal{Q}\text{ connecting }\sigma\text{ and }\zeta\}\;.
\]
If $H(\sigma)>\Gamma$, we set $\widehat{\mathcal{N}}(\sigma\,;\,\mathcal{Q})=\emptyset$.
\item For $\mathcal{P}\subseteq\mathcal{X}$ disjoint with $\mathcal{Q}$,
define
\[
\widehat{\mathcal{N}}(\mathcal{P}\,;\,\mathcal{Q})=\bigcup_{\sigma\in\mathcal{P}}\widehat{\mathcal{N}}(\sigma\,;\,\mathcal{Q})\;.
\]
\end{enumerate}
\end{defn}

With this notation, Theorem \ref{t_energy barrier} is equivalent
to $\mathcal{N}(\boxplus)\cap\mathcal{N}(\boxminus)=\emptyset$ and
$\widehat{\mathcal{N}}(\boxplus)=\widehat{\mathcal{N}}(\boxminus)$.
Since the transition must take place in the set $\widehat{\mathcal{N}}(\mathcal{S})$,
analyzing the structure of this set is crucial in the energy landscape
analysis. It will be carried out in Section \ref{sec7}.

The following lemma is useful.
\begin{lem}
\label{l_set}Suppose that $\mathcal{P}$ and $\mathcal{Q}$ are disjoint
subsets of $\mathcal{X}$. Then, it holds that
\[
\widehat{\mathcal{N}}(\mathcal{P}\cup\mathcal{Q})=\widehat{\mathcal{N}}(\mathcal{Q}\,;\,\mathcal{P})\cup\widehat{\mathcal{N}}(\mathcal{P}\,;\,\mathcal{Q})\;.
\]
\end{lem}

\begin{proof}
Since
\begin{equation}
\widehat{\mathcal{N}}(\mathcal{P}\cup\mathcal{Q})=\widehat{\mathcal{N}}(\mathcal{P})\cup\widehat{\mathcal{N}}(\mathcal{Q})\;,\label{e_set1}
\end{equation}
\[
\widehat{\mathcal{N}}(\mathcal{Q})\supset\widehat{\mathcal{N}}(\mathcal{Q}\,;\,\mathcal{P})\;,\;\text{and}\;\widehat{\mathcal{N}}(\mathcal{P})\supset\widehat{\mathcal{N}}(\mathcal{P}\,;\,\mathcal{Q})\;,
\]
it immediately follows that
\begin{equation}
\widehat{\mathcal{N}}(\mathcal{P}\cup\mathcal{Q})\supseteq\widehat{\mathcal{N}}(\mathcal{Q}\,;\,\mathcal{P})\cup\widehat{\mathcal{N}}(\mathcal{P}\,;\,\mathcal{Q})\;.\label{e_set2}
\end{equation}

Let us now prove the reversed inclusion. We now assume that there
exists $\sigma\in\mathcal{X}$ such that
\begin{align}
\sigma & \in\widehat{\mathcal{N}}(\mathcal{P}\cup\mathcal{Q})\setminus\big[\,\widehat{\mathcal{N}}(\mathcal{Q}\,;\,\mathcal{P})\cup\widehat{\mathcal{N}}(\mathcal{P}\,;\,\mathcal{Q})\,\big]\;.\label{e_set3}
\end{align}
By (\ref{e_set1}), we may assume without loss of generality that
\[
\sigma\in\widehat{\mathcal{N}}(\mathcal{P})\setminus\big[\,\widehat{\mathcal{N}}(\mathcal{Q}\,;\,\mathcal{P})\cup\widehat{\mathcal{N}}(\mathcal{P}\,;\,\mathcal{Q})\,\big]\;.
\]
Since $\sigma\in\widehat{\mathcal{N}}(\mathcal{P})\setminus\widehat{\mathcal{N}}(\mathcal{P}\,;\,\mathcal{Q})$,
we have $\sigma\notin\mathcal{P}$. Since $\sigma\in\widehat{\mathcal{N}}(\mathcal{P})$,
we can find a $\Gamma$-path connecting $\sigma$ and $\mathcal{P}$.
Let us assume that $(\omega_{t})_{t=0}^{T}$ is the shortest of all
such paths. We may assume that $\omega_{0}=\sigma$ and $\omega_{T}\in\mathcal{P}$.
\begin{itemize}
\item Suppose first that $\omega_{t}\notin\mathcal{Q}$ for all $t\in\llbracket0,\,T-1\rrbracket$.
Then the path $(\omega_{t})_{t=0}^{T}$ becomes a $\Gamma$-path in
$\mathcal{X}\setminus\mathcal{Q}$ connecting $\sigma$ and $\mathcal{P}$.
This contradicts the fact that $\sigma\notin\widehat{\mathcal{N}}(\mathcal{P}\,;\,\mathcal{Q})$.
\item Suppose next that $\omega_{t_{0}}\in\mathcal{Q}$ for some $t_{0}\in\llbracket0,\,T-1\rrbracket$.
Then, by the minimality assumption on the length of $(\omega_{t})_{t=0}^{T}$,
we must have $\omega_{t}\notin\mathcal{P}$ for all $t\in\llbracket0,\,T-1\rrbracket$.
Consequently, $(\omega_{t})_{t=0}^{t_{0}}$ becomes a path in $\mathcal{X}\setminus\mathcal{P}$
connecting $\mathcal{\sigma}$ and $\mathcal{Q}$, and hence we get
a contradiction to the fact $\sigma\notin\widehat{\mathcal{N}}(\mathcal{Q}\,;\,\mathcal{P})$.
\end{itemize}
Therefore, there is no $\sigma$ satisfying (\ref{e_set3}), and we
have proved the reversed inclusion relation of (\ref{e_set2}).
\end{proof}

\section{\label{sec6-1}Canonical Configurations and Paths}

Now, we begin to analyze the energy landscape. In this section, we
introduce the canonical configurations and paths, and then investigate
their properties. Based on these, we study the typical configurations
in the next section.

\subsection{\label{sec6.1}Canonical configurations }

\begin{figure}
\includegraphics[width=12cm]{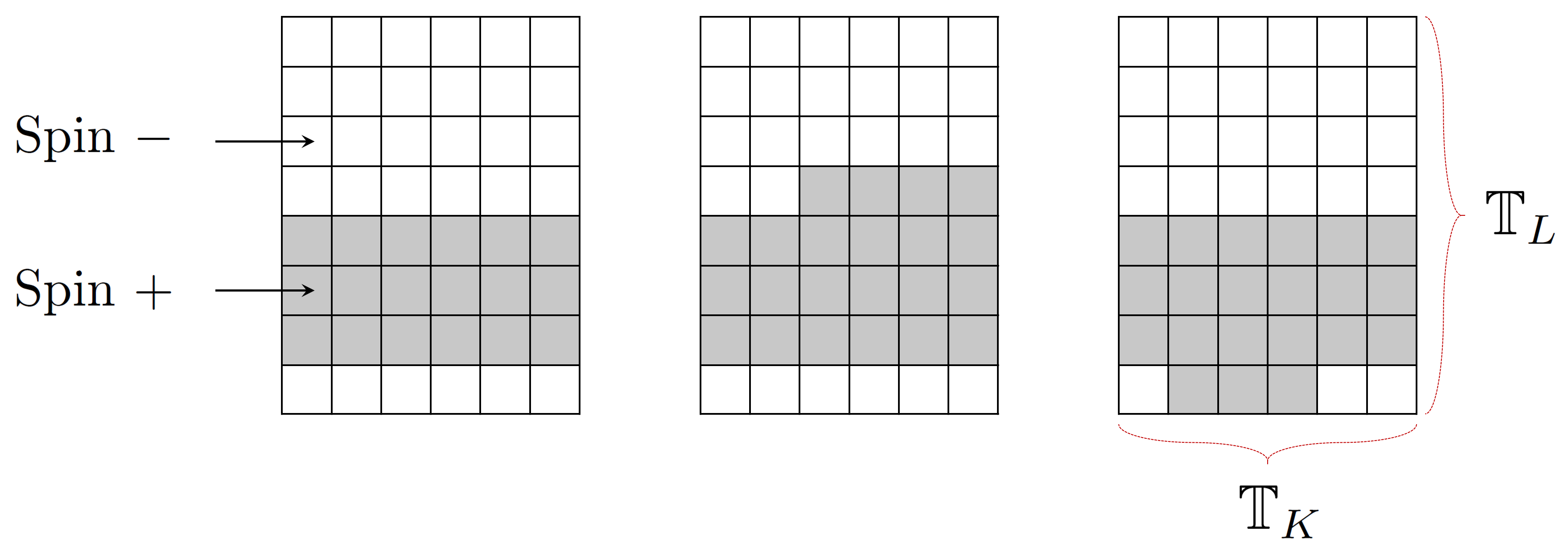}\caption{\label{fig6.1}\textbf{Canonical configurations for $(K,\,L)=(6,\,8)$.}
White and gray boxes correspond to a box with $-$ spin and $+$ spin,
respectively. Three figures represent configurations $\zeta_{2,\,3}$,
$\zeta_{2,\,3\,;\,3,\,4}^{\textup{up}}$, and $\zeta_{2,\,3\,;\,2,\,3}^{\textup{down}}$,
respectively. }
\end{figure}

\begin{defn}[Canonical configurations]
\label{d_precanreg}We refer to Figure \ref{fig6.1} for an illustration
of examples of the canonical configurations defined below. Before
defining complicated notations, we note that $k$ and $\ell$ are
used to represent elements of $\mathbb{T}_{K}$ and $\mathbb{T}_{L}$,
respectively, and $v$ and $h$ are used to denote vertical and horizontal
lengths, respectively.
\begin{itemize}
\item For $\ell\in\mathbb{T}_{L}$ and $v\in\llbracket0,\,L\rrbracket$,
denote by $\zeta_{\ell,\,v}\in\mathcal{X}$ the configuration whose
spins are $+$ on
\[
\mathbb{T}_{K}\times\{\ell+n\in\mathbb{T}_{L}:n\in\llbracket0,\,v-1\rrbracket\subseteq\mathbb{Z}\}\;.
\]
and $-$ on the remainder. Hence, we have $\zeta_{\ell,\,0}=\boxminus$
and $\zeta_{\ell,\,L}=\boxplus$ for all $\ell\in\mathbb{T}_{L}$.
For $v\in\llbracket0,\,L\rrbracket$, write
\begin{equation}
\mathcal{R}_{v}=\{\zeta_{\ell,\,v}:\ell\in\mathbb{T}_{L}\}\;.\label{e_Rv}
\end{equation}
\item For $(\ell,\,k)\in\mathbb{T}_{L}\times\mathbb{T}_{K}$ and $(v,\,h)\in\llbracket0,\,L-1\rrbracket\times\llbracket0,\,K\rrbracket$,
denote by $\zeta_{\ell,\,v\,;\,k,\,h}^{\textup{up}}\in\mathcal{X}$
the configuration whose spins are $+$ on
\[
\{x\in\Lambda:\zeta_{\ell,\,v}(x)=+\}\cup\big[\,\{k+n\in\mathbb{T}_{K}:n\in\llbracket0,\,h-1\rrbracket\subseteq\mathbb{Z}\}\times\{\ell+v\}\,\big]
\]
and $-$ on the remainder. Similarly, denote by $\zeta_{\ell,\,v\,;\,k,\,h}^{\textup{down}}\in\mathcal{X}$
whose spins are $+$ on
\[
\{x\in\Lambda:\zeta_{\ell,\,v}(x)=+\}\cup\big[\,\{k+n\in\mathbb{T}_{K}:n\in\llbracket0,\,h-1\rrbracket\subseteq\mathbb{Z}\}\times\{\ell-1\}\,\big]
\]
and $-$ on the remainder. Namely, the configuration $\zeta_{\ell,\,v\,;\,k,\,h}^{\textup{up}}$
(resp. $\zeta_{\ell,\,v\,;\,k,\,h}^{\textup{down}}$) is obtained
from $\zeta_{\ell,\,v}$ by attaching a protuberance of spin $+$
of size $h$ at the upper (resp. lower) side of the cluster of spin
$+$ of $\zeta_{\ell,\,v}$.
\item For $v\in\llbracket0,\,L-1\rrbracket$, define\textbf{
\begin{equation}
\mathcal{Q}_{v}=\bigcup_{k\in\mathbb{T}_{K}}\bigcup_{h=1}^{K-1}\{\zeta_{\ell,\,v\,;\,k,\,h}^{\textup{up}},\,\zeta_{\ell,\,v\,;\,k,\,h}^{\textup{down}}\}\;.\label{e_Qv}
\end{equation}
}Hence, $\mathcal{Q}_{v}$ consists of configurations between $\mathcal{R}_{v}$
and $\mathcal{R}_{v+1}$.
\item Finally, define
\[
\mathcal{C}=\bigcup_{v=0}^{L}\mathcal{R}_{v}\cup\bigcup_{v=0}^{L-1}\mathcal{Q}_{v}\;.
\]
In the current note, the \textit{canonical configurations} are the
configurations belonging to $\mathcal{C}$.
\end{itemize}
\end{defn}

\begin{rem}
\label{r_H2}By a direct computation, we can readily verify that $H(\sigma)\le\Gamma$
for all $\sigma\in\mathcal{C}$. In particular, we have
\[
H(\sigma)=\begin{cases}
\Gamma-2 & \text{if }\sigma\in\mathcal{R}_{v}\text{ for some }v\in\llbracket1,\,L-1\rrbracket\;,\\
\Gamma & \text{if }\sigma\in\mathcal{Q}_{v}\text{ for some }v\in\llbracket1,\,L-2\rrbracket\;.
\end{cases}
\]
\end{rem}

For the clarity of the discussion, we henceforth assume that $K<L$.
The case $K=L$ will be discussed in Section \ref{sec_K=00003DL}.
Note that the only difference for the case $K=L$ is that the configuration
obtained by rotating a canonical configuration in $\mathcal{C}$ must
play the same role, unlike the case $K<L$. This fact can be readily
taken into account in the computations, and we refer to Section \ref{sec_K=00003DL}
or \cite{Kim-Seo Potts} for further details. Note that for $K<L$,
the rows and columns play completely different roles.

\subsection{\label{sec6.2}Canonical paths}

We now explain the crucial role of canonical configurations by describing
canonical paths between $\boxminus$ and $\boxplus$ consisting of
canonical configurations. The following notation is useful.
\begin{notation}
\label{n_frakG}Suppose that $N\ge2$ is a positive integer.
\begin{itemize}
\item Denote by $\mathfrak{S}_{N}$ the collection of all connected subsets
of $\mathbb{T}_{N}$, i.e.,
\begin{equation}
\mathfrak{S}_{N}=\{P\subseteq\mathbb{T}_{N}:P=\llbracket i,\,j\rrbracket\text{ for some }i,\,j\in\mathbb{T}_{N}\text{ or }P=\emptyset\}\;.\label{e_frakG}
\end{equation}
Here, the set $\llbracket i,\,j\rrbracket\subseteq\mathbb{T}_{N}$
represents the set $\{i,\,i+1,\,\dots,\,j\}$. Note that this set
can be defined even for $j<i$. For instance, for $N=6$, the set
$\llbracket4,\,2\rrbracket$ represents $\{4,\,5,\,6,\,1,\,2\}$.
\item For two sets $P,\,P'\in\mathfrak{S}_{N}$, we write $P\prec P'$ if
$P\subseteq P'$ and $|P'|=|P|+1$.
\item A sequence $(P_{n})_{n=0}^{N}$ of sets in $\mathfrak{S}_{N}$ is
called an \textit{increasing sequence} if
\[
\emptyset=P_{0}\prec P_{1}\prec\cdots\prec P_{N}=\mathbb{T}_{N}\;.
\]
Note that, for an increasing sequence $(P_{n})_{n=0}^{N}$ in $\mathfrak{S}_{N}$,
we have that $|P_{n}|=n$ for all $n\in\llbracket0,\,N\rrbracket$.
\end{itemize}
\end{notation}

\begin{figure}
\includegraphics[width=12cm]{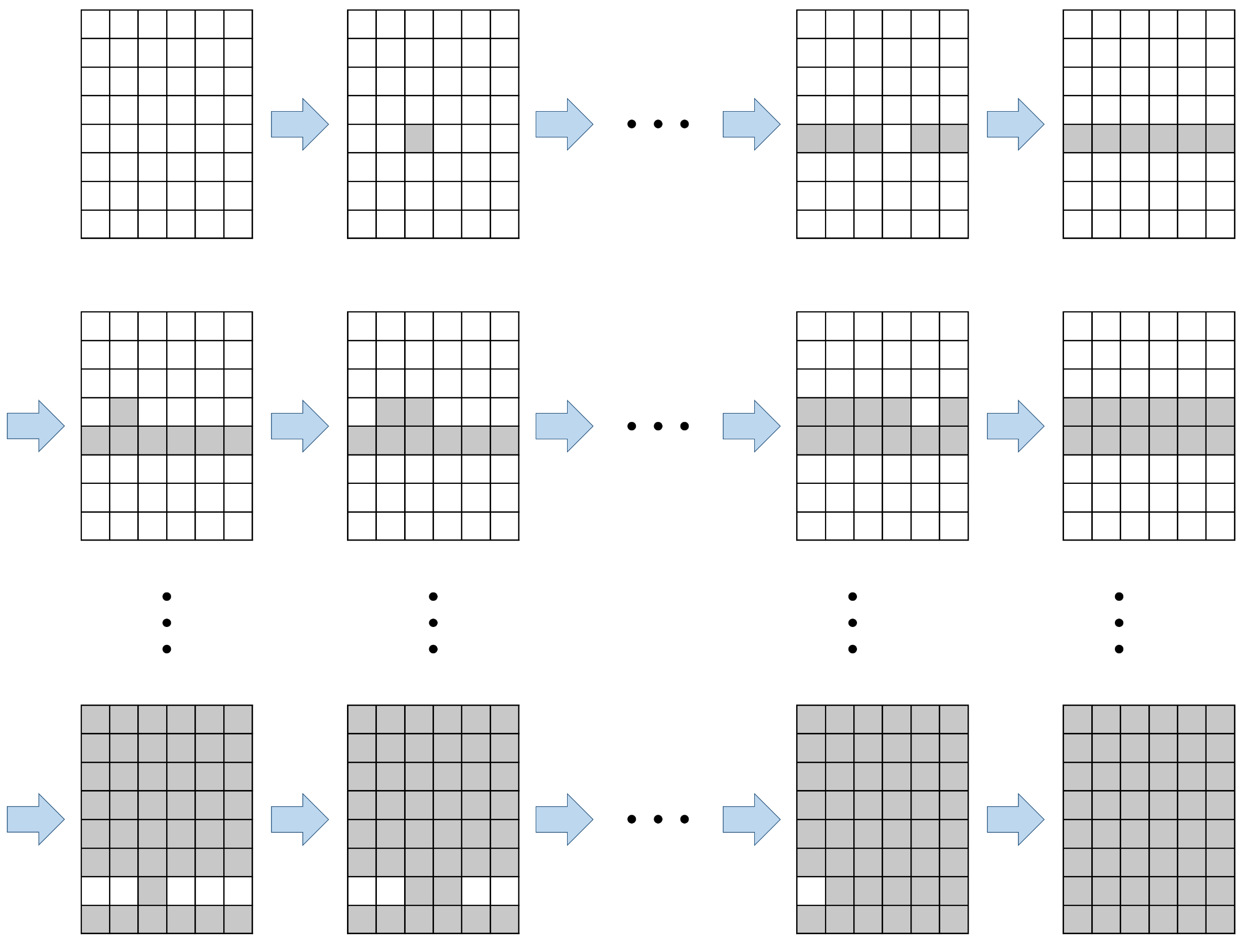}\caption{\label{fig6.2}Example of a canonical path for $(K,\,L)=(6,\,8)$.}
\end{figure}

\begin{defn}[Canonical paths]
\label{d_canpath2}We refer to Figure \ref{fig6.2} for an example
of canonical path defined below.
\begin{enumerate}
\item We first introduce a \textit{standard sequence} of subsets of $\Lambda=\mathbb{T}_{K}\times\mathbb{T}_{L}$
connecting the empty set and the full set $\Lambda$.
\begin{enumerate}
\item For $P,\,P'\in\mathfrak{S}_{L}$ with $P\prec P'$, a sequence $(A_{t})_{t=0}^{K}$
of subsets of $\Lambda$ is called a standard sequence connecting
$\mathbb{T}_{K}\times P$ and $\mathbb{T}_{K}\times P'$ if there
exists an increasing sequence $(Q_{t})_{t=0}^{K}$ in $\mathfrak{S}_{K}$
such that
\[
A_{t}=(\mathbb{T}_{K}\times P)\cup\big[\,Q_{t}\times(P'\setminus P)\,\big]\;\;\;\text{for all}\;t\in\llbracket0,\,K\rrbracket\;.
\]
\item A sequence $(A_{t})_{t=0}^{KL}$ of subsets of $\Lambda$ is called
a standard sequence connecting $\text{\ensuremath{\emptyset}}$ and
$\Lambda$ if there exists an increasing sequence $(P_{\ell})_{\ell=0}^{L}$
in $\mathfrak{S}_{L}$ such that $A_{K\ell}=\mathbb{T}_{K}\times P_{\ell}$
for all $\ell\in\llbracket0,\,L\rrbracket$, and the sub-sequence
$(A_{t})_{t=K\ell}^{K(\ell+1)}$ is a standard sequence connecting
$\mathbb{T}_{K}\times P_{\ell}$ and $\mathbb{T}_{K}\times P_{\ell+1}$
for all $\ell\in\llbracket0,\,L-1\rrbracket$.
\end{enumerate}
\item A path $(\omega_{t})_{t=0}^{KL}$ in $\mathcal{X}$ is called a \textit{canonical
path} connecting $\boxminus$ and $\boxplus$ if there exists a standard
sequence $(A_{t})_{t=0}^{KL}$ connecting $\text{\ensuremath{\emptyset}}$
and $\Lambda$ such that
\[
\omega_{t}(i,\,j)=\begin{cases}
- & \text{if }(i,\,j)\notin A_{t}\;,\\
+ & \text{if }(i,\,j)\in A_{t}\;.
\end{cases}
\]
It is easy to verify that $\omega_{0}=\boxminus$ and $\omega_{KL}=\boxplus$.
A canonical path connecting $\boxplus$ and $\boxminus$ is defined
in a similar manner. We say that a path is a canonical path if it
is a canonical path connecting either $\boxminus$ and $\boxplus$
or $\boxplus$ and $\boxminus$.
\end{enumerate}
\end{defn}

The following is an immediate consequence of the construction.
\begin{lem}
\label{l_canpath2}A canonical path consists only of canonical configurations.
In particular, for any canonical path $(\omega_{t})_{t=0}^{KL}$ connecting
$\boxminus$ and $\boxplus$, we have that
\[
\max_{t\in\llbracket0,\,KL\rrbracket}H(\omega_{t})=\Gamma\;.
\]
\end{lem}

\begin{proof}
The first assertion follows immediate from the construction. For the
second assertion, it suffices to recall Remark \ref{r_H2}.
\end{proof}
In view of the previous lemma and Theorem \ref{t_energy barrier},
a canonical path between $\boxminus$ and $\boxplus$ is an optimal
path achieving the communication height between them. We emphasize
here that the optimal transition may not always occur along this path.
Indeed, transitions from $\boxminus$ to $\mathcal{R}_{2}$ and from
$\mathcal{R}_{L-2}$ to $\boxplus$ may happen in a more complex manner,
while transitions from $\mathcal{R}_{2}$ to $\mathcal{R}_{L-2}$
should happen along a canonical path. This issue is the main topic
of the next section.

\section{\label{sec7}Typical Configurations}

The crucial notion in the energy landscape analysis between ground
states is the typical configurations defined in this section. A configuration
$\sigma$ is said to be a typical configuration if $\sigma\in\widehat{\mathcal{N}}(\mathcal{S})$.
Therefore, the typical configurations comprise all the relevant configurations
in the study of metastable transition between $\boxplus$ and $\boxminus$.

\subsection{\label{sec7.1}Typical configurations}

Let us start by defining typical configurations.

\begin{figure}
\includegraphics[width=14cm]{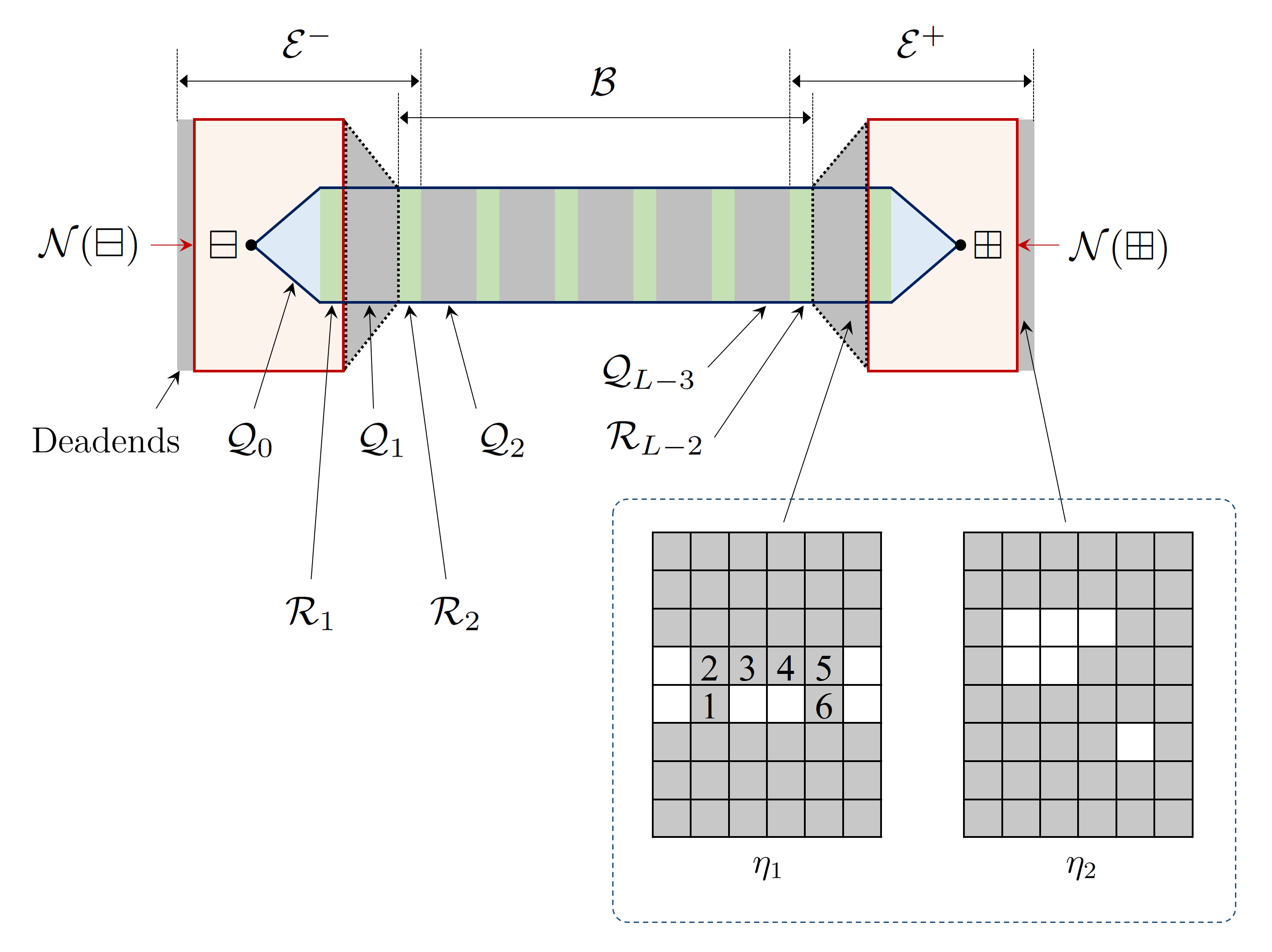}\caption{\label{fig7.1}\textbf{Structure of $\widehat{\mathcal{N}}(\mathcal{S})$
and typical configurations.} Regions consisting of configurations
with energy $\Gamma$ are colored gray. As we will verify in Proposition
\ref{p_typ2}, we can observe that
\begin{align*}
 & \widehat{\mathcal{N}}(\mathcal{S})=\mathcal{E}^{-}\cup\mathcal{E}^{+}\cup\mathcal{B}\;\;\;\;\text{and}\protect\\
 & \mathcal{E}^{-}\cap\mathcal{B}=\mathcal{R}_{2}\;\;\;\;\text{and\;\;\;\;}\mathcal{E}^{+}\cap\mathcal{B}=\mathcal{R}_{L-2}\;.
\end{align*}
The hexagonal region enclosed by the blue line denotes the set $\mathcal{C}$
of canonical configurations between $\boxminus$ and $\boxplus$.
The set $\mathcal{E}^{+}$ of edge typical configurations around $\boxplus$
consists of four regions. The first one is the neighborhood $\mathcal{N}(\boxplus)$
denoted by the red-enclosed box and the second one is $\mathcal{R}_{L-2}$.
The third one is the region consisting of configurations with energy
$\Gamma$ which are connected to $\mathcal{R}_{L-2}$ via a $\Gamma$-path
in $\mathcal{X}\setminus\mathcal{N}(\boxplus)$. An example of a configuration
belonging to this region is $\eta_{1}$. In particular, the configuration
$\eta_{1}$ is connected with a configuration in $\mathcal{R}_{L-2}$
via a $\Gamma$-path in $\mathcal{X}\setminus\mathcal{N}(\boxplus)$
which is obtained by updating six grey boxes by the order indicated
in the figure. The last region is the collection of the dead-ends
attached to $\mathcal{N}(\boxplus)$. This is a collection of configurations
with energy $\Gamma$ which are not connected to $\mathcal{R}_{L-2}$
via a $\Gamma$-path in $\mathcal{X}\setminus\mathcal{N}(\boxplus)$.
An example of a dead-end configuration is $\eta_{2}$ which has energy
$14=2\cdot6+2=2K+2$. A similar decomposition holds for $\mathcal{E}^{-}$. }
\end{figure}

\begin{defn}[Typical configurations]
\label{d_typ2}We refer to Figure \ref{fig7.1} for an illustration
of the typical configurations defined below.
\begin{itemize}
\item Define
\begin{equation}
\mathcal{B}=\bigcup_{v=2}^{L-2}\mathcal{R}_{v}\cup\bigcup_{v=2}^{L-3}\mathcal{Q}_{v}\;.\label{e_Bab2}
\end{equation}
A configuration belonging to $\mathcal{B}$ is called a \textit{bulk
typical configuration. }Then, write
\[
\mathcal{B}_{\Gamma}=\bigcup_{v=2}^{L-3}\mathcal{Q}_{v}=\{\sigma\in\mathcal{B}:H(\sigma)=\Gamma\}\;.
\]
\item Define
\begin{equation}
\mathcal{E}^{+}=\widehat{\mathcal{N}}(\boxplus\,;\,\mathcal{B}_{\Gamma})\;\;\;\text{and\;\;\;}\mathcal{E}^{-}=\widehat{\mathcal{N}}(\boxminus\,;\,\mathcal{B}_{\Gamma})\;.\label{e_EAB2}
\end{equation}
Then, we define $\mathcal{E}=\mathcal{E}^{+}\cup\mathcal{E}^{-}$.
A configuration belonging to $\mathcal{E}$ is called an \textit{edge
typical configuration. }
\end{itemize}
A configuration belonging to $\mathcal{B}\cup\mathcal{E}$ is called
a \textit{typical configuration}. Indeed, it holds that $\mathcal{B}\cup\mathcal{E}=\widehat{\mathcal{N}}(\mathcal{S})$,
and this will be verified later.
\end{defn}

Now, we explain the reason why we have decomposed typical configurations
into bulk and edge configurations. A typical transition from $\boxminus$
to $\boxplus$ of the Metropolis dynamics can be divided into three
stages. Firstly, the process passes through $\mathcal{E}^{-}$ to
arrive at $\mathcal{B}$. Then, it goes through $\mathcal{B}$ along
the canonical configurations to arrive at $\mathcal{E}^{+}$. Finally,
the process reaches at $\boxplus$ by passing through $\mathcal{E}^{+}$.
The behavior of the process at the second stage (i.e., in the bulk)
is relatively clear, and we can understand the behavior in great detail.
On the other hand, the behavior of the Metropolis dynamics on $\mathcal{E}^{-}$
and $\mathcal{E}^{+}$ is complex, and can be explained in terms of
an auxiliary Markov chain defined in Definition \ref{d_EAMc2}. We
are not able to write the constant appearing in the Eyring--Kramers
law in a simple manner because of this complex behavior of the Metropolis
dynamics in the edge typical configurations.

\subsection{\label{sec6.3}Characterization of configurations with low energy}

To investigate the typical configurations defined above, in this subsection,
we fully characterize the configurations which have energy less than
$\Gamma$. Write
\begin{equation}
\|\sigma\|_{+}=\sum_{x\in\Lambda}\mathbf{1}\{\sigma(x)=+\}\text{\;\;\;}\text{and}\;\;\;\|\sigma\|_{-}=\sum_{x\in\Lambda}\mathbf{1}\{\sigma(x)=-\}\label{e_spinnum}
\end{equation}
which denote the number of sites with spin $+$ and $-$, respectively.
\begin{prop}
\label{p_lowE}Suppose that $\sigma\in\mathcal{X}$ satisfies $H(\sigma)<\Gamma$.
Then, either (1) or (2) below must hold.
\begin{enumerate}
\item The configuration $\sigma$ belongs to $\mathcal{R}_{v}$ for some
$v\in\llbracket2,\,L-2\rrbracket$. In particular, $\mathcal{N}(\sigma)=\{\sigma\}$.
\item The configuration $\sigma$ belongs to $\mathcal{N}(\boxplus)$ or
$\mathcal{N}(\boxminus)$.
\end{enumerate}
\end{prop}

\begin{rem}
Two neighborhoods $\mathcal{N}(\boxplus)$ and $\mathcal{N}(\boxminus)$
are disjoint by Theorem \ref{t_energy barrier}.
\end{rem}

\begin{notation}
\label{n_bridge}
\begin{itemize}
\item A \textit{horizontal bridge} (resp. \textit{vertical bridge}) is a
row (resp. column), in which all spins are identical. If a bridge
consists of spin $+$ (resp. $-$), we call this bridge a $+$-bridge
(resp. $-$-bridge). Then, we denote by $B_{\pm}(\sigma)$ the number
of $\pm$-bridges in $\sigma\in\mathcal{X}$.
\item A \textit{cross} is a union of a horizontal bridge and a vertical
bridge. A cross consisting of spin $+$ (resp. $-$) is called a $+$-cross
(resp $-$-cross).
\item We denote by $r_{1},\,\dots,\,r_{L}$ the rows and $c_{1},\,\dots,\,c_{K}$
the columns of $\Lambda=\mathbb{T}_{K}\times\mathbb{T}_{L}$. For
$(v,\,h)\in\llbracket1,\,L\rrbracket\times\llbracket1,\,K\rrbracket$
and $\sigma\in\mathcal{X}$, we define
\begin{align*}
H_{r_{v}}(\sigma) & =\sum_{x,\,y\in r_{v}:\,x\sim y}\mathbf{1}\{\sigma(x)\ne\sigma(y)\}\;\text{and}\\
H_{c_{h}}(\eta) & =\sum_{x,\,y\in c_{h}:\,x\sim y}\mathbf{1}\{\sigma(x)\ne\sigma(y)\}\;,
\end{align*}
so that we can decompose the Hamiltonian in a way that
\begin{equation}
H(\sigma)=\sum_{v=1}^{L}H_{r_{v}}(\sigma)+\sum_{h=1}^{K}H_{c_{h}}(\sigma)\;.\label{e_decH2}
\end{equation}
A horizontal (resp. vertical) edge denotes an edge belonging to a
row (resp. column).
\end{itemize}
\end{notation}

The following lower bound for the Hamiltonian is a consequence of
notations and observations above.
\begin{lem}
\label{l_H2lb}It holds that
\[
H(\sigma)\ge2\,[\,K+L-B_{+}(\sigma)-B_{-}(\sigma)\,]\;.
\]
\end{lem}

\begin{proof}
The lemma follows directly from (\ref{e_decH2}) and the fact that
$H_{r_{v}}(\sigma)\ge2$ (resp. $H_{c_{h}}(\sigma)\ge2)$ if $r_{v}$
(resp. $c_{h}$) is not a bridge.
\end{proof}
We are now ready to prove Proposition \ref{p_lowE}.
\begin{proof}[Proof of Proposition \ref{p_lowE}]
Fix $\sigma\in\mathcal{X}$ with $H(\sigma)<\Gamma=2K+2$. By Lemma
\ref{l_H2lb}, we have
\[
2K+1\ge2\Big[\,K+L-B_{+}(\sigma)-B_{-}(\sigma)\,\Big]\;,
\]
and therefore $B_{+}(\sigma)+B_{-}(\sigma)\ge L$. Namely, there are
at least $L$ bridges. Let us take one of them and assume without
loss of generality that this is a $+$-bridge. Now, we consider three
cases separately.\\
\\
\textbf{(Case 1: $\sigma$ has a $+$-horizontal bridge without a
$+$-vertical one)} Since $H_{c_{h}}(\sigma)\ge2$ for all $h\in\llbracket1,\,K\rrbracket$,
we can observe from (\ref{e_decH2}) that $H_{r_{\ell}}(\sigma)=0$
for all $\ell\in\llbracket1,\,L\rrbracket$. This implies that all
rows are monochromatic, and therefore all columns are identical. Thus,
again by (\ref{e_decH2}), we get $H_{c_{k}}(\sigma)=2$ for all $k\in\llbracket1,\,K\rrbracket$,
and thus $\sigma\in\mathcal{R}_{v}$ for some $v\in\llbracket1,\,L\rrbracket$.
If $v\in\llbracket2,\,L-2\rrbracket$, then it is clear that $\mathcal{N}(\sigma)$
is a singleton since any configuration obtained from $\sigma$ by
flipping a spin has energy greater than or equal to $\Gamma$. Thus,
$\sigma$ satisfies the requirements of case (1). On the other hand,
if $v\notin\llbracket2,\,L-2\rrbracket$, we can readily observe that
$\sigma\in\mathcal{N}(\boxminus)$ or $\mathcal{N}(\boxplus)$. \\
\textbf{}\\
\textbf{(Case 2: $\sigma$ has a $+$-vertical bridge without a $+$-horizontal
one)} Since $\Delta H_{r_{v}}(\sigma)\ge2$ for all $v\in\llbracket1,\,L\rrbracket$,
we obtain from (\ref{e_decH2}) that $2K\ge2L$; hence, we obtain
a contradiction (to the assumption that $K<L$). \\
\\
\textbf{(Case 3: $\sigma$ has a $+$-cross)} Without loss of generality,
assume that $\mathbb{T}_{K}\times\{1\}$ and $\{1\}\times\mathbb{T}_{L}$
are $+$-bridges. Let us update each spin to $+$ in $\llbracket2,\,K\rrbracket\times\llbracket2,\,L\rrbracket$
in the ascending lexicographic order. The presence of spin $+$-bridges
ensures that the Hamiltonian cannot increase in the course of the
updates. Since we finally arrive at $\boxplus$, we can conclude that
\[
\Phi(\sigma,\,\boxplus)\le H(\sigma)<\Gamma\;.
\]
Thus, we have $\sigma\in\mathcal{N}(\boxplus)$.
\end{proof}

\subsection{\label{sec7.2}Properties of typical configurations}

In this subsection, we investigate the structure of typical configurations
introduced above. We start from two elementary lemmas.
\begin{lem}
\label{l_bulk1}Suppose that $\sigma\in\mathcal{B}$ and $\xi\in\mathcal{X}$
satisfy $\sigma\sim\xi$ and $H(\xi)\le\Gamma$. Then, the following
statements hold.
\begin{enumerate}
\item We have $\xi\in\mathcal{B}\cup\mathcal{Q}_{1}\cup\mathcal{Q}_{L-2}\subseteq\mathcal{C}$.
\item If $\sigma\in\mathcal{R}_{v}$ with $v\in\llbracket3,\,L-3\rrbracket$,
then $\xi\in\mathcal{B}_{\Gamma}$.
\item If $\sigma\in\mathcal{B}_{\Gamma}$, then $\xi\in\mathcal{B}$.
\end{enumerate}
\end{lem}

\begin{proof}
We consider two cases separately.
\begin{itemize}
\item \textbf{(Case 1: $\sigma\in\mathcal{R}_{v}$ for some $v\in\llbracket2,\,L-2\rrbracket$)}
Assume that $\sigma=\zeta_{\ell,\,v}$ for some $\ell\in\mathbb{T}_{L}$.
We can observe from the illustration given in Figure \ref{fig6.1}
that the only way of flipping a spin of $\sigma$ in such a way that
the resulting configuration has energy at most $\Gamma$ is either
to attach a protuberance of spin $+$ to the cluster of spin $+$
of $\sigma$ or to attach a protuberance of spin $-$ to the cluster
of spin $-$ of $\sigma$. This implies that
\[
\xi\in\{\zeta_{\ell,\,v\,;\,k,\,1}^{\textup{up}}\;,\;\;\zeta_{\ell,\,v\,;\,k,\,1}^{\textup{down}}\;,\;\;\zeta_{\ell,\,v-1;\,k,\,K-1}^{\textup{up}}\;,\;\;\zeta_{\ell+1,\,v-1;\,k,\,K-1}^{\textup{down}}:k\in\mathbb{T}_{K}\}\;.
\]
Hence, $\xi\in\mathcal{B}\cup\mathcal{Q}_{1}\cup\mathcal{Q}_{L-2}$.
This observation also implies that $\xi\in\mathcal{B}_{\Gamma}$ if
$v\in\llbracket3,\,L-3\rrbracket$, and hence part (2) is verified
here as well.
\item \textbf{(Case 2: $\sigma\in\mathcal{Q}_{v}$ for some $v\in\llbracket2,\,L-3\rrbracket$)}
Suppose that \textbf{$\sigma=\zeta_{\ell,\,v\,;\,k,\,h}^{\textup{up}}$}
for some $(k,\,\ell)\in\mathbb{T}_{K}\times\mathbb{T}_{L}$ and $h\in\llbracket1,\,K-1\rrbracket$.
In this case, we can observe that the only way of flipping a spin
of $\sigma$ without increasing the Hamiltonian is to expand or shrink
the protuberance of spin $+$ attached at $\zeta_{\ell,\,v}$, and
therefore
\[
\xi\in\{\zeta_{\ell,\,v\,;\,k,\,h-1}^{\textup{up}},\,\zeta_{\ell,\,v\,;\,k+1,\,h-1}^{\textup{up}}\}\;.
\]
Therefore, we have $\xi\in\mathcal{B}$ and hence parts (1) and (3)
are now verified. The same conclusion also holds for the case $\sigma=\zeta_{\ell,\,v\,;\,k,\,h}^{\textup{down}}$.
\end{itemize}
\end{proof}
The previous lemma implies the following result.
\begin{lem}
\label{l_bulk2}It holds that $\widehat{\mathcal{N}}(\mathcal{B}\,;\,\mathcal{C}\setminus\mathcal{B})=\mathcal{B}$.
\end{lem}

\begin{proof}
Since the energy of configurations belonging to $\mathcal{B}$ do
not exceed $\Gamma$, it follows immediately that
\[
\widehat{\mathcal{N}}(\mathcal{B}\,;\,\mathcal{C}\setminus\mathcal{B})\supset\mathcal{B}\;.
\]
Now, we claim the opposite inclusion, i.e.,
\begin{equation}
\widehat{\mathcal{N}}(\mathcal{B}\,;\,\mathcal{C}\setminus\mathcal{B})\subset\mathcal{B}\;.\label{e_bulk2}
\end{equation}
Suppose the contrary that there exists $\sigma\in\widehat{\mathcal{N}}(\mathcal{B}\,;\,\mathcal{C}\setminus\mathcal{B})$
such that $\sigma\notin\mathcal{B}$. Since $\sigma\in\widehat{\mathcal{N}}(\mathcal{B}\,;\,\mathcal{C}\setminus\mathcal{B})$,
there exists a $\Gamma$-path $(\omega_{t})_{t=0}^{T}$ in $\mathcal{X}\setminus(\mathcal{C}\setminus\mathcal{B})=(\mathcal{X}\setminus\mathcal{C})\cup\mathcal{B}$
connecting $\mathcal{B}$ and $\sigma$. Then, as $\omega_{0}\in\mathcal{B}$,
and $\omega_{T}\notin\mathcal{B}$, we can find $t_{0}\in\llbracket0,\,T-1\rrbracket$
such that $\omega_{t_{0}}\in\mathcal{B}$ and $\omega_{t_{0}+1}\notin\mathcal{B}$.
Since $(\omega_{t})_{t=0}^{T}$ is a path in $(\mathcal{X}\setminus\mathcal{C})\cup\mathcal{B}$,
we get
\[
\omega_{t_{0}+1}\in(\mathcal{X}\setminus\mathcal{B})\cap\big[\,(\mathcal{X}\setminus\mathcal{C})\cup\mathcal{B}\,\big]\subset\mathcal{X}\setminus\mathcal{C}\;.
\]
On the other hand, since $\omega_{t_{0}}\in\mathcal{B}$ we must have
$\omega_{t_{0}+1}\in\mathcal{C}$ by part (1) of Lemma \ref{l_bulk1}
and thus we have a contradiction. This proves (\ref{e_bulk2}) and
the proof is finished.
\end{proof}
Next, we prove that the two sets $\mathcal{E}^{+}$ and $\mathcal{E}^{-}$
are indeed disjoint.
\begin{prop}
\label{p_nonint}We have that $\mathcal{E}^{+}\cap\mathcal{E}^{-}=\emptyset$.
\end{prop}

\begin{proof}
Suppose the contrary that there exists a path $(\omega_{t})_{t=0}^{T}$
is a $\Gamma$-path from $\boxminus$ to $\boxplus$ in $\mathcal{X}\setminus\mathcal{B}_{\Gamma}$.
Define $u:\llbracket0,\,T\rrbracket\rightarrow\mathbb{R}$ as
\[
u(t)=B_{+}(\omega_{t})\;\;\;\;;\;t\in\llbracket0,\,T\rrbracket\;,
\]
where $B_{+}(\cdot)$ is defined in Notation \ref{n_bridge}. Then,
we have that
\begin{equation}
u(0)=0\;,\;u(T)=K+L\;,\text{ and }|u(t+1)-u(t)|\le2\text{ for all }t\in\llbracket0,\,T-1\rrbracket\;.\label{e_nonint1}
\end{equation}
Thus, the following time $t^{*}$ is well defined:
\begin{equation}
t^{*}=\min\,\{t\in\llbracket0,\,T-1\rrbracket:u(t),\,u(t+1)\ge2\}\;.\label{e_nonint2}
\end{equation}
Note that, since we need to change at least $2K-1$ spins from $\boxminus$
to get $u(t)\ge2,$ we have $t^{*}\ge2K-1$. Then, by (\ref{e_nonint1}),
we have $B_{+}(\omega_{t^{*}})=2$ or $3$. We divide the proof into
three cases as in Proposition \ref{p_lowE}.\\
\\
\textbf{(Case 1: $\omega_{t^{*}}$ has $+$-horizontal bridges without
a $+$-vertical one) }For this case, if $B_{+}(\omega_{t^{*}})=3$,
we have $B_{+}(\omega_{t^{*}-1})\ge2$ and thus we get a contradiction
to the minimality of $t^{*}$. Hence, we have $B_{+}(\omega_{t^{*}})=2$.

Since $\omega_{t^{*}}$ does have both $+$- and $-$-vertical bridges,
we get $H_{c_{h}}(\omega_{t^{*}})\ge2$ for all $h\in\llbracket1,\,K\rrbracket$.
By (\ref{e_decH2}) and the fact that $H(\omega_{t^{*}})\le\Gamma=2K+2$,
we can readily observe that $\omega_{t^{*}}\in\mathcal{R}_{2}\cup\mathcal{Q}_{2}$.
Since $\mathcal{Q}_{2}\subseteq\mathcal{B}_{\Gamma}$ and since $(\omega_{t})_{t=0}^{T}$
is a path in $\mathcal{X}\setminus\mathcal{B}_{\Gamma}$, we can conclude
that $\omega_{t^{*}}\in\mathcal{R}_{2}$. Since $H(\omega_{t^{*}+1})\le\Gamma$
and $u(t^{*}+1)\ge2$, we are forced to have $\omega_{t^{*}+1}\in\mathcal{B}_{\Gamma}$
which is a contradiction.\\
\\
\textbf{(Case 2: $\omega_{t^{*}}$ has $+$-vertical bridges without
a $+$-horizontal one)} This case is similar to \textbf{(Case 1)}.
\\
\\
\textbf{(Case 3: $\omega_{t^{*}}$ has a $+$-cross)} In this case,
$\omega_{t^{*}}$ cannot have a $-$-bridge. Thus, by (\ref{e_nonint1}),
the configuration $\omega_{t^{*}}$ has at most three bridges. Therefore,
by Lemma \ref{l_H2lb},
\[
H(\omega_{t^{*}})\ge2(K+L-3)>\Gamma\;,
\]
which contradicts the fact that $(\omega_{t})_{t=0}^{T}$ is a $\Gamma$-path.

Now, the assertion of the proposition directly follows since if $\mathcal{E}^{+}\cap\mathcal{E}^{-}\neq\emptyset$,
there must exist a $\Gamma$-path from $\boxminus$ to $\boxplus$
in $\mathcal{X}\setminus\mathcal{B}_{\Gamma}$.
\end{proof}
The previous proposition implies that any $\Gamma$-path connecting
$\boxminus$ and $\boxplus$ has to touch the set $\mathcal{B}_{\Gamma}$,
i.e., has to path through bulk typical configurations.

The next proposition concerns the relationships between bulk and edge
typical configurations.
\begin{prop}
\label{p_typ2}The following properties hold:
\begin{enumerate}
\item It holds that
\begin{equation}
\mathcal{E}^{-}\cap\mathcal{B}=\mathcal{R}_{2}\;\;\;\;\text{and}\;\;\;\;\mathcal{E}^{+}\cap\mathcal{B}=\mathcal{R}_{L-2}\;.\label{e_typ2}
\end{equation}
\item We have that $\mathcal{E}\cup\mathcal{B}=\widehat{\mathcal{N}}(\mathcal{S})$.
\end{enumerate}
\end{prop}

\begin{proof}
(1) We only prove the first one of (\ref{e_typ2}), as the second
one follows in the same manner.

First, we have $\mathcal{B}\supset\mathcal{R}_{2}$ from the definition
of $\mathcal{B}$. On the other hand, since the canonical path connecting
$\mathcal{R}_{2}$ and $\boxminus$ is a $\Gamma$-path in $\mathcal{X}\setminus\mathcal{B}_{\Gamma}$,
we also have $\mathcal{E}^{-}\supset\mathcal{R}_{2}$. Thus, we have
proved that,
\begin{equation}
\mathcal{E}^{-}\cap\mathcal{B}\supset\mathcal{R}_{2}\;.\label{e_int1}
\end{equation}
Now, we claim that the reversed inclusion also holds. To prove this
claim, we begin by observing that, since $\mathcal{B}_{\Gamma}$ and
$\mathcal{E}$ are disjoint by definition (cf. (\ref{e_EAB2})), we
can conclude that
\[
\mathcal{E}^{-}\cap\mathcal{B}\subset\mathcal{B}\setminus\mathcal{B}_{\Gamma}=\bigcup_{v\in\llbracket2,\,L-2\rrbracket}\mathcal{R}_{v}\;.
\]
For $\sigma\in\mathcal{R}_{v}$ with $v\in\llbracket3,\,L-3\rrbracket$,
we cannot have a path in $\mathcal{X}\setminus\mathcal{B}_{\Gamma}$
connecting $\boxminus$ and $\sigma$ by Lemma \ref{l_bulk1}-(2).
We therefore have $\sigma\notin\mathcal{E}^{-}$, and thus we can
conclude that
\begin{equation}
\mathcal{E}^{-}\cap\mathcal{B}\subset\mathcal{R}_{2}\cup\mathcal{R}_{L-2}\;.\label{e_int2}
\end{equation}
By the same reason with the inclusion $\mathcal{E}^{-}\supset\mathcal{R}_{2}$,
we also have $\mathcal{E}^{+}\supset\mathcal{R}_{L-2}$. Therefore,
any configuration $\sigma\in\mathcal{R}_{L-2}$ cannot belong to $\mathcal{E}^{-}$
by Proposition \ref{p_nonint}; hence, from (\ref{e_int2}), we can
deduce that
\begin{equation}
\mathcal{E}^{-}\cap\mathcal{B}\subset\mathcal{R}_{2}\;.\label{e_int3}
\end{equation}
This proves the claim and we are done.\\
\\
(2) The inclusion $\mathcal{E}\subset\widehat{\mathcal{N}}(\mathcal{S})$
is obvious from the definition of $\mathcal{E}$, and the inclusion
$\mathcal{B}\subset\widehat{\mathcal{N}}(\mathcal{S})$ also follows
immediately from the fact that any bulk typical configuration is connected
to $\boxminus$ (or $\boxplus$) via a part of a canonical path, which
is a $\Gamma$-path (cf. Remark \ref{r_H2}). Thus, we can conclude
that
\begin{equation}
\mathcal{E}\cup\mathcal{B}\subset\widehat{\mathcal{N}}(\mathcal{S})\;.\label{e_unk1}
\end{equation}
Now we prove the reversed inclusion. By Lemma \ref{l_set} with $\mathcal{P}=\mathcal{C}\setminus\mathcal{B}$
and $\mathcal{Q}=\mathcal{B}$, we get
\begin{align}
\widehat{\mathcal{N}}(\mathcal{C}) & =\widehat{\mathcal{N}}(\mathcal{B}\,;\,\mathcal{C}\setminus\mathcal{B})\cup\widehat{\mathcal{N}}(\mathcal{C}\setminus\mathcal{B}\,;\,\mathcal{B})\;.\label{unk_2}
\end{align}
By Lemma \ref{e_unk1}, we have
\begin{equation}
\widehat{\mathcal{N}}(\mathcal{B}\,;\,\mathcal{C}\setminus\mathcal{B})=\mathcal{B}\;.\label{e_unk3}
\end{equation}
Since any configuration in $\mathcal{C}\setminus\mathcal{B}$ is connected
to either $\boxminus$ or $\boxplus$ via a part of a canonical path
which is a $\Gamma$-path in $\mathcal{X}\setminus\mathcal{B}_{\Gamma}$,
we obtain
\begin{equation}
\widehat{\mathcal{N}}(\mathcal{C}\setminus\mathcal{B}\,;\,\mathcal{B})\subset\widehat{\mathcal{N}}(\mathcal{C}\setminus\mathcal{B}\,;\,\mathcal{B}_{\Gamma})\subset\widehat{\mathcal{N}}(\mathcal{S}\,;\,\mathcal{B}_{\Gamma})=\mathcal{E}\;.\label{e_unk4}
\end{equation}
By combining (\ref{unk_2}), (\ref{e_unk3}), and (\ref{e_unk4}),
we get
\[
\widehat{\mathcal{N}}(\mathcal{C})\subset\mathcal{E}\cup\mathcal{B}\;.
\]
Since $\mathcal{S}\subset\mathcal{C}$, the last inclusion implies
the opposite inclusion of (\ref{e_unk1}), and we are done.
\end{proof}

\subsection{\label{sec7.3}Characterization of edge typical configurations}

As mentioned before, edge typical configurations have far more complex
structure than bulk ones. In this subsection, we study this complex
structure in detail.

Our analysis starts with a decomposition of the form

\[
\mathcal{E}^{-}=\mathcal{O}^{-}\cup\mathcal{I}^{-}\text{\;\;\;}\text{and\;\;\;}\mathcal{E}^{+}=\mathcal{O}^{+}\cup\mathcal{I}^{+}\;,
\]
where
\[
\mathcal{O}^{\pm}=\{\sigma\in\mathcal{E}^{\pm}:H(\sigma)=\Gamma\}\;\;\;\;\text{and}\;\;\;\;\mathcal{I}^{\pm}=\{\sigma\in\mathcal{E}^{\pm}:H(\sigma)<\Gamma\}\;.
\]
Then, we analyze the structure based on this decomposition. For the
concreteness of the discussion, we focus only on $\mathcal{E}^{-}$,
as the analysis of $\mathcal{E}^{+}$ is essentially identical.

By Proposition \ref{p_lowE}, we can see that
\begin{equation}
\mathcal{I}^{-}=\mathcal{N}(\boxminus)\cup\mathcal{R}_{2}\;.\label{e_IA}
\end{equation}
We now construct a graph and a Markov chain which represent the asymptotic
behavior of the Metropolis dynamics on $\mathcal{E}^{-}$. Heuristically,
since the configurations belonging to $\mathcal{N}(\boxminus)$ are
indistinguishable in the scale $e^{\beta\Gamma}$ (as they can be
communicated by a much shorter scale), we shall identify all the configurations
in $\mathcal{N}(\boxminus)$ with $\boxminus$ and define
\begin{equation}
\overline{\mathcal{I}}^{-}=\boxminus\cup\mathcal{R}_{2}\;.\label{e_IArep}
\end{equation}
With this notation, we can write
\begin{equation}
\mathcal{E}^{-}=\mathcal{O}^{-}\cup\Big(\,\bigcup_{\sigma\in\overline{\mathcal{I}}^{-}}\mathcal{N}(\sigma)\,\Big)\;.\label{e_EAdec}
\end{equation}
Now, we define a graph structure on the vertex set $\mathscr{V^{-}}$
defined by
\begin{equation}
\mathscr{V^{-}}=\mathcal{O}^{-}\cup\overline{\mathcal{I}}^{-}\;,\label{e_VA}
\end{equation}
and define a continuous-time Markov chain on that graph.
\begin{defn}
\label{d_EAMc2}
\begin{itemize}
\item \textbf{(Graph) }We introduce a graph structure $\mathscr{G}^{-}=(\mathscr{V}^{-},\,\mathscr{E}^{-})$
where for $\sigma,\,\sigma'\in\mathscr{V^{-}}$, we say that $\{\sigma,\,\sigma'\}\in\mathscr{E}^{-}$
if and only if
\[
\begin{cases}
\sigma,\,\sigma'\in\mathcal{O}^{-}\text{ and }\sigma\sim\sigma'\text{ or}\\
\sigma\in\mathcal{O}^{-},\text{ }\sigma'\in\overline{\mathcal{I}}^{-}\text{ and }\sigma\sim\xi\text{ for some }\xi\in\mathcal{N}(\sigma')\;.
\end{cases}
\]
\item \textbf{(Markov chain)} The rate function $r^{-}:\mathscr{V}^{-}\times\mathscr{V}^{-}\rightarrow[0,\,\infty)$
is defined by, for all $\{\sigma,\,\sigma'\}\in\mathscr{E}$,
\begin{equation}
r^{\mathscr{-}}(\sigma,\,\sigma')=\begin{cases}
1 & \text{if }\sigma,\,\sigma'\in\mathcal{O}^{-}\;,\\
|\{\xi\in\mathcal{N}(\sigma):\xi\sim\sigma'\}| & \text{if }\sigma\in\overline{\mathcal{I}}^{-}\;,\;\sigma'\in\mathcal{O}^{-}\;,\\
|\{\xi\in\mathcal{N}(\sigma'):\xi\sim\sigma\}| & \text{if }\sigma\in\mathcal{O}^{-}\;,\;\sigma'\in\overline{\mathcal{I}}^{-}\;,
\end{cases}\label{e_EAMc2}
\end{equation}
and we finally set $r^{\mathscr{-}}(\sigma,\,\sigma')=0$ if $\{\sigma,\,\sigma'\}\notin\mathscr{E^{-}}$.
Then, denote by $(Z^{\mathscr{-}}(t))_{t\ge0}$ a continuous-time
Markov chain on $\mathscr{V^{-}}$ with rate $r^{\mathscr{-}}(\cdot,\,\cdot)$.
Since the rate is symmetric, the Markov chain $Z^{\mathscr{-}}(\cdot)$
is reversible with respect to the uniform distribution on $\mathscr{V}^{-}$.
\item We denote by $h_{\cdot,\,\cdot}^{-}(\cdot)$, $\mathrm{cap}^{-}(\cdot,\,\cdot)$,
$D^{-}(\cdot)$, and $\|\cdot\|_{-}$ the equilibrium potential, capacity,
Dirichlet form, and flow norm with respect to the Markov process $Z^{-}(\cdot)$,
respectively. In addition, denote by $L^{-}$ the generator of the
process $Z^{-}(\cdot)$ acting on $f:\mathscr{V}^{-}\rightarrow\mathbb{R}$
in a way that
\begin{equation}
(L^{-}f)(\sigma)=\sum_{\sigma'\in\mathscr{V}^{-}:\,\{\sigma,\,\sigma'\}\in\mathscr{E}^{-}}r^{-}(\sigma,\,\sigma')\,\{f(\sigma')-f(\sigma)\}\;.\label{e_genZA}
\end{equation}
\end{itemize}
\end{defn}

We first show that the Markov process $Z^{-}(\cdot)$ approximates
in some sense the Metropolis dynamics $\sigma_{\beta}(\cdot)$ in
$\mathcal{E}^{A}$.
\begin{prop}
\label{p_ZA2}Define a projection map $\Pi^{-}:\mathcal{E}^{-}\rightarrow\mathscr{V}^{-}$
by
\[
\Pi^{-}(\sigma)=\begin{cases}
\xi & \text{if }\sigma\in\mathcal{N}(\xi)\text{ for some }\xi\in\overline{\mathcal{I}}^{-}\;,\\
\sigma & \text{if }\sigma\in\mathcal{O}^{-}\;.
\end{cases}
\]
Then, there exists a constant $C=C(K,\,L)>0$ such that
\begin{enumerate}
\item for $\sigma,\,\sigma'\in\mathcal{O}^{-}$, we have
\begin{equation}
\Big|\,\frac{1}{2}e^{-\Gamma\beta}r^{-}(\Pi^{-}(\sigma),\,\Pi^{-}(\sigma'))-\mu_{\beta}(\sigma)\,c_{\beta}(\sigma,\,\sigma')\,\Big|\le Ce^{-(\Gamma+2)\beta}\;,\label{e_za21}
\end{equation}
\item for $\sigma\in\mathcal{O}^{-}$ and $\sigma'\in\overline{\mathcal{I}}^{-}$,
we have
\begin{equation}
\Big|\,\frac{1}{2}e^{-\Gamma\beta}r^{-}(\Pi^{-}(\sigma),\,\Pi^{-}(\sigma'))-\sum_{\xi\in\mathcal{N}(\sigma')}\mu_{\beta}(\sigma)\,c_{\beta}(\sigma,\,\xi)\,\Big|\le Ce^{-(\Gamma+2)\beta}\;.\label{e_za22}
\end{equation}
\end{enumerate}
\end{prop}

\begin{proof}
Suppose that $\sigma,\,\sigma'\in\mathcal{O}^{-}$. If $\sigma\not\sim\sigma'$,
then the left-hand side of (\ref{e_za21}) is clearly $0$. On the
other hand, if $\sigma\sim\sigma'$ so that $\{\sigma,\,\sigma'\}\in\mathscr{E}^{-}$,
then by (\ref{muprop}) and (\ref{e_EAMc2}),
\[
\Big|\,\frac{1}{2}e^{-\Gamma\beta}r^{-}(\Pi^{-}(\sigma),\,\Pi^{-}(\sigma'))-\mu_{\beta}(\sigma)\,c_{\beta}(\sigma,\,\sigma')\,\Big|=\Big|\,\frac{1}{2}e^{-\Gamma\beta}-\frac{1}{Z_{\beta}}e^{-\Gamma\beta}\,\Big|
\]
since $\mu_{\beta}(\sigma)=\mu_{\beta}(\sigma')=\frac{1}{Z_{\beta}}e^{-\Gamma\beta}$.
By (\ref{Zbest}), the right-hand side is $O(e^{-(\Gamma+2)\beta})$.
This proves part (1).

Now, we consider part (2). Let $\sigma\in\mathcal{O}^{-}$ and $\sigma'\in\overline{\mathcal{I}}^{-}$.
Similarly, we can assume $\{\sigma,\,\sigma'\}\in\mathscr{E}^{-}$
since otherwise the left-hand side of (\ref{e_za22}) is $0$. Then,
by (\ref{muprop}) and (\ref{e_EAMc2}), we can write
\begin{align*}
 & \Big|\,\frac{1}{2}e^{-\Gamma\beta}r^{-}(\Pi^{-}(\sigma),\,\Pi^{-}(\sigma'))-\sum_{\xi\in\mathcal{N}(\sigma')}\mu_{\beta}(\sigma)\,c_{\beta}(\sigma,\,\xi)\,\Big|\\
= & \Big|\,\frac{1}{2}e^{-\Gamma\beta}\,|\{\xi\in\mathcal{N}(\sigma'):\xi\sim\sigma\}|-\sum_{\xi\in\mathcal{N}(\sigma'):\,\xi\sim\sigma}\min\{\mu_{\beta}(\sigma),\,\mu_{\beta}(\xi)\}\,\Big|\\
= & |\{\xi\in\mathcal{N}(\sigma'):\xi\sim\sigma\}|\times\Big|\,\frac{1}{2}e^{-\Gamma\beta}-\frac{1}{Z_{\beta}}e^{-\Gamma\beta}\,\Big|\;,
\end{align*}
since $\min\{\mu_{\beta}(\sigma),\,\mu_{\beta}(\xi)\}=\mu_{\beta}(\sigma)$
for all $\xi\in\mathcal{N}(\sigma')$. By (\ref{Zbest}), the last
line is bounded by $KL\times O(e^{-\Gamma\beta}e^{-2\beta})=O(e^{-(\Gamma+2)\beta})$.
\end{proof}
In view of this proposition, we can assert that the equilibrium potential
$h_{\boxminus,\,\mathcal{R}_{2}}^{-}(\cdot)$ approximates the equilibrium
potential of the Metropolis dynamics in $\mathcal{E}^{-}$. For this
reason, the equilibrium potential $h_{\boxminus,\,\mathcal{R}_{2}}^{-}(\cdot)$
plays a significant role in the construction of the test function
and flow in the next sections.

Now, we are ready to define the edge constant $\mathfrak{e}$ introduced
in Section \ref{sec4.1}. Define

\begin{equation}
\mathfrak{e}=\frac{1}{|\mathscr{V}^{-}|\,\mathrm{cap}^{-}(\boxminus,\,\mathcal{R}_{2})}\;.\label{e_e2def}
\end{equation}
The appearance of $\mathrm{cap}^{-}(\boxminus,\,\mathcal{R}_{2})$
is quite natural in that the equilibrium potential $h_{\boxminus,\,\mathcal{R}_{2}}^{-}(\cdot)$
is the correct approximation of the equilibrium potential of the Metropolis
dynamics in $\mathcal{E}^{-}$. We conclude this section by showing
that the constant $\mathfrak{e}$ is small.
\begin{prop}
\label{p_e2est}We have that $\mathfrak{e}\le\frac{1}{L}$.
\end{prop}

\begin{proof}
We use the Thomson principle (cf. Theorem \ref{t_TP_rev}) to prove
the proposition. We define a test flow $\psi$ on $\mathscr{V}^{-}\times\mathscr{V}^{-}$
(with respect to the Markov process $Z^{-}(\cdot)$) as
\[
\begin{cases}
\psi(\boxminus,\,\zeta_{\ell,\,1\,;\,k,\,1}^{\textup{up}})=\frac{1}{KL} & \text{for }(k,\,\ell)\in\mathbb{T}_{K}\times\mathbb{T}_{L}\;,\\
\psi(\zeta_{\ell,\,1\,;\,k,\,h}^{\mathrm{\textup{up}}},\,\zeta_{\ell,\,1\,;\,k,\,h+1}^{\mathrm{\textup{up}}})=\frac{1}{KL} & \text{for }(k,\,\ell)\in\mathbb{T}_{K}\times\mathbb{T}_{L}\text{ and }h\in\llbracket1,\,K-1\rrbracket\;.
\end{cases}
\]
We set $\psi(\sigma,\,\sigma')=0$ for all other cases. Notice that
$\zeta_{\ell,\,1\,;\,k,\,h}^{\mathrm{\textup{up}}}\in\mathcal{O}^{-}$
for all $h\in\llbracket1,\,K-1\rrbracket$, and that $\{\boxminus,\,\zeta_{\ell,\,1\,;\,k,\,1}^{\textup{up}}\}\in\mathscr{E}^{-}$
since $\zeta_{\ell,\,1\,;\,k,\,1}^{\mathrm{\textup{up}}}\sim\zeta_{\ell,\,1}$
and $\zeta_{\ell,\,1}\in\mathcal{N}(\boxminus)$, where the latter
is readily follows from the part of a canonical path connecting $\boxminus$
and $\zeta_{\ell,\,1}$ is a $(\Gamma-2)$-path. Notice that this
is a unit flow from $\boxminus$ to $\mathcal{R}_{2}$ since
\begin{align*}
 & (\mathrm{div}\,\psi)(\boxminus)=\sum_{\ell\in\mathbb{T}_{L}}\sum_{k\in\mathbb{T}_{K}}\frac{1}{KL}=1\;,\\
 & (\mathrm{div}\,\psi)(\mathcal{R}_{2})=\sum_{\ell\in\mathbb{T}_{L}}(\mathrm{div}\,\psi)(\zeta_{\ell,\,2})=\sum_{\ell\in\mathbb{T}_{L}}\sum_{k\in\mathbb{T}_{K}}\frac{-1}{KL}=-1\;,
\end{align*}
and moreover we can readily check that
\[
(\mathrm{div}\,\psi)(\sigma)=0\text{\;\;\;for all }\sigma\in\mathscr{V}^{-}\setminus(\boxminus\cup\mathcal{R}_{2})\;.
\]
Therefore, by Theorem \ref{t_TP_rev}, we get
\begin{equation}
\mathrm{cap}^{-}(\boxminus,\,\mathcal{R}_{2})\ge\frac{1}{\|\psi\|_{-}^{2}}\;.\label{e_e2est}
\end{equation}
It remains to evaluate the flow norm $\|\psi\|_{-}^{2}$ which is
indeed equal to (since the uniform distribution is the invariant measure
for the Markov process $Z^{-}(\cdot)$)
\begin{align*}
 & \sum_{\ell\in\mathbb{T}_{L}}\sum_{k\in\mathbb{T}_{K}}\Big[\,\frac{\psi(\boxminus,\,\zeta_{\ell,\,1\,;\,k,\,1}^{\mathrm{\textup{up}}}){}^{2}}{1/|\mathscr{V}^{-}|}+\sum_{h=1}^{K-1}\frac{\psi(\zeta_{\ell,\,1\,;\,k,\,h}^{\mathrm{\textup{up}}},\,\zeta_{\ell,\,1\,;\,k,\,h+1}^{\mathrm{\textup{up}}}){}^{2}}{1/|\mathscr{V}^{-}|}\,\Big]\\
 & =LK^{2}\times\frac{|\mathscr{V}^{-}|}{K^{2}L^{2}}=\frac{|\mathscr{V}^{-}|}{L}\;.
\end{align*}
Injecting this to (\ref{e_e2est}) completes the proof.
\end{proof}
\begin{rem}
\label{r_e2est}In fact, we can verify that there exist two constants
$C_{1},\,C_{2}>0$ such that
\[
\frac{C_{1}}{KL}\le\mathfrak{e}\le\frac{C_{2}}{KL}\;.
\]
We leave this as an exercise. This can be proven with a more refined
test flow.
\end{rem}

\begin{rem}
Of course, we can also establish the results corresponding to Definition
\ref{d_EAMc2}, Propositions \ref{p_ZA2} and \ref{p_e2est} for $\mathcal{E}^{+}$
in the completely identical manner. The constant $\mathfrak{e}$ defined
for $\mathcal{E}^{+}$ should be in accordance with (\ref{e_e2def})
by the symmetry of the model.
\end{rem}

\section{\label{sec9}Upper Bound for Capacities}

In this section, we construct a test function $f_{0}:\mathcal{X}\rightarrow\mathbb{R}$
appearing in Proposition \ref{p_Capup}. For the convenience of notation,
we write
\begin{equation}
\mathfrak{h}^{\pm}(\cdot)=h_{\boxminus,\,\mathcal{R}_{2}}^{\pm}(\cdot)\label{e_shorteq}
\end{equation}
which is the equilibrium potential between $\{\boxminus\}$ and $\mathcal{R}_{2}$
with respect to the process $Z^{\pm}(\cdot)$ (cf. Definition \ref{d_EAMc2}).

\subsection{\label{sec9.1}Construction of test function}

Now, we construct a function $f_{0}:\mathcal{X}\rightarrow\mathbb{R}$.
In the end, we shall verify that this function fulfills all requirements
of the function $f_{0}$ appearing in Proposition \ref{p_Capup}.
Before defining the test function explicitly, we briefly explain the
gist of the idea. On edge typical configurations (i.e., on $\mathcal{E}^{\pm})$,
we choose $f_{0}$ as a rescale of $\mathfrak{h}^{\pm}$. This construction
mainly comes from the fact that the process $Z^{\pm}(\cdot)$ successfully
characterizes the behavior of the original process on edge typical
configurations by Proposition \ref{p_ZA2}. On the other hand, on
bulk typical configurations, we define $f$ as a rescale of the equilibrium
potential of a symmetric simple random walk on an one-dimensional
line. This is because the Metropolis dynamics behaves as an one-dimensional
random walk there thanks to the simple geometry between them.
\begin{defn}[Test function]
\label{d_testf2}We construct a test function $f_{0}:\mathcal{X}\rightarrow\mathbb{R}$
on $\mathcal{E}$, $\mathcal{B}$, and $(\mathcal{E}\cup\mathcal{B}){}^{c}=\mathcal{X}\setminus(\mathcal{E}\cup\mathcal{B})$,
separately.
\begin{enumerate}
\item \textbf{Construction of $f_{0}$ on edge typical configurations $\mathcal{E}=\mathcal{E}^{-}\cup\mathcal{E}^{+}$.}
\begin{itemize}
\item For $\sigma\in\mathcal{E}^{-}$, we recall the decomposition (\ref{e_EAdec})
of $\mathcal{E}^{-}$ and define
\begin{equation}
f_{0}(\sigma)=\begin{cases}
1-\frac{\mathfrak{e}}{\kappa}(1-\mathfrak{h}^{-}(\sigma)) & \text{if }\sigma\in\mathcal{O}^{-}\;,\\
1-\frac{\mathfrak{e}}{\kappa}(1-\mathfrak{h}^{-}(\xi)) & \text{if }\sigma\in\mathcal{N}(\xi)\text{ for some }\xi\in\overline{\mathcal{I}}^{-}\;.
\end{cases}\label{e_testf2.1}
\end{equation}
\item For $\sigma\in\mathcal{E}^{+}$, we similarly define
\begin{equation}
f_{0}(\sigma)=\begin{cases}
\frac{\mathfrak{e}}{\kappa}(1-\mathfrak{h}^{+}(\sigma)) & \text{if }\sigma\in\mathcal{O}^{+}\;,\\
\frac{\mathfrak{e}}{\kappa}(1-\mathfrak{h}^{+}(\xi)) & \text{if }\sigma\in\mathcal{N}(\xi)\text{ for some }\xi\in\overline{\mathcal{I}}^{+}\;.
\end{cases}\label{e_testf2.2}
\end{equation}
\end{itemize}
\item \textbf{Construction of $f_{0}$ on bulk typical configurations $\mathcal{B}$.}
In view of (\ref{e_Bab2}), it suffices to define this object in the
following two cases.
\begin{itemize}
\item For $\sigma\in\mathcal{R}_{v}$ with $v\in\llbracket2,\,L-2\rrbracket$,
we set
\begin{equation}
f_{0}(\sigma)=\frac{1}{\kappa}\,\Big[\,\frac{L-2-v}{L-4}\mathfrak{b}+\mathfrak{e}\,\Big]\;.\label{e_testf2.3}
\end{equation}

\item For $\sigma\in\mathcal{Q}_{v}$ with $v\in\llbracket2,\,L-3\rrbracket$,
we can write $\sigma=\zeta_{\ell,\,v\,;\,k,\,h}^{\textup{up}}$ or
$\zeta_{\ell,\,v\,;\,k,\,h}^{\text{down}}$ for some $(k,\,\ell)\in\mathbb{T}_{K}\times\mathbb{T}_{L}$
and $h\in\llbracket1,\,K-1\rrbracket$. For such $\sigma$, we set
\begin{equation}
f_{0}(\sigma)=\frac{1}{\mathfrak{\kappa}}\,\Big[\,\frac{(K+2)(L-2-v)-(h+1)}{(K+2)(L-4)}\mathfrak{b}+\mathfrak{e}\,\Big]\;.\label{e_testf2.4}
\end{equation}
\end{itemize}
\item \textbf{Construction of $f_{0}$ on the remainder $(\mathcal{E}\cup\mathcal{B})^{c}$.}
We define $f_{0}\equiv1$ on this set.
\end{enumerate}
\end{defn}

\begin{rem}
We note that $\mathcal{E}^{-}$ and $\mathcal{B}$ are not disjoint
and their intersection is $\mathcal{R}_{2}$ by Proposition \ref{p_typ2}.
However, we can easily check that our constructions of $f_{0}$ on
$\mathcal{R}_{2}$ in parts (1) and (2) of the previous definition
agree with the value $1-\mathfrak{e}/\kappa$. A similar result also
holds for $\mathcal{E}^{+}$ and $\mathcal{B}$.
\end{rem}

\subsection{\label{sec9.2}Properties of test function}

Now, we will confirm that the test function $f_{0}$ satisfies the
requirements of $f_{0}$ appearing in Proposition \ref{p_Capup}.
\begin{prop}
\label{p_testf2}The function $f_{0}$ constructed in Definition \ref{d_testf2}
belongs to $\mathfrak{C}_{1,\,0}(\{\boxplus\},\,\{\boxminus\})$ and
satisfies
\begin{equation}
\mathscr{D}_{\beta}(f_{0})=\frac{1+o_{\beta}(1)}{2\kappa}\,e^{-\Gamma\beta}\;.\label{dees}
\end{equation}
\end{prop}

\begin{proof}
For the simplicity of notation, let us write $f=f_{0}$. Since we
have $f(\boxminus)=1$ and $f(\boxplus)=0$ by part (1) of Definition
\ref{d_testf2}, we immediately have $f\in\mathfrak{C}_{1,\,0}(\{\boxplus\},\,\{\boxminus\})$.
Now, it remains to prove (\ref{dees}).

Let us divide the Dirichlet form $\mathscr{D}_{\beta}(f)$ into
\begin{equation}
\Big[\,\sum_{\{\sigma,\,\xi\}\subset(\mathcal{E}\cup\mathcal{B})^{c}}+\sum_{\sigma\in\mathcal{E}\cup\mathcal{B},\,\xi\in(\mathcal{E}\cup\mathcal{B})^{c}}+\sum_{\{\sigma,\,\xi\}\subset\mathcal{E}\cup\mathcal{B}}\,\Big]\,\mu_{\beta}(\sigma)\,c_{\beta}(\sigma,\,\xi)\,\{f(\xi)-f(\sigma)\}^{2}\;,\label{e_decDiri}
\end{equation}
where all summations are carried out for two connected configurations
$\sigma$ and $\xi$, i.e., $\sigma\sim\xi$.

The first summation is trivially $0$ by part (3) of Definition \ref{d_testf2}.
Now to consider the second summation, we recall from part (2) of Proposition
\ref{p_typ2} that $\mathcal{E}\cup\mathcal{B}=\widehat{\mathcal{N}}(\mathcal{S})$.
This implies that $H(\sigma)\le\Gamma$ and $H(\xi)\ge\Gamma+1$.
Therefore, by (\ref{muprop}), we have
\[
\mu_{\beta}(\sigma)\,c_{\beta}(\sigma,\,\xi)=\mu_{\beta}(\xi)=\frac{1}{Z_{\beta}}e^{-\beta H(\xi)}=o_{\beta}(e^{-\Gamma\beta})\;,
\]
where we implicitly used the fact that $Z_{\beta}\rightarrow2$ as
$\beta\rightarrow\infty$ at the last equality. Moreover, since $f(\sigma)\in[0,\,1]$
for all $\sigma\in\mathcal{X}$ by our construction, we can assert
that the second summation in (\ref{e_decDiri}) is $o_{\beta}(e^{-\Gamma\beta})$.

It remains to estimate the third summation of (\ref{e_decDiri}).
For $\mathcal{A}\subset\mathcal{X}$, we write
\begin{equation}
E(\mathcal{A})=\big\{\,\{\sigma,\,\xi\}\subset\mathcal{A}:\sigma\sim\xi\,\big\}\;.\label{e_edgeA}
\end{equation}
By part (1) of Proposition \ref{p_typ2}, we can decompose $E(\mathcal{E}\cup\mathcal{B})$
into
\begin{equation}
E(\mathcal{E}\cup\mathcal{B})=E(\mathcal{B})\cup E(\mathcal{E}^{-})\cup E(\mathcal{E}^{+})\;.\label{e_decE}
\end{equation}
Hence, we can further decompose the third summation of (\ref{e_decDiri})
into
\begin{equation}
\Big[\,\sum_{\{\sigma,\,\xi\}\in E(\mathcal{B})}+\sum_{\{\sigma,\,\xi\}\in E(\mathcal{E}^{-})}+\sum_{\{\sigma,\,\xi\}\in E(\mathcal{E}^{+})}\,\Big]\,\mu_{\beta}(\sigma)\,c_{\beta}(\sigma,\,\xi)\,\{f(\xi)-f(\sigma)\}^{2}\;.\label{e_mce1}
\end{equation}

Now, we compute the first summation of (\ref{e_mce1}). Decompose
\[
E(\mathcal{B})=\bigcup_{v=2}^{L-3}E(\mathcal{R}_{v}\cup\mathcal{Q}_{v}\cup\mathcal{R}_{v+1})\;,
\]
so that we can write the first summation of (\ref{e_mce1}) as
\[
\sum_{v=2}^{L-3}\,\sum_{\{\sigma,\,\xi\}\in E(\mathcal{R}_{v}\cup\mathcal{Q}_{v}\cup\mathcal{R}_{v+1})}\mu_{\beta}(\sigma)\,c_{\beta}(\eta,\,\xi)\,\{f(\xi)-f(\sigma)\}^{2}\;.
\]
This summation can be written as $\sum_{\ell\in\mathbb{T}_{L}}\sum_{k\in\mathbb{T}_{K}}$
of
\begin{align*}
 & \mu_{\beta}(\zeta_{\ell,\,v})\,c_{\beta}(\zeta_{\ell,\,v},\,\zeta_{\ell,\,v\,;\,k,\,1}^{\textup{up}})\,\{f(\zeta_{\ell,\,v\,;\,k,\,1}^{\mathrm{\textup{up}}})-f(\zeta_{\ell,\,v})\}^{2}\\
 & +\sum_{h=1}^{K-2}\mu_{\beta}(\zeta_{\ell,\,v\,;\,k,\,h}^{\textup{up}})\,c_{\beta}(\zeta_{\ell,\,v\,;\,k,\,h}^{\textup{up}},\,\zeta_{\ell,\,v\,;\,k,\,h+1}^{\textup{up}})\,\{f(\zeta_{\ell,\,v\,;\,k,\,h+1}^{\textup{up}})-f(\zeta_{\ell,\,v\,;\,k,\,h}^{\textup{up}})\}^{2}\\
 & +\sum_{h=1}^{K-2}\mu_{\beta}(\zeta_{\ell,\,v\,;\,k,\,h}^{\mathrm{\textup{up}}})\,c_{\beta}(\zeta_{\ell,\,v\,;\,k,\,h}^{\textup{up}},\,\zeta_{\ell,\,v\,;\,k-1,\,h+1}^{\mathrm{\textup{up}}})\,\{f(\zeta_{\ell,\,v\,;\,k-1,\,h+1}^{\mathrm{\textup{up}}})-f(\zeta_{\ell,\,v\,;\,k,\,h}^{\textup{up}})\}^{2}\\
 & +\mu_{\beta}(\zeta_{\ell,\,v\,;\,k,\,K-1}^{\mathrm{\mathrm{\textup{up}}}})\,c_{\beta}(\zeta_{\ell,\,v\,;\,k,\,K-1}^{\mathrm{\mathrm{\textup{up}}}},\,\zeta_{\ell,\,v+1})\,\{f(\zeta_{\ell,\,v+1})-f(\zeta_{\ell,\,v\,;\,k,\,K-1}^{\mathrm{\mathrm{\textup{up}}}})\}^{2}\;,
\end{align*}
and the same form of terms replacing up with down. By (\ref{Zbest}),
(\ref{mudef}), (\ref{muprop}), (\ref{e_testf2.3}), and (\ref{e_testf2.4}),
this equals $2\sum_{\ell\in\mathbb{T}_{L}}\sum_{k\in\mathbb{T}_{K}}$
(where $2$ is multiplied since we have to compute up/down separately)
of
\begin{align*}
 & \frac{e^{-\Gamma\beta}}{Z_{\beta}}\cdot\frac{4\mathfrak{b}^{2}}{[\kappa(K+2)(L-4)]^{2}}+\sum_{h=1}^{K-2}\frac{e^{-\Gamma\beta}}{Z_{\beta}}\cdot\frac{\mathfrak{b}^{2}}{[\kappa(K+2)(L-4)]^{2}}\\
 & +\sum_{h=1}^{K-2}\frac{e^{-\Gamma\beta}}{Z_{\beta}}\cdot\frac{\mathfrak{b}^{2}}{[\kappa(K+2)(L-4)]^{2}}+\frac{e^{-\Gamma\beta}}{Z_{\beta}}\cdot\frac{4\mathfrak{b}^{2}}{[\kappa(K+2)(L-4)]^{2}}\\
= & (1+o_{\beta}(1))\,\frac{e^{-\Gamma\beta}}{2}\,\frac{(2K+4)\mathfrak{b}^{2}}{(K+2)^{2}(L-4)^{2}\mathfrak{\kappa}^{2}}=(1+o_{\beta}(1))\,\frac{\mathfrak{b}^{2}}{(K+2)(L-4)^{2}\kappa^{2}}\,e^{-\Gamma\beta}\;.
\end{align*}
Therefore by (\ref{e_b2}), we can conclude that
\begin{align}
 & \,\sum_{\{\sigma,\,\xi\}\in E(\mathcal{B})}\mu_{\beta}(\sigma)\,c_{\beta}(\sigma,\,\xi)\,\{f(\xi)-f(\sigma)\}^{2}\nonumber \\
= & \,2\,(1+o_{\beta}(1))\sum_{v=2}^{L-3}\sum_{\ell\in\mathbb{T}_{L}}\sum_{k\in\mathbb{T}_{K}}\frac{\mathfrak{b}^{2}}{(K+2)(L-4)^{2}\kappa^{2}}\,e^{-\Gamma\beta}\label{e_mce3}\\
= & \,(1+o_{\beta}(1))\,\frac{2KL(L-4)\mathfrak{b}^{2}}{(K+2)(L-4)^{2}\kappa^{2}}\,e^{-\Gamma\beta}=\frac{\mathfrak{b}+o_{\beta}(1)}{2\kappa^{2}}\,e^{-\Gamma\beta}\;.\nonumber
\end{align}

Next, we calculate the second summation of (\ref{e_mce1}). By (\ref{e_EAdec}),
we rewrite this summation as
\begin{align*}
 & \sum_{\{\sigma_{1},\,\sigma_{2}\}\subseteq\mathcal{O}^{-}}\mu_{\beta}(\sigma_{1})\,c_{\beta}(\sigma_{1},\,\sigma_{2})\,\{f(\sigma_{2})-f(\sigma_{1})\}^{2}\\
 & +\sum_{\sigma_{1}\in\mathcal{O}^{-}}\sum_{\sigma_{2}\in\overline{\mathcal{I}}^{-}}\sum_{\xi\in\mathcal{N}(\sigma_{2})}\mu_{\beta}(\sigma_{1})\,c_{\beta}(\sigma_{1},\,\xi)\,\{f(\xi)-f(\sigma_{1})\}^{2}\;.
\end{align*}
 By Proposition \ref{p_ZA2}, this equals $1+o_{\beta}(1)$ times
\begin{equation}
\Big[\,\sum_{\{\sigma_{1},\,\sigma_{2}\}\subseteq\mathcal{O}^{-}}+\sum_{\sigma_{1}\in\mathcal{O}^{-}}\sum_{\sigma_{2}\in\overline{\mathcal{I}}^{-}}\,\Big]\,\frac{1}{2}\,e^{-\Gamma\beta}\,r^{-}(\sigma_{1},\,\sigma_{2})\,\{f(\sigma_{2})-f(\sigma_{1})\}^{2}\;.\label{e_mce4}
\end{equation}
By (\ref{e_testf2.1}), the last line becomes
\begin{align}
 & \,\frac{\mathfrak{e}^{2}}{\kappa^{2}}\,\sum_{\{\sigma_{1},\,\sigma_{2}\}\subseteq\mathscr{V}^{A}}\frac{1}{2}\,e^{-\Gamma\beta}\,r^{-}(\sigma_{1},\,\sigma_{2})\,\{\mathfrak{h}^{-}(\sigma_{2})-\mathfrak{h}^{-}(\sigma_{1})\}^{2}\nonumber \\
= & \,\frac{e^{-\Gamma\beta}\mathfrak{e}^{2}}{2\kappa^{2}}\,|\mathscr{V}^{-}|\,\mathrm{cap}^{-}(\boxminus,\,\mathcal{R}_{2})=\frac{\mathfrak{e}}{2\kappa^{2}}e^{-\Gamma\beta}\;.\label{e_mce5}
\end{align}
Therefore, we can conclude that
\begin{equation}
\sum_{\{\sigma,\,\xi\}\in E(\mathcal{E}^{-})}\mu_{\beta}(\sigma)\,c_{\beta}(\sigma,\,\xi)\,\{f(\xi)-f(\sigma)\}^{2}=\frac{\mathfrak{e}+o_{\beta}(1)}{2\mathfrak{\kappa}^{2}}e^{-\Gamma\beta}\;.\label{e_mce6-1}
\end{equation}
Similarly, we get
\begin{equation}
\sum_{\{\sigma,\,\xi\}\in E(\mathcal{E}^{+})}\mu_{\beta}(\sigma)\,c_{\beta}(\sigma,\,\xi)\,\{f(\xi)-f(\sigma)\}^{2}=\frac{\mathfrak{e}+o_{\beta}(1)}{2\mathfrak{\kappa}^{2}}e^{-\Gamma\beta}\;.\label{e_mce6}
\end{equation}
Therefore, by (\ref{e_mce1}), (\ref{e_mce3}), (\ref{e_mce6-1}),
and (\ref{e_mce6}), we conclude that the first summation of (\ref{e_decDiri})
equals
\[
\frac{\mathfrak{b}+2\mathfrak{e}+o_{\beta}(1)}{2\kappa^{2}}\,e^{-\Gamma\beta}=\frac{1+o_{\beta}(1)}{2\mathfrak{\kappa}}\,e^{-\Gamma\beta}\;,
\]
as desired.
\end{proof}
\begin{rem}
\label{r_testf2}The estimates (\ref{e_mce3}), (\ref{e_mce6-1}),
and (\ref{e_mce6}) are the reason why we term $\mathfrak{b}$ and
$\mathfrak{e}$ the bulk and edge constants, respectively.
\end{rem}

We now conclude the section with a formal proof of Proposition \ref{p_Capup}.
\begin{proof}[Proof of Proposition \ref{p_Capup} for $d=2$]
Since we have verified in the previous proposition that the function
$f_{0}$ constructed in Definition \ref{d_testf2} satisfies $f_{0}\in\mathfrak{C}(\text{\{\ensuremath{\boxminus\}},\,\{\ensuremath{\boxplus\}}})$
and the energy estimate (\ref{e_Capup}), the proof is completed.
\end{proof}

\section{\label{sec10}Lower Bound for Capacities}

In this section, we construct the test flow $\psi_{0}$ appearing
in Proposition \ref{p_Caplow}. Construction of the test flow will
be given in Section \ref{sec10.1}. Then, two properties of the test
flow appearing in (\ref{e_Caplow1}) are verified in Sections \ref{sec10.2}
and \ref{sec10.4}, respectively. Section \ref{sec10.3} is devoted
to providing some investigations of the equilibrium potential between
$\boxminus$ and $\boxplus$, which will be used in the analyses carried
out in Section \ref{sec10.4}.

\subsection{\label{sec10.1}Construction of test flow}

In this subsection, we explicitly construct a test flow $\psi_{0}$.

We explain the idea before proceeding to the construction. We again
use the convention (\ref{e_shorteq}) in this section. For the edge
typical configurations, recall that the equilibrium potential $\mathfrak{h}^{\pm}(\cdot)$
on $\mathcal{E}^{\pm}$ is the object approximating (up to some rescaling)
the equilibrium potential $h_{\boxminus,\,\boxplus}^{\beta}(\cdot)$.
Hence, we define the test flow on $\mathcal{E}^{\pm}$ as a suitable
modification of (a constant-multiple of) $\Psi_{\mathfrak{h}^{\pm}}$.
For the bulk typical configurations $\mathcal{B}$, we know the typical
behavior of the Metropolis dynamics very well, and hence we can define
$\psi_{0}$ as a simple flow from $\mathcal{R}_{2}$ to $\mathcal{R}_{L-2}$,
where the flow is constant on each edge of the transition.
\begin{defn}[Test flow]
\label{d_testfl2}In this definition, defining $\varphi(\sigma,\,\sigma')=c$
for a flow $\varphi$ implicitly implies that $\varphi(\sigma',\,\sigma)=-c$.
We now construct a flow $\psi_{0}$.
\begin{enumerate}
\item \textbf{Construction of $\psi_{0}$ on edge typical configurations
$\mathcal{E}$. }We provide an explicit construction on $\mathcal{E}^{\pm}$.
\begin{itemize}
\item If $\sigma_{1},\,\sigma_{2}\in\mathcal{O}^{\pm}$ with $\sigma_{1}\sim\sigma_{2}$,
then we set
\begin{equation}
\psi_{0}(\sigma_{1},\,\sigma_{2})=\mathfrak{e}\,r^{\pm}(\sigma_{1},\,\sigma_{2})\,[\mathfrak{h}^{\pm}(\sigma_{1})-\mathfrak{h}^{\pm}(\sigma_{2})]\;.\label{e_testfl2.1}
\end{equation}
\item If $\sigma_{1}\in\mathcal{O}^{\pm}$ and $\sigma_{2}\in\overline{\mathcal{I}}^{\pm}$,
then we set, for all $\xi\in\mathcal{N}(\sigma_{2})$ with $\xi\sim\sigma_{1}$,
\begin{equation}
\psi_{0}(\sigma_{1},\,\xi)=\frac{\mathfrak{e}\,r^{\pm}(\sigma_{1},\,\sigma_{2})\,[\mathfrak{h}^{\pm}(\sigma_{1})-\mathfrak{h}^{\pm}(\sigma_{2})]}{|\{\xi'\in\mathcal{N}(\sigma_{2}):\sigma_{1}\sim\xi'\}|}\;.\label{e_testfl2.2}
\end{equation}
\end{itemize}
\item \textbf{Construction of $\psi_{0}$ on bulk typical configurations
$\mathcal{B}$. }We need to consider the following two cases:
\begin{itemize}
\item For $(k,\,\ell)\in\mathbb{T}_{K}\times\mathbb{T}_{L}$ and $v\in\llbracket2,\,L-3\rrbracket$,
\begin{align*}
\psi_{0}(\zeta_{\ell,\,v\,;\,k,\,0}^{\textup{up}},\,\zeta_{\ell,\,v\,;\,k,\,1}^{\textup{up}}) & =\psi_{0}(\zeta_{\ell,\,v\,;\,k,\,K-1}^{\textup{up}},\,\zeta_{\ell,\,v\,;\,k,\,K}^{\mathrm{\textup{up}}})\\
=\psi_{0}(\zeta_{\ell,\,v\,;\,k,\,0}^{\textup{down}},\,\zeta_{\ell,\,v\,;\,k,\,1}^{\text{down}}) & =\psi_{0}(\zeta_{\ell,\,v\,;\,k,\,K-1}^{\mathrm{\textup{down}}},\,\zeta_{\ell,\,v\,;\,k,\,K}^{\mathrm{\textup{down}}})=\frac{2\mathfrak{b}}{(K+2)(L-4)}\;.
\end{align*}
\item For $(k,\,\ell)\in\mathbb{T}_{K}\times\mathbb{T}_{L}$, $v\in\llbracket2,\,L-3\rrbracket$,
and $h\in\llbracket1,\,K-2\rrbracket$,
\begin{align*}
 & \psi_{0}(\zeta_{\ell,\,v\,;\,k,\,h}^{\textup{up}},\,\zeta_{\ell,\,v\,;\,k,\,h+1}^{\textup{up}})=\psi_{0}(\zeta_{\ell,\,v\,;\,k,\,h}^{\textup{up}},\,\zeta_{\ell,\,v\,;\,k-1,\,h+1}^{\textup{up}})\\
= & \psi_{0}(\zeta_{\ell,\,v\,;\,k,\,h}^{\textup{down}},\,\zeta_{\ell,\,v\,;\,k,\,h+1}^{\textup{down}})=\psi_{0}(\zeta_{\ell,\,v\,;\,k,\,h}^{\textup{down}},\,\zeta_{\ell,\,v\,;\,k-1,\,h+1}^{\textup{down}})=\frac{\mathfrak{b}}{(K+2)(L-4)}\;.
\end{align*}
\end{itemize}
\item We set $\psi_{0}\equiv0$ on all the edges which are not considered
above.
\end{enumerate}
\end{defn}

\subsection{\label{sec10.2}Flow norm}

The next proposition computes the flow norm of $\psi_{0}$ to verify
the first requirement in (\ref{e_Caplow1}). In the remainder of the
current section, we write $\psi=\psi_{0}$ for the simplicity of notation.
\begin{prop}
\label{p_testfl2.1}For the flow $\psi=\psi_{0}$ constructed in Definition
\ref{d_testfl2},
\[
\|\psi\|_{\beta}^{2}=(1+o_{\beta}(1))\,2\kappa e^{\Gamma\beta}\;.
\]
\end{prop}

\begin{proof}
Since the support of $\psi$ is a subset of $\mathcal{E}\cup\mathcal{B}$,
by (\ref{e_decE}), we can write
\begin{equation}
\|\psi\|_{\beta}^{2}=\Big[\,\sum_{\{\sigma,\,\xi\}\in E(\mathcal{E}^{-})}+\sum_{\{\sigma,\,\xi\}\in E(\mathcal{E}^{+})}+\sum_{\{\sigma,\,\xi\}\in E(\mathcal{B})}\,\Big]\,\frac{\psi(\sigma,\,\xi)^{2}}{\mu_{\beta}(\sigma)\,c_{\beta}(\sigma,\,\xi)}\;.\label{fle21}
\end{equation}
By the definition of $\psi$, the first summation of (\ref{fle21})
can be written as
\[
\sum_{\{\sigma_{1},\,\sigma_{2}\}\in E(\mathcal{O}^{-})}\frac{\psi(\sigma_{1},\,\sigma_{2})^{2}}{\mu_{\beta}(\sigma_{1})\,c_{\beta}(\sigma_{1},\,\sigma_{2})}+\sum_{\sigma_{1}\in\mathcal{O}^{-}}\sum_{\sigma_{2}\in\mathcal{I}^{-}}\,\sum_{\xi\in\mathcal{N}(\sigma_{2}):\sigma_{1}\sim\xi}\frac{\psi(\sigma_{1},\,\xi)^{2}}{\mu_{\beta}(\sigma_{1})\,c_{\beta}(\sigma_{1},\,\xi)}\;.
\]
By (\ref{e_testfl2.1}), (\ref{e_testfl2.2}), and Proposition \ref{p_ZA2},
this equals $(1+o_{\beta}(1))$ times
\[
\Big[\,\sum_{\{\sigma_{1},\,\sigma_{2}\}\in E(\mathcal{O}^{-})}+\sum_{\sigma_{1}\in\mathcal{O}^{-}}\sum_{\sigma_{2}\in\mathcal{I}^{-}}\,\Big]\,\frac{2\mathfrak{e}^{2}r^{-}(\sigma_{1},\,\sigma_{2})\,\{\mathfrak{h}^{-}(\sigma_{2})-\mathfrak{h}^{-}(\sigma_{1})\}^{2}}{e^{-\Gamma\beta}}\;.
\]
By the definition of capacity, we can rewrite the last summation as
\[
2\mathfrak{e}^{2}\,e^{\Gamma\beta}\,|\mathscr{V}^{-}|\,\mathrm{cap}^{-}(\boxminus,\,\mathcal{R}_{2})=2\mathfrak{e}\,e^{\Gamma\beta}\;.
\]
Since we can apply a similar argument to the second summation of (\ref{fle21}),
we can conclude that
\begin{align}
\Big[\,\sum_{\{\sigma,\,\xi\}\in E(\mathcal{E}^{-})}+\sum_{\{\sigma,\,\xi\}\in E(\mathcal{E}^{+})}\,\Big]\,\frac{\psi(\sigma,\,\xi)^{2}}{\mu_{\beta}(\sigma)\,c_{\beta}(\sigma,\,\xi)} & =(1+o_{\beta}(1))\times2\times2\mathfrak{e}\,e^{\Gamma\beta}\nonumber \\
 & =(4\mathfrak{e}+o_{\beta}(1))\,e^{\Gamma\beta}\;.\label{fle22}
\end{align}

Now, we consider the third summation of (\ref{fle21}). By definition,
this summation is $\sum_{k\in\mathbb{T}_{K},\,\ell\in\mathbb{T}_{L}}\sum_{v=2}^{L-3}$
of
\begin{align*}
 & \Big[\,\frac{\psi(\zeta_{\ell,\,v\,;\,k,\,0}^{\textup{up}},\,\zeta_{\ell,\,v\,;\,k,\,1}^{\textup{up}})^{2}}{\mu_{\beta}(\zeta_{\ell,\,v\,;\,k,\,0}^{\textup{up}})\,c_{\beta}(\zeta_{\ell,\,v\,;\,k,\,0}^{\textup{up}},\,\zeta_{\ell,\,v\,;\,k,\,1}^{\textup{up}})}+\frac{\psi(\zeta_{\ell,\,v\,;\,k,\,K-1}^{\textup{up}},\,\zeta_{\ell,\,v\,;\,k,\,K}^{\textup{up}})^{2}}{\mu_{\beta}(\zeta_{\ell,\,v\,;\,k,\,K-1}^{\textup{up}})\,c_{\beta}(\zeta_{\ell,\,v\,;\,k,\,K-1}^{\textup{up}},\,\zeta_{\ell,\,v\,;\,k,\,K}^{\textup{up}})}\,\Big]\\
 & +\sum_{h=1}^{K-2}\Big[\,\frac{\psi(\zeta_{\ell,\,v\,;\,k,\,h}^{\mathrm{\textup{up}}},\,\zeta_{\ell,\,v\,;\,k,\,h+1}^{\mathrm{\textup{up}}})^{2}}{\mu_{\beta}(\zeta_{\ell,\,v\,;\,k,\,h}^{\mathrm{\textup{up}}})\,c_{\beta}(\zeta_{\ell,\,v\,;\,k,\,h}^{\mathrm{\textup{up}}},\,\zeta_{\ell,\,v\,;\,k,\,h+1}^{\mathrm{\textup{up}}})}+\frac{\psi(\zeta_{\ell,\,v\,;\,k,\,h}^{\mathrm{\textup{up}}},\,\zeta_{\ell,\,v\,;\,k-1,\,h+1}^{\textup{up}})^{2}}{\mu_{\beta}(\zeta_{\ell,\,v\,;\,k,\,h}^{\mathrm{\textup{up}}})\,c_{\beta}(\zeta_{\ell,\,v\,;\,k,\,h}^{\textup{up}},\,\zeta_{\ell,\,v\,;\,k-1,\,h+1}^{\textup{up}})}\,\Big]\;,
\end{align*}
and the same-form of terms can be obtained from above by replacing
up with down.

By the definition of $\psi$, (\ref{muprop}), and (\ref{Zbest}),
this expression equals $(1+o_{\beta}(1))$ times
\begin{align*}
\Big[\,\frac{32\mathfrak{b}^{2}e^{\Gamma\beta}}{(K+2)^{2}(L-4)^{2}}+\sum_{h=1}^{K-2}\frac{8\mathfrak{b}^{2}e^{\Gamma\beta}}{(K+2)^{2}(L-4)^{2}}\,\Big] & =(32+8(K-2))\,\frac{\mathfrak{b}^{2}e^{\Gamma\beta}}{(K+2)^{2}(L-4)^{2}}\\
 & =\frac{8\mathfrak{b}^{2}e^{\Gamma\beta}}{(K+2)(L-4)^{2}}\;.
\end{align*}
Hence, by the definition of $\mathfrak{b}$, the third summation of
(\ref{fle21}) equals
\begin{equation}
(1+o_{\beta}(1))\times KL(L-4)\times\frac{8\mathfrak{b}^{2}e^{\Gamma\beta}}{(K+2)(L-4)^{2}}=[2\mathfrak{b}+o_{\beta}(1)]\,e^{\Gamma\beta}\;.\label{fle23}
\end{equation}
Therefore, by (\ref{fle21}), (\ref{fle22}), and (\ref{fle23}),
we get
\[
\left\Vert \psi\right\Vert _{\beta}^{2}=2[\mathfrak{b}+2\mathfrak{e}+o_{\beta}(1)]\,e^{\Gamma\beta}=(1+o_{\beta}(1))\,2\kappa e^{\Gamma\beta}\;.
\]
This finishes the proof.
\end{proof}

\subsection{\label{sec10.3}Equilibrium potential around ground states}

It remains to verify the second requirement (\ref{e_Caplow1}) regarding
the test flow $\psi_{0}$. To this end, we first prove that the equilibrium
potential is nearly constant on the neighborhood of ground states
in this subsection. The main tool is Proposition \ref{p_eqp_bd} regarding
the estimate of the equilibrium potential.
\begin{lem}
\label{l_eqpot2}It holds that
\[
\max_{\sigma\in\mathcal{N}(\boxplus)}h_{\boxminus,\,\boxplus}^{\beta}(\sigma)=O(e^{-\beta})\;\;\;\;\text{and\;\;\;\;}\max_{\sigma\in\mathcal{N}(\boxminus)}(1-h_{\boxminus,\,\boxplus}^{\beta}(\sigma))=O(e^{-\beta})\;.
\]
\end{lem}

\begin{proof}
We prove the lemma only for the first estimate, because the second
one follows immediately from the first since $1-h_{\boxminus,\,\boxplus}^{\beta}=h_{\boxplus,\,\boxminus}^{\beta}$.

By Propositions \ref{p_eqp_bd} and \ref{prop:inc}, it holds that
\begin{equation}
h_{\boxminus,\,\boxplus}^{\beta}(\sigma)\le\frac{\mathrm{cap}_{\beta}(\sigma,\,\boxminus)}{\mathrm{cap}_{\beta}(\sigma,\,\{\boxminus,\,\boxplus\})}\le\frac{\mathrm{cap}_{\beta}(\sigma,\,\boxminus)}{\mathrm{cap}_{\beta}(\sigma,\,\boxplus)}\;.\label{e_eqpot2.1}
\end{equation}
Now, we estimate $\mathrm{cap}_{\beta}(\sigma,\,\boxplus)$ and $\mathrm{cap}_{\beta}(\sigma,\,\boxminus)$
separately.

We first give a lower bound of $\mathrm{cap}_{\beta}(\sigma,\,\boxplus)$
via the Thomson principle (Theorem \ref{t_TP_rev}). As $\sigma\in\mathcal{N}(\boxplus)$,
there exists a $(\Gamma-1)$-path $(\omega_{t})_{t=0}^{T}$ connecting
$\boxplus$ and $\sigma$, where $T$ is bounded by a constant depending
only on $K$ and $L$. We define a test flow $\phi$ on $\mathcal{X}$
by
\[
\phi(\omega_{t},\,\omega_{t+1})=-\phi(\omega_{t+1},\,\omega_{t})=1\text{ for }t\in\llbracket0,\,T-1\rrbracket\;,
\]
and $\phi=0$ otherwise. This construction implies that $\phi$ is
a unit flow from $\{\boxplus\}$ to $\{\sigma\}$. Since $(\omega_{t})_{t=0}^{T}$
is a $(\Gamma-1)$-path, by Proposition \ref{p_inv} and (\ref{muprop}),
\[
\mu_{\beta}(\omega_{t})\,c_{\beta}(\omega_{t},\,\omega_{t+1})=\min\,\{\mu_{\beta}(\omega_{t}),\,\mu_{\beta}(\omega_{t+1})\}\le\frac{1+o_{\beta}(1)}{2}\,e^{-(\Gamma-1)\beta}\;.
\]
Therefore, we obtain
\[
\|\phi\|_{\beta}^{2}=\sum_{t=0}^{T-1}\frac{\phi(\omega_{t},\,\omega_{t+1})^{2}}{\mu_{\beta}(\omega_{t})\,c_{\beta}(\omega_{t},\,\omega_{t+1})}\le\sum_{t=0}^{T-1}\frac{q+o_{\beta}(1)}{e^{-(\Gamma-1)\beta}}\le Ce^{(\Gamma-1)\beta}\;.
\]
Hence, by Theorem \ref{t_TP_rev},
\begin{equation}
\mathrm{cap}_{\beta}(\sigma,\,\boxplus)\ge\frac{1}{\|\phi\|_{\beta}^{2}}\ge\frac{1}{C}\,e^{-(\Gamma-1)\beta}\;.\label{e_eqpot2.2}
\end{equation}

Next, we establish an upper bound for $\mathrm{cap}_{\beta}(\sigma,\,\boxminus)$.
To this end, we first observe from our construction of $f_{0}$ (cf.
Definition \ref{d_testf2}) that $f_{0}\in\mathfrak{C}_{1,\,0}(\mathcal{N}(\boxminus),\,\mathcal{N}(\boxplus))$.
Therefore, by the symmetry of capacities (cf. (\ref{e_cap_synm})),
the monotonicity of capacities (cf. Proposition \ref{prop:inc}),
and the Dirichlet principle (cf. Theorem \ref{t_DP_rev}), we have
\begin{equation}
\mathrm{cap}_{\beta}(\sigma,\,\boxminus)=\mathrm{cap}_{\beta}(\boxminus,\,\sigma)\le\mathrm{cap}_{\beta}(\mathcal{N}(\boxminus),\,\mathcal{N}(\boxplus))\le\mathscr{D}_{\beta}(f_{0})\le C\,e^{-\Gamma\beta}\label{e_eqpot2.3}
\end{equation}
for some constant $C>0$, where the last bound follows from Proposition
\ref{p_testf2}.

The proof is completed by (\ref{e_eqpot2.1}), (\ref{e_eqpot2.2}),
and (\ref{e_eqpot2.3}).
\end{proof}

\subsection{\label{sec10.4}Divergence of test flow}

Now, we investigate the divergence of the test flow $\psi_{0}$ constructed
in Definition \ref{d_testfl2}. For simplicity, we again write $\psi=\psi_{0}$
throughout the current subsection. We first check that this flow is
divergence-free on bulk typical configurations.
\begin{lem}
\label{l_fl2bulk1}We have \textbf{$(\mathrm{div}\,\psi)(\sigma)=0$}
for all $\sigma\in\mathcal{B}\setminus\mathcal{E}$.
\end{lem}

\begin{proof}
Let us fix $\sigma\in\mathcal{B}\setminus\mathcal{E}$.

If $\sigma=\zeta_{\ell,\,v}$ for some $\ell\in\mathbb{T}_{L}$ and
$v\in\llbracket3,\,L-3\rrbracket$, then we can write $(\mathrm{div}\,\psi)(\sigma)$
as
\[
\sum_{k\in\mathbb{T}_{K}}[\psi(\sigma,\,\zeta_{\ell,\,v\,;\,k,\,1}^{\textup{up}})+\psi(\sigma,\,\zeta_{\ell,\,v\,;\,k,\,1}^{\textup{down}})+\psi(\sigma,\,\zeta_{\ell,\,v-1\,;\,k,\,K-1}^{\textup{up}})+\psi(\sigma,\,\zeta_{\ell+1,\,v-1\,;\,k,\,K-1}^{\textup{down}})]\;.
\]
By recalling Definition \ref{d_testfl2}, this summation is equal
to
\[
\sum_{k\in\mathbb{T}_{K}}\Big[\,\frac{2\mathfrak{b}}{(K+2)(L-4)}+\frac{2\mathfrak{b}}{(K+2)(L-4)}-\frac{2\mathfrak{b}}{(K+2)(L-4)}-\frac{2\mathfrak{b}}{(K+2)(L-4)}\,\Big]=0\;.
\]

If $\sigma=\zeta_{\ell,\,v\,;\,k,\,h}^{+}$ for some $(k,\,\ell)\in\mathbb{T}_{K}\times\mathbb{T}_{L}$,
$v\in\llbracket2,\,L-3\rrbracket$, and $h\in\llbracket1,\,K-1\rrbracket$,
then we can write $(\mathrm{div}\,\psi)(\sigma)$ as
\begin{align*}
 & \phi(\sigma,\,\zeta_{\ell,\,v\,;\,k,\,h+1}^{\textup{up}})+\phi(\sigma,\,\zeta_{\ell,\,v\,;\,k-1,\,h+1}^{\textup{up}})+\phi(\sigma,\,\zeta_{\ell,\,v\,;\,k,\,h-1}^{\textup{up}})+\phi(\sigma,\,\zeta_{\ell,\,v\,;\,k+1,\,h-1}^{\textup{up}})\\
 & =\frac{\mathfrak{b}}{(K+2)(L-4)}+\frac{\mathfrak{b}}{(K+2)(L-4)}-\frac{\mathfrak{b}}{(K+2)(L-4)}-\frac{\mathfrak{b}}{(K+2)(L-4)}=0\;.
\end{align*}
The cases $\sigma=\zeta_{\ell,\,v\,;\,k,\,h}^{\textup{up}}$ and $\zeta_{\ell,\,v\,;\,k,\,h}^{\textup{down}}$
can be handled in the same manner. This concludes the proof.
\end{proof}
Next, we show that $\psi$ is divergence-free on $\mathcal{R}_{2}$
and $\mathcal{R}_{L-2}$.
\begin{lem}
\label{l_fl2bulk2}It holds that \textbf{$(\mathrm{div}\,\psi)(\sigma)=0$}
for all $\sigma\in\mathcal{R}_{2}\cup\mathcal{R}_{L-2}$.
\end{lem}

\begin{proof}
We only consider the divergence on $\mathcal{R}_{2}$, since the proof
for $\mathcal{R}_{L-2}$ is identical. Recall the generator $L^{-}$
from (\ref{e_genZA}). Then, since the uniform measure on $\mathscr{V}^{-}$
is the invariant measure for the process $Z^{-}(\cdot)$, by the expression
(\ref{e_l241}) of capacity,
\begin{equation}
\mathrm{cap}^{-}(\boxminus,\,\mathcal{R}_{2})=-\frac{1}{|\mathscr{V}^{-}|}\,\sum_{\sigma\in\mathcal{R}_{2}}\,\sum_{\xi\in\mathscr{V}^{-}\setminus\mathcal{R}_{2}}r^{-}(\sigma,\,\xi)\,\{\mathfrak{h}^{-}(\sigma)-\mathfrak{h}^{-}(\xi)\}\;.\label{ee31}
\end{equation}
On the other hand, by the definition of $\psi$, we can write
\begin{equation}
\sum_{\sigma\in\mathcal{R}_{2}}\,\sum_{\xi\in\mathcal{E}^{-}}\psi(\sigma,\,\xi)=\mathfrak{e}\sum_{\sigma\in\mathcal{R}_{2}}\,\sum_{\xi\in\mathscr{V}^{-}\setminus\mathcal{R}_{2}}r^{-}(\sigma,\,\xi)\,\{\mathfrak{h}^{-}(\sigma)-\mathfrak{h}^{-}(\xi)\}\;.\label{ee32}
\end{equation}
By (\ref{ee31}) and (\ref{ee32}), we get
\begin{equation}
\sum_{\sigma\in\mathcal{R}_{2}}\,\sum_{\xi\in\mathcal{E}^{-}}\psi(\sigma,\,\xi)=-\mathfrak{e}\,|\mathscr{V}^{-}|\,\mathrm{cap}^{-}(\boxminus,\,\mathcal{R}_{2})=-1\;,\label{e_fl2bulk2.1}
\end{equation}
where the second identity follows from the definition of $\mathfrak{e}$.
On the other hand, by the definition of $\psi$,
\begin{equation}
\sum_{\sigma\in\mathcal{R}_{2}}\,\sum_{\xi\in\mathcal{B}}\psi(\sigma,\,\xi)=L\times2K\times\frac{2\mathfrak{b}}{(K+2)(L-4)}=1\;.\label{e_fl2bulk2.2}
\end{equation}
By adding (\ref{e_fl2bulk2.1}) and (\ref{e_fl2bulk2.2}), we obtain
\[
\sum_{\sigma\in\mathcal{R}_{2}}(\textup{div}\,\psi)(\sigma)=0\;.
\]
Since $\textup{div}\,\psi$ is a constant function on $\mathcal{R}_{2}$
by symmetry, we can conclude that $(\textup{div}\,\psi)(\sigma)=0$
for all $\sigma\in\mathcal{R}_{2}$.
\end{proof}
Next, we show that the flow $\psi$ is divergence-free on $\mathcal{O}^{-}$
and $\mathcal{O}^{+}$.
\begin{lem}
\label{l_fl2edge}We have \textbf{$(\mathrm{div}\,\psi)(\sigma)=0$}
for all $\sigma\in\mathcal{O}^{-}$ and $\sigma\in\mathcal{O}^{+}$.
\end{lem}

\begin{proof}
We only consider the case $\sigma\in\mathcal{O}^{-}$ since the case
$\sigma\in\mathcal{O}^{+}$ can be handled in the same manner. By
the definition of $\psi$, we can write
\begin{equation}
(\mathrm{div}\,\psi)(\sigma)=\mathfrak{e}\,(L^{-}\mathfrak{h}^{-})(\sigma)\;.\label{e_fl2edge}
\end{equation}
By (\ref{e_eq_prop}), we can conclude from (\ref{e_fl2edge}) that
$(\mathrm{div}\,\psi)(\sigma)=0$ for all $\sigma\in\mathscr{V}^{A}\setminus(\boxminus\cup\mathcal{R}_{2})$.
It suffices to observe from (\ref{e_IArep}) and (\ref{e_EAdec})
that $\mathscr{V}^{-}\setminus(\boxminus\cup\mathcal{R}_{2})=\mathcal{O}^{-}$.
\end{proof}
We can conclude from Lemmas \ref{l_fl2bulk1}, \ref{l_fl2bulk2},
and \ref{l_fl2edge} that the flow $\psi$ is almost a divergence-free
flow.
\begin{prop}
\label{p_fl2}The flow $\psi$ is divergence-free on $\mathcal{X}\setminus\mathcal{N}(\mathcal{S})$.
\end{prop}

\begin{proof}
We can decompose the set $\mathcal{X}\setminus\mathcal{N}(\mathcal{S})$
as
\[
(\mathcal{B}\setminus\mathcal{E})\cup\mathcal{R}_{2}\cup\mathcal{R}_{L-2}\cup\mathcal{O}^{-}\cup\mathcal{O}^{+}\cup\big(\,\mathcal{X}\setminus(\mathcal{E}\cup\mathcal{B})\,\big)\;.
\]
Since it follows immediately from the definition that $(\textup{div}\,\psi)\equiv0$
on $\mathcal{X}\setminus(\mathcal{E}\cup\mathcal{B})$, we can conclude
the proof from Lemmas \ref{l_fl2bulk1}, \ref{l_fl2bulk2}, and \ref{l_fl2edge}.
\end{proof}
Now, we are ready to prove the second requirement of (\ref{e_Caplow1}).
\begin{prop}
\label{p_testfl2.2}We have that
\[
\sum_{\sigma\in\mathcal{X}}h_{\boxminus,\,\boxplus}^{\beta}(\sigma)\,(\mathrm{div}\,\psi)(\sigma)=1+o_{\beta}(1)\;.
\]
\end{prop}

\begin{proof}
In view of Proposition \ref{p_fl2}, it suffices to prove that
\begin{equation}
\sum_{\sigma\in\mathcal{N}(\boxminus)}h_{\boxminus,\,\boxplus}^{\beta}(\sigma)\,(\mathrm{div}\,\psi)(\sigma)=1+o_{\beta}(1)\;\;\;\text{and\;\;\;}\sum_{\sigma\in\mathcal{N}(\boxplus)}h_{\boxminus,\,\boxplus}^{\beta}(\sigma)\,(\mathrm{div}\,\psi)(\sigma)=o_{\beta}(1)\;.\label{e_flr1}
\end{equation}
We focus only on the former, since the proof for the latter is essentially
identical. By Lemma \ref{l_eqpot2}, we have
\begin{align*}
\sum_{\sigma\in\mathcal{N}(\boxminus)}h_{\boxminus,\,\boxplus}^{\beta}(\sigma)\,(\mathrm{div}\,\psi)(\sigma) & =(1+o_{\beta}(1))\sum_{\sigma\in\mathcal{N}(\boxminus)}\,\sum_{\xi:\,\xi\sim\sigma}\psi(\sigma,\,\xi)\\
 & =(1+o_{\beta}(1))\sum_{\sigma\in\mathcal{N}(\boxminus)}(\textup{div}\,\psi)(\sigma)\;.
\end{align*}
Note that, in the previous computation, we implicitly used the fact
that neither $\psi$ nor $\mathcal{N}(\boxminus)$ depends on $\beta$.

Now, we claim that
\begin{equation}
\sum_{\sigma\in\mathcal{N}(\boxminus)}(\textup{div}\,\psi)(\sigma)=1\;.\label{e_flr3}
\end{equation}
By (\ref{e_testfl2.2}), we can rewrite the left-hand side of the
previous identity as
\[
\sum_{\sigma\in\mathcal{N}(\boxminus)}\,\sum_{\xi\in\mathcal{O}^{-}:\,\xi\sim\sigma}\psi(\sigma,\,\xi)=-\sum_{\sigma\in\mathcal{N}(\boxminus)}\,\sum_{\xi\in\mathcal{O}^{-}:\,\xi\sim\sigma}\frac{\mathfrak{e}\,r^{-}(\xi,\,\boxminus)\,[\mathfrak{h}^{-}(\xi)-\mathfrak{h}^{-}(\boxminus)]}{|\{\xi'\in\mathcal{N}(\boxminus):\xi\sim\xi'\}|}\;.
\]
Since $r^{-}(\cdot,\,\cdot)$ is symmetric, we can further rewrite
as
\[
-\mathfrak{e}\sum_{\xi\in\mathcal{O}^{-}:\{\xi,\,\boxminus\}\in\mathscr{E}^{A}}r^{-}(\boxminus,\,\xi)\,[\mathfrak{h}^{-}(\xi)-\mathfrak{h}^{-}(\boxminus)]\;.
\]
By (\ref{e_l241}) and the definition of $\mathfrak{e}$, the last
display equals
\[
\mathfrak{e}\,|\mathscr{V}^{-}|\,\mathrm{cap}^{-}(\boxminus,\,\mathcal{R}_{2})=1\;.
\]
This completes the proof for the first estimate of (\ref{e_flr1})
and concludes the proof.
\end{proof}
We conclude this section with the proof of Proposition \ref{p_Caplow}.\textbf{ }
\begin{proof}[Proof of Proposition \ref{p_Caplow}]
Let $\psi_{0}$ be the test flow defined in Definition \ref{d_testfl2}.
Then, the two properties appearing in (\ref{e_Caplow1}) for $\psi_{0}$
have been verified in Propositions \ref{p_testfl2.1} and \ref{p_testfl2.2}.
This completes the proof of Proposition \ref{p_Caplow}.
\end{proof}

\section{\label{sec_K=00003DL}Comments on Case $K=L$}

Now, we suppose that $K=L$. We define $\theta:\mathbb{T}_{K}^{2}\rightarrow\mathbb{T}_{K}^{2}$
as
\[
\theta(k,\,\ell)=(\ell,\,k)\;\;\;\;;\;(k,\,\ell)\in\mathbb{T}_{K}^{2}\;.
\]
Then, define an operator $\Theta:\mathcal{X}\rightarrow\mathcal{X}$
as, for $\sigma\in\mathcal{X}$,
\[
\Theta(\sigma)(x)=\sigma(\theta(x))\;\;\;\;;\;x\in\mathbb{T}_{K}^{2}\;.
\]
Then, the collection of canonical configurations should be $\mathcal{C}\cup\Theta(\mathcal{C})$.
Similarly, the definitions of bulk typical configurations and edge
typical configurations should be extended to $\mathcal{B}\cup\Theta(\mathcal{B})$
and $\mathcal{E}\cup\Theta(\mathcal{E})$, respectively. With these
new definitions of canonical and typical configurations, we can perform
similar computations to prove the Eyring--Kramers law.

\newpage

\part{\label{pt4}Condensing Zero-range Processes}

In this third part of the lecture note, we consider a class of interacting
particle systems known as the zero-range processes. The particles
comprising this model are sticky and therefore tend to condensed at
a site. The movements of this condensate are the metastable behavior
of this model. To precisely understand the successive movements of
the condensate, we use the Markov chain model reduction technique
in the context of the metastability to analyze this model. According
to the general methodology known as the martingale approach developed
in \cite{B-L TM,B-L TM2,B-L MG}, the proof of the Markov chain model
reduction for metastable Markov processes is largely based on the
potential theory.

This connection between the Markov chain model reduction and the potential
theory is relatively clear if the underlying model is reversible.
On the other hand, if the model is non-reversible, not only the estimates
of the capacity but also deriving the Markov chain model reduction
from such estimates are complicated.

In this part, we will try to explain the general method for carrying
out these tasks as clearly as possible. We will use the generalized
Dirichlet and Thomson principles for the non-reversible Markov processes
(cf. Theorem \ref{t_GTP}) to derive sharp estimates of capacities
between metastable sets, and then use a robust method developed in
\cite{L-Seo NR1} to derive the Markov chain model reduction from
there.

We note that the current part is largely based on the article \cite{Seo NRZRP}.

\section{Zero-range processes\label{sec2}}

In this section, we introduce a class of zero-range processes exhibiting
the condensation phenomenon.

\subsubsection*{Underlying random walk}

A zero-range process is a system of interacting particles. We start
by explaining the dynamics of the underlying particles comprising
the zero-range process. Let $\kappa\ge2$ be an integer and denote
by
\[
S=\mathbb{T}_{\kappa}=\mathbb{Z}/\kappa\mathbb{Z}
\]
the cycle of length $\kappa$. Denote by $\mathbb{X}(\cdot)$ the
continuous-time Markov process on $S$ with rate
\[
r(x,\,y)=\begin{cases}
p & \text{if }y=x+1\;,\\
1-p & \text{if }y=x-1\;,\\
0 & \text{otherwise\;.}
\end{cases}
\]
We note that the addition and subtraction in $\mathbb{T}_{\kappa}$
are always carried out modulo $\kappa$. We denote by $L_{\mathbb{X}}$
and $D_{\mathbb{X}}$ the generator and Dirichlet form associated
with the process $\mathbb{X}(\cdot)$. We note that the potential
theory of the process $\mathbb{X}(\cdot)$ has been analyzed in Exercise
\ref{ex26}. We denote by $\mathbf{P}_{x}$, $x\in S$, the law of
the underlying Markov process $\mathbb{X}(\cdot)$ starting from a
site $x\in S$.

\subsubsection*{Zero-range processes}

The zero-range process is defined as an interacting system of $N$
particles, where particles basically follow the law of the process
$\mathbb{X}(\cdot)$ defined above, but interact through the zero-range
interaction explained below.

Let $a:\mathbb{N}\rightarrow\mathbb{R}$ and $g:\mathbb{N}\rightarrow\mathbb{R}$
(with the convention $\mathbb{N}=\{0,\,1,\,2,\,\cdots\}$) be functions
defined by
\begin{equation}
a(n)=\begin{cases}
1 & \mbox{if }n=0\;,\\
n^{\alpha} & \mbox{if }n\ge1\;,
\end{cases}\label{def_a}
\end{equation}
and
\begin{equation}
g(n)=\begin{cases}
0 & \mbox{if }n=0\;,\\
a(n)/a(n-1) & \mbox{if }n\ge1\;,
\end{cases}\label{def_g}
\end{equation}
where the parameter $\alpha$ stands for the stickiness of constituent
particles. We assume that $\alpha>1$ in this note. We will discuss
this assumption for $\alpha$ in Remark \ref{rem:alpha}.

For $N\in\mathbb{N}$, define $\mathcal{\mathcal{H}}_{N}\subset\mathbb{N}^{S}$
as the space of configurations on $S$ with $N$ particles:
\[
\mathcal{\mathcal{H}}_{N}=\Bigl\{\,\eta=(\eta_{x})_{x\in S}\in\mathbb{N}^{S}:\sum_{x\in S}\eta_{x}=N\,\Bigr\}\;.
\]
Here, $\eta\in\mathbb{N}^{S}$ represents the entire set of particle
configurations on $S$ and $\eta_{x}$, $x\in S$, represents the
number of particles at $x$.

Now we are ready to define the zero-range process. For $N\in\mathbb{N}$,
the zero-range process $\{\eta_{N}(t):t\ge0\}$ consisting of $N$
particles is defined as a continuous-time Markov process on $\mathcal{H}_{N}$
associated with the generator
\[
(\mathscr{L}_{N}\mathbf{f})(\eta)=\sum_{x,\,y\in S}g(\eta_{x})r(x,\,y)(\mathbf{f}(\sigma^{x,\,y}\eta)-\mathbf{f}(\eta))\;\;\;\;;\;\eta\in\mathcal{H}_{N}\;,
\]
for $\mathbf{f}:\mathcal{H}_{N}\rightarrow\mathbb{R}$, where $\sigma^{x,\,y}\eta\in\mathcal{H}_{N}$
represents the configuration obtained from $\eta$ by sending a particle
at site $x$ to $y$ (if possible), that is, $\sigma^{x,\,y}\eta=\eta$
if $\eta_{x}=0$, and
\[
(\sigma^{x,\,y}\eta)_{z}=\begin{cases}
\eta_{z}-1 & \mbox{if }z=x\\
\eta_{z}+1 & \mbox{if }z=y\\
\eta_{z} & \mbox{otherwise\;,}
\end{cases}
\]
if $\eta_{x}\ge1$. Of course, we have $\sigma^{x,\,x}\eta=\eta$
for all $x\in S$ and $\eta\in\mathcal{H}_{N}.$ For $\eta\in\mathcal{H}_{N}$,
denote by $\mathbb{P}_{\eta}^{N}$ the law of the zero-range process
$\eta_{N}(\cdot)$ starting from $\eta$, and denote by $\mathbf{\mathbb{E}}_{\eta}^{N}$
the corresponding expectation.
\begin{notation}
A function on $\mathcal{H}_{N}$ will always be denoted by bold font
such as $\mathbf{f}$ or $\mathbf{g}$ to distinguish such functions
from functions on $S$.
\end{notation}

Heuristically, under the zero-range dynamics defined above, one of
the particles at site $x$ jumps to site $y$ at a rate $g(\eta_{x})r(x,\,y)$.
We can observe two important features of the dynamics at this point.
Firstly, since the rate $g(\eta_{x})r(x,\,y)$ is independent of $\eta_{z}$,
$z\neq x$, we can observe that each particle interacts only with
the particles at the same site through the function $g(\cdot)$. This
is the reason that this interacting particle system is called a zero-range
process.

Secondly, in view of (\ref{def_a}) and (\ref{def_g}), this jump
rate $g(\eta_{x})r(x,\,y)$ decreases as $\eta_{x}$($\ge2$) becomes
larger. Namely, a particle is deactivated as there are more particles
grouped together with that particle. For this reason, we can observe
that particles of the zero-range process are sticky. This sticky interaction
eventually causes the condensation of particles as defined in the
next section.
\begin{xca}
\label{ex121}Prove that the zero-range process defined above is irreducible.
\end{xca}

\subsubsection*{Invariant measure }

For $\eta\in\mathcal{H}_{N}$, let us write
\begin{equation}
a(\eta)=\prod_{x\in S}a(\eta_{x})\;.\label{rec1}
\end{equation}
By Exercise \ref{ex121}, the zero-range process has a unique invariant
measure. One can readily verify that this invariant measure $\mu_{N}(\cdot)$
on $\mathcal{H}_{N}$ is given by
\begin{equation}
\mu_{N}(\eta)=\frac{N^{\alpha}}{Z_{N}}\frac{1}{a(\eta)}\;\;\;\;;\;\eta\in\mathcal{H}_{N}\;,\label{inv_zrp}
\end{equation}
where $Z_{N}$ is the partition function turning $\mu_{N}$ into a
probability measure, i.e.,
\[
Z_{N}=N^{\alpha}\sum_{\eta\in\mathcal{H}_{N}}\frac{1}{a(\eta)}\;.
\]

\begin{xca}
\label{ex_zepinv}
\begin{enumerate}
\item Prove that $\mu_{N}(\cdot)$ is the invariant measure for the zero-range
process $\eta_{N}(\cdot)$.
\item Prove that the zero-range process $\eta_{N}(\cdot)$ is reversible
if and only if $p=1/2$.
\end{enumerate}
\end{xca}

Define

\[
\Gamma_{\alpha}=\sum_{n=0}^{\infty}\frac{1}{a(n)}=1+\sum_{n=1}^{\infty}\frac{1}{n^{\alpha}}<\infty\;,
\]
 where the last inequality holds since we have assumed that $\alpha>1$.
Then, define
\begin{equation}
Z=\kappa\Gamma_{\alpha}^{\kappa-1}\;.\label{e20}
\end{equation}
The following proposition explains the appearance of the somewhat
unnecessary $N^{\alpha}$ term at (\ref{inv_zrp}).
\begin{prop}
\label{p_zn}We have that
\[
\lim_{N\rightarrow\infty}Z_{N}=Z\;.
\]
\end{prop}

Since our primary concern is the connection between the potential
theory and the metastability of the zero-range processes, we shall
not prove all the detailed properties of the zero-range processes.
Instead, we refer to \cite{B-L ZRP} for the proof. For this proposition,
we refer to \cite[Proposition 2.1]{B-L ZRP} for the proof.

\subsubsection*{Dirichlet form}

We write $\mathscr{D}_{N}(\mathbf{f})$, $\mathbf{f}:\mathcal{H}_{N}\rightarrow\mathbb{R}$,
the Dirichlet form associated with the zero-range process $\eta_{N}(\cdot)$,
i.e.,
\[
\mathscr{D}_{N}(\mathbf{f})=\left\langle \mathbf{f},\,-\mathscr{L}_{N}\mathbf{f}\right\rangle _{\mu_{N}}\;.
\]
By summation by parts, we can rewrite this Dirichlet form as
\[
\mathscr{D}_{N}(\mathbf{f})=\frac{1}{2}\sum_{x\in S}\sum_{y\in S}\mu_{N}(\eta)\,g(\eta_{x})\,r(x,\,y)\left[\mathbf{f}(\sigma^{x,\,y}\eta)-\mathbf{f}(\eta)\right]^{2}\;.
\]

\subsubsection*{Equilibrium potentials and capacities}

In the investigation of the zero-range process, both the potential
theories of the underlying random walk $\mathbb{X}(\cdot)$ and of
the zero-range process $\eta_{N}(\cdot)$ are important. Hence, in
order to avoid confusion, we have to carefully define potential theoretical
notions for these processes.
\begin{itemize}
\item Denote by $\tau_{A}$ and $\tau_{\mathcal{A}}$ the hitting times
of the sets $A\subset S$ and $\mathcal{A}\subset\mathcal{H}_{N}$,
respectively. In this part, the subsets of $S$ will be denoted by
plain capital letters, while the subsets of $\mathcal{H}_{N}$ are
denoted by calligraphic capital letters.
\item For two disjoint and non-empty sets $A$ and $B$ of $S$, we denote
by $h_{A,\,B}:S\rightarrow[0,\,1]$ and $\textup{cap}_{\mathbb{X}}(A,\,B)$
the equilibrium potential and the capacity with respect to the underlying
process $\mathbb{X}(\cdot)$, respectively:
\begin{align*}
h_{A,\,B}(x) & :=\mathbf{P}_{x}[\tau_{A}<\tau_{B}]\;\text{\;\;\;;}\;x\in S\;,\text{ and}\\
\textup{cap}_{\mathbb{X}}(A,\,B) & :=D_{\mathbb{X}}(h_{A,\,B})\;.
\end{align*}
 For two disjoint and non-empty sets $\mathcal{A}$ and $\mathcal{B}$
of $\mathcal{H}_{N}$, we denote by $\mathbf{h}_{\mathcal{A},\,\mathcal{B}}:\mathcal{H}_{N}\rightarrow[0,\,1]$
and $\textup{cap}_{N}(\mathcal{A},\,\mathcal{B})$ the equilibrium
potential and the capacity with respect to the zero-range processes
$\eta_{N}(\cdot)$, respectively:
\begin{align*}
\mathbf{h}_{\mathcal{A},\,\mathcal{B}}(\eta) & =\mathbf{h}_{\mathcal{A},\,\mathcal{B}}^{N}(\eta):=\mathbb{P}_{\eta}^{N}\left[\tau_{\mathcal{A}}<\mathcal{\tau_{\mathcal{B}}}\right]\;\;;\;\eta\in\mathcal{H}_{N}\;,\text{ and }\\
 & \textup{cap}_{N}(\mathcal{A},\,\mathcal{B}):=\mathscr{D}_{N}(\mathbf{h}_{\mathcal{A},\,\mathcal{B}})\;.
\end{align*}
\end{itemize}

\subsubsection*{Adjoint and symmetrized processes}

Define the adjoint rate $r^{\dagger}(\cdot,\,\cdot)$and the symmetrized
rate $r^{s}(\cdot,\,\cdot)$ as
\[
r^{\dagger}(x,\,y)=\begin{cases}
1-p & \text{if }y=x+1\;,\\
p & \text{if }y=x-1\;,\\
0 & \text{otherwise\;,}
\end{cases}\;\;\;\;\text{and\;\;\;\;}r^{s}(x,\,y)=\begin{cases}
1/2 & \text{if }y=x+1\;,\\
1/2 & \text{if }y=x-1\;,\\
0 & \text{otherwise\;,}
\end{cases}
\]
so that
\[
r^{s}(x,\,y)=\frac{1}{2}\left[r(x,\,y)+r^{\dagger}(x,\,y)\right]\;.
\]
Denote by $(\mathbb{X}^{\dagger}(t))_{t\ge0}$ and $(\mathbb{X}^{s}(t))_{t\ge0}$
the Markov processes on $S$ with rate $r^{\dagger}(\cdot,\,\cdot)$
and $r^{s}(\cdot,\,\cdot)$, respectively.
\begin{xca}
Prove that $\mathbb{X}^{\dagger}(\cdot)$ and $\mathbb{X}^{s}(\cdot)$
are the adjoint and symmetrized processes, respectively, of the underlying
process $\mathbb{X}(\cdot)$.
\end{xca}

We write $L_{\mathbb{X}}^{\dagger}$ and $L_{\mathbb{X}}^{s}$ for
the generators of the processes $\mathbb{X}^{\dagger}(\cdot)$ and
$\mathbb{X}^{s}(\cdot)$, respectively. In addition, we write $h_{A,\,B}^{\dagger}(\cdot)$
for the equilibrium potential with respect to the process $\mathbb{X}^{\dagger}(\cdot)$.

Next we define two generators $\mathscr{L}_{N}^{\dagger}$ and $\mathscr{L}_{N}^{s}$
acting on $\mathbf{f}:\mathcal{H}_{N}\rightarrow\mathbb{R}$ as
\begin{align*}
 & (\mathscr{L}_{N}^{\dagger}\mathbf{f})(\eta)=\sum_{x,y\in S}g(\eta_{x})r^{\dagger}(x,\,y)\,(\mathbf{f}(\sigma^{x,\,y}\eta)-\mathbf{f}(\eta))\;\text{and}\\
 & (\mathscr{L}_{N}^{s}\mathbf{f})(\eta)=\sum_{x,y\in S}g(\eta_{x})r^{s}(x,\,y)\,(\mathbf{f}(\sigma^{x,\,y}\eta)-\mathbf{f}(\eta))\;,
\end{align*}
respectively. Denote by $(\eta_{N}^{\dagger}(t))_{t\ge0}$ and $(\eta_{N}^{s}(t))_{t\ge0}$
the continuous-time Markov processes on $\mathcal{H}_{N}$ generated
by $\mathscr{L}_{N}^{\dagger}$ and $\mathscr{L}_{N}^{s}$, respectively.
\begin{xca}
Prove that $\eta_{N}^{\dagger}(\cdot)$ and $\eta_{N}^{s}(\cdot)$
are the adjoint and symmetrized processes, respectively, of the zero-range
process $\eta_{N}(\cdot)$.
\end{xca}

We write $\mathbf{h}_{\mathcal{A},\,\mathcal{B}}^{\dagger}(\cdot)$
for the equilibrium potential with respect to the adjoint process
$\eta_{N}^{\dagger}(\cdot)$. We also write $\textup{cap}_{N}^{s}(\mathcal{A},\,\mathcal{B})$
for the capacity with respect to the symmetrized process $\eta_{N}^{s}(\cdot)$.\textbf{ }

\section{Condensation Phenomenon}

\textcolor{black}{In this section, we explain the condensation phenomena
of the zero-range processes defined in the previous section.}

\subsubsection*{Metastable valleys }

We first define an auxiliary sequences to concretely define the metastable
sets of the zero-range processes. For two sequences $(a_{N})_{N\in\mathbb{N}},\,(b_{N})_{N\in\mathbb{N}}$
of positive real numbers, the notation $a_{N}\ll b_{N}$ implies that
\[
\lim_{N\rightarrow\infty}\frac{b_{N}}{a_{N}}=\infty\;.
\]
Let $(\ell_{N})_{N\in\mathbb{N}}$ be sequences of positive integer
such that
\begin{equation}
1\ll\ell_{N}\ll N^{(1+\alpha)/(1+(\kappa-1)\alpha)}\;.\label{cond_ell}
\end{equation}
We explain later the reason for imposing this complicated upper bound
for $\ell_{N}$.

For each $x\in S$, the \textit{metastable valley }or\textit{ metastable
set} $\mathcal{E}_{N}^{x}\subset\mathcal{H}_{N}$ is defined as the
set of configurations such that all but at most $\ell_{N}\ll N$ (by
(\ref{cond_ell}) since $\kappa\ge2$) particles are condensed at
site $x$:
\[
\mathcal{E}_{N}^{x}=\left\{ \eta\in\mathcal{H}_{N}:\eta_{x}\ge N-\ell_{N}\right\} \;.
\]

Define
\begin{equation}
\mathcal{E}_{N}=\bigcup_{x\in S}\mathcal{E}_{N}^{x}\;\;\;\;\text{and\;\;\;\;}\Delta_{N}=\mathcal{H}_{N}\setminus\mathcal{E}_{N}\;.\label{endn}
\end{equation}

\subsubsection*{Condensation of particles}

The following theorem shows that the zero-range process defined above
exhibit a phenomenon known as condensation of particles.
\begin{thm}
\label{t01}It holds that
\[
\lim_{N\rightarrow\infty}\mu_{N}(\mathcal{E}_{N}^{x})=\frac{1}{\kappa}\text{ for all }x\in S\;.
\]
Therefore, the invariant measure $\mu_{N}(\cdot)$ is concentrated
on the metastable sets defined above in the sense that
\[
\lim_{N\rightarrow\infty}\mu_{N}(\mathcal{E}_{N})=1\;\;\text{\;\;and\;\;\;\;}\lim_{N\rightarrow\infty}\mu_{N}(\Delta_{N})=0\;.
\]
\end{thm}

\begin{proof}
We refer to \cite[display (3.2)]{B-L ZRP} for a proof.
\end{proof}
\begin{rem}
This theorem holds for any sequence $(\ell_{N})_{N\in\mathbb{N}}$
satisfying $1\ll\ell_{N}\ll N$. The condition $\ell_{N}\ll N^{(1+\alpha)/(1+(\kappa-1)\alpha)}$
appearing in (\ref{cond_ell}) is used only in the investigation of
the metastable behavior explained in the next section (cf. conditions
(H1) and (H3) introduced later in (\ref{H1}) and (\ref{H3}), respectively).
\end{rem}

This theorem assert that, with dominating probability (as $N$ gets
larger), almost all particles are condensed at a single site. This
phenomenon is called a condensation of particles. Hence, if the zero-range
process starts from any configuration, it will eventually form a condensate
at a certain site. Subsequently, this condensate will move around
sites of $S$. Such movements of the condensate, which are often referred
to as the inter-valley dynamics, are the metastable behavior of the
zero-range processes and are our main concern that will be discussed
in the next section.

\section{\label{sec32-1}Markov Chain Model Reduction}

In this section, we introduce the main results regarding the analysis
of the metastable behavior of the zero-range process, and then outline
a general framework regarding the Markov chain model reduction of
the metastable behavior that can be applied to the current model.
This general framework is called the martingale approach, which is
developed in \cite{B-L TM,B-L TM2,B-L MG} and then enhanced in \cite{L-L-M}.

\subsection{Order process }

In this section, we introduce the so-called order process which represents
the inter-valley dynamics and hence plays a significant role in the
Markov chain model reduction. All the definitions introduced in the
current section can be made for a general class of Markov processes,
but we define them only in the context of the zero-range processes
for the convenience of the discussion.

\subsubsection*{Trace process}

The trace process of the zero-range process $\eta_{N}(\cdot)$ on
the set $\mathcal{E}_{N}$ (cf. (\ref{endn})) is defined as
\[
T^{\mathcal{E}_{N}}(t)=\int_{0}^{t}\mathbf{1}\left\{ \eta_{N}(s)\in\mathcal{E}_{N}\right\} ds\;\;\;\;;\;t\ge0\ .
\]
This random time represents the total amount of time for which the
zero-range process stays in $\mathcal{E}_{N}$ up to time $t$. We
denote by $S^{\mathcal{E}_{N}}:[0,\,\infty)\rightarrow[0,\,\infty)$
the generalized inverse of the non-decreasing function $T^{\mathcal{E}_{N}}(\cdot)$,
i.e.,
\[
S^{\mathcal{E}_{N}}(t)=\sup\left\{ s\ge0:T^{\mathcal{E}_{N}}(s)\le t\right\} \;\;\;\;;\;t\ge0\;.
\]
The trace process $(\eta_{N}^{\mathcal{E}_{N}}(t))_{t\ge0}$ of the
zero-range process $\eta_{N}(\cdot)$ on the set $\mathcal{E}_{N}$
is defined by
\[
\eta_{N}^{\mathcal{E}_{N}}(t)=\eta_{N}(S^{\mathcal{E}_{N}}(t))\;\;\;\;;\;t\ge0\;.
\]
By carefully looking at the definition, one can observe that the trajectory
of $\eta_{N}^{\mathcal{E}_{N}}(\cdot)$ is obtained from that of the
zero-range process $\eta_{N}(\cdot)$ by removing the excursion of
$\eta_{N}(\cdot)$ on the set $\Delta_{N}$ (cf. (\ref{endn})). This
is the reason that the process $\eta_{N}^{\mathcal{E}_{N}}(\cdot)$
is called the trace process of $\eta_{N}(\cdot)$ on $\mathcal{E}_{N}$.
\begin{xca}
(The answers to the following questions can be found in \cite{B-L TM,B-L TM2})
\begin{enumerate}
\item Prove that $\eta_{N}^{\mathcal{E}_{N}}(\cdot)$ is indeed an irreducible
continuous-time Markov process on $\mathcal{E}_{N}$.
\item Prove that the invariant measure of the process $\eta_{N}^{\mathcal{E}_{N}}(\cdot)$
is the conditioned measure $\mu_{N}^{\mathcal{E}_{N}}(\cdot)$ of
$\mu_{N}(\cdot)$ on $\mathcal{E}_{N}$, i.e.,
\[
\mu_{N}^{\mathcal{E}_{N}}(\eta)=\frac{\mu_{N}(\eta)}{\mu_{N}(\mathcal{E}_{N})}\;\;\;\;;\;\eta\in\mathcal{E}_{N}\;.
\]
\item Prove that the Markov process $\eta_{N}^{\mathcal{E}_{N}}(\cdot)$
is reversible if $\eta_{N}(\cdot)$ is reversible. Is the converse
true?
\end{enumerate}
\end{xca}

\subsubsection*{Order process}

We note that the trace process $\eta_{N}^{\mathcal{E}_{N}}(\cdot)$
includes all the information about the behavior of the process $\eta_{N}(\cdot)$
on $\mathcal{E}_{N}$. However, in view of the metastable behavior,
we are only concerned with the inter-valley dynamics and are not interested
in the exact location in a metastable valley within which the zero-range
process is staying. Hence, the order process is defined as the process
obtained from the trace process by discarding this information.

More precisely, we define a projection function $\Psi:\mathcal{E}_{N}\rightarrow S$
as
\[
\Psi(\eta)=\sum_{x\in S}x\cdot\mathbf{1}\{x\in\mathcal{E}_{N}^{x}\}
\]
and then define the\textit{ order process} as
\[
Y_{N}(t)=\Psi(\eta_{N}^{\mathcal{E}_{N}}(N^{1+\alpha}t))\;\;\;\;;\;t\ge0\;.
\]

To explain the meaning of the order process, we first consider the
projected trace process
\[
W_{N}(t)=\Psi(\eta_{N}^{\mathcal{E}_{N}}(t))\;\;\;\;;\;t\ge0\;.
\]
This process $W_{N}(t)$ indicates the label of the valley at which
the trace process $\eta_{N}^{\mathcal{E}_{N}}(t)$ is staying. Hence,
this process captures all the relevant information regarding the inter-valley
dynamics of the process $\eta_{N}(\cdot)$ on $\mathcal{E}_{N}$.
We defined $Y_{N}(t)$ as a speeded-up version of this process, namely,
\[
Y_{N}(t)=W_{N}(N^{1+\alpha}t)
\]
since we observe the transitions between metastable valleys in the
time scale of $N^{1+\alpha}$.

It takes a long time to move a condensate from one site to another
since the particles are sticky and hence tend to keep the condensate.
We can also notice that this transition time scale $N^{1+\alpha}$
is increasing in $\alpha$. This is a natural result since the parameter
$\alpha$ corresponds to the stickiness of the particles.

\subsection{Markov chain model reduction via convergence of order process}

\subsubsection*{Markov chain model reduction}

We note here that the order process may not be a Markov process. However,
one can usually prove that, in the metastable situation, the order
process converges to a certain limiting Markov process $Y(\cdot)$
on $S$. Heuristically, this is mainly because the process entering
a metastable valley will spend long enough time inside the valley
to forget the entering location. This is indeed the case for the zero-range
process, and the following is the main theorem regarding the Markov
chain model reduction. We remark that the limiting Markov process
$Y(\cdot)$ for the zero-range process is defined in the next paragraph.
\begin{thm}
\label{t_zrpmain}The following hold:
\begin{enumerate}
\item Suppose that $\eta_{N}(0)\in\mathcal{E}_{N}^{x}$ for all $N\in\mathbb{N}$
for some $x\in S$. Then, the law of the order process $Y_{N}(\cdot)$
converges to the law of limiting Markov process $Y(\cdot)$ starting
at $x$.
\item For all $T>0$, it holds that
\[
\lim_{N\rightarrow\infty}\sup_{\eta\in\mathcal{E}_{N}}\mathbb{E}_{\eta}^{N}\left[\int_{0}^{T}\mathbf{1}\{\eta_{N}(N^{1+\alpha}t)\in\Delta_{N}\}dt\right]=0\;.
\]
\end{enumerate}
\end{thm}

If the zero-range process $\eta_{N}(\cdot)$ spends non-negligible
amount of time at $\Delta_{N}=\mathcal{H}_{N}\setminus\mathcal{E}_{N}$,
then the trace process $\eta_{N}^{\mathcal{E}_{N}}(\cdot)$, which
is obtained by turning off the clock when the zero-range process $\eta_{N}(\cdot)$
stays at $\Delta_{N}$, discards too much information regarding the
inter-valley dynamics of the zero-range process. Part (2) of the previous
theorem implies that, in the scale $N^{1+\alpha}$, the zero-range
process does not spend meaningful amount of time at $\Delta_{N}$
and hence the trace process is indeed a good approximation of $\eta_{N}(\cdot)$
in view of the inter-valley dynamics. This gives authority to part
(1) which asserts that the inter-valley dynamics of the trace process
(and hence the zero-range process by part (2)) is approximated by
the limiting Markov process $Y(\cdot)$. So far, we have explained
a general way to derive the Markov chain model reduction via convergence
of order process.

We discuss the strategy to prove Theorem \ref{t_zrpmain} in Section
\ref{sec61-1}.

\subsubsection*{Limiting Markov process}

We next define the limiting Markov process $Y(\cdot)$ for the zero-range
process. Define a constant by
\begin{equation}
I_{\alpha}=\int_{0}^{1}u^{\alpha}(1-u)^{\alpha}du\;.\label{ialpha}
\end{equation}
Define $a:S\times S\rightarrow[0,\,\infty)$ by
\begin{equation}
a(x,\,y)=\frac{\kappa}{\Gamma_{\alpha}I_{\alpha}}\textup{cap}_{\mathbb{X}}(x,\,y)\;\;\;\;;\;x,\,y\in S\;,\label{eax}
\end{equation}
where we remind here that the notation $x,\,y\in S$ implies that
$x$ and $y$ are different. Note that the capacity $\textup{cap}_{\mathbb{X}}(x,\,y)$
has been computed in Exercise \ref{ex26}. Now the limiting Markov
process $(Y(t))_{t\ge0}$ is define as a continuous-time Markov process
on $S$ with rate $a(\cdot,\,\cdot)$. We denote by $\mathbf{Q}_{x}$
the law of process $Y(\cdot)$ starting at $x\in S$.

Since $a(x,\,y)>0$ for all $x,\,y>0$, the irreducibility is clear
for the process $Y(\cdot)$. Denote by $\nu(\cdot)$ the uniform measure
on $S$:
\begin{equation}
\nu(x)=\frac{1}{\kappa}\;\;\;\;;\;x\in S\;.\label{mux}
\end{equation}

\begin{xca}
Prove that that the unique invariant measure of the irreducible Markov
process $Y(\cdot)$ is $\nu(\cdot)$ and furthermore, that the process
$Y(\cdot)$ is reversible. (Hint: use (\ref{e_cap_synm}))
\end{xca}

A remarkable fact here is that the limiting Markov process is always
reversible, while the underlying zero-range process is not, especially
for $p\neq1/2$.

\subsection{Markov chain model reduction via convergence of marginal distributions}

An alternative way of describing the Markov chain model reduction
was developed in \cite{L-L-M}. This methodology does not discard
the excursions of the zero-range process on $\Delta_{N}$ (and hence
does not use the trace and order processes) but proves the convergence
result with a weaker notion of convergence, namely the convergence
of finite dimensional distributions. This is the nature of the problem;
without removing noisy excursions at $\Delta_{N}$, we cannot expect
the convergence in path space with the usual mode of convergence.
We refer to \cite{B-L TM} for more detail. Instead, the soft topology
introduced in \cite{Lan-soft} can be used to prove the convergence.

To explain this alternative method, let us define a projection function
$\widehat{\Psi}:\mathcal{H}_{N}\rightarrow S\cup\{\mathfrak{0}\}$
as
\[
\widehat{\Psi}(\eta)=\begin{cases}
x & \mbox{if }x\in\mathcal{E}_{N}^{x}\;,\\
\mathfrak{0} & \mbox{if }x\in\Delta_{N}\;.
\end{cases}
\]
Then, define a process $(\widehat{Y}_{N}(t))_{t\ge0}$ as
\[
\widehat{Y}_{N}(t)=\widehat{\Psi}(\eta_{N}(N^{1+\alpha}t))\;\;\;\;;\;t\ge0
\]
Then, the process $\widehat{Y}_{N}(\cdot)$ is a process on $\widehat{S}=S\cup\{\mathfrak{0}\}$
and may not be a Markov process. We note that the order process $Y_{N}(\cdot)$
is a trace process of $\widehat{Y}_{N}(\cdot)$ on the set $S$.

Define an extended limiting process $(\widehat{Y}(t))_{t\ge0}$ on
$\widehat{S}$ as a continuous-time Markov process with jump rate
\[
\widehat{a}(x,\,y)=\begin{cases}
a(x,\,y) & \text{if }x,\,y\in S\;,\\
0 & \text{otherwise\;.}
\end{cases}
\]
Hence, $\mathfrak{0}$ is merely a cemetery point of the Markov process
$\widehat{Y}(\cdot)$. Denote by $\widehat{\mathbf{Q}}_{x}$, $x\in\widehat{S}$,
the law of process $\widehat{Y}(\cdot)$ that starts at $x$.
\begin{xca}
Prove that the measure $\widehat{\nu}(\cdot)$ on $\widehat{S}$ defined
by
\[
\widehat{\nu}(x)=\begin{cases}
\nu(x) & \text{if }x\in S\;,\\
0 & \text{otherwise\;.}
\end{cases}
\]
is an invariant measure of the Markov process $\widehat{Y}(\cdot)$.
\end{xca}

The following is the second way of establishing a Markov chain model
reduction of the metastable behavior developed in \cite{L-L-M}.
\begin{thm}
\label{t_zrp2}For all $x\in S$ and for all $(\eta_{N})_{N\in\mathbb{N}}$
such that $\eta_{N}\in\mathcal{E}_{N}^{x}$ for all $N$, the finite
dimensional distributions of the process $\widehat{Y}_{N}(\cdot)$
under $\mathbb{P}_{\eta_{N}}^{N}$ converges to that of the law $\widehat{\mathbf{Q}}_{x}$,
as $N$ tends to infinity.
\end{thm}

The proof of this theorem is close to that of Theorem \ref{t_zrpmain}
and will be explained in the next subsection.

\subsection{\label{sec61-1}Martingale approach}

In this section, we explain the general principle developed in \cite{B-L TM,B-L TM2,B-L MG,L-L-M}.
This principle, which is now called the martingale approach to the
metastability reduces the proof of Theorems \ref{t_zrpmain} and \ref{t_zrp2}
to the verification of several sufficient conditions.

To explain the general principle in the context of zero-range process,
we now explain several essential notions.
\begin{itemize}
\item Recall that $\eta_{N}^{\mathcal{E}_{N}}(\cdot)$ is a Markov process
on $\mathcal{E}_{N}$. Denote by $j_{N}:\mathcal{E}_{N}\times\mathcal{E}_{N}\rightarrow[0,\,\infty)$
the jump rate of the process $\eta_{N}^{\mathcal{E}_{N}}(\cdot)$.
\item For $x,\,y\in S$, the \textit{mean jump rate} between two valleys
$\mathcal{E}_{N}^{x}$ and $\mathcal{E}_{N}^{y}$ is defined by
\[
r_{N}(x,\,y)=\frac{1}{\mu_{N}(\mathcal{E}_{N}^{x})}\sum_{\eta\in\mathcal{E}_{N}^{x}}\sum_{\zeta\in\mathcal{E}_{N}^{y}}\mu_{N}(\eta)\,j_{N}(\eta,\,\zeta)\;.
\]
\item For each $x\in S$, let $\xi_{N}^{x}\in\mathcal{H}_{N}$ be the configuration
such that all particles are located at site $x$.
\item For $x\in S$, write $\breve{\mathcal{E}}_{N}^{x}=\mathcal{E}_{N}\setminus\mathcal{E}_{N}^{x}$.
\item For $x,\,y\in S$, write $\breve{\mathcal{E}}_{N}^{x,\,y}=\mathcal{E}_{N}\setminus\left(\mathcal{E}_{N}^{x}\cup\mathcal{E}_{N}^{y}\right)$.
\end{itemize}
Now we introduce several sufficient conditions for the Markov chain
model reduction.
\begin{itemize}
\item Condition \textbf{(H0)}: For all $x,\,y\in S$,
\begin{equation}
\lim_{N\rightarrow\infty}N^{1+\alpha}\,r_{N}(x,\,y)=a(x,\,y)\;.\label{H0}
\end{equation}
Hence, the mean jump rate between two valleys $\mathcal{E}_{N}^{x}$
and $\mathcal{E}_{N}^{y}$ is approximately $a(x,\,y)/N^{1+\alpha}$.
This is the reason that we accelerated the process by a factor of
$N^{1+\alpha}$ in the definition of the order process. This accurate
estimate of the mean jump rate is the crucial and most difficult step
in the proof of Theorems \ref{t_zrpmain} and \ref{t_zrp2}.
\item Condition \textbf{(H1)}: For all $x\in S$,
\begin{equation}
\lim_{N\rightarrow\infty}\sup_{\eta,\,\zeta\in\mathcal{E}_{N}^{x}}\frac{\textup{cap}_{N}(\mathcal{E}_{N}^{x},\,\breve{\mathcal{E}}_{N}^{x})}{\textup{cap}_{N}(\eta,\,\zeta)}=0\;.\label{H1}
\end{equation}
This condition implies that, for any $\eta,\,\zeta\in\mathcal{E}_{N}^{x}$,
the process starting at $\eta\in\mathcal{E}_{N}^{x}$ hits the configuration
$\zeta\in\mathcal{E}_{N}^{x}$ before hitting the set $\breve{\mathcal{E}}_{N}^{x}$,
i.e., before arriving at one of other valleys, with dominating probability.
We term this phenomenon a visiting property. We discuss this further
in Remark \ref{rem_zrp1}.
\begin{xca}
Prove the last assertion. (Hint: Proposition \ref{p_eqp_bd})
\end{xca}

\item Condition \textbf{(H2)}: For all $x\in S$,
\begin{equation}
\lim_{N\rightarrow\infty}\frac{\mu_{N}(\Delta_{N})}{\mu_{N}(\mathcal{E}_{N}^{x})}=0\;.\label{H2}
\end{equation}
This condition implies that the set $\Delta_{N}$ is negligible compared
to $\mathcal{E}_{N}^{x}$ with respect to the invariant measure. \textit{We
emphasize that this condition is a direct consequence of Theorem \ref{t01}. }
\item Condition \textbf{(H3)}: For all $x\in S$,
\begin{equation}
\lim_{\delta\rightarrow0}\limsup_{N\rightarrow\infty}\sup_{\eta\in\mathcal{E}_{N}^{x}}\sup_{2\delta\le s\le3\delta}\mathbb{P}_{\eta}^{N}\left[\eta_{N}(N^{1+\alpha}\,s)\in\Delta_{N}\right]=0\;.\label{H3}
\end{equation}
This implies that, if the zero-range process starts from a valley
it will still be in the same valley after a short time. Note that
we cannot replace $\sup_{2\delta\le s\le3\delta}$ with $\sup_{0\le s\le\delta}$,
since if the process starts at the boundary of $\mathcal{E}_{N}^{x}$,
then with a non-negligible probability it leaves the valley within
a few steps. This condition (\ref{H3}) implies that, even after such
an escape from the valley, the process returns to the valley immediately.
\end{itemize}
The next theorem is a consequence of \cite[Theorem 2.1]{B-L TM2}
and \cite[Proposition 1.1]{L-L-M},
\begin{thm}
Suppose that conditions (H0), (H1), and (H2) hold. Then, Theorem \ref{t_zrpmain}
holds. Moreover, if condition (H3) additionally holds, then Theorem
\ref{t_zrp2} also holds.
\end{thm}

Therefore, to prove Theorems \ref{t_zrpmain} and \ref{t_zrp2}, it
suffices to verify the conditions (H0), (H1), (H2), and (H3):
\begin{itemize}
\item Condition (H0) will be proven in Proposition \ref{p_H0}. This is
the most difficult part of the current problem. The main components
of the proof are the estimate of capacities, the sector condition,
and the argument developed in \cite{L-Seo NR1} based on collapsed
processes.
\item Conditions (H1) and (H3) are proven based on Propositions \ref{p_H1}
and \ref{p_H3}, respectively, based on the estimate of capacities.
\item As we have mentioned above, the condition (H2) is a consequence of
Theorem \ref{t01}.
\end{itemize}
\begin{rem}
\label{rem:alpha}In fact, Theorem \ref{t01} holds only for $\alpha\ge1$
(for the critical case $\alpha=1$, we should be more careful about
the selection of $\ell_{N}$, see \cite{L-M-Seo1}) and hence the
metastable behavior must be studied for all $\alpha\ge1$. Below is
the history of the research on this problem in chronological order:
\begin{enumerate}
\item Beltran and Landim \cite{B-L ZRP} first analyzed the reversible case
$p=1/2$ with $\alpha>1$.
\item Landim \cite{L-ZRP} analyzed the totally asymmetric case $p=1$ with
$\alpha>3$.
\item Seo \cite{Seo NRZRP} analyzed the general case $p\in[0,\,1]$ with
$\alpha>2$.
\item Landim, Marcondes and Seo \cite{L-M-Seo1,L-M-Seo2} analyzed the critical
case $\alpha=1$ with $p=1/2$.
\end{enumerate}
We note that the articles \cite{B-L ZRP,L-M-Seo1,L-M-Seo2,Seo NRZRP}
considered a more general case, i.e., the particle system on any finite
set consisting of  any underlying random walk $\mathbb{X}(\cdot)$.
The articles \cite{B-L ZRP,L-M-Seo1,L-M-Seo2} assumed the reversibility
of the zero-range process. Moreover, \cite{L-M-Seo1,L-M-Seo2} assumed
that the invariant measure for the underlying random walk $\mathbb{X}(\cdot)$
is the uniform measure on $S$. We also emphasize here that \cite{L-ZRP}
is the first rigorous quantitative analysis of the metastable behavior
of a non-reversible Markov process.
\end{rem}

\begin{rem}
The current part of this lecture note is mainly derived from article
\cite{Seo NRZRP}. With a more refined argument, we are able to weaken
the assumption $\alpha>2$ of \cite{Seo NRZRP} to $\alpha>1$. The
critical case $\alpha=1$ for the non-reversible case is largely unknown
at this moment. We discuss in the next remark the difficulty of the
critical case.
\end{rem}

\begin{rem}
\label{rem_zrp1}If $\ell_{N}$ is too large, then there are too many
configurations inside the valley and hence the visiting property explained
in condition (H1) may not hold. In fact, the upper bound of $\ell_{N}$
given in (\ref{cond_ell}) is imposed to verify condition (H1). For
the critical case $\alpha=1$, in order to ensure that (H1) is in
force, we have to take $\ell_{N}$ so small that the metastable valley
$\mathcal{E}_{N}^{x}$ with such $\ell_{N}$ violates Theorem \ref{t01}
(i.e., the condition (H2)).\textit{ In conclusion, the critical zero-range
process cannot satisfy two condition (H1) and (H2) simultaneously,
no matter what value we give to $\ell_{N}$.} This is the reason that
the critical case cannot be handled with the martingale approach described
here. Recently, \cite{L-M-Seo2} developed a new approach based on
the analysis of the solution of certain form of resolvent equations
and used this approach to investigate the metastable behavior of critical
case with $p=1/2$.
\end{rem}

\subsection{Outlook of the remainder of Part \ref{pt4}}

In the remainder of the note, we verify conditions (H0), (H1), and
(H3).
\begin{itemize}
\item In Section \ref{sec16}, we explain and prove the capacity estimates
between valleys. The proof is based on the generalized Dirichlet and
Thomson principles and hence we need to construct the test functions
and flows.
\item In Section \ref{sec17}, we prove condition (H0).
\item In Section \ref{sec18}, we prove conditions (H1) and (H3).
\end{itemize}

\section{\label{sec16}Estimate of Capacities }

In this section, we provide, up to the construction of test objects,
the estimate of the capacity between metastable valleys based on generalized
Dirichlet and Thomson principles.

\subsubsection*{Main result}

For $f:S\rightarrow\mathbb{R}$, the generator of the limiting Markov
process $Y(\cdot)$ on $S$ can be written as
\begin{equation}
(\mathfrak{L}_{Y}f)(x)=\sum_{y\in S\setminus\{x\}}\,\frac{\kappa\,\textup{cap}_{X}(x,\,y)}{\Gamma_{\alpha}\,I_{\alpha}}\left[f(y)-f(x)\right]\;\;\;\;\;;\;x\in S\;.\label{hu3}
\end{equation}
As we have mentioned before, the invariant measure for $Y(\cdot)$
is the uniform measure $\nu(\cdot)$ on $S$, i.e.,
\[
\nu(x)=\frac{1}{\kappa}\;\;\text{for all }\;x\in S\;.
\]
Therefore, the Dirichlet form with respect to the process $Y(\cdot)$
acting on $f:S\rightarrow\mathbb{R}$ such a way that
\[
\mathfrak{D}_{Y}(f)=\sum_{x\in S}\nu(x)\,f(x)\left[-(\mathfrak{L}_{Y}\,f)(x)\right]=\frac{1}{2}\sum_{x\in S}\sum_{y\in S}\frac{\textup{cap}_{X}(x,\,y)}{\Gamma_{\alpha}\,I_{\alpha}}\left[f(y)-f(x)\right]^{2}\;
\]
Recall that $\mathbf{Q}_{x}$ denote the law of the process $Y(\cdot)$
starting from $x\in S$. For two disjoint non-empty sets $A$ and
$B$ of $S$, the\textbf{ }equilibrium potential and capacity between
$A$ and $B$ with respect to the process $Y(\cdot)$ are defined
by
\begin{equation}
\mathfrak{h}_{A,\,B}(x)=\mathbf{\mathbf{Q}}_{x}(\tau_{A}<\tau_{B})\;\;\text{for}\;x\in S_{\star}\;\;\text{and\;\;}\textup{cap}_{Y}(A,\,B)=\mathfrak{D}_{Y}(\mathfrak{h}_{A,\,B})\;,\label{hu4}
\end{equation}
respectively.

For a non-empty set $A\subseteq S$, we write
\[
\mathcal{E}_{N}(A)=\bigcup_{x\in A}\mathcal{E}_{N}^{x}\;.
\]
The following theorem is the main capacity estimate for the zero-range
processes
\begin{thm}
\label{t_zrpcap}For disjoint, non-empty subsets $A,\,B$ of $S$,
we have that
\[
\lim_{N\rightarrow\infty}N^{1+\alpha}\,\textup{cap}_{N}(\mathcal{E}_{N}(A),\,\mathcal{E}_{N}(B))=\textup{cap}_{Y}(A,\,B)\;.
\]
\end{thm}

In addition, if $(A,\,B)$ is a partition of $S$, that is, $A\cup B=S$,
the equilibrium potential $\mathfrak{h}_{A,\,B}(\cdot)$ becomes the
indicator function on $A$, and hence by (\ref{hu4}) we immediately
obtain the following result as a corollary of the previous theorem.
\begin{cor}
\label{c02}Suppose that two disjoint, non-empty subsets $A,\,B$
of $S$ satisfy $A\cup B=S$. Then,
\[
\lim_{N\rightarrow\infty}N^{1+\alpha}\,\textup{cap}_{N}(\mathcal{E}_{N}(A),\,\mathcal{E}_{N}(B))=\frac{1}{\Gamma_{\alpha}I_{\alpha}}\,\sum_{x\in A}\sum_{y\in B}\textup{cap}_{X}(x,\,y)\;.
\]
\end{cor}

Now we discuss how we can prove Theorem \ref{t_zrpcap}.

\subsubsection*{Strategy to prove Theorem \ref{t_zrpcap}}

Let us now turn to the proof of Theorem \ref{t_zrpcap}, which is
based on the generalized Dirichlet and Thomson principles (cf. Theorem
\ref{t_genDTP_nonrev}). We explain how we can apply these principles
in the context of the zero-range processes.

We start from the test functions and flows. Let us first introduce
a new parameter $\epsilon>0$ denoting small numbers. The parameter
$\epsilon$ will be sent to $0$ in the end (after sending $N$ to
$\infty$).
\begin{rem}
Henceforth, all constants are assumed to depend only on $p,\,\kappa$,
$\alpha$, and $\epsilon$ and are independent of $N$. Furthermore,
we write $a(N,\,\epsilon)=o_{N}(1)$ and $b(N,\,\epsilon)=o_{\epsilon}(1)$
if
\begin{align*}
 & \lim_{N\rightarrow\infty}a(N,\,\epsilon)=0\text{ for all }\epsilon>0\;\;\text{ and }\\
 & \lim_{\epsilon\rightarrow0}\sup_{N\in\mathbb{N}}b(N,\,\epsilon)=0\;,
\end{align*}
respectively. The dependencies of the constant and the $o_{N}(1)$
term on the parameter $\epsilon$ do not incur any problem, as we
always send $N$ to infinity first before sending $\epsilon$ to $0$.
\end{rem}

Throughout the remainder of the current section, let us fix two disjoint
non-empty subsets $A$ and $B$ of $S$.

In \cite[Section 7]{Seo NRZRP}, for sufficiently large $N\in\mathbb{N}$,
two functions
\[
\mathbf{V}_{A,\,B}=\mathbf{V}_{A,\,B}^{N,\,\epsilon}:\mathcal{H}_{N}\rightarrow\mathbb{R}\;\;\;\;\text{and\;\;\;\;}\mathbf{V}_{A,\,B}^{\dagger}=\mathbf{V}_{A,\,B}^{\dagger,N,\,\epsilon}:\mathcal{H}_{N}\rightarrow\mathbb{R}
\]
approximating the equilibrium potentials $\mathbf{h}_{\mathcal{E}_{N}(A),\mathcal{\,E}_{N}(B)}$
and $\mathbf{h}_{\mathcal{E}_{N}(A),\,\mathcal{E}_{N}(B)}^{\dagger}$,
respectively, are constructed. It is also verified there that these
functions enjoy the following properties.
\begin{prop}
\label{p_zrpfn}For all small enough $\epsilon$ and large enough
$N$, two functions $\mathbf{V}_{A,\,B}$ and $\mathbf{V}_{A,\,B}^{\dagger}$
satisfy the following properties:
\begin{enumerate}
\item It hold that $\mathbf{V}_{A,\,B},\,\mathbf{V}_{A,\,B}^{\dagger}\in\mathfrak{C}_{1,\,0}(\mathcal{E}_{N}(A),\,\mathcal{E}_{N}(B))$.
Moreover, for all $x\in S\setminus\{A,\,B\}$, it holds that
\[
\mathbf{V}_{A,\,B}(\eta)=\mathbf{V}_{A,\,B}^{\dagger}(\eta)=\mathfrak{h}_{A,\,B}(x)\text{ \;\;\;for all }\eta\in\mathcal{E}_{N}^{x}\;.
\]
\item It holds that
\[
N^{1+\alpha}\mathscr{D}_{N}(\mathbf{V}_{A,\,B}),\,N^{1+\alpha}\mathscr{D}_{N}(\mathbf{V}_{A,\,B}^{\dagger})\le[1+o_{N}(1)+o_{\epsilon}(1)]\,\textup{cap}_{Y}(A,\,B)\;.
\]
\end{enumerate}
\end{prop}

We next construct test flows approximating $\Phi_{\mathbf{h}_{\mathcal{E}_{N}(A),\mathcal{\,E}_{N}(B)}}^{*}$
and $\Phi_{\mathbf{h}_{\mathcal{E}_{N}(A),\mathcal{\,E}_{N}(B)}^{\dagger}}$
(cf. (\ref{flow})). The natural candidates are $\Phi_{\mathbf{V}_{A,\,B}}^{*}$
and $\Phi_{\mathbf{V}_{A,\,B}^{\dagger}}$. However, the divergences
of these flows are larger than required along the saddle tube between
metastable sets (cf. \cite[Section 7.2]{Seo NRZRP}) and hence we
need to perform a local surgery to cancel these divergences out without
impacting approximating features of the flows $\Phi_{\mathbf{V}_{A,\,B}}^{*}$
and $\Phi_{\mathbf{V}_{A,\,B}^{\dagger}}$. This procedure is the
most complicated part in the analysis of the zero-range process. The
consequences of this correction procedure can be summarized as follows.
\begin{prop}
\label{p_zrpflow}For all small enough $\epsilon$ and large enough
$N$, there exist flows
\[
\Phi_{A,\,B}=\Phi_{A,\,B}^{N,\,\epsilon}\in\mathfrak{F}_{N}\;\;\;\;\text{and\;\;\;\;}\Phi_{A,\,B}^{\dagger}=\Phi_{A,\,B}^{\dagger,\,N,\,\epsilon}\in\mathfrak{F}_{N}
\]
enjoying the following properties.
\begin{enumerate}
\item The flows $\Phi_{A,\,B}$ and $\Phi_{A,\,B}^{\dagger}$ approximate
$\Phi_{\mathbf{V}_{A,\,B}}^{*}$ and $\Phi_{\mathbf{V}_{A,\,B}^{\dagger}}$
in the sense that
\begin{align*}
\left\Vert \Phi_{A,\,B}-\Phi_{\mathbf{V}_{A,\,B}}^{*}\right\Vert ^{2} & =[o_{N}(1)+o_{\epsilon}(1)]\,N^{-(1+\alpha)}\;\;\;\text{and}\\
\left\Vert \Phi_{A,\,B}^{\dagger}-\Phi_{\mathbf{V}_{A,\,B}^{\dagger}}\right\Vert ^{2} & =[o_{N}(1)+o_{\epsilon}(1)]\,N^{-(1+\alpha)}\;.
\end{align*}
\item The divergence of $\Phi_{A,\,B}$ is negligible on $\Delta_{N}$ in
the sense that
\[
\sum_{\eta\in\Delta_{N}}\left|(\textup{div }\Phi_{A,\,B})(\eta)\right|=o_{N}(1)\,N^{-(1+\alpha)}\;.
\]
\item The divergence of $\Phi_{A,\,B}$ is negligible on $\mathcal{E}_{N}^{x}$,
$x\in S\setminus(A\cup B)$, in the sense that
\begin{align}
 & (\textup{div }\Phi_{A,\,B})(\mathcal{E}_{N}^{x})=o_{N}(1)\,N^{-(1+\alpha)}\;\;\text{and}\label{e4120}\\
 & \sum_{\eta\in\mathcal{E}_{N}^{x}}\mathbf{\mathbf{h}}_{\mathcal{E}_{N}(A),\,\mathcal{E}_{N}(B)}(\eta)\,(\textup{div }\Phi_{A,\,B})(\eta)=o_{N}(1)\,N^{-(1+\alpha)}\;.\label{e412}
\end{align}
\item The divergence of $\Phi_{A,\,B}$ satisfies
\begin{align*}
 & (\textup{div }\Phi_{A,\,B})(\mathcal{E}_{N}(A))=[1+o_{N}(1)]\,N^{-(1+\alpha)}\,\textup{cap}_{Y}(A,\,B)\;\;\mbox{and}\\
 & (\textup{div }\Phi_{A,\,B})(\mathcal{E}_{N}(B))=-[1+o_{N}(1)]\,N^{-(1+\alpha)}\,\textup{cap}_{Y}(A,\,B)\;.
\end{align*}
\end{enumerate}
The flow $\Phi_{A,\,B}^{\dagger}$ also satisfies properties (2),
(3), and (4).
\end{prop}

The proof of this proposition is given in \cite[Section 8]{Seo NRZRP}.

\textit{Since the proofs of Proposition \ref{p_zrpfn} and \ref{p_zrpflow}
are too technical and hence are not suitable as contents of a lecture
note, we refer to the interested readers to the article \cite{Seo NRZRP}.
Instead, we will now focus on how we can prove the Markov chain model
reduction based on this constructions. }

By (2), (3), and (4) of the previous proposition, we have the following
estimate that enables the application of the generalized Dirichlet
and Thomson principles.
\begin{lem}
\label{lem66}We have that
\begin{align}
\sum_{\eta\in\mathcal{H}_{N}}\mathbf{\mathbf{h}}_{\mathcal{E}_{N}(A),\,\mathcal{E}(B)}(\eta)\,(\textup{div }\Phi_{A,\,B})(\eta) & =[1+o_{N}(1)]\,N^{-(1+\alpha)}\,\textup{cap}_{Y}(A,\,B)\;\;\text{and}\label{eq:eh1}\\
\sum_{\eta\in\mathcal{H}_{N}}\mathbf{\mathbf{h}}_{\mathcal{E}_{N}(A),\,\mathcal{E}(B)}(\eta)\,(\textup{div }\Phi_{A,\,B}^{\dagger})(\eta) & =[1+o_{N}(1)]\,N^{-(1+\alpha)}\,\textup{cap}_{Y}(A,\,B)\;.\label{eh2}
\end{align}
\end{lem}

\begin{proof}
We only consider (\ref{eq:eh1}) since the proof of (\ref{eh2}) is
identical. The summation on the left-hand side of (\ref{eq:eh1})
can be divided into
\begin{equation}
\sum_{\eta\in\mathcal{E}_{N}(A)}+\sum_{\eta\in\mathcal{E}_{N}(B)}+\sum_{x\notin A\cup B}\sum_{\eta\in\mathcal{E}_{N}^{x}}+\sum_{\eta\in\Delta_{N}}\;.\label{e_decsm}
\end{equation}
Since $\mathbf{\mathbf{h}}_{\mathcal{E}_{N}(A),\,\mathcal{E}(B)}\equiv1$
on $\mathcal{E}_{N}(A)$, by part (4) of Proposition \ref{p_zrpflow},
the first summation is equal to
\[
[1+o_{N}(1)]\,N^{-(1+\alpha)}\,\textup{cap}_{Y}(A,\,B)\;.
\]
Since $\mathbf{\mathbf{h}}_{\mathcal{E}_{N}(A),\,\mathcal{E}(B)}\equiv0$
and $\mathcal{E}_{N}(B)$, the second summation in (\ref{e_decsm})
is trivially $0$. The third summation is $o_{N}(1)\,N^{-(1+\alpha)}$
by the second estimate of (3) of Proposition \ref{p_zrpflow}. Finally,
as $|\mathbf{\mathbf{h}}_{\mathcal{E}_{N}(A),\,\mathcal{E}(B)}|\le1$,
the last summation is $o_{N}(1)\,N^{-(1+\alpha)}$ by (2) of Proposition
\ref{p_zrpflow}.
\end{proof}
Now by accepting Proposition \ref{p_zrpfn} and \ref{p_zrpflow},
we can complete the proof of Theorem \ref{t_zrpcap}.
\begin{proof}[Proof of Theorem \ref{t_zrpcap}]
Inspired by the optimizer of Theorem \ref{t_genDTP_nonrev}-(1),
let us take
\begin{equation}
\mathbf{f}=\frac{\mathbf{V}_{A,\,B}+\mathbf{V}_{A,\,B}^{\dagger}}{2}\;\;\;\;\mbox{and}\;\;\;\;\phi=\frac{\Phi_{A,\,B}^{\dagger}-\Phi_{A,\,B}}{2}\;.\label{e_test1}
\end{equation}
Note that $\mathbf{f}\in\mathfrak{C}_{1,\,0}(\mathcal{E}_{N}(A),\,\mathcal{E}_{N}(B))$
by (1) of Proposition \ref{p_zrpfn}. Thus, by the generalized Dirichlet
principle (i.e., Theorem \ref{t_genDTP_nonrev}-(1)) and Lemma \ref{lem66},
we can write
\begin{equation}
\textup{cap}_{N}(\mathcal{E}_{N}(A),\,\mathcal{E}_{N}(B))\le\left\Vert \Phi_{\mathbf{f}}-\phi\right\Vert ^{2}+o_{N}(1)\,N^{-(1+\alpha)}\;,\label{e630}
\end{equation}
where $\Vert\cdot\Vert=\Vert\cdot\Vert_{\mathfrak{F}}$ denotes the
flow norm with respect to the zero-range process $\eta_{N}(\cdot)$.

Let us write
\begin{equation}
\Phi_{A,\,B}=\Phi_{\mathbf{V}_{A,\,B}}^{*}+\Theta_{N}\;\;\mbox{and}\;\;\Phi_{A,\,B}^{\dagger}=\Phi_{\mathbf{V}_{A,\,B}^{\dagger}}+\Theta_{N}^{\dagger}\;,\label{ec2}
\end{equation}
so that we have
\begin{align}
\Phi_{\mathbf{f}}-\phi & =\Phi_{(\mathbf{V}_{A,\,B}+\mathbf{V}_{A,\,B}^{\dagger})/2}-\frac{\Phi_{\mathbf{V}_{A,\,B}^{\dagger}}-\Phi_{\mathbf{V}_{A,\,B}}^{*}}{2}+\frac{\Theta_{N}-\Theta_{N}^{\dagger}}{2}\;.\nonumber \\
 & =\Psi_{\mathbf{V}_{A,\,B}}+\frac{\Theta_{N}-\Theta_{N}^{\dagger}}{2}\label{ec33}
\end{align}
By (2) of Proposition \ref{p_zrpfn}, it holds that
\begin{equation}
\left\Vert \Psi_{\mathbf{V}_{A,\,B}}\right\Vert ^{2}=\mathscr{D}_{N}(\mathbf{V}_{A,B})\le[1+o_{N}(1)+o_{\epsilon}(1)]\,N^{-(1+\alpha)}\,\textup{cap}_{Y}(A,\,B)\;.\label{ec331}
\end{equation}
On the other hand, by (2) of Proposition \ref{p_zrpflow} and definition
(\ref{ec2}), it holds that
\begin{equation}
\left\Vert \frac{\Theta_{N}-\Theta_{N}^{\dagger}}{2}\right\Vert ^{2}=[o_{N}(1)+o_{\epsilon}(1)]\,N^{-(1+\alpha)}\label{ec332}
\end{equation}
Therefore, (\ref{ec33}), (\ref{ec331}), (\ref{ec332}), and the
triangle inequality, we can conclude that
\begin{equation}
\left\Vert \Phi_{\mathbf{f}}-\phi\right\Vert ^{2}\le[1+o_{N}(1)+o_{\epsilon}(1)]\,N^{-(1+\alpha)}\,\textup{cap}_{Y}(A,\,B)\;.\label{e64up}
\end{equation}
Inserting this into (\ref{e630}), we obtain the following upper bound
of the capacity.
\begin{equation}
\textup{cap}_{N}(\mathcal{E}_{N}(A),\,\mathcal{E}_{N}(B))\le[1+o_{N}(1)+o_{\epsilon}(1)]\,N^{-(1+\alpha)}\,\textup{cap}_{Y}(A,\,B)\;.\label{e64u}
\end{equation}

Now we use the generalized Thomson principle to obtain the lower bound.
Based on the optimizer of Theorem \ref{t_genDTP_nonrev}-(2) and our
guess of the asymptotic limit of capacity $\textup{cap}_{N}(\mathcal{E}_{N}(A),\,\mathcal{\mathcal{E}}_{N}(B))$,
we take

To this end, let
\begin{equation}
\mathbf{g}=\frac{\mathbf{V}_{A,\,B}^{\dagger}-\mathbf{V}_{A,\,B}}{2\,N^{-(1+\alpha)}\,\textup{cap}_{Y}(A,\,B)}\;\;\;\;\mbox{and}\;\;\;\;\psi=\frac{\Phi_{A,\,B}^{\dagger}+\Phi_{A,\,B}}{2\,N^{-(1+\alpha)}\,\textup{cap}_{Y}(A,\,B)}\;.\label{e_test2}
\end{equation}
By (1) of Proposition \ref{p_zrpfn}, we have $\mathbf{g}\in\mathfrak{C}_{0,\,0}(\mathcal{E}_{N}(A),\,\mathcal{E}_{N}(B))$.
Moreover, by Lemma \ref{lem66}, it holds that
\[
\sum_{\eta\in\mathcal{H}_{N}}\mathbf{\mathbf{h}}_{\mathcal{E}_{N}(A),\,\mathcal{E}(B)}(\eta)\,(\textup{div }\psi_{A,\,B})(\eta)=1+o_{N}(1)+o_{\epsilon}(1)\;.
\]
Therefore, by the generalized Thomson principle (i.e., Theorem \ref{t_genDTP_nonrev}-(2)),
we can conclude that
\begin{equation}
\textup{cap}_{N}(\mathcal{E}_{N}(A),\,\mathcal{E}_{N}(B))\ge\frac{1+o_{N}(1)+o_{\epsilon}(1)}{\left\Vert \Phi_{\mathbf{g}}-\psi\right\Vert ^{2}}\;.\label{e644}
\end{equation}
Now it remains to compute the flow norm $\left\Vert \Phi_{\mathbf{g}}-\psi\right\Vert ^{2}$.
To this end, using (\ref{ec2}), let us write
\[
\Phi_{\mathbf{g}}-\psi=-\frac{1}{N^{-(1+\alpha)}\,\textup{cap}_{Y}(A,\,B)}\left[\Psi_{\mathbf{V}_{A,\,B}}+\frac{\Theta_{N}+\Theta_{N}^{\dagger}}{2}\right]\;.
\]
Then, by similar computations as in the upper bound. we can conclude
that
\begin{equation}
\Vert\Phi_{\mathbf{g}}-\psi\Vert^{2}\le\frac{1+o_{N}(1)+o_{\epsilon}(1)}{N^{-(1+\alpha)}\,\textup{cap}_{Y}(A,\,B)}\;.\label{e645}
\end{equation}
Combining (\ref{e644}) and (\ref{e645}), we can finally obtain the
lower bound on the capacity:
\begin{equation}
\textup{cap}_{N}(\mathcal{E}_{N}(A),\,\mathcal{E}_{N}(B))\ge[1+o_{N}(1)+o_{\epsilon}(1)]\,N^{-(1+\alpha)}\,\textup{cap}_{Y}(A,\,B)\;.\label{e64lo}
\end{equation}
By the upper bound (\ref{e64u}) and lower bound (\ref{e64lo}), we
can conclude that
\begin{align*}
[1+o_{\epsilon}(1)]\,\textup{cap}_{Y}(A,\,B)\le & \liminf_{N\rightarrow\infty}N^{1+\alpha}\textup{cap}_{N}(\mathcal{E}_{N}(A),\,\mathcal{E}_{N}(B))\\
\le & \limsup_{N\rightarrow\infty}N^{1+\alpha}\textup{cap}_{N}(\mathcal{E}_{N}(A),\,\mathcal{E}_{N}(B))\le[1+o_{\epsilon}(1)]\,\textup{cap}_{Y}(A,\,B)\;,
\end{align*}
where the error terms $o_{\epsilon}(1)$ are now dependent only on
$\epsilon$. Since the two terms in the middle are independent of
$\epsilon$, by sending $\epsilon$ to $0$, we can complete the proof.
\end{proof}
From the previous proof, the estimate obtained in Proposition (\ref{p_zrpfn})
can be strengthened as follows.
\begin{cor}
\label{cor66}We have that
\[
\mathscr{D}_{N}(\mathbf{V}_{A,\,B})=\left(1+o_{N}(1)+o_{\epsilon}(1)\right)N^{-(1+\alpha)}\,\textup{cap}_{Y}(A,\,B)\;.
\]
\end{cor}

\section{\label{sec17}Estimate of Mean Jump Rates}

In this section, we verify (in Proposition \ref{p_H0}) the condition
(H0) for the zero-range process by estimating the mean jump rate $r_{N}(x,\,y)$
for $x,\,y\in S$.

For the reversible case, we can readily reduce the estimate of the
mean-jump rate to that of the capacity between valleys. More precisely,
it has been verified in \cite[Lemma 6.8]{B-L TM} that, for the reversible
case, that is, the case $p=1/2$, the mean jump rate satisfies the
following expression
\begin{equation}
r_{N}(x,\,y)=\frac{1}{2}\left[\textup{cap}_{N}(\mathcal{E}_{N}^{x},\,\breve{\mathcal{E}}_{N}^{x})+\textup{cap}_{N}(\mathcal{E}_{N}^{y},\,\breve{\mathcal{E}}_{N}^{y})-\textup{cap}_{N}(\mathcal{E}_{N}^{x}\cup\mathcal{E}_{N}^{y},\,\breve{\mathcal{E}}_{N}^{x,\,y})\right]\label{er1}
\end{equation}
for all $x,\,y\in S$. Hence, the estimate of the mean jump rate is
a direct consequence of Theorem \ref{t_zrpcap}.

Unfortunately, a the relationship (\ref{er1}) is no longer valid
in the non-reversible case and the estimation of the mean jump rate
$r_{N}(x,\,y)$ becomes a more challenging task. The general strategy
for this task has been developed in \cite[Section 8]{L-Seo NR1}.
The following is a summary of this strategy.
\begin{enumerate}
\item Define the mean holding rate by
\[
\lambda_{N}(x)=\sum_{y\in S\setminus\{x\}}r_{N}(x,\,y)\;.
\]
Then, in \cite[display (A.8)]{B-L TM2}, it has been verified that
the holding rate $\lambda_{N}(x)$ satisfies
\begin{equation}
\lambda_{N}(x)=\frac{\textup{cap}_{N}(\mathcal{E}_{N}^{x},\,\breve{\mathcal{E}}_{N}^{x})}{\mu_{N}(\mathcal{E}_{N}^{x})}\;.\label{e610}
\end{equation}
Therefore, by estimating the capacity $\textup{cap}_{N}(\mathcal{E}_{N}^{x},\,\breve{\mathcal{E}}_{N}^{x})$
and applying Theorem \ref{t01}, we can obtain the sharp asymptotics
of $\lambda_{N}(x)$.
\item The second step is to compute the sharp asymptotics of $r_{N}(x,\,y)/\lambda_{N}(x)$
using the collapsed process introduced in Section \ref{sec_col}.
More precisely, we fix $x\in S$, and we consider a process $\overline{\eta}_{N}(\cdot)$
which is the collapsed process obtained by collapsing the metastable
set $\mathcal{E}_{N}^{x}$ into a single point $\mathfrak{e}$. Denote
by $\overline{\mathbb{P}}_{\mathfrak{e}}^{N}$ the law of this collapsed
process starting from $\mathfrak{e}$. Then it has been proven in
\cite[Proposition 4.2]{B-L TM2} that
\begin{equation}
\frac{r_{N}(x,\,y)}{\lambda_{N}(x)}=\overline{\mathbb{P}}_{\mathfrak{e}}^{N}\bigl[\tau_{\mathcal{E}_{N}^{y}}<\tau_{\breve{\mathcal{E}}_{N}^{x,\,y}}\bigr]\;.\label{e611}
\end{equation}
Surprisingly, we can estimate the right-hand side based on the capacity
estimate for the collapsed process along with the sector condition
of the zero-range process which will be verified in Section \ref{sec_sector}.
\item Since we can obtain an estimate of $\lambda_{N}(x)$ and $r_{N}(x,\,y)/\lambda_{N}(x)$
by (\ref{e610}) and (\ref{e611}), we can finally obtain an estimate
of the mean jump rate $r_{N}(x,\,y)$. This argument is rigorously
explained in Section \ref{sec17}.
\end{enumerate}
In order to focus only on the effectiveness of potential theoretic
computations, we will not attempt to prove (\ref{e610}) and (\ref{e611})
in the current note; we refer to \cite{B-L TM,B-L TM2}. Instead,
we shall directly apply this strategy to verify the condition (H0)
for the zero-range processes. We note that we again assume Propositions
\ref{p_zrpfn} and \ref{p_zrpflow} (and hence all the results obtained
in previous sections) throughout this section.

\subsection{\label{sec_sector}Sector condition}

In this section, we prove the sector condition (cf. Definition \ref{sector})
for the zero-range process. This sector condition is one of the essential
ingredients of the method developed in \cite{L-Seo NR1} which will
be applied to the zero-range process in this section.
\begin{prop}
\label{psector}There exists a constant $C_{0}>0$ such that for all
$\mathbf{f},\,\mathbf{g}:\mathcal{H}_{N}\rightarrow\mathbb{R}$, we
have
\[
\left\langle \mathbf{g},\,-\mathcal{\mathscr{L}}_{N}\mathbf{f}\right\rangle _{\mu_{N}}^{2}\le C_{0}\mathscr{\,D}_{N}(\mathbf{f})\mathscr{\,D}_{N}(\mathbf{g})\;.
\]
\end{prop}

For $u\in S$, denote by $\omega^{u}=(\omega_{x}^{u})_{x\in S}\in\mathcal{H}_{1}$
the configuration with one particle at site $u$, namely,
\[
\omega_{x}^{u}=\begin{cases}
1 & \mbox{if }x=u\\
0 & \mbox{otherwise.}
\end{cases}
\]
Therefore, for $u\in S$ and $\eta\in\mathcal{H}_{N}$, the configuration
$\eta+\omega^{u}\in\mathcal{H}_{N+1}$ is the one obtained from $\eta$
by adding a particle from site $u$. Similarly, the configuration
$\eta-\omega^{u}\in\mathcal{H}_{N-1}$ is the one obtained from $\eta$
by removing a particle at site $u$, provided that $\eta_{u}\ge1$.

With this notation, we can observe the following convenient identity:
for $u\in S$ and for $\eta\in\mathcal{H}_{N}$ with $\eta_{u}\ge1$,
\begin{equation}
\mu_{N}(\eta)\,g(\eta_{u})=a_{N}\,\mu_{N-1}(\eta-\omega^{u})\;,\label{u1}
\end{equation}
where $a_{N}$ is defined by
\[
a_{N}=\frac{N^{\alpha}\,Z_{N-1}}{(N-1)^{\alpha}\,Z_{N}}\;.
\]
By Proposition \ref{p_zn}, it follows immediately that
\begin{equation}
\lim_{N\rightarrow\infty}a_{N}=1\;.\label{u2}
\end{equation}

\begin{proof}[Proof of Proposition \ref{psector}]
Fix $\mathbf{f},\,\mathbf{g}:\mathcal{H}_{N}\rightarrow\mathbb{R}$.
By (\ref{u1}) the change of variable $\eta-\omega^{x}=\zeta$, we
can write
\begin{align}
 & \mathscr{D}_{N}(\mathbf{f})\nonumber \\
= & \frac{1}{2}\sum_{\eta\in\mathcal{H}_{N}}\,\sum_{x\in S}\sum_{y\in S}\mu_{N}(\eta)\,g(\eta_{x})\,r(x,\,y)\left[\mathbf{f}(\eta)-\mathbf{f}(\sigma^{x,\,y}\eta)\right]^{2}\label{edr}\\
= & \frac{a_{N}}{2}\sum_{\zeta\in\mathcal{H}_{N-1}}\sum_{x\in S}\sum_{y\in S}\mu_{N-1}(\zeta)\,r(x,\,y)\left[\mathbf{f}(\zeta+\omega^{x})-\mathbf{f}(\zeta+\omega^{y})\right]^{2}\;.\nonumber
\end{align}
By a similar computation,

\begin{equation}
\begin{aligned} & \left\langle \mathbf{g},\,-\mathcal{\mathscr{L}}_{N}\mathbf{f}\right\rangle _{\mu_{N}}\\
= & \sum_{\eta\in\mathcal{H}_{N}}\,\sum_{x\in S}\sum_{y\in S}\mu_{N}(\eta)\,g(\eta_{x})\,r(x,\,y)\left[\mathbf{f}(\eta)-\mathbf{f}(\sigma^{x,\,y}\eta)\right]\mathbf{g}(\eta)\\
= & a_{N}\,\sum_{\zeta\in\mathcal{H}_{N-1}}\sum_{x\in S}\sum_{y\in S}\mu_{N-1}(\zeta)\,r(x,\,y)\left[\mathbf{f}(\zeta+\omega^{x})-\mathbf{f}(\zeta+\omega^{y})\right]\mathbf{g}(\zeta+\omega^{x})\;.
\end{aligned}
\label{ep41}
\end{equation}

For $\zeta\in\mathcal{H}_{N-1}$, write
\begin{equation}
\overline{\mathbf{g}}(\zeta)=\frac{1}{\kappa}\sum_{z\in S}\mathbf{g}(\zeta+\omega^{z})\;.\label{avg}
\end{equation}
Since we obviously have
\begin{align*}
\sum_{x,\,y\in S}r(x,\,y)\left[\mathbf{f}(\zeta+\omega^{x})-\mathbf{f}(\zeta+\omega^{y})\right] & =0\;,
\end{align*}
we can deduce from (\ref{ep41}) that

\begin{equation}
\begin{aligned} & \left\langle \mathbf{g},\,-\mathcal{\mathscr{L}}_{N}\mathbf{f}\right\rangle _{\mu_{N}}\\
= & a_{N}\sum_{\zeta\in\mathcal{H}_{N-1}}\sum_{x\in S}\sum_{y\in S}\mu_{N-1}(\zeta)\,r(x,\,y)\left[\mathbf{f}(\zeta+\omega^{x})-\mathbf{f}(\zeta+\omega^{y})\right]\left[\mathbf{g}(\zeta+\omega^{x})-\overline{\mathbf{g}}(\zeta)\right]\\
\le & \,\frac{a_{N}}{2}\sum_{\zeta\in\mathcal{H}_{N-1}}\sum_{x\in S}\sum_{y\in S}\mu_{N-1}(\zeta)r(x,\,y)\left(\left[\mathbf{f}(\zeta+\omega^{x})-\mathbf{f}(\zeta+\omega^{y})\right]^{2}+\left[\mathbf{g}(\zeta+\omega^{x})-\overline{\mathbf{g}}(\zeta)\right]^{2}\right)\\
= & \mathscr{D}_{N}(\mathbf{f})+\frac{a_{N}}{2}\sum_{\zeta\in\mathcal{H}_{N-1}}\sum_{x\in S}\mu_{N-1}(\zeta)\left[\mathbf{g}(\zeta+\omega^{x})-\overline{\mathbf{g}}(\zeta)\right]^{2}\;,
\end{aligned}
\label{ee1p}
\end{equation}
where the last line follows from (\ref{edr}) and the fact that $\sum_{y\in S}r(x,\,y)=1$.

Then, by (\ref{avg}), we can write
\begin{equation}
\sum_{x\in S}\left[\mathbf{g}(\zeta+\omega^{x})-\overline{\mathbf{g}}(\zeta)\right]^{2}=\frac{1}{\kappa}\sum_{u,\,v\in S}\left[\mathbf{g}(\zeta+\omega^{u})-\mathbf{g}(\zeta+\omega^{v})\right]^{2}\;.\label{ee1}
\end{equation}
by the Cauchy--Schwarz inequality,
\begin{align*}
\left[\mathbf{g}(\zeta+\omega^{u})-\mathbf{g}(\zeta+\omega^{v})\right]^{2} & =\left[\sum_{x=u}^{v-1}\mathbf{g}(\zeta+\omega^{x+1})-\mathbf{g}(\zeta+\omega^{x})\right]^{2}\\
 & \le\left[\sum_{x\in S}\left|\mathbf{g}(\zeta+\omega^{x+1})-\mathbf{g}(\zeta+\omega^{x})\right|\right]^{2}\\
 & \le\kappa\sum_{x\in S}[\mathbf{g}(\zeta+\omega^{x+1})-\mathbf{g}(\zeta+\omega^{x})]^{2}\;.
\end{align*}
Inserting this into (\ref{ee1}) yields that
\[
\sum_{x\in S}\left[\mathbf{g}(\zeta+\omega^{x})-\overline{\mathbf{g}}(\zeta)\right]^{2}\le\sum_{x\in S}[\mathbf{g}(\zeta+\omega^{x+1})-\mathbf{g}(\zeta+\omega^{x})]^{2}\;.
\]
Therefore, we have
\begin{align*}
 & \frac{a_{N}}{2}\sum_{\zeta\in\mathcal{H}_{N-1}}\sum_{x\in S}\mu_{N-1}(\zeta)\left[\mathbf{g}(\zeta+\omega^{x})-\overline{\mathbf{g}}(\zeta)\right]^{2}\\
\le & \frac{a_{N}}{2}\sum_{\zeta\in\mathcal{H}_{N-1}}\sum_{x\in S}\mu_{N-1}(\zeta)[\mathbf{g}(\zeta+\omega^{x+1})-\mathbf{g}(\zeta+\omega^{x})]^{2}\\
= & \frac{a_{N}}{2p}\sum_{\zeta\in\mathcal{H}_{N-1}}\sum_{x\in S}\mu_{N-1}(\zeta)r(x,\,x+1)[\mathbf{g}(\zeta+\omega^{x+1})-\mathbf{g}(\zeta+\omega^{x})]^{2}\\
\le & \frac{a_{N}}{2p}\sum_{\zeta\in\mathcal{H}_{N-1}}\sum_{x\in S}\sum_{y\in S}\mu_{N-1}(\zeta)r(x,\,y)[\mathbf{g}(\zeta+\omega^{y})-\mathbf{g}(\zeta+\omega^{x})]^{2}=\frac{1}{p}\mathscr{D}_{N}(\mathbf{g})\;,
\end{align*}
where the last line follows from (\ref{edr}).

Finally, inserting this into (\ref{ee1p}), we get
\[
\left\langle \mathbf{g},\,-\mathcal{\mathscr{L}}_{N}\mathbf{f}\right\rangle _{\mu_{N}}\le\mathcal{\mathscr{D}}_{N}(\mathbf{f})+\frac{1}{p}\mathcal{\mathscr{D}}_{N}(\mathbf{g})\;.
\]
By Remark \ref{rmk_sector}, we are done.
\end{proof}
Henceforth, the constant $C_{0}$ is always used to denote the constant
appearing in Proposition \ref{psector}. The following corollary is
now immediate from the above proposition and Propositions \ref{p_capcomp1}
and \ref{p_capcomp2}. We recall that $\textup{cap}_{N}^{s}(\cdot,\,\cdot)$
denotes the capacity with respect to the symmetrized process.
\begin{cor}
\label{csector}For any two disjoint, non-empty subsets $\mathcal{A},\,\mathcal{B}$
of $\mathcal{H}_{N}$,
\[
\textup{cap}_{N}^{s}(\mathcal{A},\,\mathcal{B})\le\textup{cap}_{N}(\mathcal{A},\,\mathcal{B})\le C_{0}\textup{cap}_{N}^{s}(\mathcal{A},\,\mathcal{B})\;.
\]
\end{cor}

\subsection{\label{sec64-1}Capacity estimates for collapsed processes}

Another essential ingredient of the method of \cite{L-Seo NR1} is
the sharp estimate of capacity with respect to the collapsed processes.
In this subsection, we explain this ingredient. In the remainder of
the current section, we will fix $x_{0}\in S$.

\subsubsection*{Definition of collapsed processes}

We first define collapsed processes and then explain the notation
regarding the collapsed process in terms of the zero-range processes.

Let $\overline{\mathcal{H}}_{N}=(\mathcal{H}_{N}\setminus\mathcal{E}_{N}^{x_{0}})\cup\{\mathfrak{e}\}$
be the set obtained from $\mathcal{H}_{N}$ by collapsing the metastable
set $\mathcal{E}_{N}^{x_{0}}$ into a single point $\mathfrak{e}$.
Denote by $(\overline{\eta}_{N}(t))_{t\ge0}$ the collapsed process
on $\overline{\mathcal{H}}_{N}$ which is obtained from $\eta_{N}(\cdot)$
by collapsing the set $\mathcal{E}_{N}^{x_{0}}$ to $\mathfrak{e}$.
Let $\overline{\mu}_{N}(\cdot)$ be a measure on $\overline{\mathcal{H}}_{N}$
defined by
\[
\begin{cases}
\overline{\mu}_{N}(\eta)=\mu_{N}(\eta) & \text{if }\eta\in\mathcal{H}_{N}\setminus\mathcal{E}_{N}^{x_{0}}\;,\\
\overline{\mu}_{N}(\mathfrak{e})=\mu_{N}(\mathcal{E}_{N}^{x_{0}})\;.
\end{cases}
\]
Then, by Exercise \ref{ex_collinv}, we get the following lemma.
\begin{lem}
The Markov chain $\overline{\eta}_{N}(\cdot)$ is irreducible on $\overline{\mathcal{H}}_{N}$,
and its unique invariant measure is $\overline{\mu}_{N}(\cdot)$.
\end{lem}

We now redefine the notation regarding the collapsed process in terms
of the zero-range process
\begin{itemize}
\item We denote by $\mathscr{\overline{L}}_{N}$ the generator of the collapsed
chain $\overline{\eta}_{N}(\cdot)$, and let $\mathscr{\overline{L}}_{N}^{\dagger}$
and $\overline{\mathscr{L}}_{N}^{\,s}$ denote the adjoint generator
and the symmetrized generator of $\mathscr{\overline{L}}_{N}$, respectively
(in the space $L^{2}(\overline{\mu}_{N})$). The continuous-time Markov
processes on $\overline{\mathcal{H}}_{N}$ generated by $\mathscr{\overline{L}}_{N}^{\dagger}$
and $\overline{\mathscr{L}}_{N}^{\,s}$ are denoted by $\overline{\eta}_{N}^{\dagger}(\cdot)$
and $\overline{\eta}_{N}^{s}(\cdot)$, respectively.
\item Let $\overline{\mathscr{D}}_{N}(\cdot)$ be the Dirichlet form associated
with the generator $\mathscr{\overline{L}}_{N}$.
\item Denote by $\overline{\mathbb{P}}_{\eta}^{N}$, $\eta\in\overline{\mathcal{H}}_{N}$,
the law of process $\overline{\eta}_{N}(\cdot)$ starting from $\eta$.
\item We denote by $\overline{\mathfrak{F}}_{N}$ the space of flow associated
with the collapsed process $\overline{\eta}_{N}(\cdot)$. The inner
product and flow norm associated with this flow structure will be
denoted by $\left\langle \cdot,\,\cdot\right\rangle _{\overline{\mathfrak{F}}_{N}}$
and $\Vert\cdot\Vert_{\overline{\mathfrak{F}}_{N}}$, respectively.
\item For each flow $\phi\in\mathfrak{F}_{N}$, we denote by $\overline{\phi}\in\overline{\mathfrak{F}}_{N}$
the collapsed flow in the sense of (\ref{eq:coll_flow}).
\item For each $\mathbf{f}:\mathcal{H}_{N}\rightarrow\mathbb{R}$ which
is constant over $\mathcal{E}_{N}^{x_{0}}$, we denote by $\overline{\mathbf{f}}:\overline{\mathcal{H}}_{N}\rightarrow\mathbb{R}$
the collapsed function in the sense of (\ref{col_f}).
\item For $\mathbf{f}:\overline{\mathcal{H}}_{N}\rightarrow\mathbb{R}$
we define flows $\overline{\Phi}_{\mathbf{f}}$, $\overline{\Phi}_{\mathbf{f}}^{*}$
and $\overline{\Psi}_{\mathbf{f}}$ as in (\ref{eq: col  phi_f})-(\ref{eq: col psi_f}).
\item For two disjoint non-empty subsets $\mathcal{A}$ and $\mathcal{B}$
of $\overline{\mathcal{H}}_{N}$, we denote by $\overline{\mathbf{h}}_{\mathcal{A},\mathcal{\,B}}$
and $\overline{\textup{cap}}_{N}(\mathcal{A},\,\mathcal{B})$ the
equilibrium potential and capacity between $\mathcal{A}$ and $\mathcal{B}$
with respect to the collapsed process $\overline{\eta}_{N}(\cdot)$.
In addition, we write $\overline{\textup{cap}}_{N}^{s}(\mathcal{A},\,\mathcal{B})$
for the capacity between $\mathcal{A}$ and $\mathcal{B}$ with respect
to process $\overline{\eta}_{N}^{s}(\cdot)$.
\end{itemize}
\begin{rem}
Notice that $\overline{\mathbf{h}}_{\mathcal{A},\mathcal{\,B}}$ and
$\overline{\mathbf{h}_{\mathcal{A},\,\mathcal{B}}}$ are different
objects. Since the equilibrium potential $\mathbf{h}_{\mathcal{A},\,\mathcal{B}}$
may not be constant on $\mathcal{E}_{N}^{x_{0}}$, we may not be able
to define the collapsed function $\overline{\mathbf{h}_{\mathcal{A},\,\mathcal{B}}}$.
\end{rem}

By Lemma \ref{lem:col_sector} and Proposition \ref{psector}, we
get the following proposition where $C_{0}$ is the constant appearing
in Proposition \ref{psector}
\begin{prop}
\label{prop_colsec}The collapsed process $\overline{\eta}_{N}(\cdot)$
satisfies a sector condition with constant $C_{0}$. Hence, for any
two disjoint non-empty subsets $\mathcal{A},\,\mathcal{B}$ of $\overline{\mathcal{H}}_{N}$,
it holds that
\[
\overline{\textup{cap}}_{N}^{s}(\mathcal{A},\,\mathcal{B})\le\textup{\ensuremath{\overline{\textup{cap}}}}_{N}(\mathcal{A},\,\mathcal{B})\le C_{0}\,\textup{\ensuremath{\overline{\textup{cap}}}}_{N}^{s}(\mathcal{A},\,\mathcal{B})\;.
\]
\end{prop}

\subsubsection*{Capacity estimates}

The following lemma, which is a direct consequence of Lemma \ref{lem:coll_CAP}
asserts that we are able to reduce the computation of capacity with
respect to the collapsed process to that of the original zero-range
process when one of the sets involved is $\{\mathfrak{e}\}$.
\begin{lem}
\label{lem613}For all non-empty subsets $\mathcal{A}$ of $\mathcal{H}_{N}\setminus\mathcal{E}_{N}^{x_{0}}$,
\[
\textup{\ensuremath{\overline{\textup{cap}}}}_{N}(\mathcal{A},\,\mathfrak{e})=\textup{cap}_{N}(\mathcal{A},\,\mathcal{E}_{N}^{x_{0}})\;.
\]
\end{lem}

In view of Exercise \ref{ex47}, the following estimate is not a simple
consequence of Theorem \ref{t_zrpcap} (or Corollary \ref{c02}).
We need an independent proof.
\begin{prop}
\label{p613}For two disjoint and non-empty subsets $A$ and $B$
of $S\setminus\{x_{0}\}$ satisfying $A\cup B=S\setminus\{x_{0}\}$,
it holds that
\[
\textup{\ensuremath{\overline{\textup{cap}}}}_{N}(\mathcal{E}_{N}(A),\,\mathcal{E}_{N}(B))=[1+o_{N}(1)+o_{\epsilon}(1)]\,N^{-(1+\alpha)}\,\textup{cap}_{Y}(A,\,B)\;.
\]
\end{prop}

The proof of this proposition will be given in next subsection.

\subsection{Capacity estimates for collapsed processes}

We now prove Proposition \ref{p613} by several steps. Throughout
this subsection, we fix two disjoint and non-empty subsets $A,\,B$
satisfying the condition of Proposition \ref{p613}. Recall the test
functions $\mathbf{V}_{A,\,B}$ and $\mathbf{V}_{A,\,B}^{\dagger}$
from Proposition \ref{p_zrpfn} and the test flows $\Phi_{A,\,B}$
and $\Phi_{A,\,B}^{\dagger}$ from Proposition \ref{p_zrpflow}. Since
$\mathbf{V}_{A,\,B}$ and $\mathbf{V}_{A,\,B}^{\dagger}$ are constant
on $\mathcal{E}_{N}^{x_{0}},$we can collapse them; let us write $\overline{\mathbf{V}}_{A,\,B}=\overline{\mathbf{V}_{A,\,B}}$
and $\overline{\mathbf{V}}_{A,\,B}^{\dagger}=\overline{\mathbf{V}_{A,\,B}^{\dagger}}$.
Note that, by Proposition \ref{p_zrpfn} we have that
\[
\overline{\mathbf{V}}_{A,\,B}(\mathfrak{e})=\mathbf{\overline{V}}_{A,\,B}^{\dagger}(\mathfrak{e})=\mathfrak{h}_{A,\,B}(x_{0})\;.
\]

\begin{lem}
\label{lem615} It holds that
\[
\bigl\Vert\overline{\Psi}_{\overline{\mathbf{V}}_{A,\,B}}\bigr\Vert_{\overline{\mathcal{\mathfrak{F}}}_{N}}^{2}=[1+o_{N}(1)+o_{\epsilon}(1)]\,N^{-(1+\alpha)}\textup{\,cap}_{Y}(A,\,B)\;.
\]
\end{lem}

\begin{proof}
By Exercise \ref{ex:eqcon}\textbf{ }and Lemma \ref{lem:coll_Phi_f},
we obtain
\[
\bigl\Vert\overline{\Psi}_{\overline{\mathbf{V}}_{A,\,B}}\bigr\Vert_{\mathfrak{\overline{F}}_{N}}^{2}=\bigl\Vert\overline{\Psi_{\mathbf{V}_{A,\,B}}}\bigr\Vert_{\mathfrak{\overline{F}}_{N}}^{2}=\bigl\Vert\Psi_{\mathbf{V}_{A,\,B}}\bigr\Vert^{2}\;.
\]
It is now enough to invoke Corollary \ref{cor66} to complete the
proof.

Let $\overline{\Phi}_{A,\,B}=\overline{\Phi_{A,\,B}}$ and $\overline{\Phi}_{A,\,B}^{\dagger}=\overline{\Phi_{A,\,B}^{\dagger}}$
be the collapsed flow of $\Phi_{A,\,B}$ of $\Phi_{A,\,B}^{\dagger}$,
respectively.
\end{proof}
\begin{lem}
\label{lem617}It holds that
\begin{align}
\sum_{\eta\in\overline{\mathcal{H}}_{N}}\overline{\mathbf{\mathbf{h}}}_{\mathcal{E}_{N}(A),\,\mathcal{E}_{N}(B)}(\eta)(\textup{div }\overline{\Phi}_{A,\,B})(\eta) & =[1+o_{N}(1)]\,N^{-(1+\alpha)}\,\textup{cap}_{Y}(A,\,B)\;\;\text{and}\label{eq:eh1-1}\\
\sum_{\eta\in\overline{\mathcal{H}}_{N}}\overline{\mathbf{\mathbf{h}}}_{\mathcal{E}_{N}(A),\,\mathcal{E}_{N}(B)}(\eta)(\textup{div }\overline{\Phi}_{A,\,B}^{\dagger})(\eta) & =[1+o_{N}(1)]\,N^{-(1+\alpha)}\,\textup{cap}_{Y}(A,\,B)\;.\label{eh2-1}
\end{align}
\end{lem}

\begin{proof}
It suffices to prove (\ref{eq:eh1-1}) as the proof of (\ref{eh2-1})
is essentially the same. In view of Lemma \ref{lem66}, it suffices
to check
\begin{equation}
\overline{\mathbf{\mathbf{h}}}_{\mathcal{E}_{N}(A),\,\mathcal{E}_{N}(B)}(\mathfrak{e})(\textup{div }\overline{\Phi}_{A,\,B})(\mathfrak{e})-\sum_{\eta\in\mathcal{E}_{N}^{x_{0}}}\mathbf{h}_{\mathcal{E}_{N}(A),\,\mathcal{E}_{N}(B)}(\eta)(\textup{div }\Phi_{A,\,B})(\eta)=o_{N}(1)N^{-(1+\alpha)}\;.\label{e4123}
\end{equation}
By (\ref{eq:coll_div}) and (\ref{e4120}), we have
\[
(\textup{div }\overline{\Phi}_{A,\,B})(\mathfrak{e})=(\textup{div }\Phi_{A,\,B})(\mathcal{E}_{N}^{x_{0}})=o_{N}(1)N^{-(1+\alpha)}\;.
\]
Thus, the first term at the left-hand side of (\ref{e4123}) is $o_{N}(1)N^{-(1+\alpha)}$.
On the other hand, the second term is $o_{N}(1)N^{-(1+\alpha)}$ by
(\ref{e412}). Hence, we have (\ref{e4123}).
\end{proof}
Now we are ready to prove Proposition \ref{p613} by using generalized
Dirichlet and Thomson principles.
\begin{proof}[Proof of Proposition \ref{p613}]
The proof is similar to that of Theorem \ref{t_zrpcap}. We begin
by recalling the functions $\mathbf{f},\,\mathbf{g}$ and the flows
$\phi,\,\psi$ from (\ref{e_test1}) and (\ref{e_test2}). Then, by
the definition of the collapsing procedure, it is obvious that
\[
\overline{\mathbf{f}}\in\mathfrak{C}_{1,\,0}(\mathcal{E}_{N}(A),\,\mathcal{\mathcal{E}}_{N}(B))\;\;\;\;\text{and\;\;}\;\;\overline{\mathbf{g}}\in\mathfrak{C}_{0,\,0}(\mathcal{E}_{N}(A),\,\mathcal{\mathcal{E}}_{N}(B))\;.
\]
Since we can write
\begin{equation}
\overline{\Phi}_{\,\mathbf{\overline{f}}}-\overline{\phi}=\overline{\Psi}_{\mathbf{\overline{V}}_{A,\,B}}-\frac{\overline{\Theta}_{N}^{\dagger}-\overline{\Theta}_{N}}{2}\;,\label{e6131}
\end{equation}
where $\overline{\Theta}_{N}$ and $\overline{\Theta}_{N}^{\dagger}$
are the collapsed flows of $\Theta_{N}$ and $\Theta_{N}^{\dagger}$
defined in (\ref{ec2}), respectively. By part (1) of Proposition
\ref{p_zrpflow} and Lemma \ref{lem:coll_flow}, we have
\begin{equation}
\bigl\Vert\overline{\Theta}_{N}\bigr\Vert_{\overline{\mathfrak{F}}_{N}}^{2}=\left(o_{N}(1)+o_{\epsilon}(1)\right)N^{-(1+\alpha)}\;\;\mbox{and}\;\;\bigl\Vert\overline{\Theta}_{N}^{\dagger}\bigr\Vert_{\overline{\mathfrak{F}}_{N}}^{2}=\left(o_{N}(1)+o_{\epsilon}(1)\right)N^{-(1+\alpha)}\;.\label{e6132}
\end{equation}
Thus, by Theorem \ref{t_genDTP_nonrev}-(1), Lemma \ref{lem615},
and Lemma \ref{lem617}, we get the following upper bound:
\begin{equation}
\textup{\ensuremath{\overline{\textup{cap}}}}_{N}(\mathcal{E}_{N}(A),\,\mathcal{E}_{N}(B))\le[1+o_{N}(1)+o_{\epsilon}(1)]\,N^{-(1+\alpha)}\,\textup{cap}_{Y}(A,\,B)\;.\label{epp1}
\end{equation}

For the opposite inequality, we can repeat the same arguments with
test function $\overline{\mathbf{g}}$ and test flow $\overline{\psi}$
to deduce
\begin{equation}
\textup{\ensuremath{\overline{\textup{cap}}}}_{N}(\mathcal{E}_{N}(A),\,\mathcal{E}_{N}(B))\ge\left(1+o_{N}(1)+o_{\epsilon}(1)\right)N^{-(1+\alpha)}\,\textup{cap}_{Y}(A,\,B)\;.\label{epp2}
\end{equation}
By (\ref{epp1}) and (\ref{epp2}), the proof is completed.
\end{proof}
\begin{xca}
Prove (\ref{epp2}) by using the generalized Thomson principle.
\end{xca}

\begin{xca}
In fact, the condition $A\cup B=S\setminus\{x_{0}\}$ in Proposition
\ref{p613} is redundant. We imposed this condition only because we
do not need a general result without this restriction. Prove the general
result without this restriction.
\end{xca}

\subsection{\label{sec66-1}Estimate of mean jump rate}

Now we are ready to estimate the mean jump rate. In view of (\ref{e611}),
to obtain the sharp asymptotics of the mean jump rate $r_{N}(x_{0},\,y)$
for $x_{0},\,y\in S$, the crucial object to be estimated is the probability
$\overline{\mathbb{P}}_{\mathfrak{e}}^{N}[\tau_{\mathcal{E}_{N}^{y}}<\tau_{\breve{\mathcal{E}}_{N}^{x_{0},\,y}}]$.
This estimate follows from the following proposition.
\begin{prop}
\label{p618}For two disjoint and non-empty subsets $A$ and $B$
of $S\setminus\{x_{0}\}$ satisfying $A\cup B=S\setminus\{x_{0}\}$,
we have that
\[
\lim_{N\rightarrow\infty}\overline{\mathbb{P}}_{\mathfrak{e}}^{N}\left[\tau_{\mathcal{E}_{N}(A)}<\tau_{\mathcal{E}_{N}(B)}\right]=\mathfrak{h}_{A,\,B}(x_{0})\;.
\]
\end{prop}

\begin{proof}
The proof relies on Propositions \ref{prop_colsec}, \ref{p613} and
Lemma \ref{lem615}. Recall the equilibrium potential $\overline{\mathbf{h}}_{\mathcal{E}_{N}(A),\mathcal{\,E}_{N}(B)}$
between $\mathcal{E}_{N}(A)$ and $\mathcal{E}_{N}(B)$, with respect
to the collapsed chain $\overline{\eta}_{N}(\cdot)$. Then, by Proposition
\ref{p613},
\begin{equation}
\left\Vert \overline{\Psi}_{\overline{\mathbf{h}}_{\mathcal{E}_{N}(A),\,\mathcal{E}_{N}(B)}}\right\Vert _{\overline{\mathfrak{F}}_{N}}^{2}=\overline{\textup{cap}}_{N}(\mathcal{E}_{N}(A),\mathcal{\,E}_{N}(B))=[1+o_{N}(1)+o_{\epsilon}(1)]\,N^{-(1+\alpha)}\textup{\,cap}_{Y}(A,\,B)\;.\label{ek1}
\end{equation}
By Lemma \ref{lem615},
\begin{equation}
\left\Vert \,\overline{\Psi}_{\overline{\mathbf{V}}_{A,\,B}}\right\Vert _{\overline{\mathfrak{F}}_{N}}^{2}=[1+o_{N}(1)+o_{\epsilon}(1)]\,N^{-(1+\alpha)}\textup{\,cap}_{Y}(A,\,B)\;.\label{ek2}
\end{equation}
By (\ref{e6131}), (\ref{e6132}), (\ref{ek1}), and the Cauchy-Schwarz
inequality, we get
\begin{equation}
\begin{aligned} & \left\langle \overline{\Psi}_{\,\overline{\mathbf{V}}_{A,\,B}},\,\overline{\Psi}_{\overline{\mathbf{h}}_{\mathcal{E}_{N}(A),\,\mathcal{E}_{N}(B)}}\right\rangle _{\overline{\mathfrak{F}}_{N}}\\
 & =\left\langle \overline{\Phi}_{\mathbf{\overline{f}}}-\overline{\phi},\,\overline{\Psi}_{\overline{\mathbf{h}}_{\mathcal{E}_{N}(A),\,\mathcal{E}_{N}(B)}}\right\rangle _{\overline{\mathfrak{F}}_{N}}+\left(o_{N}(1)+o_{\epsilon}(1)\right)N^{-(1+\alpha)}\;,
\end{aligned}
\label{e4511}
\end{equation}
where $\overline{\mathbf{f}}$ and $\overline{\phi}$ are the objects
defined in the proof of Proposition \ref{p613}. By the same computation
as in (\ref{e451}), we can write
\begin{equation}
\begin{aligned} & \left\langle \overline{\Phi}_{\mathbf{\,\overline{f}}}-\overline{\phi},\,\overline{\Psi}_{\overline{\mathbf{h}}_{\mathcal{E}_{N}(A),\,\mathcal{E}_{N}(B)}}\right\rangle _{\overline{\mathfrak{F}}_{N}}\\
 & =\overline{\textup{cap}}_{N}(\mathcal{E}_{N}(A),\,\mathcal{E}_{N}(B))-\sum_{\eta\in\mathcal{\overline{H}}_{N}\setminus\mathcal{E}_{N}(A\cup B)}\overline{\mathbf{h}}_{\mathcal{E}_{N}(A),\,\mathcal{E}_{N}(B)}(\eta)\,(\textup{div }\overline{\phi})(\eta)
\end{aligned}
\label{e4512}
\end{equation}
Thus, by combining (\ref{e4511}), (\ref{e4512}) and Proposition
\ref{p613}, we get ,
\begin{equation}
\left\langle \overline{\Psi}_{\overline{\mathbf{V}}_{A,\,B}},\,\overline{\Psi}_{\overline{\mathbf{h}}_{\mathcal{E}_{N}(A),\,\mathcal{E}_{N}(B)}}\right\rangle _{\overline{\mathfrak{F}}_{N}}=\left(1+o_{N}(1)+o_{\epsilon}(1)\right)N^{-(1+\alpha)}\,\textup{cap}_{Y}(A,\,B)\label{ek3}
\end{equation}

Define $\mathbf{u}=\overline{\mathbf{h}}_{\mathcal{E}_{N}(A),\,\mathcal{E}_{N}(B)}-\overline{\mathbf{V}}_{A,\,B}$.
Then, by (\ref{ek1}), (\ref{ek2}), and (\ref{ek3}) we get
\begin{equation}
\begin{aligned}\bigl\Vert\overline{\Psi}_{\mathbf{u}}\bigr\Vert_{\overline{\mathfrak{F}}_{N}}^{2} & =\bigl\Vert\,\overline{\Psi}_{\,\overline{\mathbf{h}}_{\mathcal{E}_{N}(A),\,\mathcal{E}_{N}(B)}}\bigr\Vert_{\overline{\mathfrak{F}}_{N}}^{2}+\bigl\Vert\,\overline{\Psi}_{\,\overline{\mathbf{V}}_{A,\,B}}\bigr\Vert_{\overline{\mathfrak{F}}_{N}}^{2}-2\left\langle \overline{\Psi}_{\,\overline{\mathbf{V}}_{A,\,B}},\,\overline{\Psi}_{\,\overline{\mathbf{h}}_{\mathcal{E}_{N}(A),\mathcal{\,E}_{N}(B)}}\right\rangle _{\overline{\mathfrak{F}}_{N}}\\
 & =\left(o_{N}(1)+o_{\epsilon}(1)\right)N^{-(1+\alpha)}\;.
\end{aligned}
\label{ekk1}
\end{equation}
As $\mathbf{u}(\mathfrak{e})=\overline{\mathbf{h}}_{\mathcal{E}_{N}(A),\,\mathcal{E}_{N}(B)}(\mathfrak{e})-\mathfrak{h}_{A,\,B}(x_{0})$
and $\mathbf{u}(\eta)=0$ for all $\eta\in\mathcal{E}_{N}(A\cup B)$,
we can write
\[
\mathbf{u}=\left(\,\overline{\mathbf{h}}_{\mathcal{E}_{N}(A),\,\mathcal{E}_{N}(B)}(\mathfrak{e})-\mathfrak{h}_{A,\,B}(x)\right)\mathbf{u}_{0}
\]
for some $\mathbf{u}_{0}\in\mathfrak{C}_{1,\,0}(\{\mathfrak{o}\},\,\mathcal{E}_{N}(A\cup B))$.
With this notation, we can write
\begin{equation}
\bigl\Vert\overline{\Psi}_{\mathbf{u}}\bigr\Vert_{\overline{\mathfrak{F}}_{N}}^{2}=\overline{\mathscr{D}}_{N}(\mathbf{u})=\left(\,\overline{\mathbf{h}}_{\mathcal{E}_{N}(A),\,\mathcal{E}_{N}(B)}(\mathfrak{e})-\mathfrak{h}_{A,\,B}(x)\right)^{2}\,\overline{\mathscr{D}}_{N}(\mathbf{u}_{0})\;.\label{ekk2}
\end{equation}
By the Dirichlet principle for reversible dynamics (cf. Theorem \ref{t_DP_rev})
and the sector condition for the collapsed process (cf. Proposition
\ref{prop_colsec}), we have that
\begin{equation}
\overline{\mathscr{D}}_{N}(\mathbf{u}_{0})\ge\overline{\textup{cap}}_{N}^{\,s}(\mathfrak{e},\,\mathcal{E}_{N}(A\cup B))\ge C_{0}^{-1}\,\overline{\textup{cap}}_{N}(\mathfrak{e},\,\mathcal{E}_{N}(A\cup B))\;.\label{ekk21}
\end{equation}
By Lemma \ref{lem613} and Theorem \ref{t_zrpcap},
\begin{align}
\overline{\textup{cap}}_{N}(\mathfrak{e},\,\mathcal{E}_{N}(A\cup B)) & =\textup{cap}_{N}(\mathcal{E}_{N}^{x},\,\mathcal{E}_{N}(A\cup B))\nonumber \\
 & =[1+o_{N}(1)+o_{\epsilon}(1)]\,N^{-(1+\alpha)}\textup{\,cap}_{Y}(x,\,A\cup B)\;.\label{ekk22}
\end{align}
By (\ref{ekk21}) and (\ref{ekk22}), we can conclude that
\[
\overline{\mathscr{D}}_{N}(\mathbf{u}_{0})\ge C\,[1+o_{N}(1)+o_{\epsilon}(1)]\,N^{-(1+\alpha)}
\]
for some constant $C>0$. Inserting this and (\ref{ekk1}) into (\ref{ekk2}),
we get
\[
\left[\,\overline{\mathbf{h}}_{\mathcal{E}_{N}(A),\,\mathcal{E}_{N}(B)}(\mathfrak{e})-\mathfrak{h}_{A,\,B}(x)\right]^{2}\le o_{N}(1)+o_{\epsilon}(1)\;.
\]
By taking $\limsup_{N\rightarrow\infty}$ and then $\limsup_{\epsilon\rightarrow0}$,
we get
\[
\limsup_{N\rightarrow\infty}\left|\overline{\mathbf{h}}_{\mathcal{E}_{N}(A),\,\mathcal{E}_{N}(B)}(\mathfrak{e})-\mathfrak{h}_{A,\,B}(x)\right|=0
\]
and we are done.
\end{proof}
Now we are ready to verify condition (H0) for the zero-range process.
\begin{prop}
\label{p_H0}The condition (H0) holds for the zero-range processes.
In other words, for all $x,\,y\in S$,
\[
\lim_{N\rightarrow\infty}N^{1+\alpha}\,r_{N}(x,\,y)=a(x,\,y)\;.
\]
\end{prop}

\begin{proof}
By (\ref{e610}), Theorem \ref{t01}, and Corollary \ref{c02}, we
get
\begin{equation}
\lambda_{N}(x)=\frac{\textup{cap}_{N}(\mathcal{E}_{N}^{x},\,\breve{\mathcal{E}}_{N}^{x})}{\mu(\mathcal{E}_{N}^{x})}=\left(1+o_{N}(1)\right)N^{-(1+\alpha)}\,\frac{\kappa}{\Gamma_{\alpha}I_{\alpha}}\sum_{y\in S\setminus\{x\}}\textup{cap}_{X}(x,\,y)\;.\label{e6191}
\end{equation}
Recall from (\ref{hu4}) the definition of $\mathfrak{h}_{y,\,S\setminus\{x,\,y\}}$.
Write
\[
\tau=\inf\left\{ t\ge0:Y(t)\neq Y(0)\right\} \;.
\]
Then, one can observe that
\[
\mathfrak{h}_{y,\,S\setminus\{x,\,y\}}(x)=\mathbf{Q}_{x}\left(Y(\tau)=y\right)=\,\frac{\textup{cap}_{X}(x,\,y)}{\sum_{y\in S\setminus\{x\}}\textup{cap}_{X}(x,\,y)}\;.
\]
Thus, by (\ref{e611}) and Proposition \ref{p618}, we get
\begin{equation}
\frac{r_{N}(x,\,y)}{\lambda_{N}(x)}=\left(1+o_{N}(1)\right)\mathfrak{h}_{y,\,S\setminus\{x,\,y\}}(x)=\left(1+o_{N}(1)\right)\,\frac{\textup{cap}_{X}(x,\,y)}{\sum_{y\in S\setminus\{x\}}\textup{cap}_{X}(x,\,y)}\;.\label{e6192}
\end{equation}
We can complete the proof by multiplying (\ref{e6191}) and (\ref{e6192}).
\end{proof}

\section{\label{sec18}Conditions (H1) and (H3)}

Since we have verified conditions (H0) and (H2), it now remains to
verify conditions (H1) and (H3). Verification of these conditions
also use the capacity estimate obtained in Theorem \ref{t_zrpcap}
and the sector condition obtained in Proposition \ref{prop_colsec}.
We again assume the results obtained in Section \ref{sec16}.

We first prove the following lemma.
\begin{lem}
\label{lem_Cap}For any $x\in S$, there exists a constant $C$ such
that
\[
\inf_{\eta,\,\zeta\in\mathcal{E}_{N}^{x}}\textup{cap}_{N}(\eta,\,\zeta)\ge\frac{C}{\ell_{N}^{\alpha(\kappa-1)+1}}
\]
\end{lem}

\begin{proof}
We fix $x\in S$ and $\eta,\,\zeta\in\mathcal{E}_{N}^{x}$. We first
find a lower bound for $\textup{cap}_{N}^{s}(\eta,\,\zeta)$. For
$\xi,\,\xi'\in\mathcal{H}_{N}$, we denote by $R_{N}(\xi,\,\xi')$
the jump rate of the symmetrized zero-range process from $\xi$ to
$\xi'$:
\[
R_{N}(\xi,\,\xi')=\sum_{x\in S}\sum_{y\in S}g(\xi_{x})r^{s}(x,\,y)\mathbf{1}\{\xi'=\sigma^{x,\,y}\xi\}\;,
\]
where $r^{s}(x,\,y)=\frac{1}{2}\mathbf{1}\{|x-y|=1\}$. Take a path
$(\omega_{t})_{t=0}^{T}$ in $\mathcal{E}_{N}^{x}$ connecting $\eta$
and $\zeta$ in the sense that $\omega_{t}\in\mathcal{E}_{N}^{x}$
for all $t\in\llbracket0,\,T\rrbracket$ and moreover satisfies
\[
\omega_{0}=\eta\;,\quad\;\omega_{T}=\zeta\;\quad\text{and }\;\;R_{N}(\omega_{t},\,\omega_{t+1})>0\text{ for all }t\in\llbracket0,\,T-1\rrbracket\;.
\]
The existence of such a path with $T\le C\ell_{N}$ where $C$ is
a constant that only depends on $\kappa$ is obvious. Define a flow
$\phi\in\mathfrak{F}_{N}$ by
\[
\phi(\xi,\,\xi')=\begin{cases}
1 & \text{if }(\xi,\,\xi')=(\omega_{t},\,\omega_{t+1})\text{ for some }t\in\llbracket0,\,T-1\rrbracket\;,\\
-1 & \text{if }(\xi,\,\xi')=(\omega_{t+1},\,\omega_{t})\text{ for some }t\in\llbracket0,\,T-1\rrbracket\;,\\
0 & \text{otherwise.}
\end{cases}
\]
Then,
\begin{equation}
\Vert\phi\Vert_{\mathfrak{F}_{N}}^{2}=\sum_{t=0}^{T-1}\frac{1}{\mu_{N}(\omega_{t})R_{N}(\omega_{t},\,\omega_{t+1})}\;.\label{ebb33}
\end{equation}
Since $g(k)\ge1$ for all $k\ge1$, if $\omega_{t+1}=\sigma^{x,\,y}\omega_{t}$
for some $x,\,y\in S$ with $|y-x|=1$,
\begin{equation}
\mu_{N}(\omega_{t})R_{N}(\omega_{t},\,\omega_{t+1})=\frac{N^{\alpha}}{Z_{N}}\frac{1}{a(\omega_{t})}\times\frac{1}{2}g((\omega_{t})_{x})\ge C\frac{N^{\alpha}}{N^{\alpha}\ell_{N}^{\alpha(\kappa-1)}}=C\frac{1}{\ell_{N}^{\alpha(\kappa-1)}}\;,\label{ebb34}
\end{equation}
where we use a trivial bound
\[
a(\xi)=a(\xi_{x})\prod_{y\in S\setminus\{x\}}a(\xi_{y})\le N^{\alpha}\ell_{N}^{\alpha(\kappa-1)}\;\;\;\text{for all }\xi\in\mathcal{E}_{N}^{x}
\]
and Proposition \ref{p_zn} at the inequality of (\ref{ebb34}). Inserting
(\ref{ebb34}) and the bound $T\le C\ell_{N}$ to (\ref{ebb33}),
we get
\[
\Vert\phi\Vert_{\mathfrak{F}_{N}}^{2}\le C\ell_{N}\times\ell_{N}^{\alpha(\kappa-1)}=\text{\ensuremath{C\ell_{N}^{\alpha(\kappa-1)+1}}}
\]
Since $\phi$ is the unit flow from $\{\eta\}$ to $\{\zeta\}$, by
the Thomson principle for the reversible Markov process (cf. Theorem
\ref{t_TP_rev}),
\[
\textup{cap}_{N}^{s}(\eta,\,\zeta)\ge\frac{1}{\Vert\phi\Vert_{\mathfrak{F}_{N}}^{2}}\ge\frac{C}{\ell_{N}^{\alpha(\kappa-1)+1}}\;.
\]
Now the proof of lemma is completed by Corollary \ref{csector}.
\end{proof}
\begin{xca}
In the previous proof, prove the existence of a path $(\omega_{t})_{t=0}^{T}$
in $\mathcal{E}_{N}^{x}$ connecting $\eta$ and $\zeta$ with $T\le C\ell_{N}$
for some constant $C$ depending only on $\kappa$.
\end{xca}

Now we verify condition (H1).
\begin{prop}
\label{p_H1}The condition (H1) holds for the zero-range processes.
\end{prop}

\begin{proof}
Fix $x\in S$. For $\eta,\,\zeta\in\mathcal{E}_{N}^{x}$, by Theorem
\ref{t_zrpcap} and Lemma \ref{lem_Cap}, there exists $C>0$ such
that
\begin{equation}
\frac{\textup{cap}_{N}(\mathcal{E}_{N}^{x},\,\breve{\mathcal{E}}_{N}^{x})}{\textup{cap}_{N}(\eta,\,\zeta)}\le C\frac{\ell_{N}^{\alpha(\kappa-1)+1}}{N^{1+\alpha}}=o_{N}(1)\label{H1-1}
\end{equation}
where the last equality follows from the condition (\ref{cond_ell})
on $\ell_{N}$.
\end{proof}
At this moment, we shall check that the condition (H3) is in force
for the zero-range processes.
\begin{prop}
\label{p_H3}The condition (H3) holds for the zero-range processes.
\end{prop}

\begin{proof}
Fix $x\in S$. Recall that $\xi_{N}^{x}\in\mathcal{E}_{N}^{x}$ represent
a configuration such that all the particles are located at site $x$.
By \cite[Lemma 3.4]{L-L-M}, it suffices to verify that
\begin{align}
 & \lim_{N\rightarrow\infty}\sup_{\eta\in\mathcal{E}_{N}^{x}}\mathbb{P}_{\eta}^{N}[\tau_{\xi_{N}^{x}}>N^{1+\alpha}\delta]=0\;\;\text{for all \ensuremath{\delta>0\;,\;\;\text{and}}}\label{fd1}\\
 & \lim_{\delta\rightarrow0}\limsup_{N\rightarrow\infty}\sup_{\delta<t<3\delta}\mathbb{P}_{\xi_{N}^{x}}^{N}[\eta_{N}(N^{1+\alpha}t)\in\Delta_{N}]=0\;.\label{fd2}
\end{align}
For (\ref{fd1}), by the Markov inequality and (\ref{exphit2}), we
have
\begin{equation}
\mathbb{P}_{\eta}^{N}[\tau_{\xi_{N}^{x}}>N^{1+\alpha}\delta]\le\frac{1}{N^{1+\alpha}\delta}\mathbb{E}_{\eta}^{N}\left[\tau_{\xi_{N}^{x}}\right]\le\frac{1}{N^{1+\alpha}\,\delta}\frac{1}{\textup{cap}_{N}(\eta,\,\xi_{N}^{x})}\;,\label{fd3}
\end{equation}
where at the second inequality we use the trivial bound $h_{\eta,\,\xi_{N}^{x}}\le1$.
By Lemma \ref{lem_Cap},
\[
\mathbb{P}_{\eta}^{N}[\tau_{\xi_{N}^{x}}>N^{1+\alpha}\delta]\le\frac{C}{\delta}\frac{\ell_{N}^{\alpha(\kappa-1)+1}}{N^{1+\alpha}}\;.
\]
The proof of (\ref{fd1}) now follows from the condition (\ref{cond_ell})
on $\ell_{N}$.

For (\ref{fd2}), note first from the definition of $\mu_{N}$ that
we have $\mu_{N}(\xi_{N}^{x})=Z_{N}^{-1}$. Hence, for $t>0$, since
$\mu_{N}$ is the invariant measure,
\[
\mathbb{P}_{\xi_{N}^{x}}^{N}\left[\eta_{N}(N^{1+\alpha}\,t)\in\Delta_{N}\right]\le\frac{\mathbb{P}_{\mu_{N}}^{N}\left[\eta_{N}(N^{1+\alpha}\,t)\in\Delta_{N}\right]}{\mu_{N}(\xi_{N}^{x})}=\frac{\mu_{N}(\Delta_{N})}{\mu_{N}(\xi_{N}^{x})}=Z_{N}\,\mu_{N}(\Delta_{N})\;.
\]
Hence, (\ref{fd2}) follows directly from Proposition \ref{p_zn}
and Theorem \ref{t01}.
\end{proof}

\end{document}